\documentclass[10pt]{amsart}
\usepackage{latexsym,amssymb,color,amsmath,amsthm,graphicx,listings,comment,wrapfig}
\usepackage[section]{placeins}
   
\definecolor{red1}{rgb}{1,0.9,0.9} \definecolor{blue1}{rgb}{0.9,0.9,1} \definecolor{green1}{rgb}{0.9,1,0.9} 
\definecolor{yellow1}{rgb}{1,1,0.9} \definecolor{yellow2}{rgb}{1,1,0.8}

\def\question#1{ \vspace{2mm} \begin{center} \fcolorbox{green1}{green1}{ \parbox{11.2cm}{{\bf Question:} #1}} \vspace{2mm} \end{center} }
\def\conjecture#1{ \vspace{2mm} \begin{center} \fcolorbox{green1}{green1}{ \parbox{11.2cm}{{\bf Conjecture:} #1}} \vspace{2mm} \end{center} }
\def\resultremark#1{ \vspace{2mm} \begin{center} \fcolorbox{yellow1}{yellow1}{ \parbox{11.2cm}{{\bf Remark:} #1}} \vspace{2mm} \end{center} }

\title{Some experiments in number theory}
\author{Oliver Knill}
\date{Jun 19, 2016}
\address{Department of Mathematics \\ Harvard University \\ Cambridge, MA, 02138, USA }

\begin{document}
\maketitle

\begin{abstract}
Following footsteps of Gauss, Euler, Riemann, Hurwitz, Smith, Hardy, Littlewood, Hedlund, Khinchin and
Chebyshev, we visit some topics in elementary number theory. For matrices defined by Gaussian primes we 
observe a circular spectral law for the eigenvalues. We experiment then with various
Goldbach conjectures for Gaussian primes, Eisenstein primes, Hurwitz primes or Octavian primes. These conjectures
relate with Landau or Bunyakovsky or Andrica type conjectures for rational primes.
The Landau problem asking whether infinitely many predecessors of primes are square 
is also related to a determinant problem for the prime matrices under consideration. 
Some of these matrices are adjacency matrices of bipartite graphs. Their Euler characteristics in turn
is related to the prime counting function. When doing statistics of Gaussian primes on rows,
we detect a sign of correlations: rows of even distance for example look asymptotically correlated.
The expectation values of prime densities were conjectured to converge by Hardy-Littlewood almost 100 years ago. 
We probe the convergence to these constants, following early experimenters. After factoring out the dihedral symmetry
of Gaussian primes, they are bijectively related to the standard primes but the sequence of angles appears
random. A similar story happens for Eisenstein primes. Gaussian or Eisenstein primes have now a unique angle
attached to them. We also look at the eigenvalue distribution of greatest common divisor matrices whose explicitly known 
determinants are given number theoretically by Jacobi totient functions and where 
unexplained spiral patterns can appear in the spectrum. Related are a class of graphs for which the
vertex degree density is related to the Euler summatory totient function. 
We then apply cellular automata maps on prime configurations. Examples are
Conway's life and moat-detecting cellular automata which we ran on Gaussian primes.  
Related to prime twin conjectures and more general pattern conjectures for Gaussian primes is the question 
whether "life" exists arbitrary far away from the origin, even if is primitive life in form of
a blinker obtained from a prime twin. Most questions about Gaussian primes can be asked for Hurwitz
primes inside the quaternions, for which the zeta function is just shifted. There is a Goldbach statement
for quaternions: we see experimentally that every Lipschitz integer with entries larger than 1 is a sum of two Hurwitz 
primes with positive entries and every Hurwitz prime with entries larger than $3$ is a sum of a Hurwitz and
Lipschitz prime. For Eisenstein primes, we see that all but finitely many Eisenstein integers with coordinates 
larger than 2 can be written as a sum of two Eisenstein primes with positive coordinates. We also predict that every Eisenstein
integer is the sum of two Eisenstein primes without any further assumption. For coordinates larger than 1,
there are two curious ghost examples. For Octonions, we see that there are arbitrary large Gravesian integer 
with entries larger than 1 which are not the sum of 
two Kleinian primes with positive coordinates but we ask whether every Octavian integers
larger than some constant $K$ is a sum of two Octavian primes with positive coordinates. 
Finally we look at some spectra of almost periodic pseudo 
random matrices defined by Diophantine irrational rotations, where fractal spectral phenomena occur. 
The matrix is the real part of a van der Monde matrix whose determinant has relations with the curlicue problem
in complex analysis or the theory of partitions of integers. Diophantine properties allow to estimate 
the growth rate of the determinants of these complex matrices if the rotation number is the golden mean. 
\end{abstract}

\section{Introduction}

In this medley of experiments, we pick up some number theoretical themes.
Our approach is mostly elementary and experimental and sometimes in linear algebra, 
graph theoretical or statistics context. The few mathematical remarks which appear 
in this text all have quick derivations. 
The topics have emerged in the last couple of years while teaching 
linear algebra courses or a course on ``teaching math with a historical perspective" in the
``math for teaching" program at the Harvard extension school.
New phenomena came up especially while writing exam problems or computer algebra projects. 
Many of the featured experiments remain unexplained. We take the 
opportunity to illuminate the material also from a historical perspective 
focusing on mathematicians like Gauss, Euler, Goldbach, Riemann, Hurwitz, Smith,
Hardy and Littlewood and Chebyshev, who were all interested in multiple
fields of Mathematics like number theory and analysis.  \\

The history of mathematics illustrates that mathematical explorations initially are often 
experimental by nature: the Pythagorean theorem was first explored experimentally without proof,
as writings on Clay tablets show \cite{BabylonianTexts}. Descartes discovered both the Euler polyhedral 
formula as well as the Goldbach conjecture by doing experiments, Fermat discovered the two square theorem
experimentally. Conjectures by Fermat, Euler Gauss as well as Hardy 
and Littlewood came from experiments. The computations of Gauss counting primes
led to the prime number theorem. Also his first work on quadratic reciprocity was 
experimental at first, before proofs got available, and Gauss himself got a few proofs. 
Many results in random matrix theory were exploratory at the beginning, 
including work of Wigner, Ginibre or Girko, who found the circular law
mentioned below. Some open problems are constantly probed experimentally, like finding more structure in the
roots of the zeta function by pushing the limits of Goldbach, investigating the statistics of frequencies of 
prime twins or searching for even perfect integers or perfect Euler bricks.
The situation of Hardy and Littlewood is remarkable as it might have been the first time that 
``pure math" research mathematicians started to get assisted by 
collaborators who helped doing the computations, before electronic 
computers became available. Even today, in a time where ``experimental mathematics" has become
its own field \cite{BBG,CrandallPomerance} and groups of research 
mathematicians like ``polymath" collaborate, including experimental mathematicians 
who write and run code. Experiments feed the intuition and also help to fill computer assisted parts. 
The story of unsolved problems in number theory is closely linked to computations and
experiments.  \\

Primes in division algebras like Gaussian or Hurwitz or Octavian primes in Complex, Quaternionic or
Octonion spaces give plenty of opportunity to make experiments. We got dragged into this also by reading
\cite{MazurStein} who proposed an exercise in Section 5, which totally mixed up our other summer plans.
We certainly scratch only the surface while exploring this a bit. We take the opportunity to include seemingly 
remote topics like greatest common divisor matrices or matrices obtained from rotations using Diophantine angles.
We will see for example, again following Hardy and Littlewood, that Goldbach problems in division algebras
can be ported guided by calculus problems related to the circle method. There is plenty of opportunity for new questions. 
An example of a question which seems never have been asked is the existence of ``Gaussian prime life" 
arbitrary far from the origin or whether a Prime twin theorem holds 
for Hurwitz or Octavian integers. Evidence that Goldbach holds for Hurwitz primes comes from 
relating the toughest boundary cases with Landau type problem which have quantitative 
generalizations given by Hardy-Liouville statistical laws which we test ourself up to $2^{37}$. 
In the Gaussian prime case already, we have full assurance that Gaussian Goldbach is hard, as it would 
imply an open Landau problem. Connections with other mathematical fields like topology comes
in as prime numbers define classes of graphs. One can look for example at all the positive integers
as the vertex set of a graph and connect two numbers $a,b$ if $a+ib$ is a Gaussian prime. 
Similarly, one can take the set of Gaussian integers 
$a+ib,c+id$ in the first quadrant of the complex plane and connect two, if either $a+ib+cj+dk$ 
is a Lipschitz prime or $(a+ib+cj+dk)/2$ is a Hurwitz prime. \\

Primes in division algebras form an even grander arena for explorations than the 
traditional primes. It could even be an El Dorado for early explorations or in education. 
For Gaussian primes, there were attempts \cite{Eric} 
in the 1960'ies to implement them into secondary education curricula.
The question why quaternions have disappeared from mainstream calculus is interesting.
\cite{ConwaySmith} give as a reason the success of notations put forward by Gibbs \cite{GibbsWilson} 
with precursor notes distributed as early as 1881. 
Coauthor Wilson in that Gibbs textbook acknowledges that quaternions were useful in helping the text.
While quaternions are maybe a bit in the background in calculus frameworks, they are unbeaten in 
elegance in number theory: Gaussian primes produce the most natural frame work for proving some Diophantine problems
and Hurwitz proof of the Lagrange four square theorem using quaternions can not be beaten in simplicity. 
This theorem implies that unlike for Gaussian primes, where half of the rational 
primes remain prime, all rational primes decay into Quaternionic primes. \\

The richness of the topic can be illustrated by the presence of open questions for Gaussian,
Hurwitz and Octonion integers. Especially in the context of primes. One can ask
questions about the {\bf existence} of geometric patterns in the set of Hurwitz integers, 
like the prime twin problem for Hurwitz integers, additive number theoretical questions of 
Goldbach type in division algebras as done here, percolation problems of topological type like
moat type problems in Hurwitz integers,
as well as growth problems which deal with asymptotic behavior.
We start also to look at correlation, determinant and trace type questions for Gaussian primes. \\

An other twist comes in when looking at primes as the input for dynamical systems.
We can apply cellular automata to prime constellations for example. 
Cellular automata are dynamical systems, which emerged from topological dynamics
are in some sense partial differential equations, where not only time
and space is discrete, but where also the target space of the functions, here the alphabet $A$ is 
discrete. They are useful in computer graphics, in algorithms for edge detection, for seeing 
morphological features, for smoothing or sharpening operations. Almost all image processing filters are 
cellular automata acting on the color frames of the lattice of pixels.
On can apply such maps to configurations like Gaussian primes. 
Dynamical systems of a different kind enter when constructing matrices whose
entries are $A_{nm} = \cos(k m \alpha + m \beta)$ where $\alpha,\beta$ are Diophantine. These are real parts of 
van der Monde matrices $B_{nm}$. The spectra are unexplained. We are able to give a bound on the growth of
the determinant of the $B$ matrices if $\alpha$ is the golden mean. \\
 
{\bf Acknowledgements}: Some computations in this report were run on the Odyssey 
cluster supported by the FAS Division of Science, Research Computing Group 
at Harvard University. 

\section{Gaussian primes}   

{\bf Gaussian primes} are the irreducible elements in the ring of 
{\bf Gaussian integers}. The {\bf arithmetic norm} $N(z)$ of a Gaussian integer $z=a+ib$ is defined as $a^2+b^2$.
If $|z|$ denotes the {\bf absolute value} of the complex number $z$ then $N(z)=|z|^2$ is the norm. 
Integers of norm $1$ are called {\bf units}. The four unit elements $\{1,i,-1,-i\}$ form a cyclic multiplicative 
subgroup of the ring $Z[i]$. A number in $Z[i]$ is called a {\bf Gaussian prime} if it 
is not the product of two other numbers with smaller norm and norm larger than $1$.
Examples of Gaussian primes are $2+i$ or $3$. Examples of non-primes are $1+3i$ or $5$ as $(1+3i)=(1+i)(1+2i)$
and $5=(2+i)(2-i)$. Also $0$ and the units are not primes by definition. 
Since a unique factorization theorem holds in the ring $Z[i]$ which is 
part of a {\bf division algebra} satisfying $|zw|=|z| |w|$, one can produce a list of Gaussian
primes easily using the {\bf baby prime test}: a Gaussian integer $z$ is prime, if and only if for all 
$w$ satisfying $1<|w| \leq \sqrt{|z|}$, the Gaussian integer $z$ is not a multiple of $w$. Fortunately,
there is a much faster way to give a list of all Gaussian primes: one can just
take the list of traditional primes $\{2,3,5,7,11,\dots \}$ in $\mathbb{N}$
called {\bf rational primes} and associate to each a
{\bf Gaussian quadruple} $z,iz,-z,-iz$ or an {\bf octuple} 
$z,iz,-z,iz,\overline{z},i\overline{z},-\overline{z},i\overline{z}$
of Gaussian primes. The Gaussian primes can be seen as a ``cover" of the rational primes: take a 
Gaussian prime $z$. Draw 
\begin{wrapfigure}{l}{4.1cm} \begin{center}
\includegraphics[width=4cm]{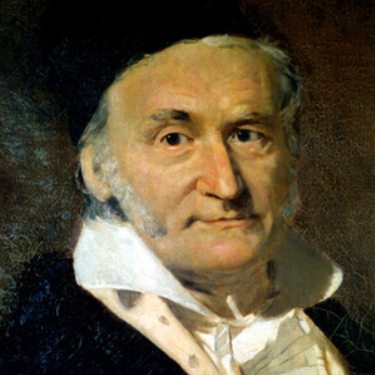} 
\end{center} \end{wrapfigure}
the circle centered at $0$ passing through a prime $p$. The radius $|p|$ is either 
$\sqrt{2}$ or a $(4k+3)$-prime or then the square root of a $(4k+1)$-prime. 
In the first two cases, there are $4$ other Gaussian primes $u z$
on the circle where $u$ is one of the four {\bf units} $\{1,i,-1,-i\} = \{ |z|=1 \} \cap Z[i]$. 
In the last case, if $z=a+ib$ with $a \neq b$, there are $8$ primes $u z$ and $u z'$, where $u$
is a unit and $z'=b+ia = i \overline{z}$. {\bf Unique factorization} holds in $Z[i]$ modulo units. 
The quotient $X=\mathbb{C}/A$ of the complex plane by the dihedral group $A$ generated by multiplication by units and the 
conjugation $z \to \overline{z}$ is a cone on which the usual primes live. Due to the fact that for the primes on the 
axes have {\bf prime length} and the primes off the axes have {\bf prime norm} $N(z)=|z|^2$, the order is shuffled:
the $(4k+1)$-primes $p$ are located at height $\sqrt{p}$ while the $(4k+3)$-primes $p$ are located at height $p$. 
This will enrich slightly the Gauss zeta function by multiplying it with a beta function as Dirichlet already knew.
The cone $X$ is an example of what topologists call an {\bf orbifold}. This space $\mathbb{C}/A$ is the 
quotient of the manifold $\mathbb{C}$ by the equivalence relation defined by a finite group action. 
A geometer would say that the complex plane as a {\bf ramified cover} over the space $X$. 
The rational primes $2,3,11 \dots$ are special in that their `
ramification index" is different and $2$ is very special as they are also invariant under conjugation. 
Topological jargon has entered number theory not accidentally: some number theorists were
also geometers like Hurwitz, Gauss and Euler. There are plenty of accessible texts explaining the structure 
of Gaussian primes well \cite{StillwellNumber,Conrad,ConwaySmith,Schlackow}.
Gaussian integers are not without applications: for example in the study of eigenfunctions of the 
Laplacian on the torus \cite{Bourgain2014} or to study discrete velocity models 
for the Boltzmann equation \cite{Fainsilber}. Engineers
use it already for watermarks or codes \cite{Bouyuklieva},
RSA \cite{KovalVerkhovsky} or zero-knowledge identifications \cite{Valluri}. 
Any question about primes in the integers $\mathbb{Z}$ can also be asked in the integral domain
$\mathbb{Z}[i]$ of Gaussian integers \cite{Gauss1831}.  \\

{\bf Remarks.} \\
{\bf 1)} According to Dickson \cite{dicksonI}, 
the Fermat two square theorem was noticed empirically first by Albert Girard, 
who was also the first to introduce the now common notations $\sin,\cos,\tan$ for trigonometric
functions as well as the recurrence for the Fibonacci sequence \cite{dicksonI} page 393. 
Since Fermat announced a proof of his theorem to Mersenne on December 25, that two-square theorem 
is also called the {\bf Christmas theorem}. \\
{\bf 2)} There are various definitions for orbifold. We use the version telling that
the quotient $M/A$ of a finite group $A$ acting on a smooth manifold $M$ by smooth diffeomorphisms,
such that for every $a \in A$, the set of fixed points is a smooth submanifold of $M$. In the Gaussian
integer case, the group $A$ is the dihedral group $D_4$ on the complex plane generated by rotations 
$z \to i z$ and reflections $z \to \overline{z}$. All elements clearly have fixed point
sets which qualify. They are either the entire manifold, or a line or then a single point.  \\
{\bf 3)} The covering story also applies to the ring $Z[i]$ of Gaussian integers. After dividing out
the dihedral group $A$ generated by the multiplication by units and
the conjugation $z \to \overline{z}$, we obtain a structure in which the Gaussian primes
correspond bijectively to the usual primes. But the larger dihedral symmetry also includes
{\bf complex conjugation}. The nature of an odd rational prime $p$ can
be read off from the {\bf Jacobi symbol} $(-1)^{(p-1)/2} = \left( \frac{-1}{p} \right) = (-1|p)$
which is in the case of primes $p$ also called the {\bf Legendre symbol}.
Each prime of the form $4k+1$ is by the {\bf Fermat two square theorem}
of the form $p=a^2+b^2$ and so the product of two Gaussian primes
$(a+ib),(a-ib)$, the Gaussian primes on the axes are the $4k+3$
rational primes. Away from the axes, we have also the Gaussian primes $1+i$ and its conjugates. \\ 
{\bf 4)} A mathematical message of this section is that one can put a natural order
structure the Gaussian primes after a suitable identification. It motivates questions like the following:
define the sequence of numbers $\theta(n)={\rm arg}(p_n) - \pi/8 \in (-\pi/8,\pi/8)$,
where $p_n$ is the $n$'th prime of the form $4k+1$. As usual in topological dynamics, the half sequence 
$\{ \theta(n) \}_{n \in \mathbb N}$ 
defines a compact topological space $X$, the hull in the compact product space $[-\pi/8,\pi/8]^\mathbb{Z}$ given 
by all the accumulation points. The shift $T$ is now  homeomorphism of this space and notions like topological
entropy are defined. One can for example ask about the structure of the {\bf invariant measures} of $T$ etc.  
A natural conjecture is that the limiting system for these {\bf prime angles} is a Bernoulli system,
meaning that on the hull, the random variable $X_k(\omega) = \theta_k(\omega)$ are all independent identically distributed, with
uniform distribution on $[-\pi/8,\pi/8]$. Especially, {\bf central limit theorems} and {\bf laws of iterated logarithms}
would hold. There are similar easy looking problems, like the system given by the integer digits of $\pi$ or $\sqrt{2}$
to some base, where a Bernoulli statement is believed to be true but not accessible to a proof yet. \\
{\bf 5)} One way to see that a number is composite if there are two fundamentally different ways
to write it as a sum of two squares $n=a^2+b^2=c^2+d^2$ is due to Euler who pointed out 
$n=( (a-c)^2 + (b-d)^2 )*( (a+c)^2 + (b-d)^2)/(4 (b-d)^2 )$ which leads to a factorization. 
The relation with factorization could indicates that finding the angle 
$\theta = {\rm arg}(a+ib)$ of a prime $p=a^2+b^2$ could be hard if only $p=4k+1$ is known. We are not
aware of a faster determination of the angle $\theta(p)$ than searching in $O(\sqrt{p})$ time. 
One can ask for example whether it is possible to find the angle $\theta$ in $O(p^{1/4})$ computation
steps from $p$ or even in $O(\log(p))$ time.  \\

We can now look at the random walk when using the angles $\theta(n)$:
 
\question{
How fast does $S(n) = \sum_{k=1}^n \theta(p(n))$ grow, where $p(n)$ is the $n$'th prime
of the form $4k+1$.  }

Our experiments indicate that this random walk does not behave differently than a usual
random variable with uniform distribution in $[-\pi/8,\pi/8]$. It is known that the
distribution of $\theta$ is the uniform distribution. (References are given in \cite{Kubilyus}).
We made experiments with correlations of neighboring $\theta$:
the sequences $Y(n)=X(2n)$ and $Z(n)=X(2n+1)$ for example appear to be decorrelated, indicating
that it is difficult to predict the angle of prime $k$ if the angle of prime $k-1$ is known. \\
Any way, we are not aware of a {\bf fast way} to compute the unique $(a,b)$ with $0,b<a<\pi/4$
satisfying $a^2+b^2=p$ if $p$ is a $4k+1$ prime. \\

\resultremark{ On the quotient space $Z[i]/D_4$, the Gaussian primes 
admit a natural total order and one-to-one bijection to the rational primes. }

Of course, since both sets are rational, there existed also other bijections. The point
is that we can relate them naturally. While this does not appear to be useful, it is a source for
new questions like about the distribution of the angle. 

\begin{figure}[!htpb]
\scalebox{0.80}{\includegraphics{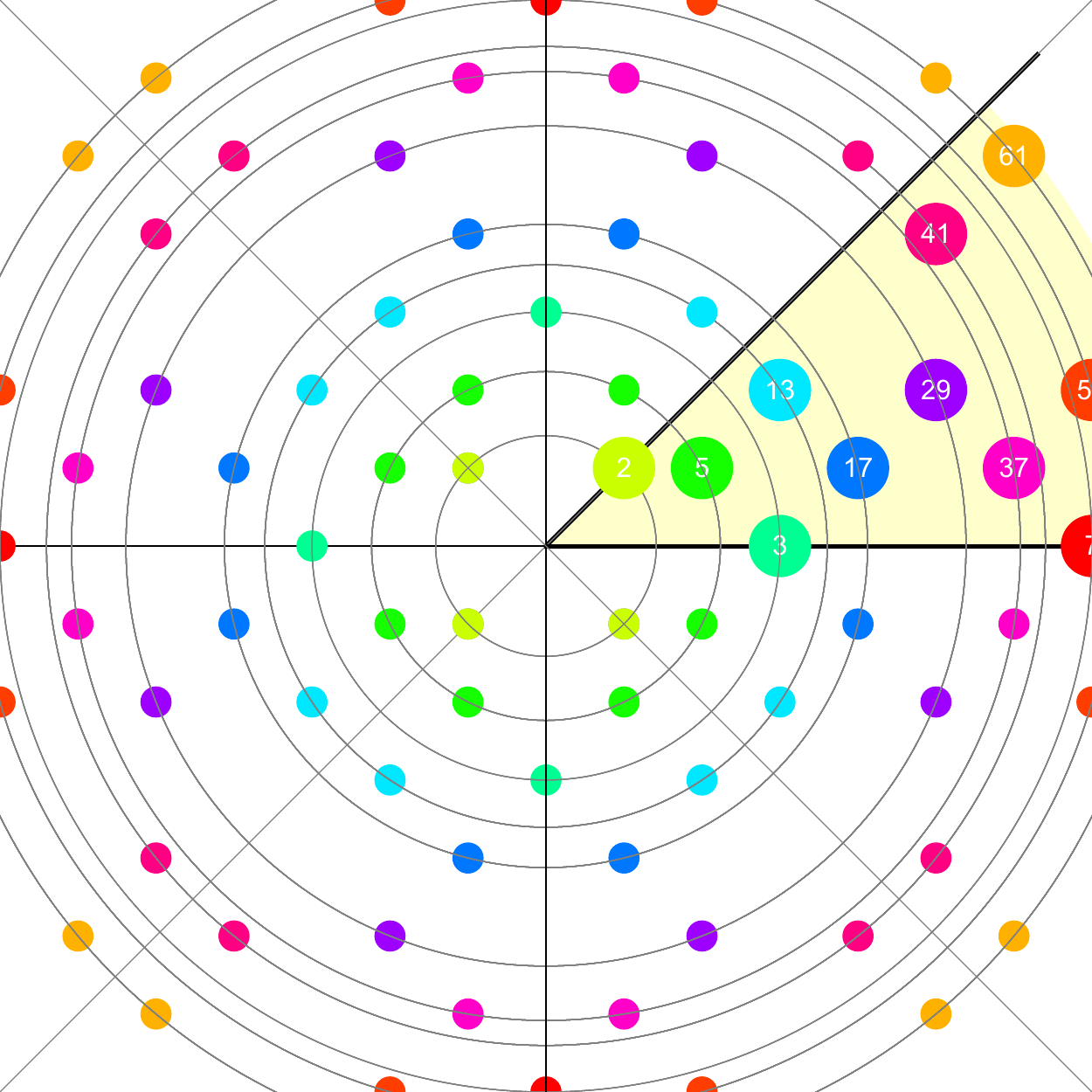}}
\caption{
Gaussian primes are either located on the real or imaginary axes, given by primes
of the form $4k+3$, or then lattice points on circles of radius $r$,
where $p=r^2$ is prime. In other words, the circles on which Gaussian primes are located
are in one to one correspondence to the standard rational primes. Only the order is
shuffled. The second smallest circle corresponds to the primes $\pm 2 \pm i$ and
so to the prime $5$. The third then to the prime $3$.
Each of the primes on the symmetry lines come with multiplicity
$4$, the rest with multiplicity $8$. The picture shows the first octant sector, which
is the fundamental domain for the action of the dihedral group $A=D_4$ on the complex plane $\mathbb{C}$. 
}
\label{circles}
\end{figure}

\begin{figure}[!htpb]
\scalebox{1.0}{\includegraphics{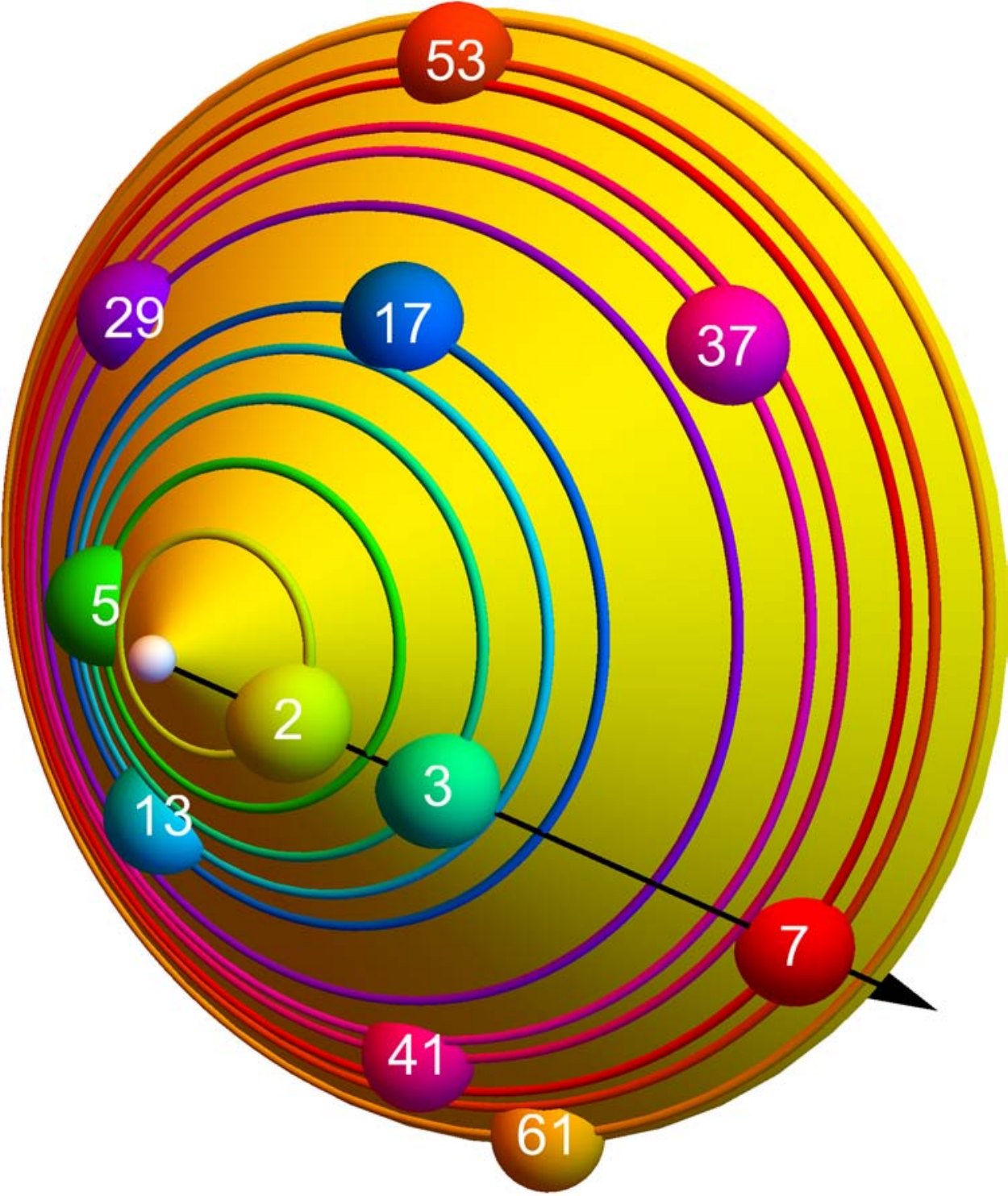}}
\caption{ 
When drawn on the cone orbifold $\mathbb{C}/A$, where $A$ is the dihedral symmetry group,
every traditional rational prime can be matched in a one-to-one manner 
with an equivalence class of Gaussian primes. The rational primes of the form $4k+3$ are 
on the real axes where the cone is glued. The angles $\theta(n)$ at which 
prime $p(n)$ is located, appear to pretty random. All correlation tests done so far 
indicate this. }
\label{circles}
\end{figure}

\begin{figure}[!htpb]
\scalebox{0.8}{\includegraphics{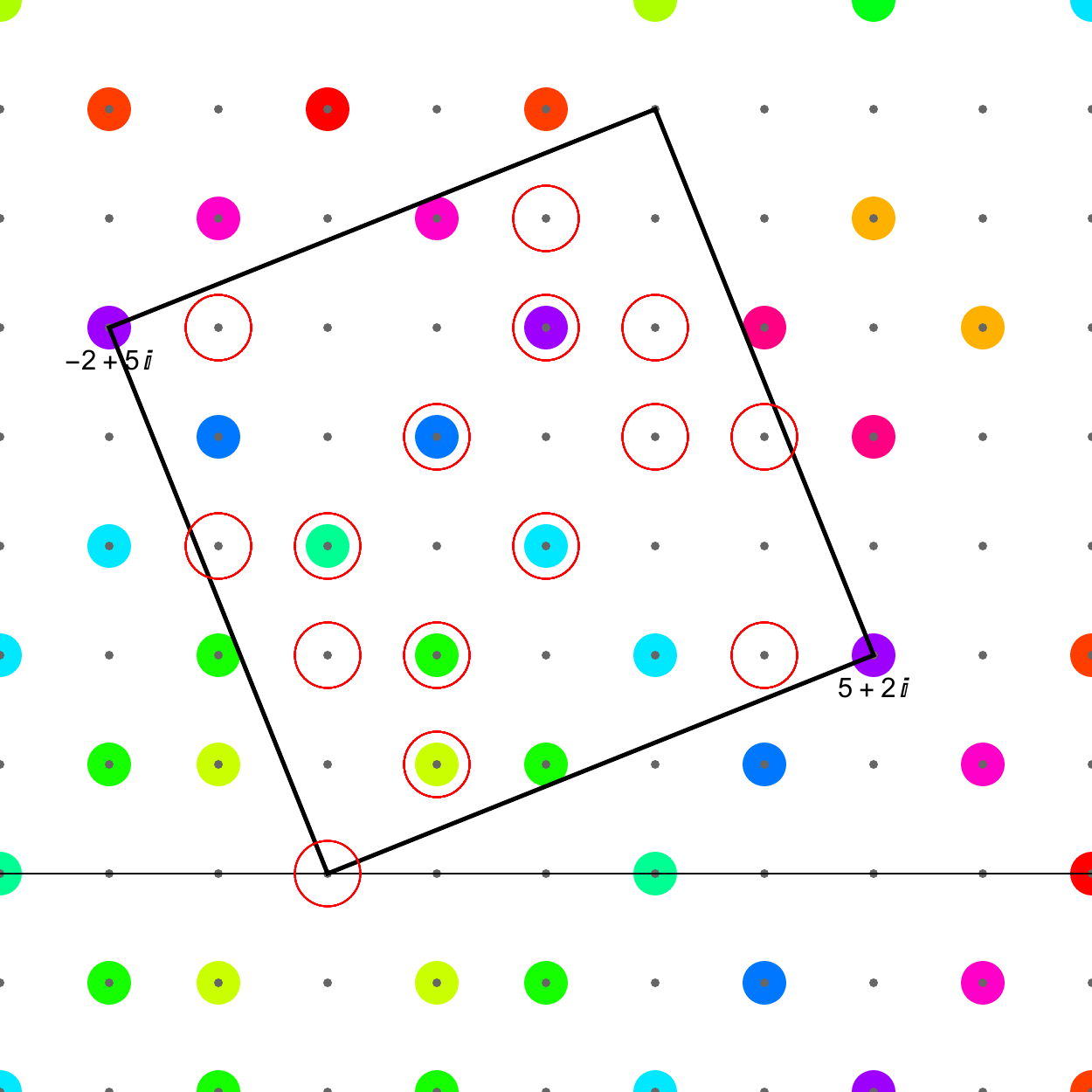}}
\caption{
We see the region $\mathbb{C}$ modulo a prime $p=5+2i$ and unit multiples.
A number $z$ mod $p$ is in the square. There are $N(p)$ equivalence 
classes. For odd $N(p)$, half of them are quadratic residues. 
The quadratic reciprocity result which Gauss proved already knew for
Gaussian primes tells which primes are quadratic residues. 
}
\label{circles}
\end{figure}

\begin{figure}[!htpb]
\scalebox{0.50}{\includegraphics{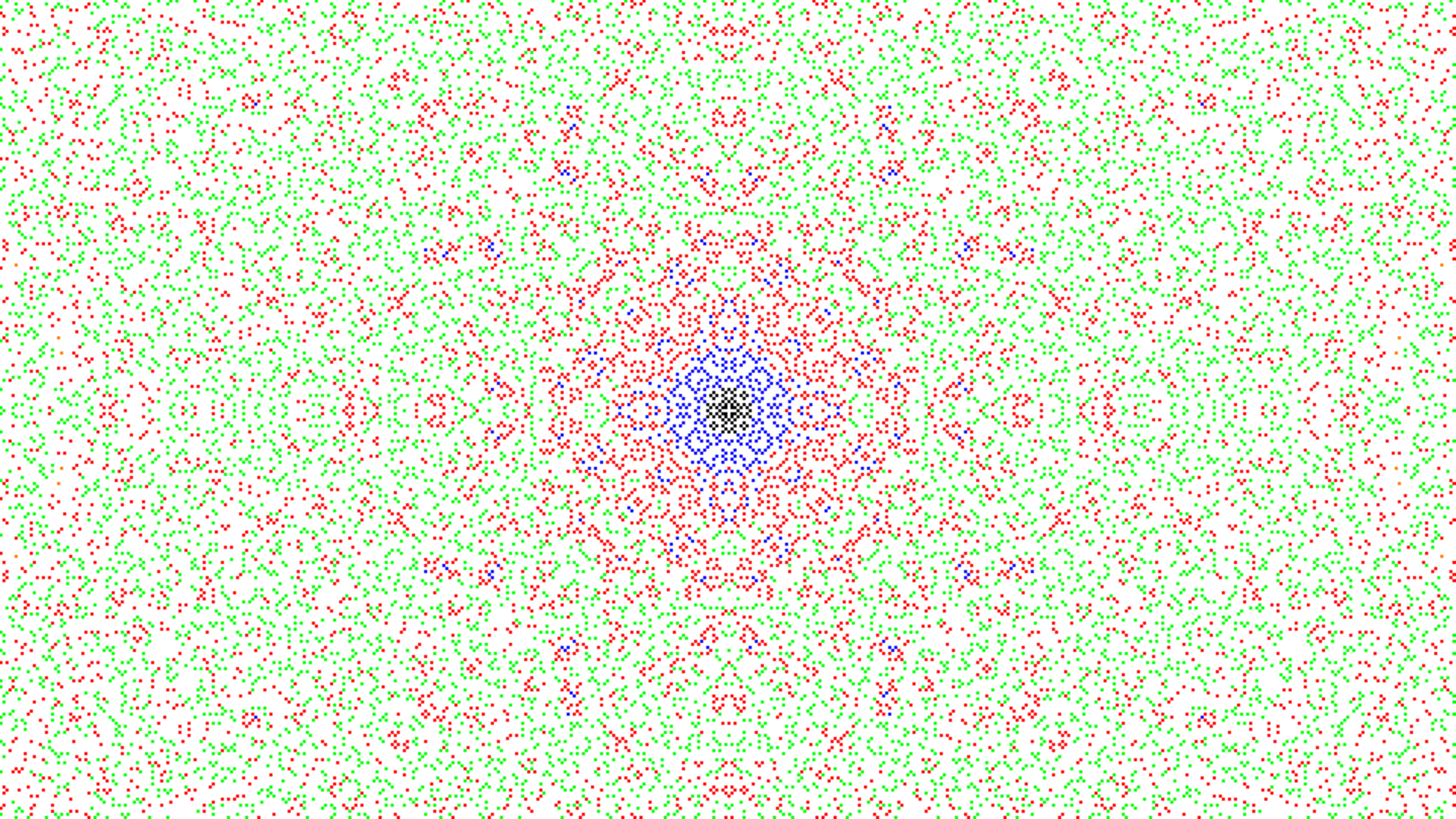}}
\caption{
A larger part of the Gaussian primes. 
}
\label{circles}
\end{figure}

\begin{figure}[!htpb]
\scalebox{1.0}{\includegraphics{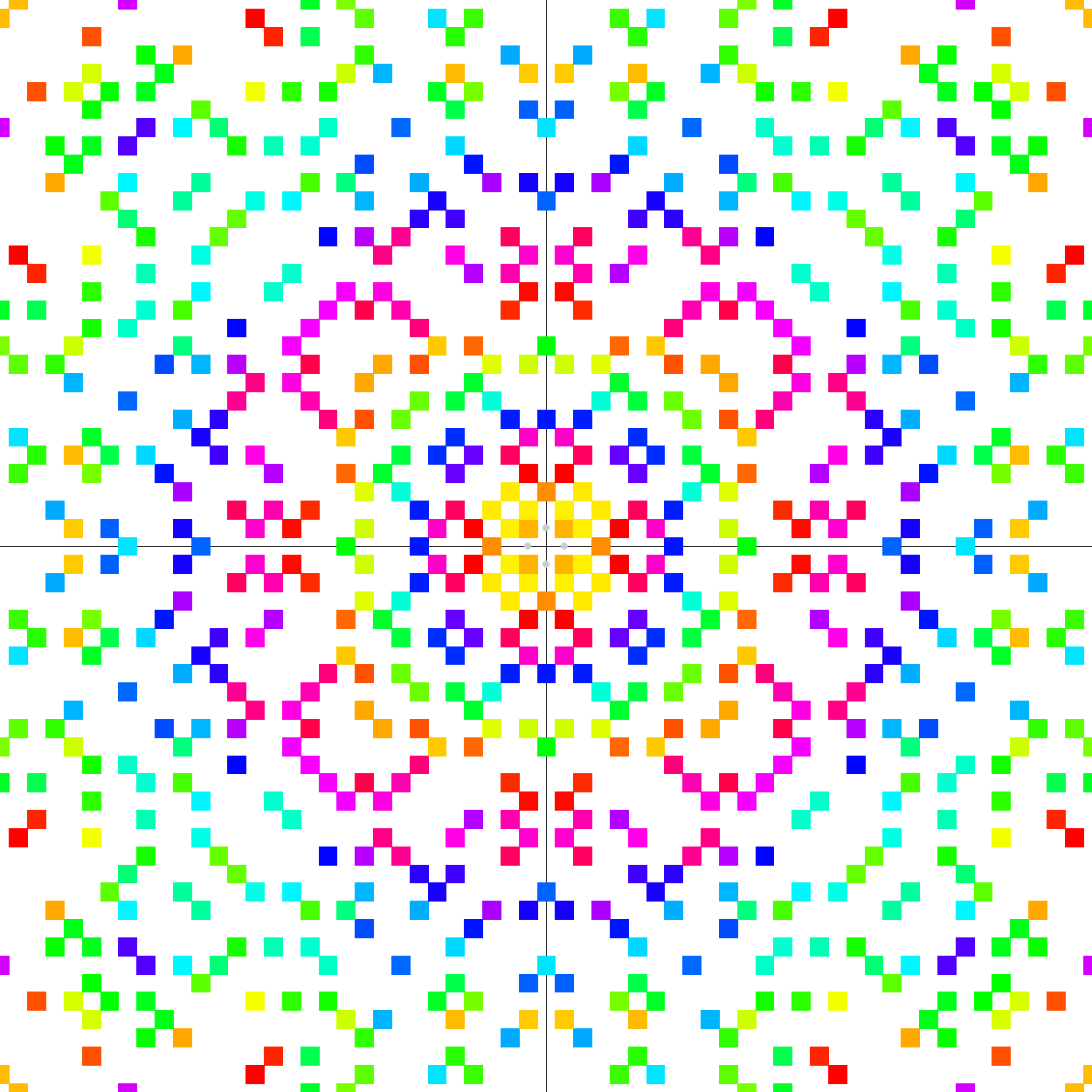}}
\caption{
Gaussian primes have more structure near the origin which 
is amplified by the dihedral symmetry. 
According to a ``Fortune" article of June 1958 
(which we unfortunately can not get hold of),
the Gaussian prime pattern has been used for a tea cloth. 
}
\end{figure}

\begin{figure}[!htpb]
\scalebox{0.2}{\includegraphics{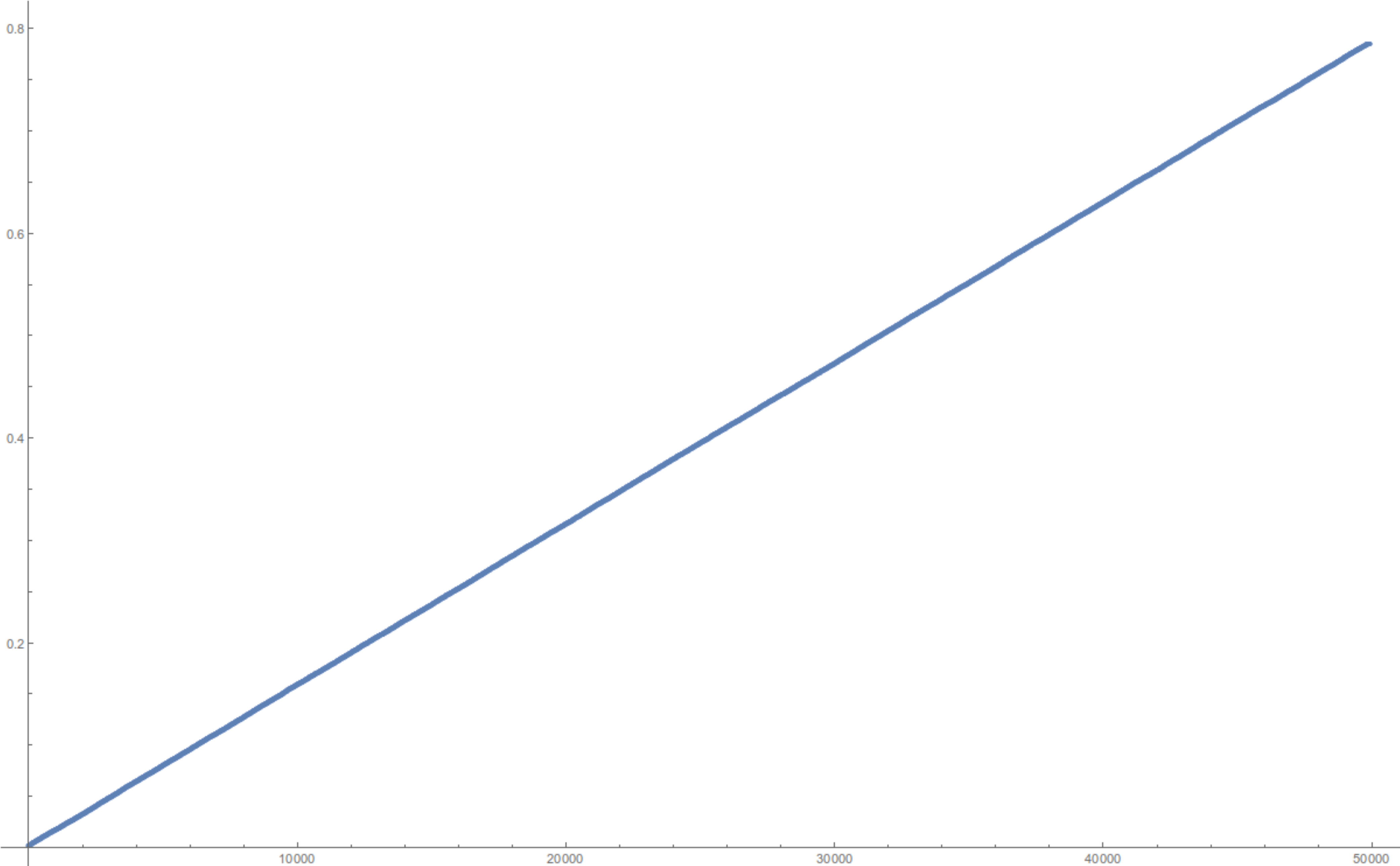}}
\scalebox{0.2}{\includegraphics{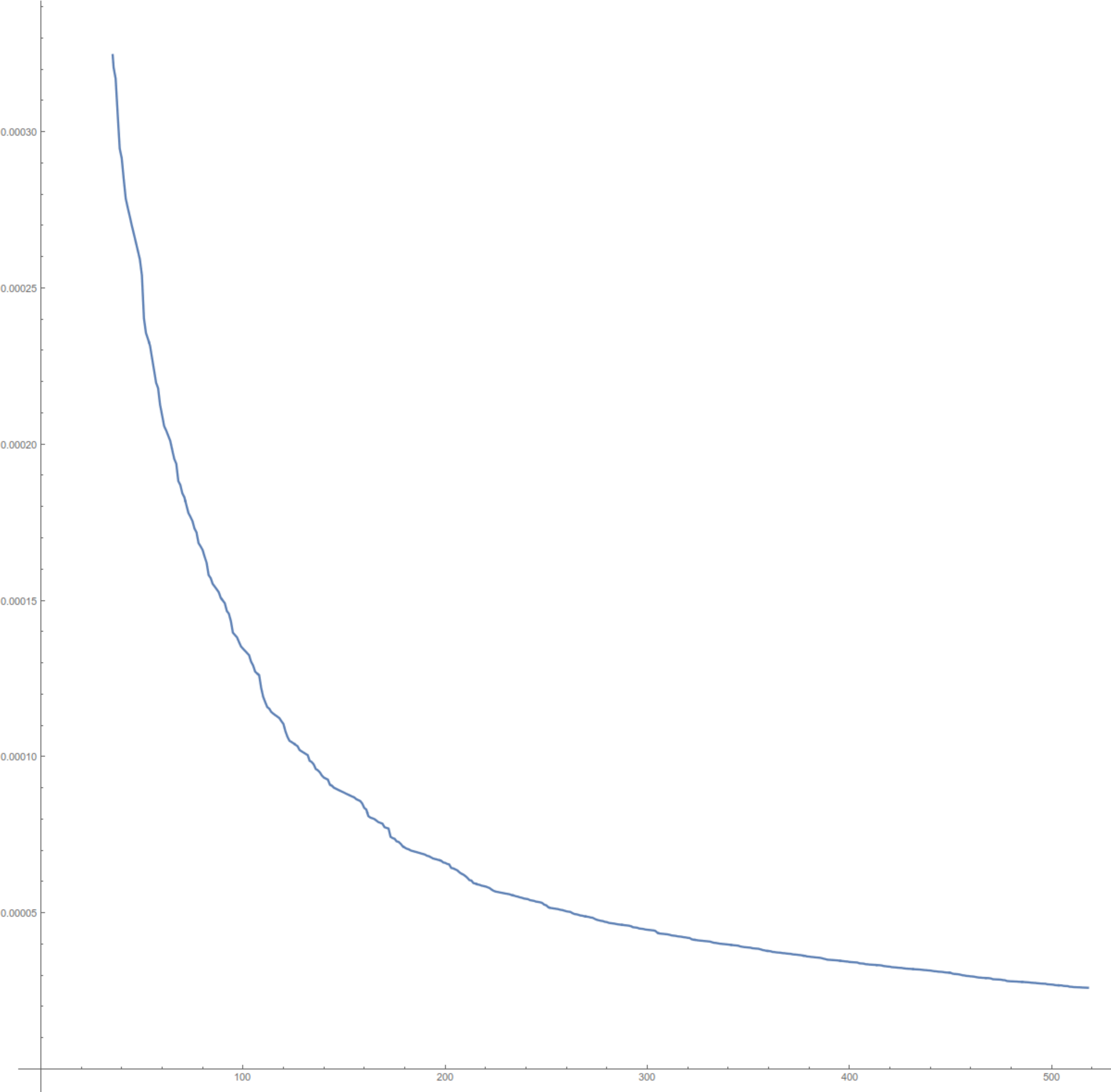}}
\caption{
The argument ${\rm arg}(z)$ over primes becomes uniform: no angular
direction is favored. The function $p \to {\rm arg}(P)$, where $P$
is the unique Gaussian prime with argument $\theta(p) \in [0,\pi/4)$ belonging to $p$ 
is a random variable which is pretty random for the primes $p$ of the form $4k+1$. 
The second figure shows the correlation between the vectors 
$X(n) = (\theta(P(1), \dots \theta(P(2k+1)), \dots \theta(P(2n+1))$ 
and $Y(n) = (\theta(P(2)), \dots, \theta(P(2k)), \dots \theta(P(2n)))$ as a function 
of $n$.
}
\end{figure}

{\bf Remarks:} \\
{\bf 1)} Some results have been pushed over from rational primes to Gaussian primes.
There is for example an analogue of Dirichlet's theorem on arithmetic progressions:
for an arbitrary finite set in $\mathbb{Z}[i]$, there exist infinitely many
$a \in \mathbb{Z}[i]$ and $r \in \mathbb{R}^+$ such that $a+r \sum_{f \in F} v_f$
is a Gaussian prime \cite{Tao2006}. An other example is that the 
density of the {\bf prime quotients} $p/q$ in $R^+$ generalizes to the proven statement
that the Gaussian prime quotients $p/q$ are dense in the complex plane \cite{Garcia}. 
Patterns are explored in \cite{JordanRabung76,GaussianZoo}. \\
{\bf 2)} The connection between Gaussian primes and rational primes has been 
used by Gauss as a tool for {\bf quadratic reciprocity}. 
Also Dirichlet realized the connection. The reason is that 
the multiplicative part of $Z[i]/p$ which has size $N(p)-1$ 
can be generated as $a^k \; {\rm mod} \;  p$ as in the rational case. Now define 
the birational Jacobi symbol $(a|p)$ as $a^{(N(p)-1)/2}$ mod $p$. If this
is $-1$, then there is a solution to $x^2=a$ {\rm mod} $p$ by Fermat's little
theorem. \\
{\bf 3)} While the bijective map $p \to a+ib$ with $a>0, b \geq 0, a<b$ from rational primes
to equivalence classes of Gaussian primes $\mathbb{C}/A$ has an inverse $a+ib \to N(a+ib)$
which can be computed with no effort, the map itself needs some searching. Since most
cases are where $a$ is very small, this is usually found fast, but its not elegant. \\
{\bf 4)} While the identification of Gaussian primes with rational primes after taking
equivalence classes is nice, there is still much more structure in the Gaussian primes.
We will see that when formulating a Goldbach conjecture, when applying cellular automata
map on them or by looking at the Gaussian Riemann zeta function, which is an example of 
a {\bf Dedekind zeta function}. It leads to interesting structures like the {\bf complex M\"obius
function} which is defined as in the rational case but for which the corresponding Mertens
function fluctuates more. \\
{\bf 5)} There is a ``one line" topological proof of the two square theorem \cite{Zagier90}. But as
it relies on a fixed point theorem, it is not constructive. Lets briefly sketch it: it uses the compact surface
with boundary $x^2+4yz=13, x \geq 0, y \geq 0, z \geq 0$ on which there are two 
involutions $T(x,y,z) = (x,z,y)$ and $S(x,y,z)$ which is $(x+2z,z,y-x-z)$ if $x<y-z$
and $2y-x,y,x-y+z)$ if $y-x<x<2y$ and $(x-2y,x-y+z,y)$ if $2y<x$. The map $T$ has exactly one fixed point
$(1,1,k)$ if $p=4k+1$. Since the fixed points $T$ and the number of fixed points of $S$ are the same
modulo $2$, also $S$ must have a fixed point of the form $(x,y,y)$. This is topological as the argument
about fixed points is essentially a Riemann-Hurwitz argument relating the difference 
$\chi(X)-2\chi(X/T)$ as a sum of ramification indices of fixed points.  \\
{\bf 6)} Also the Minkowski proof using geometry of numbers or the Jacobi proof using theta functions
$\theta(x) = \sum_n x^{n^2}$ is not constructive. To the later:  as 
$\theta(x)^2$ is the generating function of the number $2 a(n)$ counting the way to write a 
number $n$ as a sum of two squares and this is equal to $1+4 \sum_{n \geq 0} x^{4n+1}/(1-x^{4n+1}) - x^{4n+3}/(1-x^{4n+3})$,
this is $a(n) = 4(d_1(n)-d_3(n))$, as $2d_k(n)$ counting twice the number of divisors larger than 1 of $n$ 
which are congruent to $k$ modulo $4$ has the generating function $\sum_n x^{4n+k}/(1-x^{4n+k})$.  
For a number $n=p$ which is a $4k+1$ prime, this means $a(n)=4$. Also this is not constructive. 
Also the proof using the solution $x^2=-1$ modulo $p$ which is possible
as $-1$ is a quadratic residue using finding``square roots" of $-1$ to reduce the problem. 

\section{Lipschitz and Hurwitz primes}

Next to the field of {\bf complex numbers} comes the skew field of {\bf quaternions}. 
Quaternions are numbers of the form $a+bi+cj+dk=(a,b,c,d)$, where
$i,j,k$ are symbols satisfying $i^2=j^2=k^2=ijk=-1$. Discovered by Hamilton in 1843 
they introduced simultaneously the {\bf inner product} and {\bf cross product} to vector
calculus because $(0,v_1,v_2,v_3) \cdot (0,w_1,w_2,w_3) = (-v \cdot w, v \times w)$, even so
the later structures were only promoted by Gibbs. Taking the product of three such quaternions
produces the {\bf triple product} $(0,v_1,v_2,v_3) \cdot (0,v_1,v_2,v_3) \cdot (0,w_1,w_2,w_3)
=(u \cdot (v \times w), u \times (v \times w))$. We see that both the {\bf triple scalar product} as
well as the {\bf triple cross product} are present. The textbook \cite{GibbsWilson} is the key
to \begin{wrapfigure}{l}{4.1cm} 
\includegraphics[width=4cm]{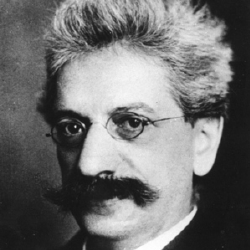} 
\end{wrapfigure}
see how experts in quaternions used their knowledge to introduce simpler and now popular notions and 
by doing so erase the quaternions from calculus textbooks. 

\footnote{
The quaternions remind of the Middlesex canal which connected the Merrimack river with Boston:
it was used to build the railroads, which in turn made the canal obsolete.
Similarly, the success and applicability of the dot and cross product, which were historically
built through quaternions led to their disappearance in vector calculus \cite{ConwaySmith}
}

Quaternions are useful in {\bf computer vision} to realize a rotation in space around a 
unit vector $u$ by an angle $\theta$ by building the two quaternions $s=v_1 i + v_2 j + v_3 k$ and 
the $r=\cos(\theta/2) + \sin(\theta/2)(u_1 i + u_2 j + u_3 k)$, then get the rotated vector $Rv$
from the spacial components of the quaternion $r \cdot s \cdot r^*$. 

Number theory in the ring of quaternions is a bit more strange, as quaternion multiplication 
is {\bf no more commutative}. Furthermore, the ring of integer quaternions is no 
more a {\bf unique factorization domain}. Even the Euclidean algorithm comes short. 
Hurwitz \cite{Hurwitz1919} found a way out of it and 
realized that one can get a {\bf Euclidean domain} when including half units. 
Quaternions with "integer" coordinates are called {\bf Lipschitz integers} 
$a+i b + j c + k d$, where $(a,b,c,d) \in \mathbb{Z}^4$. The others 
which must be of the form $(a,b,c,d) \in \mathbb{Z}^4+ 1/2(1,1,1,1)$ to satisfy 
$N(z) \in \mathbb{N}$ are called the {\bf Hurwitz integers}.
In this text, we want to use these two names to distinguish between the two distinct 
classes of integers and call their union simply {\bf quaternion integers}. 
The quaternion integers form a non-commutative ring of the quaternion division ring. 

\begin{figure}[!htpb]
\scalebox{0.50}{\includegraphics{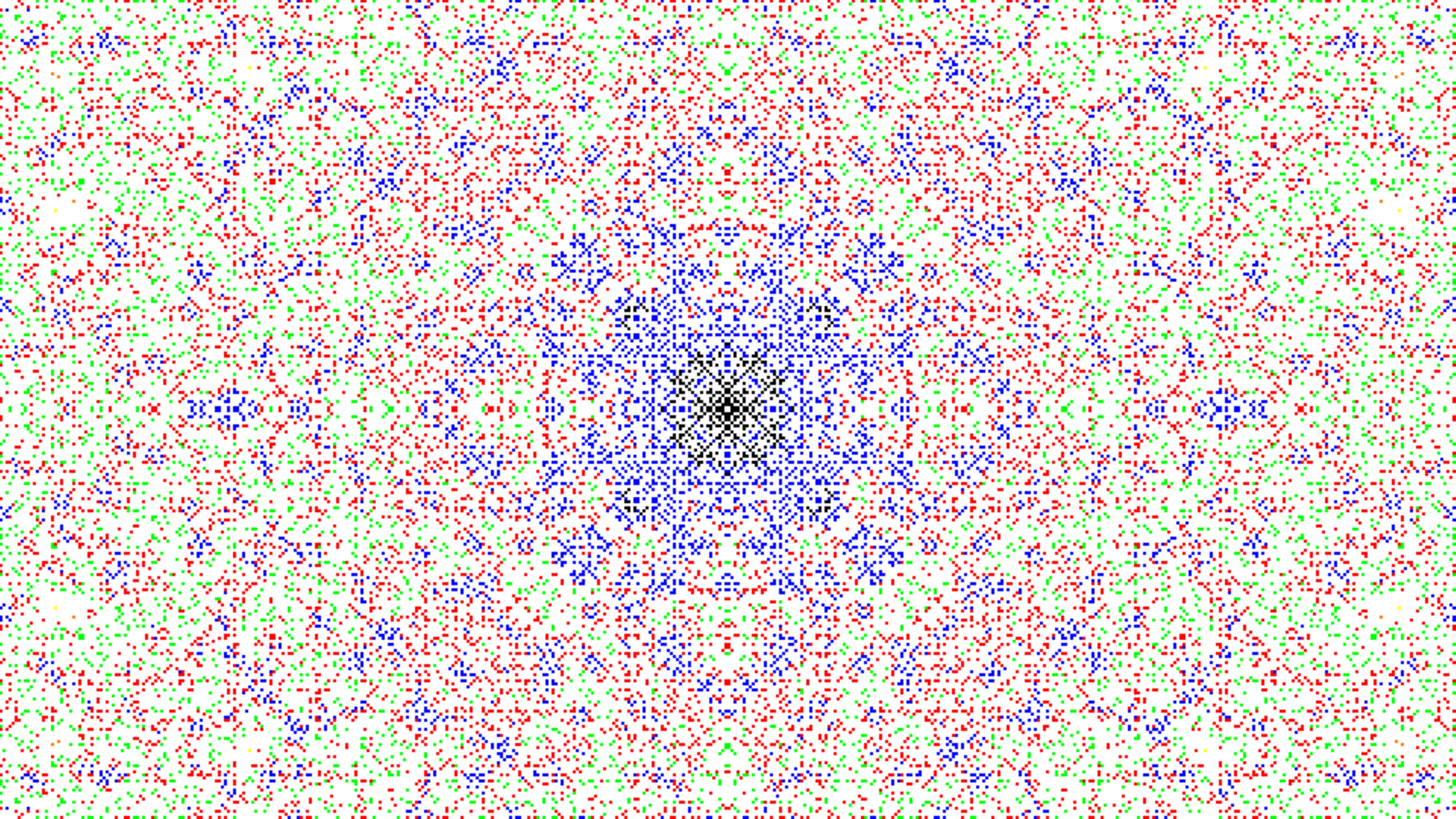}}
\caption{
A slice through the Hurwitz primes: it is the two-dimensional 
plane $(1,1,2x+1,2y+1)/4$. The pixel color $(x,y)$ depends on 
the number of neighbors in that plane. 
}
\label{circles}
\end{figure}

\begin{figure}[!htpb]
\scalebox{0.8}{\includegraphics{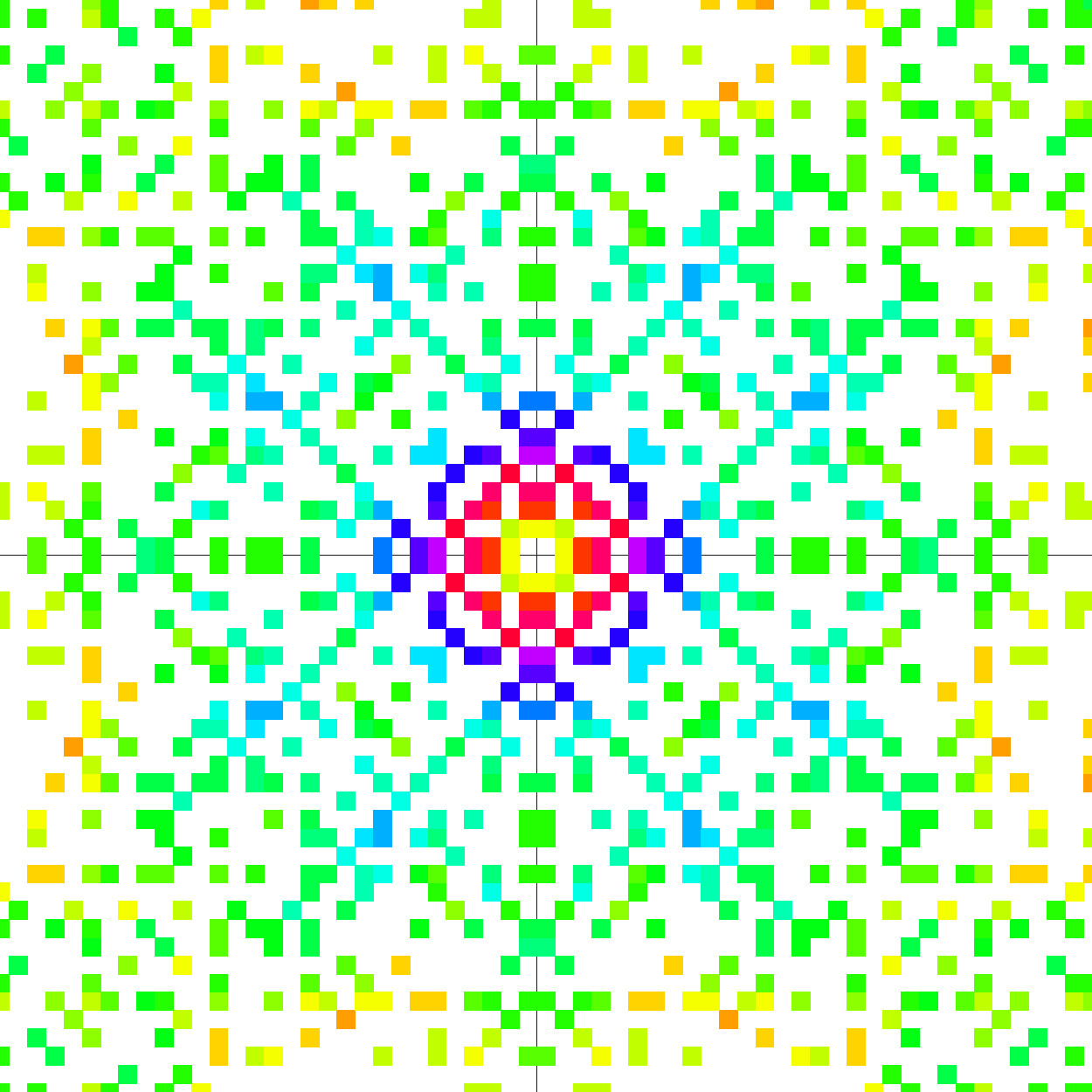}}
\label{circles}
\end{figure}

There are $24$ {\bf units} in the set of quaternion integers.
Eight of them are Lipschitz integers given by permutations of
$(\pm 1,0,0,0)$ and $16$ are Hurwitz integers of the form $(\pm 1,\pm 1, \pm 1, \pm 1)/2$. 
The units build a subgroup of the quaternion group. It goes under the name 
{\bf binary tetrahedral group} $U=SL(2,3)$. As quaternions can be written as unitary $2 \times 2$
matrices, $U$ is a discrete finite subgroup of the unitary group $SU(2)$. 
It can be identified as the semi-direct product of the {\bf quaternion sub group} 
$Q$ built by the $8$ Lipschitz units and the cyclic group $Z_3$, generated by conjugation 
$i \to j \to k \to i$. The group $U$ is a finitely presented group with relations 
$(ab)^2=a^3=b^3=1$ and generators $a=(1+i+j+k)/2$ and $b=(1+i+j-k)/2$. All cyclic subgroups have order
$2,3$ or $6$. The element $-a$ for example generates a cyclic group of order $6$.  \\

There is an other group $V$ with $24$ elements acting on the ring of Hurwitz integers. It is the 
group $S_4$ of permutations of the $4$ basis elements. Both $U$ and $V$ are non-abelian but solvable. 
Group theorists have counted that there are exactly $15$ finite groups of order $24$ up to isomorphisms 
and the two just mentioned groups belong to the interesting ones. One can see $V$ also as a subgroup of the 
{\bf group of automorphisms} of $U$. The set of {\bf inner automorphisms} $x \to e^{-1} x e$ of $U$ 
has order $12$ and is isomorphic to the {\bf alternating group} of four elements. 
But it is not the group of even permutations of $\{1,i,j,k\}$ 
which has the same order: while the rotation $i \to j \to k \to i$ is inner, the rotation $1 \to i \to j \to 1$
is of course not, because $1$ is mapped into itself by any inner automorphism.  \\

Before we move on, lets just remember that we had a similar situation for Gaussian integers. There was
also the group of units $U$ with four elements and the group of automorphisms $V$. 
commutativity has prevented the existence of non-trivial inner automorphisms and conjugation was the only 
outer automorphism, as the units formed the group $U=Z_4$ which has $V=Z_2$ as automorphism group. 
When identifying integers using both groups $U,V$, we factor out the dihedral group $D_4$ leading
to the bijection between rational primes and {\bf Gaussian prime classes}. 

\begin{figure}[!htpb]
\scalebox{0.20}{\includegraphics{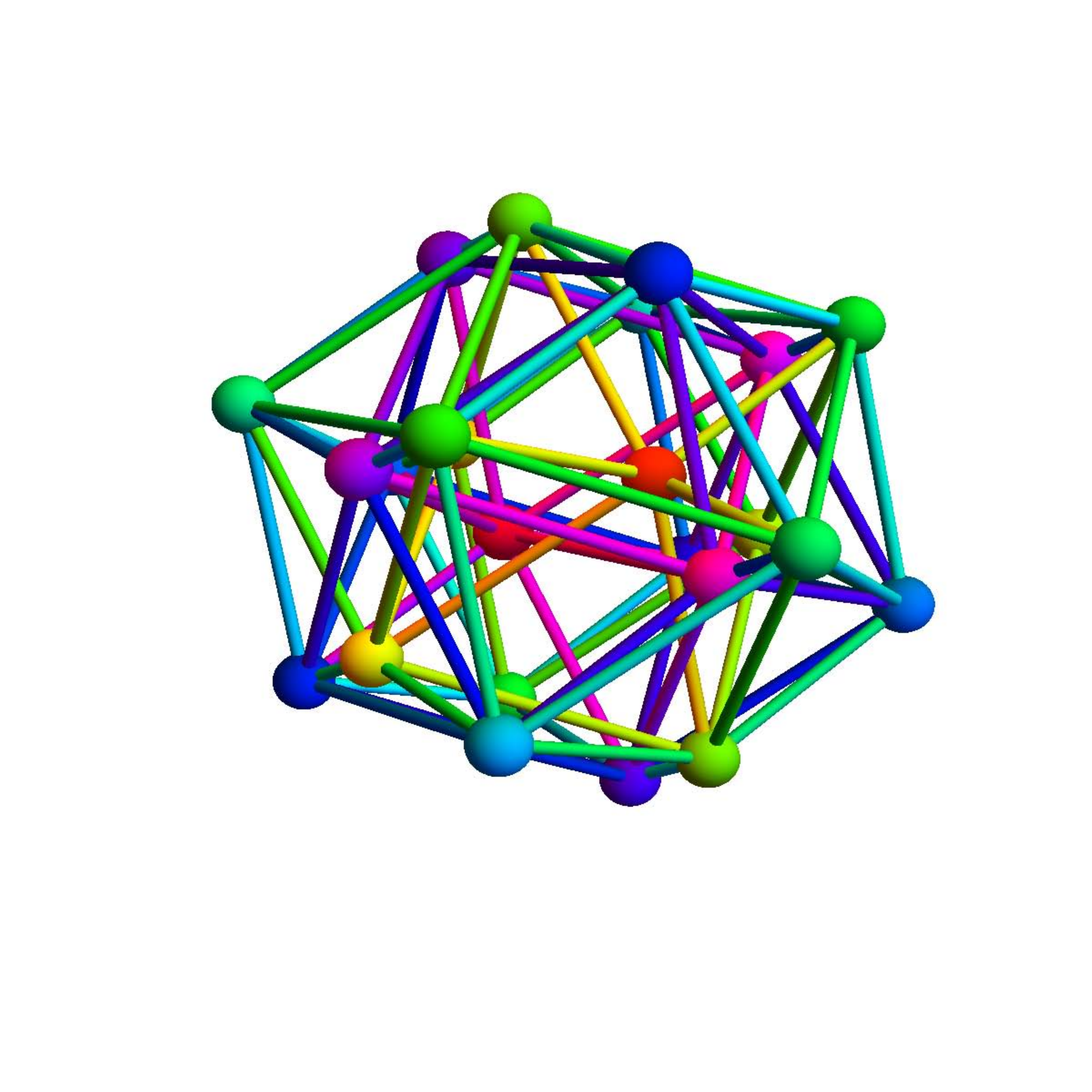}}
\caption{
The units of the Gaussian integers is the {\bf binary
tetrahedral group}. It is a discrete subgroup of $SU(2)$
with 24 elements. It contains the $8$ Lipschitz units located
on the 4 coordinate axes  $\pm 1, \pm i, \pm j, \pm k$ and additionally
the 16 units of the form $(a,b,c,d)/2$ with $a,b,c,d \in Z_2$. The
members of $U$ can be arranged as the vertices of the {\bf 24 cell},
one of the six {\bf Platonic solids} of dimension $3$ realized as regular
4-polytop located on the 3-sphere in four dimensional space. 
}
\label{circles}
\end{figure}

\begin{figure}[!htpb]
\scalebox{0.2}{\includegraphics{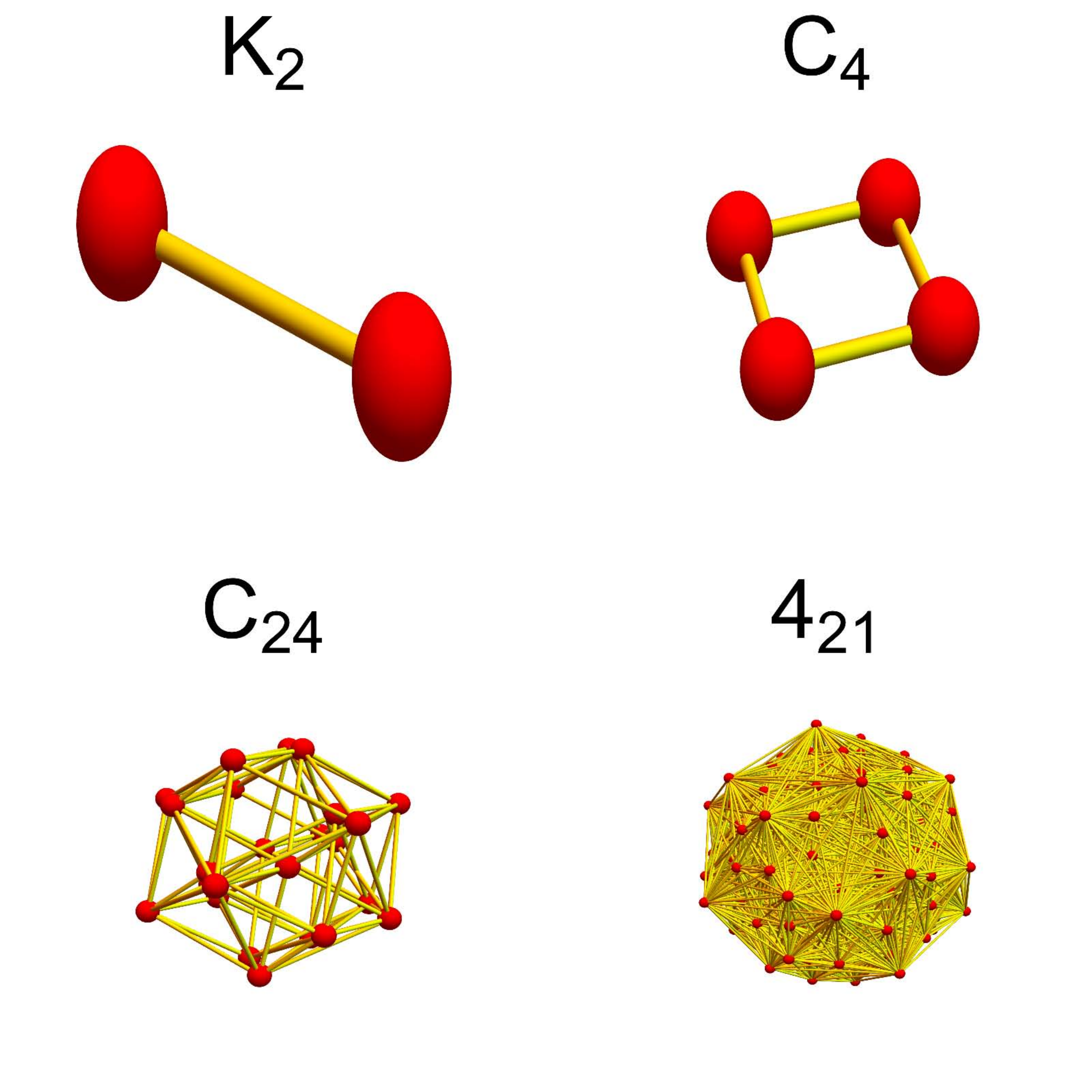}}
\caption{
These are illustrations of the units in the four normed division algebras
$\mathbb{R},\mathbb{C},\mathbb{H},\mathbb{O}$. We see the $K_2$ for $\{ \pm 1\} \subset \mathbb{R}$,
$C_4 = \{1,i,-1,-i\} \subset \mathbb{C}$, the 24 cell in $\mathbb{H}$ and the Gosset polytop
in $\mathbb{Q}$ generating the $E_8$ lattice. The all are known to generate the densest sphere packings,
except for $\mathbb{C}$, where one has to look at the lattice generated by the 
Eisenstein integers instead.
\label{division}
}
\end{figure}

Returning back to the quaternions, we have two groups $U,V$ acting on the ring $X$ of quaternion integers. 
It is so possible to do identifications with either of them and form $X/U$ or $X/V$. We can also factor out
both of them. Note that as $U$ and $V$ are both subgroups of the symmetry group $S_{24}$ of all permutations on 
$U$, one can look at the subgroup of order $48$ generated by both. \\

The {\bf irreducible elements} in in the ring of Hurwitz integers are the {\bf quaternion primes}
or simply {\bf primes} if the context is clear. 
They are known to be Lipschitz or Hurwitz integers $z$ with prime norm $p=N(z)=z \overline{z}$. 
Geometrically, they are the integer or half integer lattice points on the sphere of radius 
$\sqrt{p}$ in $\mathbb{R}^4$. To render terminology simpler, we call the integer coordinates ones the 
{\bf Lipschitz primes} and the others called {\bf Hurwitz primes}. Since by the {\bf Lagrange's 
four square theorem}, every positive integer can be written as a sum of four squares, 
there are no quaternion primes on the four coordinate axes. We can factor an integer $p +0 i + 0 j + 0 k$ 
as $(a,b,c,d) \cdot (a,-b,-c,-d)=a^2+b^2+c^3+d^3=p$. Actually, Hurwitz gave a proof of the Lagrange four
square theorem using quaternions. Hurwitz proved also \cite{Hurwitz1919}:

\resultremark{(Hurwitz)
If $p$ is a rational prime, then the number of quaternion
primes modulo $U$ sitting above $p$ is exactly $p+1$. }

The result holds only for odd $p$ as for $p=2$, there there is only the prime $(1,1,0,0)$. 
Note that $2=(1+i)(1-i)= (1+j)(1-j)=(1+k)(1-k)$ are factorizations of $2$. 
For $p=3$,  the $p+1=4$ prime classes are $i+j+k,(1+3i+j+k)/2, (1+3j+i+k)/2, (1+3k+i+j)/2$.   \\

Let us look at the Hurwitz integers on the sphere of radius $\sqrt{p}$
and factor out first only the automorphism group $S_4$ which commutes the vector entries
of a Hurwitz integer $(a,b,c,d)$. The equivalence classes are Hurwitz integers of the form
$(a,b,c,d)$ with $a,b,c,d \geq 0$ and $a \leq b \leq c \leq d$. We call such a
Hurwitz integer {\bf positively ordered}. \\

In the case $p=2$, the only representative is $(0,0,1,1)$, for $p=3$, the classes are represented by the 
primes $(0,1,1,1)$ and $(1,1,1,3)/2$. As an other example, for $p=13$, the positively ordered  
Hurwitz primes lying above $p$ are $(0,0,2,3)$,$(1,1,1,7)/2$, $(1,1,5,5)/2$,$(1,2,2,2)$, and 
$(3,3,3,5)/2$.  We see that unlike for Gaussian integers, there are now in general 
several equivalence classes of primes sitting above each prime.  \\

Now, these classes can be divided into equivalence classes using the group $U$: 
two positively ordered primes $p,q$ are {\bf equivalent} if either $p/q$ is unit
of $p/\overline{q}$ is a unit. 

\begin{figure}[!htpb]
\scalebox{0.15}{\includegraphics{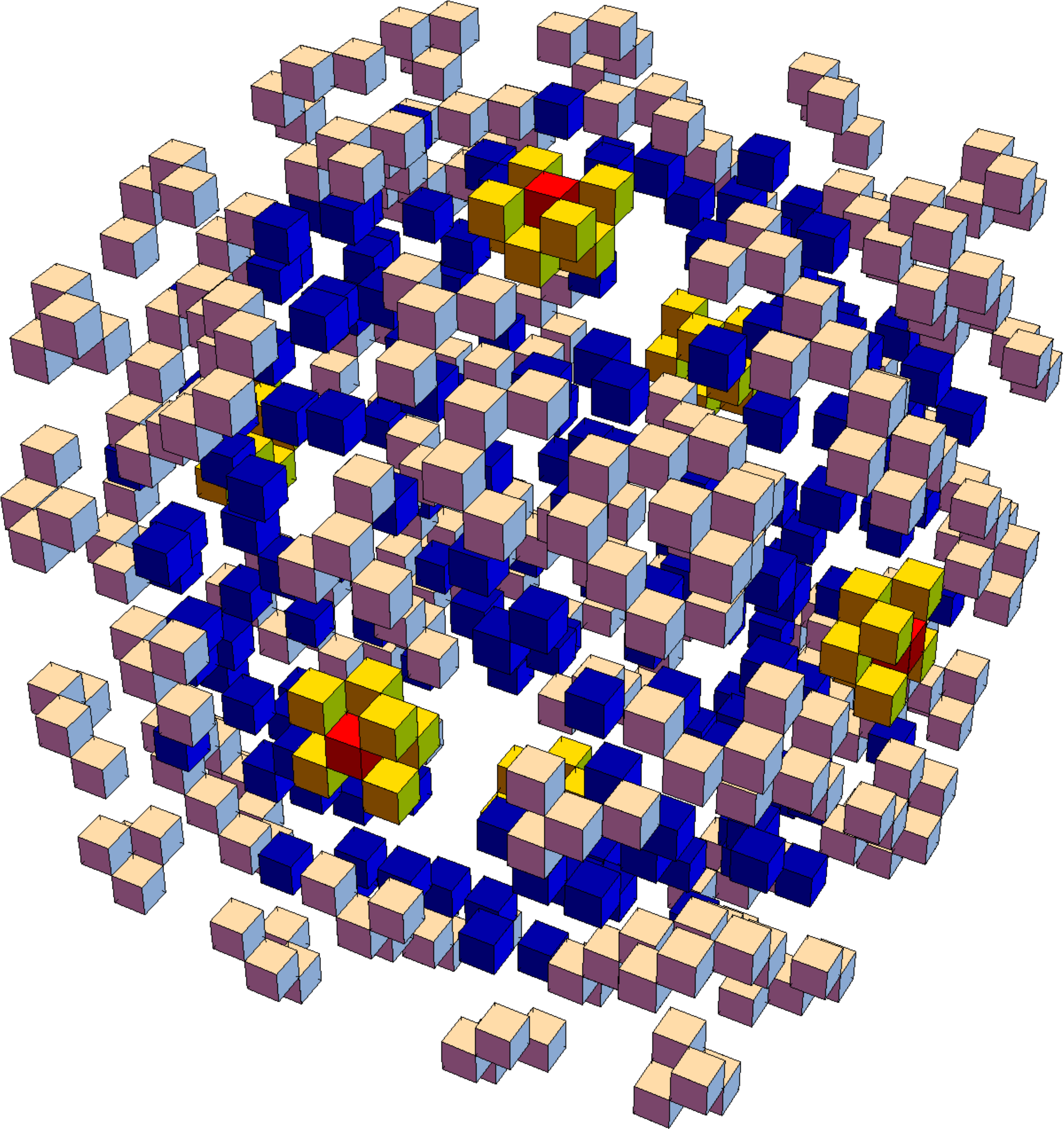}}
\caption{
An illustration of a slice through all Lipschitz primes $(a,x,y,z)$ for $a=701$.
They are given by the integer vectors $(a,x,y,z)$ for which $a^2+x^2+y^2+z^2$ is
a rational prime.
}
\label{circles}
\end{figure}

\resultremark{
The orbits of the group $U$ acting on these positively ordered primes belonging
to an odd rational prime $p$ all have length $2$ or $3$. 
}

\begin{figure}[!htpb]
\scalebox{0.15}{\includegraphics{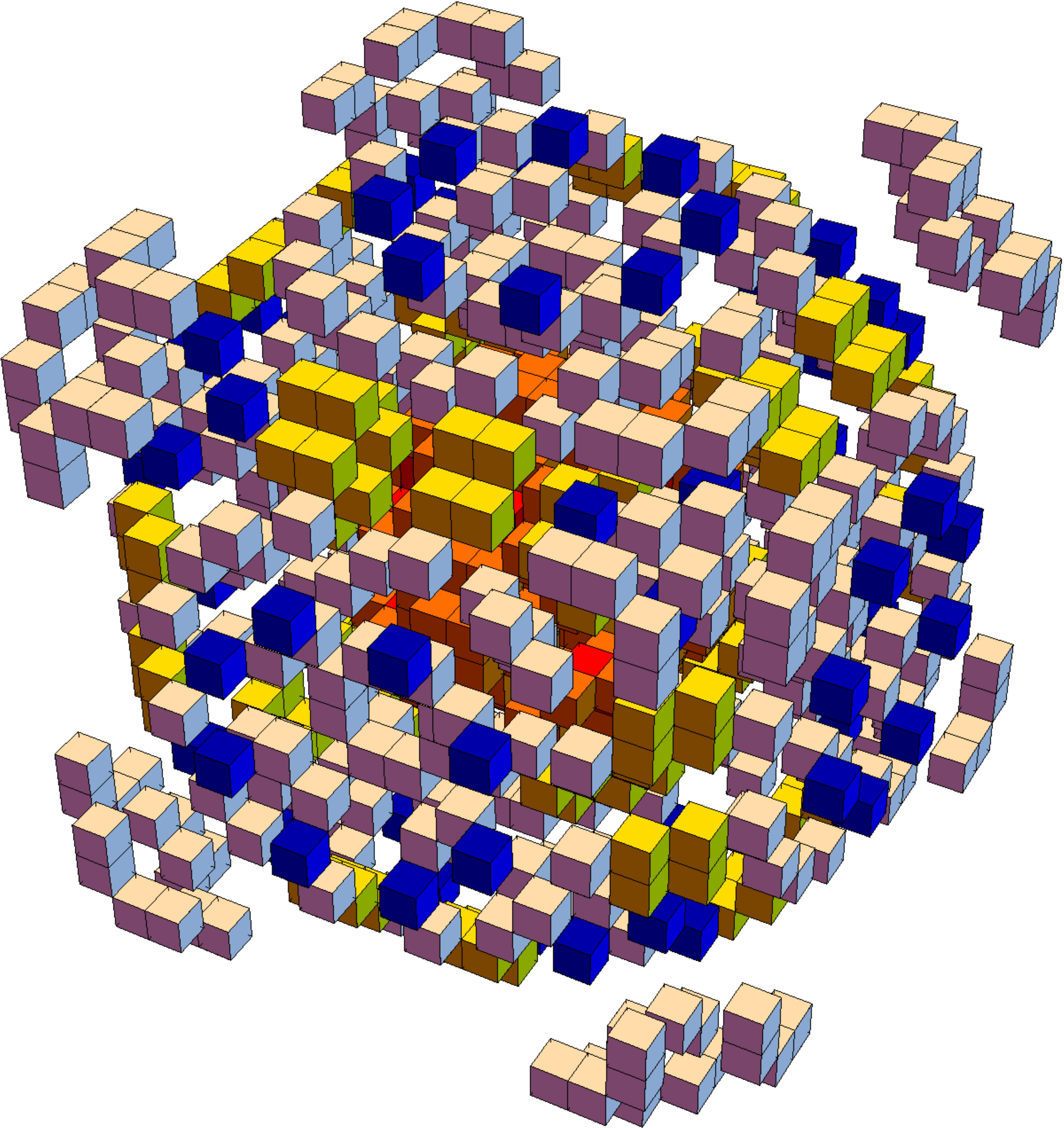}}
\caption{
An illustration of a slice through all Hurwitz primes $(1+2a,1+2x,1+2y,1+2z)/2$ for $a=3001$.
They are given by the integer vectors $(a,x,y,z)+(1,1,1,1)/2$ for which the sum of the
squares is a rational prime.
}
\label{circles}
\end{figure}

For example, there are exactly two positively ordered units. They are 
$\{ (1,0,0,0), (1,1,1,1)/2 \}$. The positively ordered
Lipschitz primes with coefficients $\leq 2$ are
$\{ (0,0,1,2),(0,1,1,1),(1,1,1,2),(1,2,2,2) \}$ with norm $5,3,7,13$.
The positively ordered Hurwitz primes with coefficients $\leq 2$ are
$\{ (1,1,1,3)/2,(1,1,3,3)/2,(1,3,3,3)/2 \}$ with norm $3,5,7$. \\

One can visualize the primes by intersecting them with two dimensional 
planes while the plane moves further and further away from the origin.

\begin{figure}[!htpb]
\scalebox{0.2}{\includegraphics{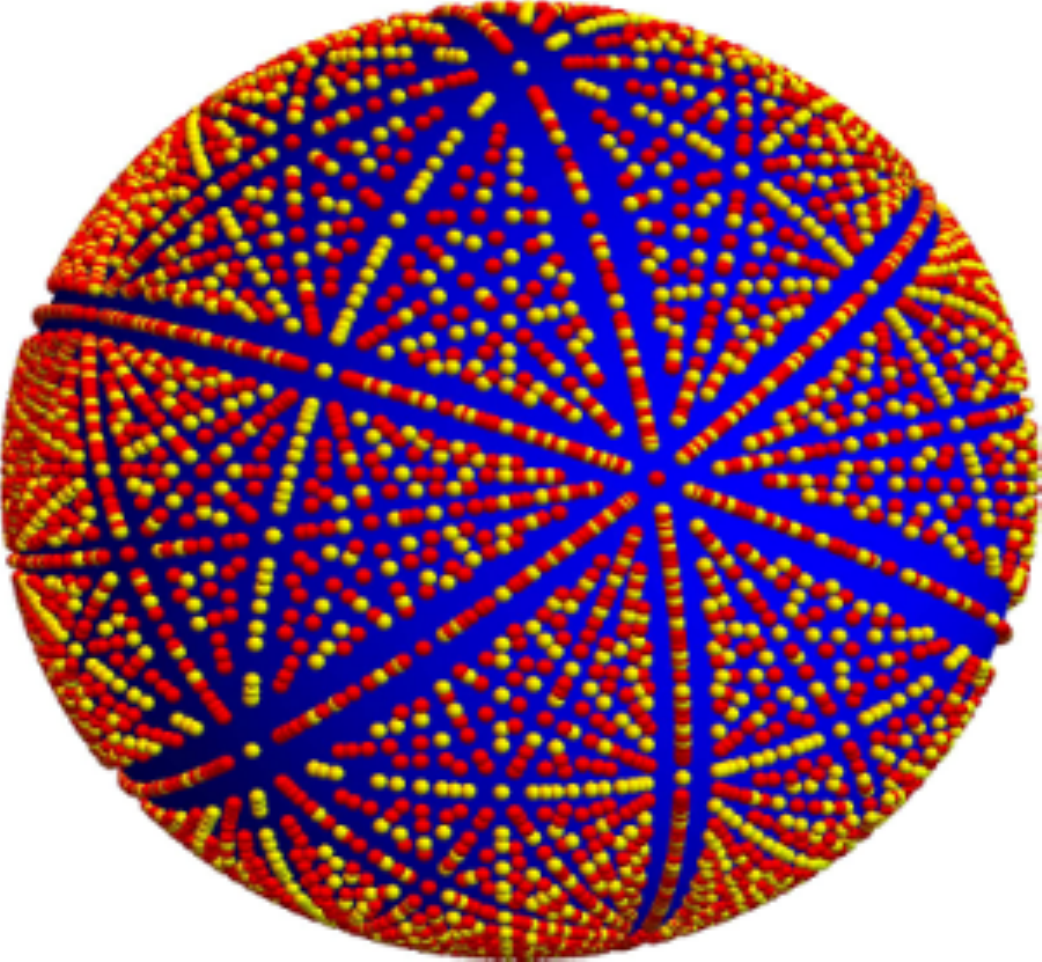}}
\caption{
All the Hurwitz (red) and Lipschitz primes (yellow) $(a,b,c,d)$ with
$|a|,|b|,|c|,|d|<M=8$, projected onto the first three coordinates and then projected onto the
unit sphere. 
}
\label{circles}
\end{figure}

\begin{figure}[!htpb]
\scalebox{0.3}{\includegraphics{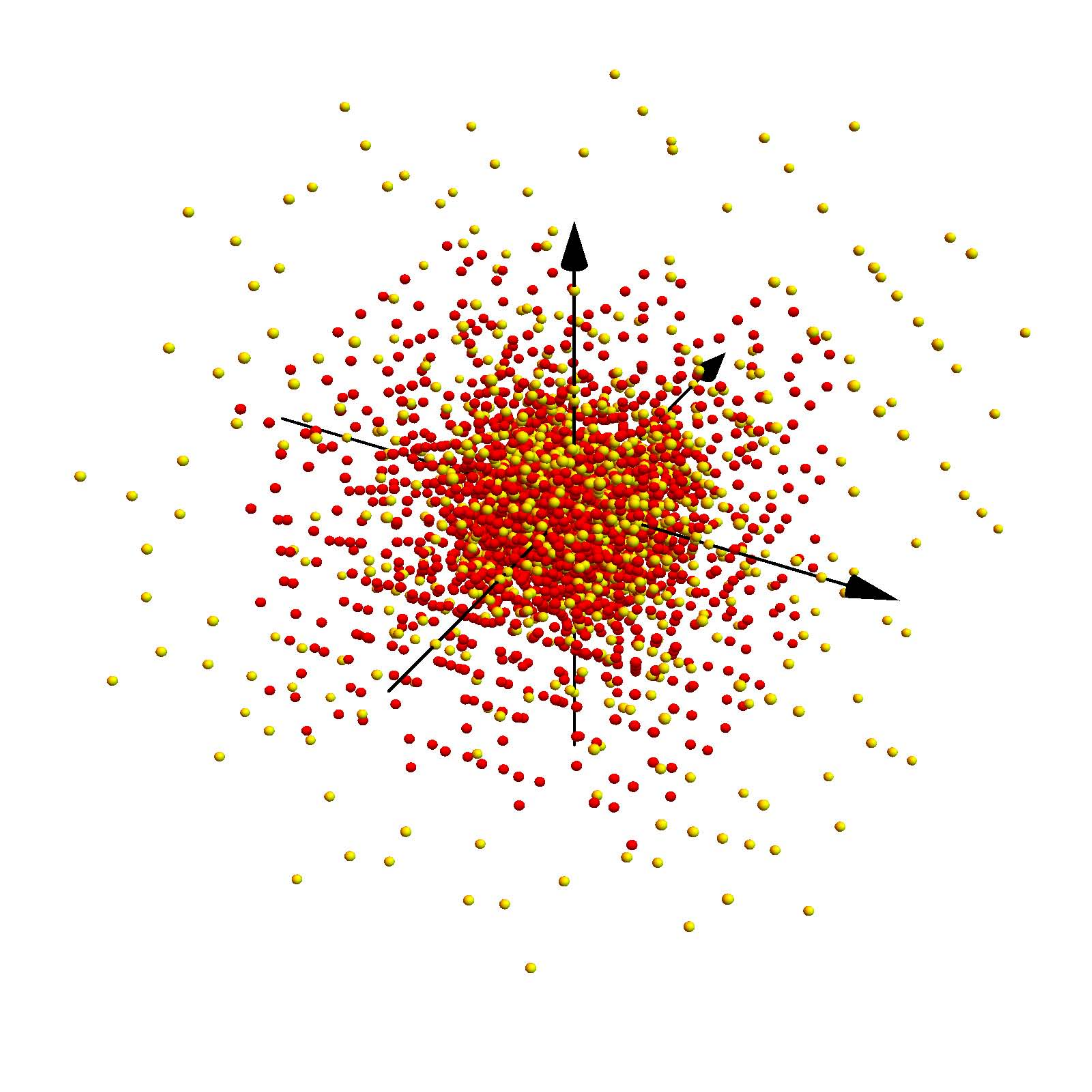}}
\caption{
A stereographic projection of the Hurwitz prime sphere from $\mathbb{H}$
to $\mathbb{R}^3$ but for the cube radius $M=5$.
}
\label{circles}
\end{figure}

\section{Goldbach for Gaussian primes}  

In this section, we formulate a Goldbach conjecture for Gaussian integers.
It is not the first version of the kind:
Mitsui got a version for integer rings in number fields, \cite{Mitsui}, 
a second one is given \cite{HolbenJordan} for Gaussian primes. 
Holben and Jordan use an angle condition: their statement is that
every even Gaussian integer of norm larger than $2$
can be written as a sum of two Gaussian integers $a,b$ such that 
${\rm arg}(a/z),{\rm arg}(b/z) < \pi/4$.
\begin{wrapfigure}{l}{4.1cm} \begin{center}
\includegraphics[width=4cm]{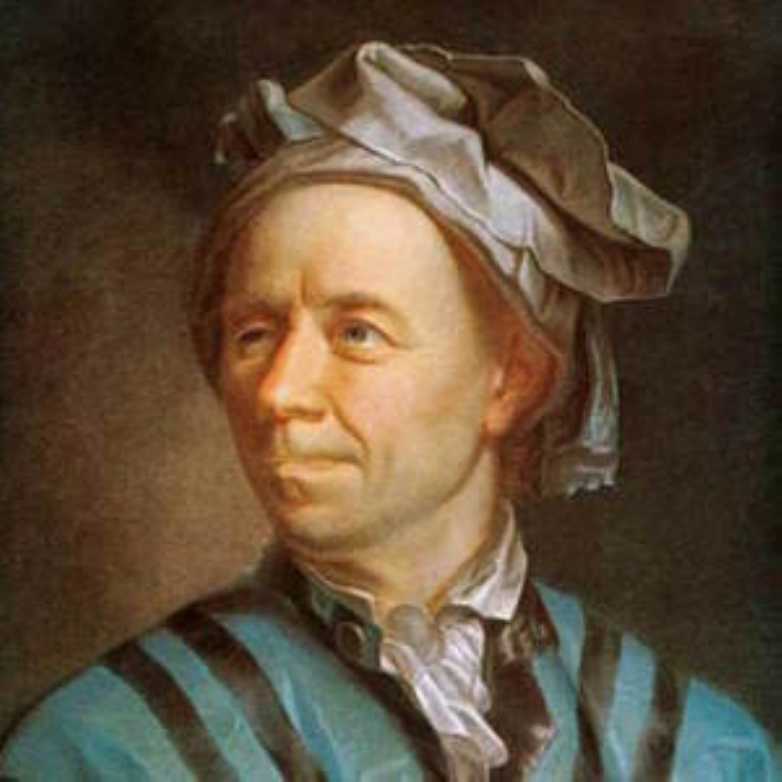}
\end{center} \end{wrapfigure}
Mitsui assumes only a positivity condition for real conjugates and so no positivity condition in the Gaussian case.
He just asks every even $z$ to be representable as a sum of two primes. \\

As Goldbach in $\mathbb{Z}$ deals with the sum of positive primes, we use primes in the
{\bf open first quadrant} of the complex plane: 
$$   Q = \{ a+i b \; | \; a>0,b>0 \} \; . $$ 
A Gaussian integer $z$ is called {\bf even} if it is of the form $z (1+i)$, where $z$ is a 
Gaussian integer. The even Gaussian integers $z$ are of the form $z=a+ib$ with $a+b$ even. Equivalently, they 
are even if $N(z)=a^2+b^2$ is even. The {\bf Goldbach conjecture for Gaussian primes} is:

\conjecture{Every {\bf even} Gaussian integer $z=a+ib$ with $a>1,b>1$ is the sum of two Gaussian primes in $Q$.}

In other word, for every Gaussian integer $z=(1+i) w$ in $Q$ we have a sum $z=p+q$, where $p,q$ are Gaussian primes
in $Q$. It implies the following {\bf ternary Gaussian prime conjecture}:

\conjecture{Every Gaussian integer $a+ib$ in $Q$ satisfying $a>2,b>2$ is the sum of three Gaussian primes in $Q$.}

Details are important. The formulation uses the {\bf open} first quadrant, which does not include points on 
the real and positive axes. It uses a {\bf definite number of summands} and not ``at least" two or 
three summands as $0$ is not counted as a prime. We also do not mean that the primes are required to 
be different. The even Gaussian integer $2+2i$ for example can only be written as the sum of two identical primes 
$p=1+i$ and $q=1+i$. It could well be that for large $|z|$, the additional requirement that the summands
have to be different could be imposed additionally. One could weaken the conjecture also by allowing the
summands to be in the {\bf closed first quadrant} 
$\overline{Q} = \{ a+ib \; | \; a \geq 0, b \geq 0 \}$,
but there is strong enough evidence for the open version. The open version given here also has the advantage
that it relates to the difficult Landau problem, so that it is not in danger of being ``trivial". As any
conjecture it is in danger to be wrong however. 

\begin{figure}[!htpb]
\scalebox{0.5}{\includegraphics{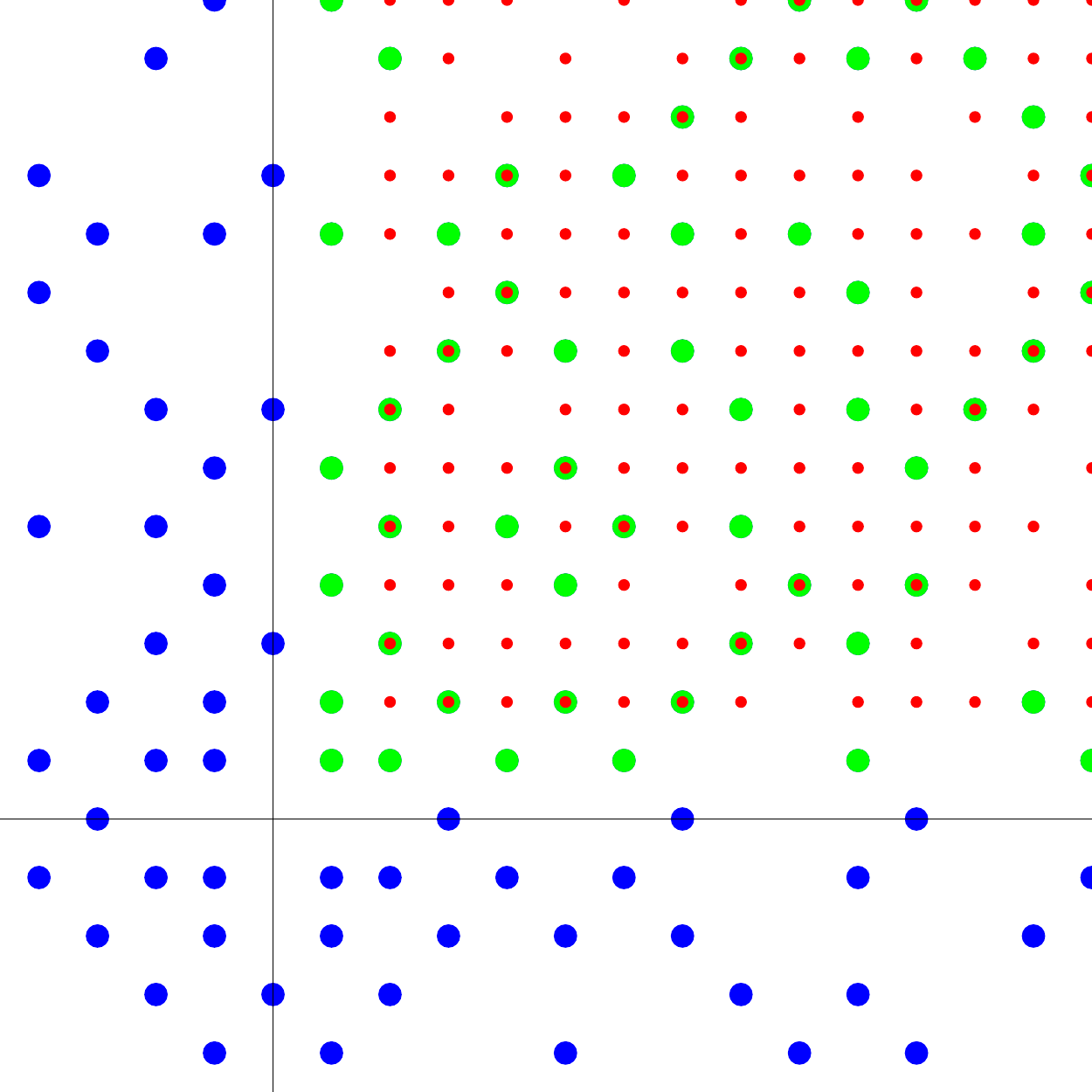}}
\caption{
The Gaussian Goldbach conjecture.
}
\label{circles}
\end{figure}

\begin{figure}[!htpb]
\scalebox{0.5}{\includegraphics{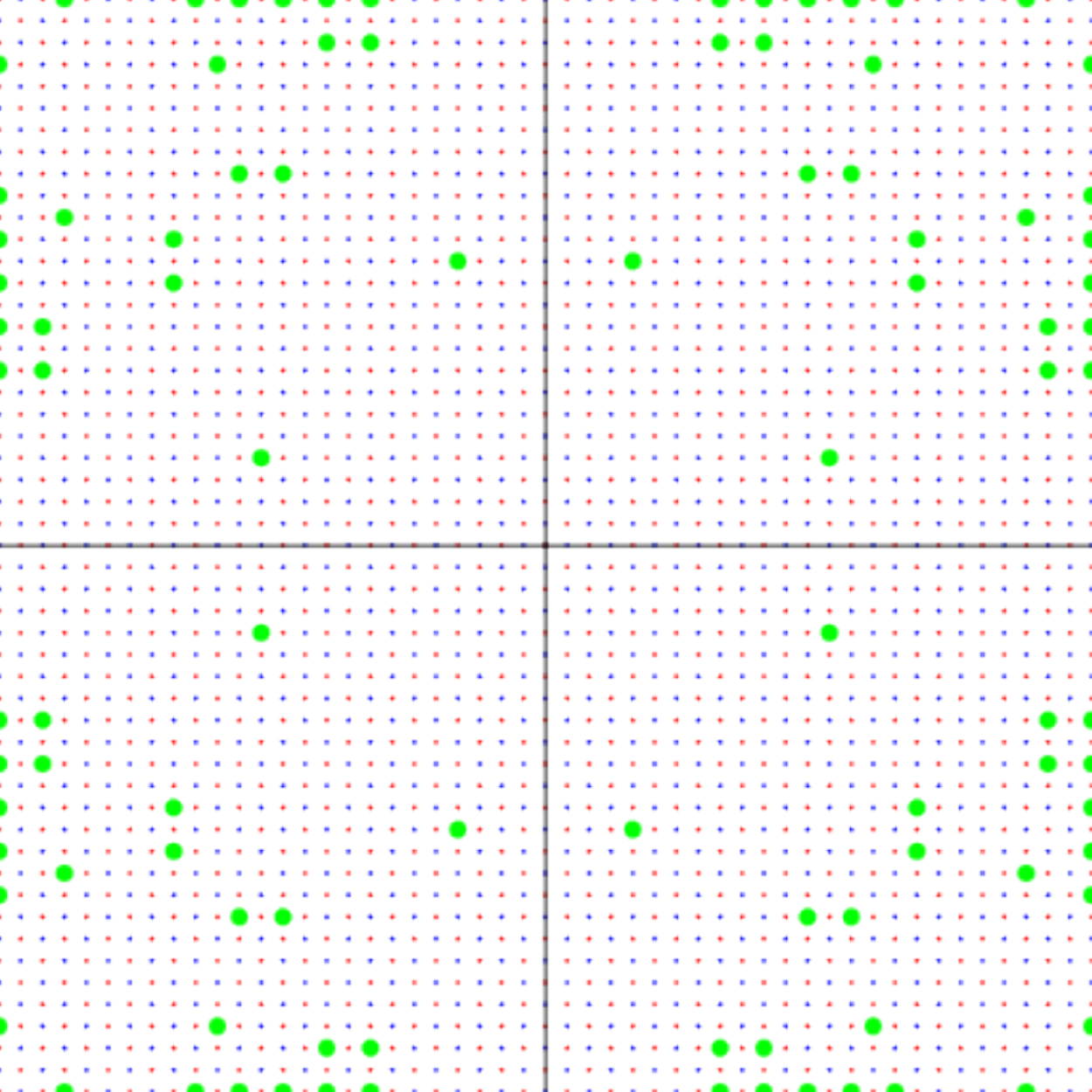}}
\caption{
Which Gaussian integers are not the sum of two Gaussian primes? The smallest
example is $4+13i$. A rather weak Gaussian Goldbach conjecture predicts that
all such examples are odd. 
}
\label{circles}
\end{figure}

The statement is natural because it can be reformulated for 
Taylor or Fourier coefficients of powers $f^2$ and $f^3$, where a
{\bf Gaussian Goldbach function} given by
$$  f=\sum_{a+i b \in P \cap Q} x^a y^b/(a! b!) $$ 
sums over the set of Gaussian primes in $Q$. For $x=e^{i \phi}/2, y=e^{i \theta}/2$, 
one gets smooth double periodic functions $f(\phi,\theta)$. These are Gaussian Goldbach functions
on the two-dimensional torus. The factorial terms are an example but assure that the functions are smooth. 
As for rational primes, the Goldbach problem for Gaussian primes
is now a {\bf calculus problem} which lies at the basis of the {\bf Hardy-Littlewood circle method}.
A naive hope is to find coefficients $g_{ab}$ with support on the Gaussian primes
such that $f=\sum g_{ab} x^a y^b$ is expressible as a finite sum of known functions,
allowing to check the conjecture by computing Taylor or Fourier
coefficients of $f^2$: in the Taylor picture, the Gaussian Goldbach conjecture. 

\resultremark{
The Goldbach conjecture follows from the statement that $(f^2)_{kl} \neq 0$ for every pair $k,l>1$
for which $k+l$ is even. Alternatively, it follows from $\hat{f^2}_{k,l} \neq 0$ for every $k,l>1$
with $k+l$ even. 
}

There are interesting subproblems of the Goldbach problem for Gaussian integers. One can restrict
for example to the {\bf diagonal} and ask whether every Gaussian integer $k+ik$ with $k>1$ can be 
written as a sum of two primes with positive entries. Here is the {\bf diagonal Gaussian Goldbach}
conjecture: 

\conjecture{Every Gaussian diagonal integer $k+ki$ with $k>1$ can be written as a sum of
two Gaussian diagonal primes $a+ib, c+id$ in $Q$.}

While Gaussian Goldbach would imply the diagonal version, the point of a diagonal version is
that it appears to be much easier so that there is some chance that it could be solved with some
already available techniques. As the map from a $(4k+1)$-prime $p=a+ib$ to $a+b$ 
with $b<a,a>0,b>0$ is defined, the question is related to the existence of {\bf large gaps}
on the class of $4k+1$ primes. If gaps near $x$ are of size smaller than $O(\sqrt{x})$
then the result follows. The classical {\bf Andrica conjecture} asks whether the gap size is always smaller than 
$2 \sqrt{x}+1$. The diagonal Goldbach subquestion appears easier as there are in general many more
ways to write $k+ki$ as a sum of two Gaussian primes. Also, the general belief is that
the prime gaps are much smaller: the {\bf Cramer conjecture} for example 
predicts that gaps are not bigger than $C \log(x)^2$. And the same
can be expected for gaps in the set of primes of the order $4k+1$.  \\

An other possibility is to weaken the conjecture and for a Schnirelman type result:

\conjecture{
There exists a constant $M$, such that every Gaussian integer $z = a+ib$ with $a \geq M,b \geq M$ 
in $Q$ can be written as a sum of maximal $M$ Gaussian primes in $Q$. }

This statement looks more approachable as the corresponding statement for rational 
integers has shown. Schnirelman proved the rational case using the concept of {\bf Schnirelman density}. 
But it still would imply the Landau problem. And the boundary case
is the Achilles heel as it could be that the Schnirelman density is zero at the boundary. 

\begin{figure}[!htpb]
\scalebox{0.22}{\includegraphics{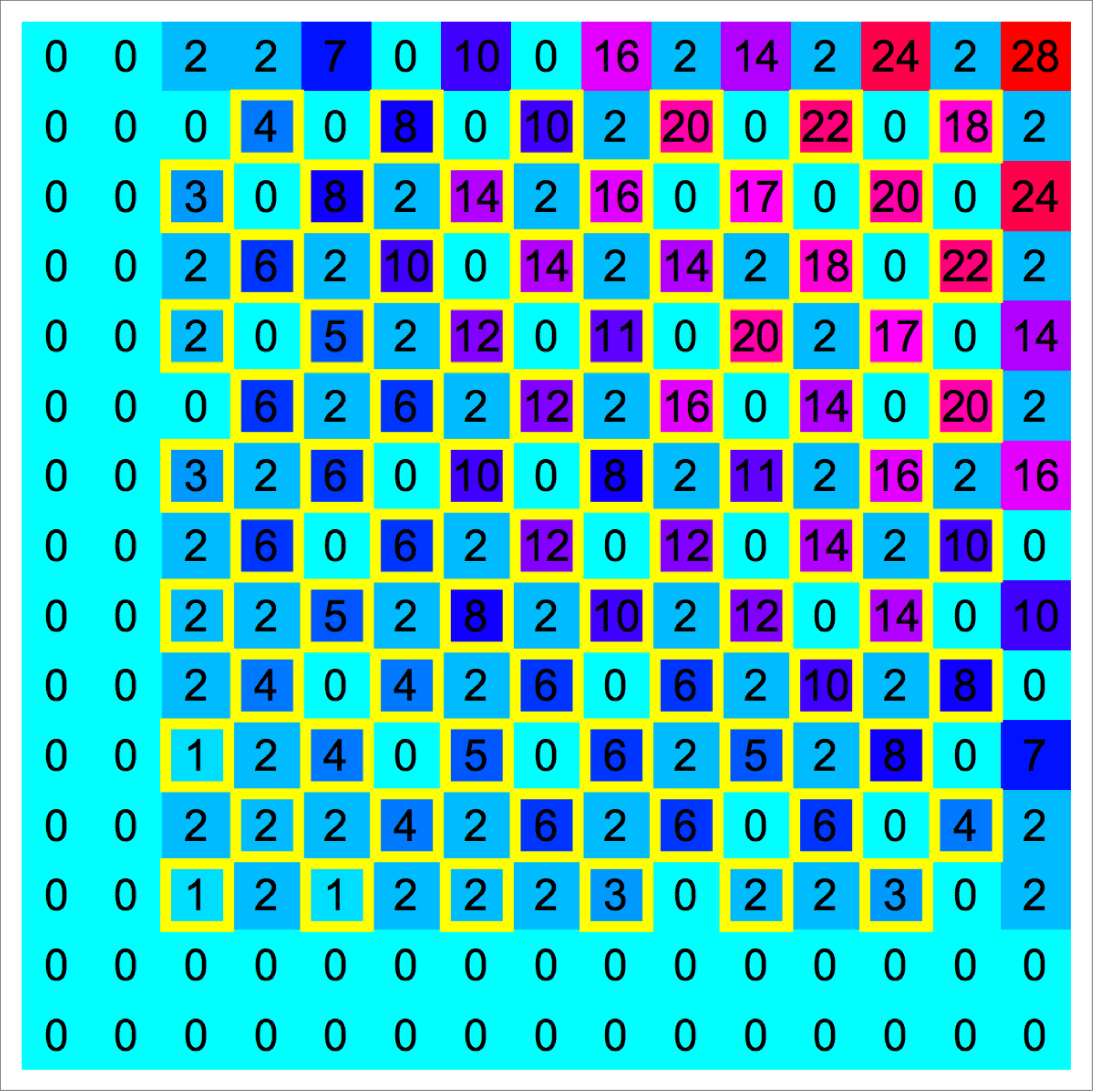}}
\caption{
The {\bf Gaussian Goldbach problem} for the open first quadrant $Q={\mathbb Z}^+[i]$ claims that
every even $a+ib \in {\mathbb Z}^+[i]$ with $a,b \geq 2$ can be written as a sum of 
exactly two primes in ${\mathbb Z}^+[i]$.
The matrix entries in the picture show how many times a Gaussian integer can be written as a sum
of two primes. A first challenge is to establish the claim for the diagonal. 
The hardest appears the boundary row, where the result implies Landau's open problem
about the infinitude of primes of the form $n^2+1$ because the lowest row entries
$z=c+2i$ must be the sum of two primes of the form $a+i, b+i$ which both must have
prime  norm $a^2+1,b^2+1$.
Algebraically, the Goldbach conjecture is a claim about the Goldbach function
$f=\sum_{a+ib \in P^+} x^a y^b$. The claim is that the 
even derivatives $(f^2)_{kl}$ are positive if $k,l \geq 2$. 
\label{goldbachmatrix}
}
\end{figure}

\begin{figure}[!htpb]
\scalebox{0.2}{\includegraphics{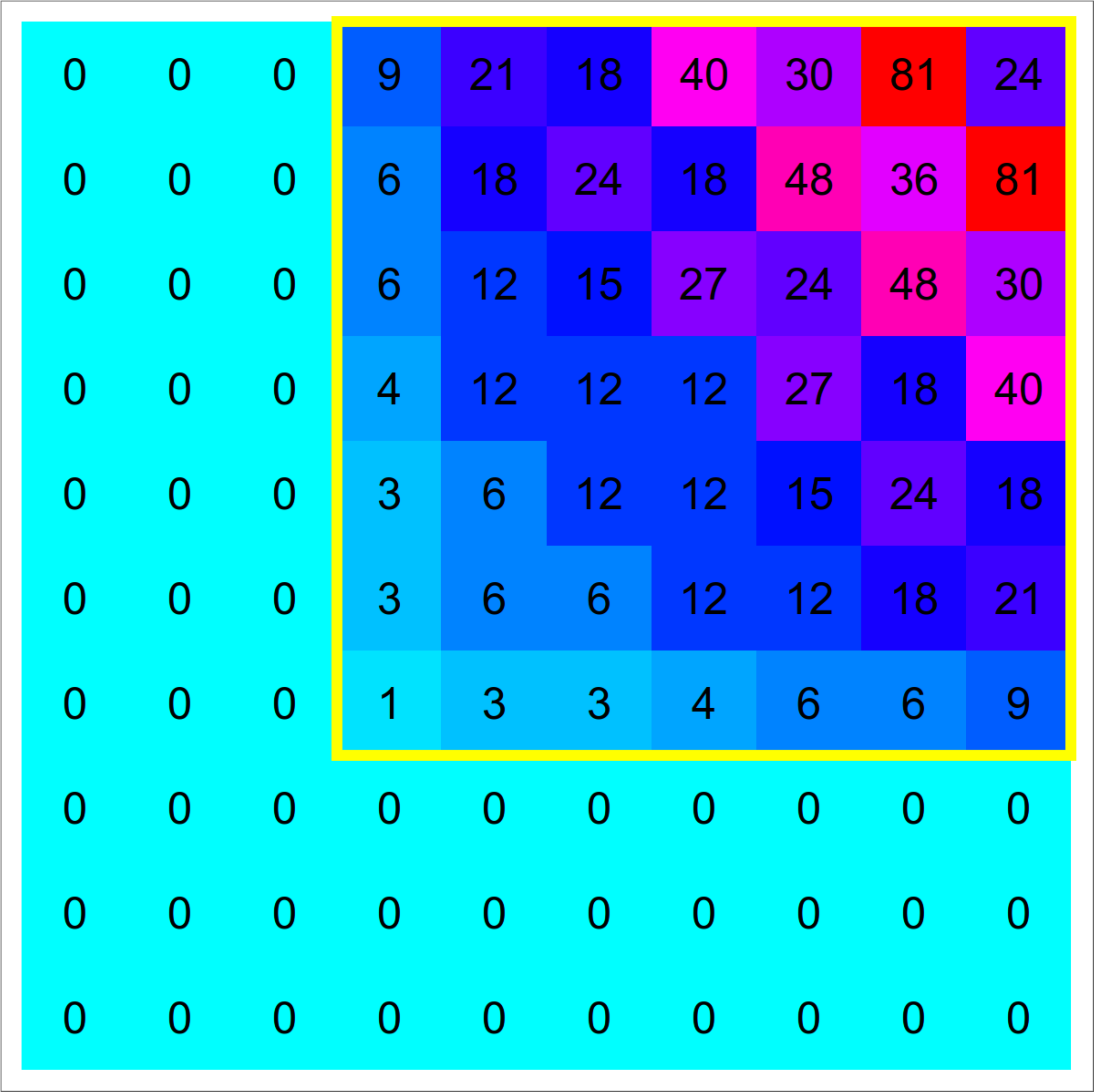}}
\caption{
The {\bf ternary Gauss Goldbach problem} claims that every Gaussian integer $a+ib \in {\mathbb Z}^+[i]$
with $a,b > 2$ is the sum of exactly three Gaussian primes in ${\mathbb Z}^+[i]$.
Algebraically, it tells that coefficients of $f^3$ are positive if $a \geq 3, b \geq 3$. 
\label{ternarymatrix}
}
\end{figure}

\begin{figure}[!htpb]
\scalebox{0.13}{\includegraphics{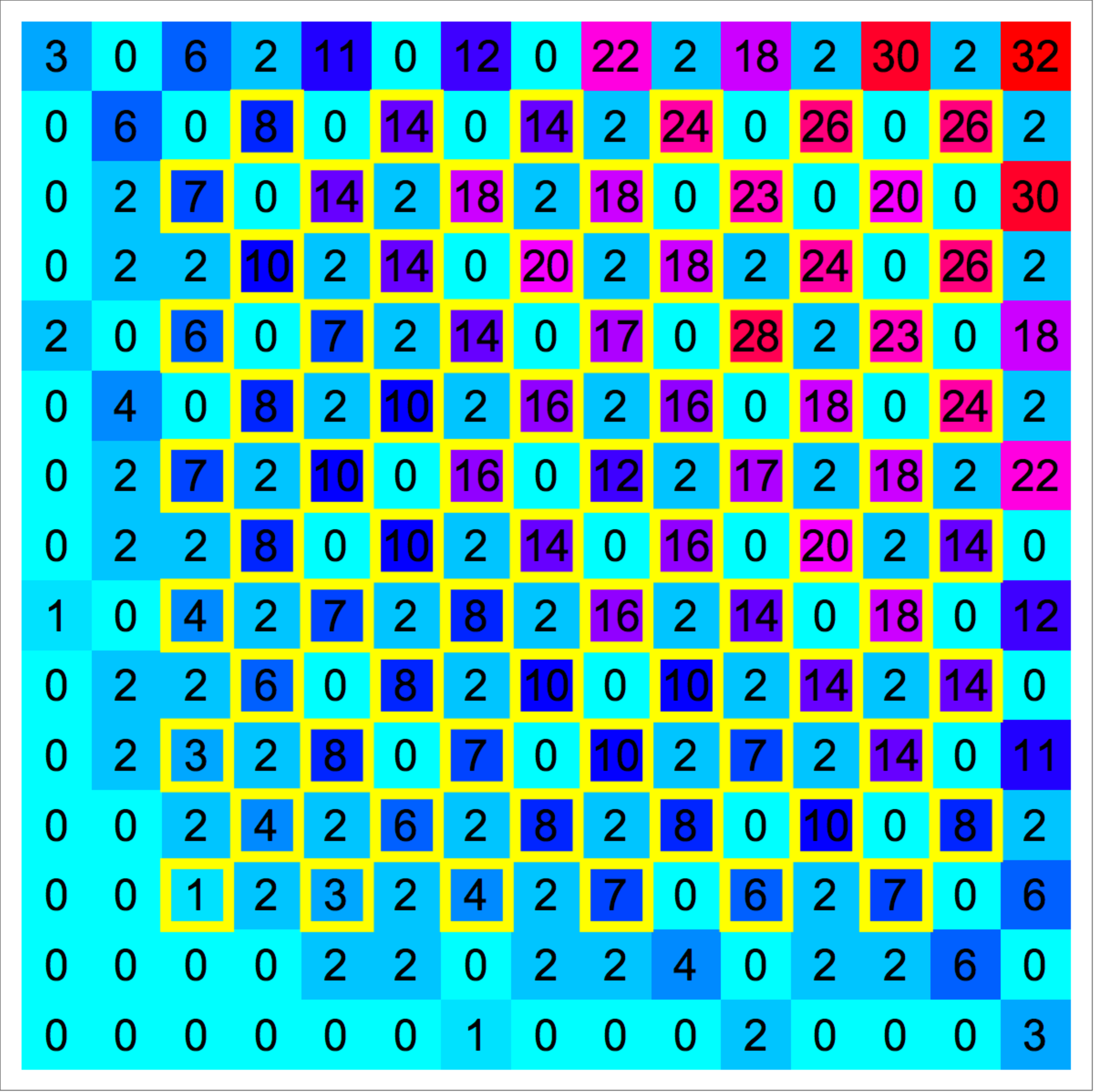}}
\scalebox{0.13}{\includegraphics{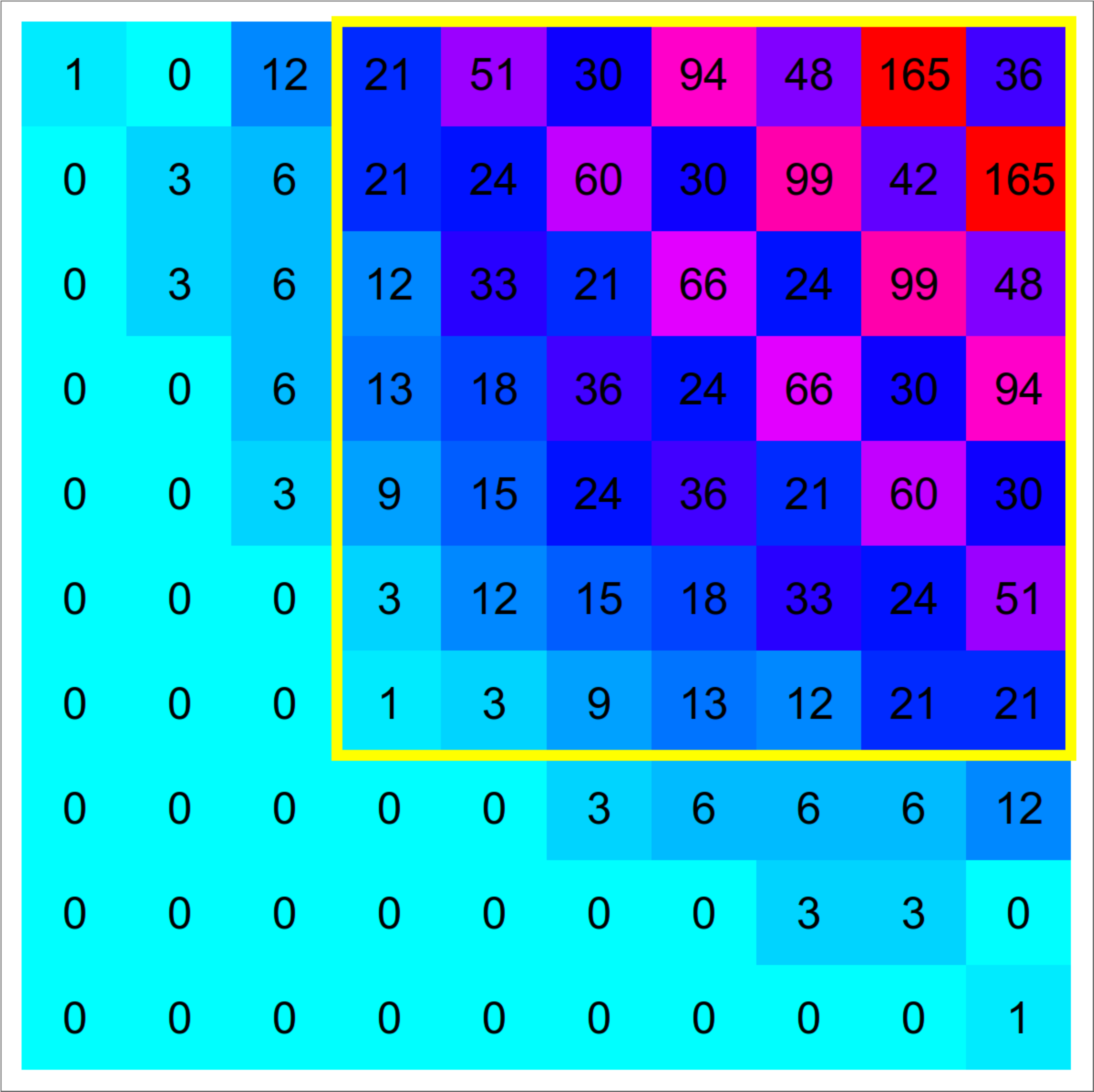}}
\caption{
A weaker {\bf Goldbach problem} for Gaussian integers is that
for every even $n+m$ and $n,m \geq 2$ (now including the boundary of $Q$),
the prime $z=n+i m$ can be written as a sum of two Gaussian primes $z=p+q$.
\label{weak goldbachmatrix}
}
\end{figure}

The ternary Goldbach problem is implied by the Gauss Goldbach problem and also
weaker. Again, establishing the diagonal case appears to be the easiest:
write every integer $k \geq 3$ as a sum $k=a+b+c=d+e+f$,
where $a^2+d^2,b^2+e^2,c^2+f^2$ are prime.
Establishing the result for the third row $k+3i$ again would imply the open
Landau problem about the infinitude of primes of the form $n^2+1$. So, also the
ternary Goldbach problem appears to be hard.

\begin{figure}[!htpb]
\scalebox{0.12}{\includegraphics{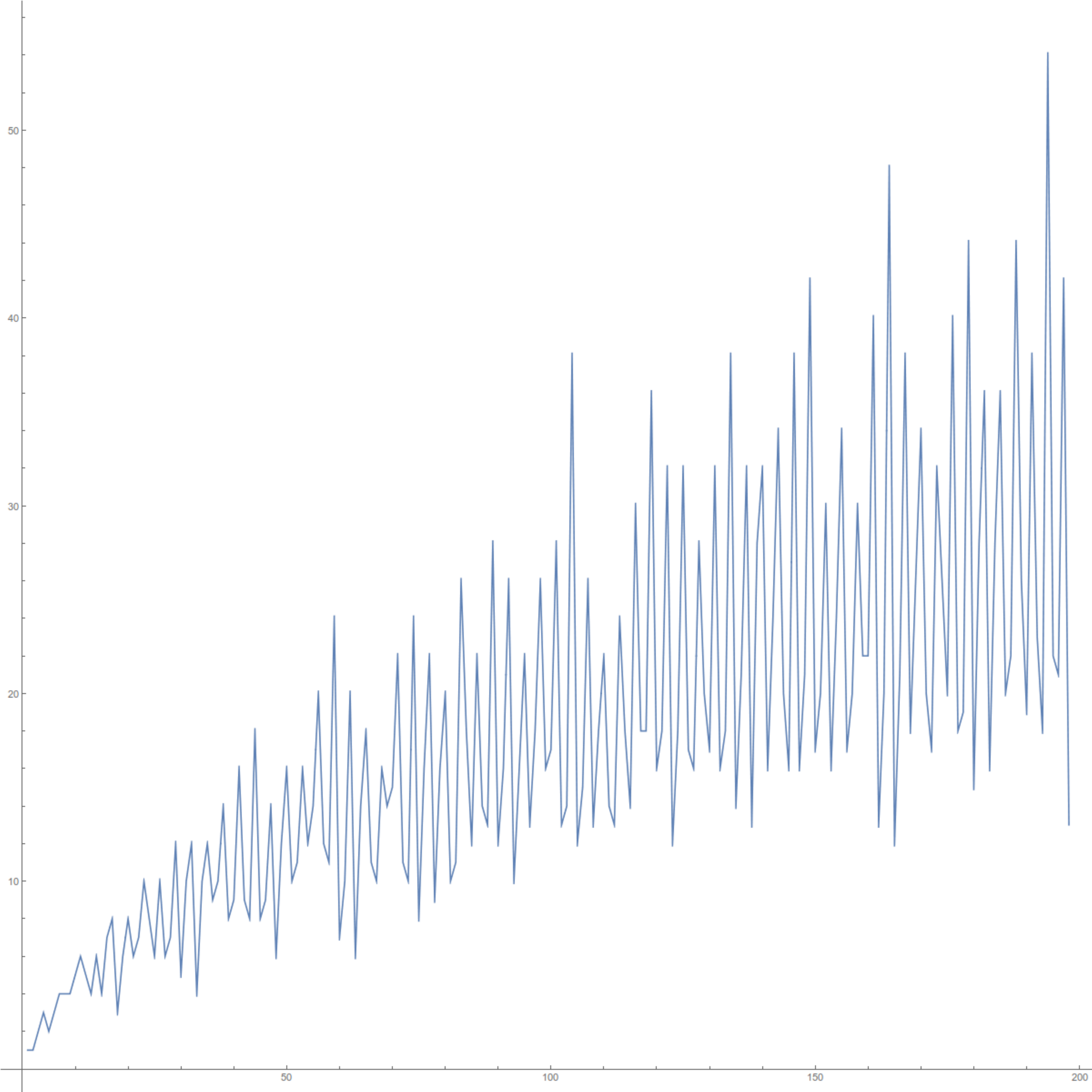}}
\scalebox{0.12}{\includegraphics{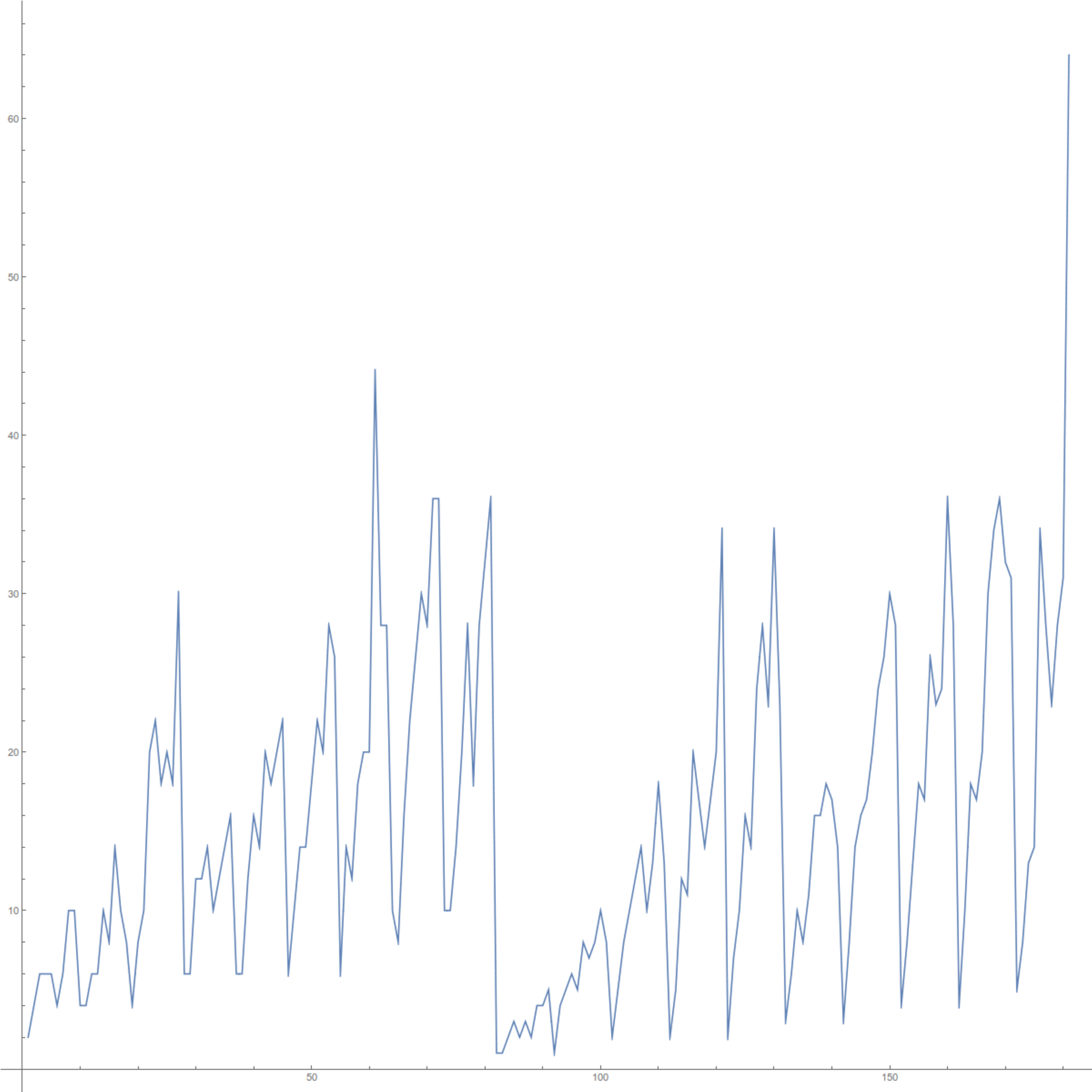}}
\caption{
The coefficients of $g(x) = f(x)^2$ with
$f(x)=\sum_{p} x^p$ form the {\bf Goldbach comet}
for primes. The coefficient $x^n$ tells in how many ways one
can write $n$ as a sum of two primes.
The picture to the right shows the comet in the case of Gaussian
primes, where the coefficient of $x^n y^m$ tells in how many times the
Gaussian integer $n+i m$ can be written as a sum of two Gaussian primes.
The Goldbach conjecture for Gaussian primes claims that every
Gaussian integer $2n+2m i$ with $n,m>0$ is the sum of two Gaussian primes
with nonnegative coordinates.
}
\end{figure}

What about the following statement?

\conjecture{
Every even Gaussian integer is the sum of two Gaussian primes. 
}

One can not leave away the evenness condition. Already the real number $29$ is not a sum of two
Gaussian primes. It is actually the second largest integer which can not be written as a sum or
difference of rational primes. This statement looks like the type of statement which Mitsui looked at in
number fields. \\

{\bf Remarks.} \\
{\bf 1)} The Goldbach statement also implies that there are are infinitely many rational primes of the form
$n^2+4$. Proof: assume there are only finitely many. If every even Gaussian integer of the form $n + 3i$
could be written as a sum of two primes as such a summation would require one to be of the form
$p+i$, the other $q+2i$. Assume the set of primes of the form $q+2i$ is finite.
It would imply that every even integer $w$ of the form $n+3i$ has a prime $w$ within some
fixed distance $R+3$. Define $Z$ as the product of all Gaussian integers $z=a+ib$ with
$|z| \leq R$. This number has no prime in a neighborhood $B(Z,R+3)$ because any $P=Z+A$ has
the factor $A$. Therefore, also $Z+3i$ has no prime in a neighborhood of radius $R$.

{\bf 2)} Lets look again at \cite{HolbenJordan} which the formulate as ``conjecture E": 
if $n$ is a Gaussian integer with $n \overline n = N(n)>2$, then it can be written as a 
sum $n=p+q$ of two primes for which the angles between $n$ and $p$ as well as $n$ and $q$ are 
both $\leq \pi/4$. \cite{HolbenJordan} can strengthen their statement differently: 
their ``conjecture F" claims that for $N(n)>10$, one can write $n$ as a sum of two primes 
$p,q$ for which the angles between $n$ and $p$ and $n$ and $q$ are both $\leq \pi/6$.
The angle could be reduced, but not arbitrarily since $4$ is then not the sum of two Gaussian 
primes. Still, the {\bf Gaussian Goldbach conjectures with angle conditions} is different from the
open quadrant formulation given here. Our condition is easier to reformulate 
algebraically and also can be stated almost identically in all division algebras. \\

{\bf 3)} \cite{Mitsui} formulates a conjecture in general number fields. It does not use a cone restriction.
Gaussian and Eisenstein primes are special cases. In the Gaussian case, it states that every even Gaussian integer is the 
sum of two Gaussian primes. The evenness condition is necessary. The smallest Gaussian integer which 
is not the sum of two Gaussian primes is $4+13i$. The Holben-Jordan conjecture implies the Mitsui statement.
In the Eisenstein case, we see
that every Eisenstein integer is the sum of two Eisenstein primes without evenness condition. 
The question makes sense also in the $\mathbb{Z}$: is every even integer the sum of two {\bf signed primes},
where the set of {\bf signed primes} is $\{ \dots, -7,-5,-3,-2,2,3,5,7, \dots \}$. The smallest number
which is not the sum of two signed primes is $23$. All even numbers seem to be the sum of two signed primes!
We are not aware of a proof of this {\bf signed Goldbach statement}. It might not be of interest,
since Goldbach implies it. \\

{\bf 4)} The Goldbach conjecture is not the only statement which involves the 
additive structure and primes (which inherently rely on the multiplicative 
structure of the ring): any additive function $f(zw)=f(z)+f(w)$ which satisfies
$f(p+1)=0$ for all Gaussian primes is $0$ \cite{MehtaViswandham}. \\

{\bf 5)} We do not know for sure whether the angle or open quadrant statement are not
equivalent but believe the open quadrant statement harder (due to the Landau boundary
problem). Both have their advantages: the angle version can be strengthened or weakened
by changing the angle. The quadrant version of Goldbach is algebraically natural, especially 
when reformulating using Taylor series where no negative powers can occur. \\

{\bf 6)} The work \cite{Kubilyus} extends Rademacher, Hecke and Vinogradov
to show that the number of primes in a sector ${\rm arg}(z) \in {\alpha,\beta},
|z| \leq r$ is $(\beta-\alpha)/(2\pi) \int_2^r dx/\log(x)$ with
an error of the form $O(r e^{-c \sqrt{\log(r)}})$ meaning that
in any radial sector one sees the distribution of the prime number theorem.
The prime twin and Goldbach problem has been stated in \cite{HolbenJordan} using some
angle conditions. Holben-Jordan also define {\bf Gaussian prime twins} as Gaussian
prime pairs of distance $\sqrt{2}$. (We believe the prime twin problem for Gaussian
primes must have been asked before but we don't find an earlier source than 
Holben-Jordan from 1968.)
The twin prime conjecture was also featured in Hilbert's problem 8 and is part
of the list of Landau's problems. 

\section{Goldbach for Hurwitz primes}

Most number theoretical question can also be asked for Hurwitz integers: one can in particular look at the 
{\bf twin prime} or {\bf Goldbach problems}. Also the {\bf zeta function} $\sum_z 1/N(z)$ can be considered for
Hurwitz primes. It will turn out to be a scaled and shifted Riemann zeta function because it 
can be written as $\sum_p (p+1)/p^s$ by Hurwitz's result. But lets first look at Goldbach for Hurwitz. 
It might more likely to be true as we are in higher dimensions but numerical experiments already become
costly. It would not be a complete surprise if there was large counter example. \\

Lets look first a bit at the well known history of the one-dimensional problem \cite{Guy}.
The Goldbach conjecture for rational primes was proposed in 1742 by 
Christian Goldbach in a letter to Leonard Euler \cite{HelmutKoch,GoldbachBuch}. 
Goldbach's original statement was that every 
\begin{wrapfigure}{l}{4.1cm} \begin{center}
\includegraphics[width=4cm]{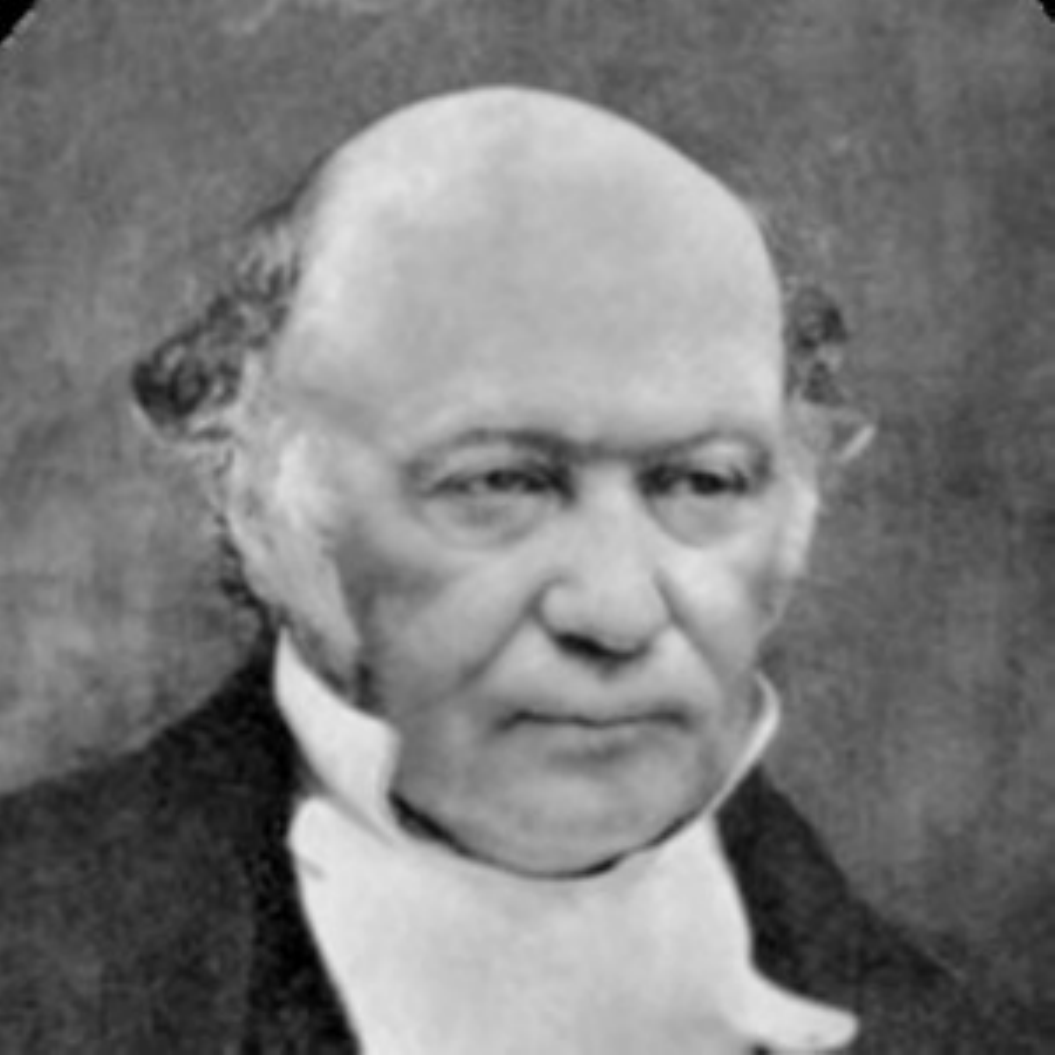}
\end{center} \end{wrapfigure}
integer $n \geq 5$ is a sum of three primes, which Euler reformulated as every even function 
$\geq 4$ being a sum of $2$ primes. Goldbach's footnote in the handwritten letter is 
well readable. Some transcribe it as $>2$ rather than $\geq 5$ which would imply Goldbach considered $1$
as a possible summand in "aggregatum trium numerorum primorum" as custom at the time \cite{Pintz}.
The statement that every odd integer larger than 4 can be written as a sum of three primes has
first been formulated by Waring in 1770 \cite{Pintz}. 
The book \cite{IwaniecKowalski} transcribes on page 443 with $\geq 5$. More on the correspondence
can be found in \cite{HelmutKoch}. 
As Erd\"os liked to point out, Descartes already earlier voiced a similar 
conjecture \cite{Schechter,Vaughan} so that the Goldbach formulation sometimes is also called 
{\bf Descartes Conjecture} \cite{Pintz}. Dickson \cite{dicksonI} mentions on page 421 that Descartes
formulated the conjecture in the form that every even number is the sum of 1,2 or 3 primes and
that Waring conjectured in 1770 that every odd number is either a prime or the sum of three primes.
Many have done experiments. Even Georg Cantor checked Goldbach up to 1000 \cite{dicksonI} p. 422. 
Much progress has been done.
Landmark results in theory were Hardy-Littlewood \cite{HardyLittlewood1923} with the circle method,
the Lev Schnirelemans theorem \cite{Hua1982}
using density and Ivan Vinogradov's theorem \cite{Vinogradov} using trigonometric sums.
The genesis paper for the circle method is in \cite{Vaughanbook} attributed to the paper of Hardy with
Srinivasa Ramanujan in 1918 \cite{HardyRamanujan}.
Chen's theorem \cite{ChenJinrun}  tells that any sufficiently large even $n$ is the sum of a prime
and a semi prime. Taos theorem tells that any odd number larger than $1$ is the sum of at most
$5$ primes. Helfgott announced the unconditional solution of the ternary Goldbach problem
\cite{Helfgott2014}. \\ 

For now, lets look at the quaternion case, where we are not aware even of experiments about
Goldbach have been done. Recall that there are two type of primes, the {\bf Lipschitz primes}, which are of 
the form $(a,b,c,d)$ with integers $(a,b,c.d)$. The rest of the {\bf quaternion primes}. 
They are of the form $(a+1/2,b+1/2,c+1/2,d+1/2)$ with integers $(a,b,c,d)$.
To make the formulations easier, lets call the primes in set of Hurwitz integers just "primes" or 
"quaternion primes" and call the quaternion primes which are not Lipschitz primes the 
{\bf Hurwitz primes}. Hurwitz himself called primes "Primquaternion" which is the German 
expression for prime quaternion. \\

Define $Q$ as the set of Hurwitz integers for which all coordinates are positive. 

\conjecture{Every Lipschitz integer quaternion with entries $>1$ is the sum of two Hurwitz primes in $Q$.}

We believe that it will be difficult to prove if it holds up to be true. We also hope for

\question{Every Hurwitz integer quaternion with entries $>2$ is the sum of a Hurwitz and Lipschitz primes in $Q$.}

The Hurwitz integer $(3/2,3/2,3/2,3/2)$ is not the sum of a Hurwitz and Lipschitz prime because the only
decomposition would be $(1,1,1,1) + (1,1,1,1)/2$ but both are not prime. Together:

\question{Every integer quaternion with entries larger than $2$ is the sum of two quaternion primes in $Q$.}

Here are some computations showing in how many ways a Lipschitz integer can be written as a sum of two
Hurwitz primes. Since we are in a 4 dimensional lattice, we fix the first two coordinates $a,b$, 
then build the matrix $G(a,b)$ for which $G_{cd}(a,b)$ tells in how many ways one can write 
the Lipschitz prime $(a,b,c,d)$ as a sum of two Hurwitz primes:

\begin{tiny}
$$
G(1,1)=\left[ \begin{array}{ccccc}
 0 & 0 & 1 & 2 & 3 \\
 0 & 2 & 4 & 4 & 2 \\
 1 & 4 & 3 & 4 & 6 \\
 2 & 4 & 4 & 8 & 6 \\
 3 & 2 & 6 & 6 & 5 \\
\end{array} \right], 
G(1,2)=\left[ \begin{array}{ccccc}
 0 & 2 & 4 & 4 & 2 \\
 2 & 6 & 6 & 6 & 10 \\
 4 & 6 & 8 & 10 & 10 \\
 4 & 6 & 10 & 10 & 10 \\
 2 & 10 & 10 & 10 & 10 \\
\end{array} \right] $$
\end{tiny}

\begin{tiny}
$$ G(1,3)=\left[ \begin{array}{ccccc}
 1 & 4 & 3 & 4 & 6 \\
 4 & 6 & 8 & 10 & 10 \\
 3 & 8 & 13 & 12 & 11 \\
 4 & 10 & 12 & 10 & 14 \\
 6 & 10 & 11 & 14 & 20 \\
\end{array} \right],
G(1,4)=\left[ \begin{array}{ccccc}
 2 & 4 & 4 & 8 & 6 \\
 4 & 6 & 10 & 10 & 10 \\
 4 & 10 & 12 & 10 & 14 \\
 8 & 10 & 10 & 14 & 20 \\
 6 & 10 & 14 & 20 & 14 \\
\end{array} \right],
G(1,5)=\left[ \begin{array}{ccccc}
 3 & 2 & 6 & 6 & 5 \\
 2 & 10 & 10 & 10 & 10 \\
 6 & 10 & 11 & 14 & 20 \\
 6 & 10 & 14 & 20 & 14 \\
 5 & 10 & 20 & 14 & 39 \\
\end{array} \right]  \; . $$
\end{tiny}

\begin{tiny}
$$
G(2,1)=\left[
\begin{array}{ccccc}
 2 & 6 & 6 & 6 & 10 \\
 6 & 14 & 14 & 16 & 16 \\
 6 & 14 & 20 & 24 & 24 \\
 6 & 16 & 24 & 22 & 24 \\
 10 & 16 & 24 & 24 & 34 \\
\end{array} \right],
G(2,2)=\left[
\begin{array}{ccccc}
 4 & 6 & 8 & 10 & 10 \\
 6 & 14 & 20 & 24 & 24 \\
 8 & 20 & 20 & 26 & 36 \\
 10 & 24 & 26 & 34 & 38 \\
 10 & 24 & 36 & 38 & 42 \\
\end{array} \right] $$
\end{tiny}

\begin{tiny}
$$
G(2,3)=\left[
\begin{array}{ccccc}
 4 & 6 & 10 & 10 & 10 \\
 6 & 16 & 24 & 22 & 24 \\
 10 & 24 & 26 & 34 & 38 \\
 10 & 22 & 34 & 38 & 36 \\
 10 & 24 & 38 & 36 & 56 \\
\end{array} \right],
G(2,4)=\left[
\begin{array}{ccccc}
 2 & 10 & 10 & 10 & 10 \\
 10 & 16 & 24 & 24 & 34 \\
 10 & 24 & 36 & 38 & 42 \\
 10 & 24 & 38 & 36 & 56 \\
 10 & 34 & 42 & 56 & 44 \\
\end{array} \right] $$
\end{tiny}

\begin{tiny}
$$ 
G(3,3)=\left[
\begin{array}{ccccc}
 3 & 8 & 13 & 12 & 11 \\
 8 & 20 & 20 & 26 & 36 \\
 13 & 20 & 46 & 34 & 33 \\
 12 & 26 & 34 & 40 & 46 \\
 11 & 36 & 33 & 46 & 61 \\
\end{array}
\right],
G(3,4)=\left[
\begin{array}{ccccc}
 4 & 10 & 12 & 10 & 14 \\
 10 & 24 & 26 & 34 & 38 \\
 12 & 26 & 34 & 40 & 46 \\
 10 & 34 & 40 & 52 & 62 \\
 14 & 38 & 46 & 62 & 72 \\
\end{array}
\right]
G(3,5)=\left[
\begin{array}{ccccc}
 6 & 10 & 11 & 14 & 20 \\
 10 & 24 & 36 & 38 & 42 \\
 11 & 36 & 33 & 46 & 61 \\
 14 & 38 & 46 & 62 & 72 \\
 20 & 42 & 61 & 72 & 96 \\
\end{array}
\right] $$
\end{tiny}

\begin{tiny}
$$ G(4,4) = \left[
\begin{array}{ccccc}
 8 & 10 & 10 & 14 & 20 \\
 10 & 22 & 34 & 38 & 36 \\
 10 & 34 & 40 & 52 & 62 \\
 14 & 38 & 52 & 62 & 74 \\
 20 & 36 & 62 & 74 & 86 \\
\end{array}
\right],
G(4,5)=\left[
\begin{array}{ccccc}
 6 & 10 & 14 & 20 & 14 \\
 10 & 24 & 38 & 36 & 56 \\
 14 & 38 & 46 & 62 & 72 \\
 20 & 36 & 62 & 74 & 86 \\
 14 & 56 & 72 & 86 & 104 \\
\end{array}
\right]
G(5,5) = \left[
\begin{array}{ccccc}
 6 & 10 & 14 & 20 & 14 \\
 10 & 24 & 38 & 36 & 56 \\
 14 & 38 & 46 & 62 & 72 \\
 20 & 36 & 62 & 74 & 86 \\
 14 & 56 & 72 & 86 & 104 \\
\end{array}
\right]
$$ 
\end{tiny}

For example, here are all the 14 summands of the Lipschitz integer $z=(2,3,2,2)$ 
as a sum of two Hurwitz integers $z=p+q$. The rational primes $N(p),N(q)$ can vary. \\

\begin{center}
\begin{tabular}{|l|l|l|l|} \hline
2p            & 2q            & N(p) & N(q) \\ \hline
(1, 1, 1, 3)  & ( 3, 5, 3, 1) & 3  & 11 \\
(1, 1, 3, 1)  & ( 3, 5, 1, 3) & 3  & 11 \\
(1, 3, 1, 3)  & ( 3, 3, 3, 1) & 5  & 7  \\
(1, 3, 3, 1)  & ( 3, 3, 1, 3) & 5  & 7  \\
(1, 3, 3, 3)  & ( 3, 3, 1, 1) & 7  & 5  \\
(1, 5, 1, 1)  & ( 3, 1, 3, 3) & 7  & 7  \\
(1, 5, 3, 3)  & ( 3, 1, 1, 1) & 11 & 3  \\
(3, 1, 1, 1)  & ( 1, 5, 3, 3) & 3  & 11 \\
(3, 1, 3, 3)  & ( 1, 5, 1, 1) & 7  & 7  \\
(3, 3, 1, 1)  & ( 1, 3, 3, 3) & 5  & 7  \\
(3, 3, 1, 3)  & ( 1, 3, 3, 1) & 7  & 5  \\
(3, 3, 3, 1)  & ( 1, 3, 1, 3) & 7  & 5  \\
(3, 5, 1, 3)  & ( 1, 1, 3, 1) & 11 & 3  \\
(3, 5, 3, 1)  & ( 1, 1, 1, 3) & 11 & 3  \\ \hline
\end{tabular}
\end{center} 

Already for $z=(2,2,2,2)$, the smallest allowable integer for the conjecture, there are 14 summands.
They are of the form $(3,1,1,1)/2 + (1,3,3,3)/2$ (8 cases) or $(1,1,3,3)/2+(3,3,1,1)/2$  (6 cases).  \\

Let us look at the special case, when the Lipschitz integer is $z=(2,2,2,n)$. 
The two primes summing up to it must then have the form of one of the following cases modulo permutations: 
$p = (1,1,1,x)/2,  q=(3,3,3,n-x)/2$,  or then 
$p = (1,1,3,x)/2,  q=(3,3,1,2n-x)/2$ for an unknown odd integer $x$. 
In the first case, we need simultaneously to have $(3+x^2)/4$  and $(27+(2n-x)^2)/4$ to be rational primes. 
In the second case, we need simultaneously to have $(11+x^2)/4$ and $(19+(2n-x)^2)/4$ to be prime. Since $x$
needs to be odd for $p$ to be a Hurwitz prime, we can 
write $x=2k+1$. Now $(3+x^2)/4=1+k+k^2$ and $(27+(2n-x)^2)/4 = 7-(n-k)+(n-k)^2$.
In the second case $3+k+k^2$ and $5-(n-k)+(2n-k)^2$. We see:
If the Hurwitz Goldbach conjecture holds, then for any $n$, there exists $k<n$ for which 
both $3+k+k^2$ and $7+k+k^2-n-2kn+n^2$ are prime or for which both $1+k+k^2$ and 
$5+k+k^2-n-2kn+n^2$ are prime. Now, if Hurwitz-Goldbach is true and if there existed only finitely 
many primes of the form $3+k+k^2$ and $1+k+k^2$, then for all $m=n-k$ large enough, $7-m+m^2$ and $5-m+m^2$ 
would always have to be prime. This is obviously not true if $m$ is a multiple of $7$ or $5$. We see that 
Goldbach implies a special case of the Bunyakovsky conjecture, which like Landau's problem is 
likely not so easy to prove:

\resultremark{  
If the Hurwitz Goldbach conjecture is true, then one of the sequences 
$1+k+k^2$ or $3+k+k^2$ contains infinitely many primes.}

Let us quickly verify the {\bf Bunyakovsky conditions} which must be checked for the {\bf Bunyakovsky conjecture}. 
First, $\phi_3(k) = 1+k+k^2$ is already the cyclotomic polynomial and also $k^2+k+3$ satisfies the conditions
of the Bunyakovsky conjecture: the maximal coefficient of the polynomial $f$ is positive, the coefficients have
no common divisor and there is a pair of integers $n,m$ such that $f(n),f(m)$ have no common divisor.
Hardy-Littlewood analogue density results have been formulated in \cite{BatemanHorn}.
The set of $k$ for which $k^2+k+1$ is prime is the sequence $A002384$ in \cite{A002384}.

\begin{figure}[!htpb]
\scalebox{0.4}{\includegraphics{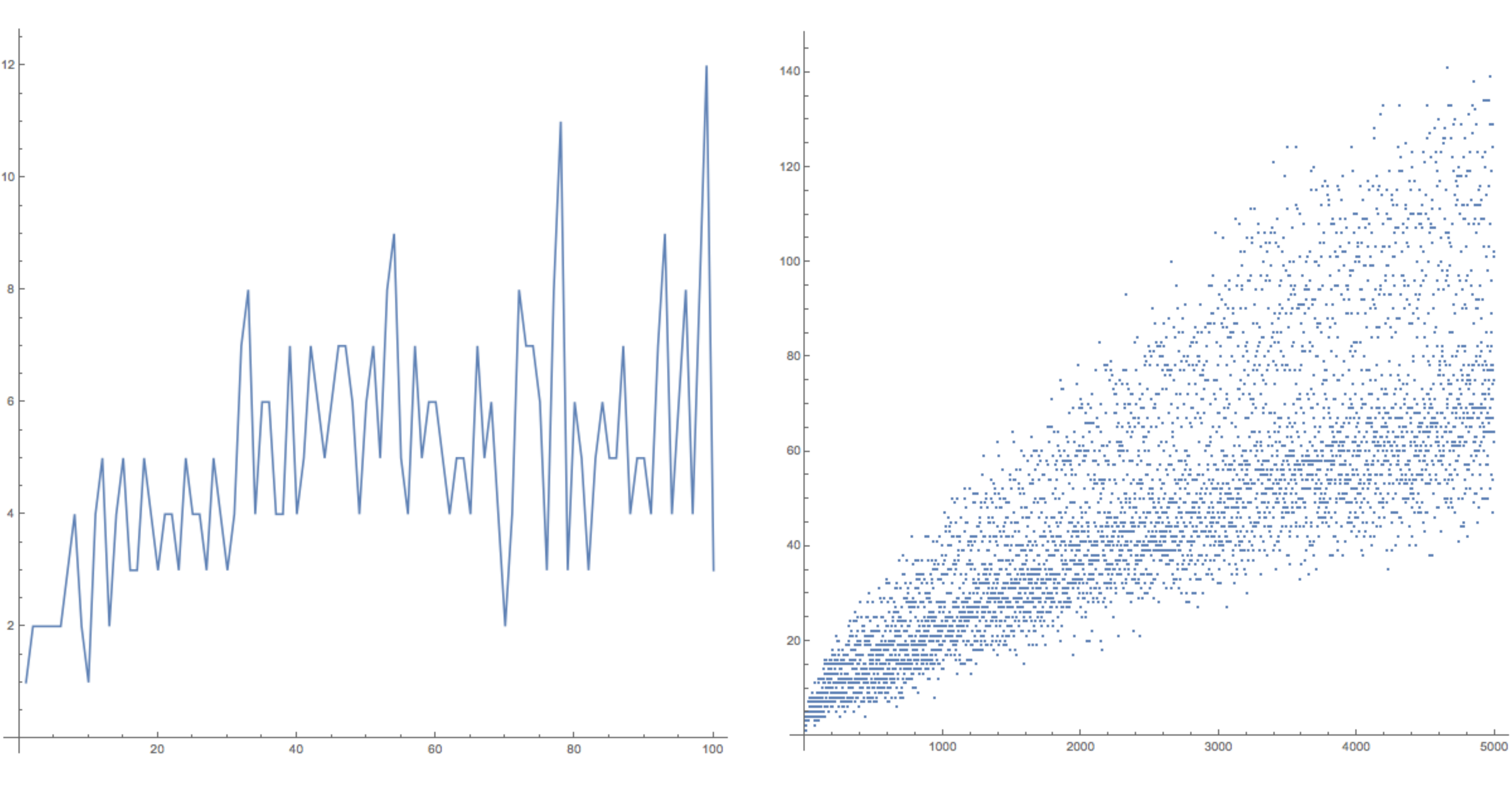}}
\caption{
The {\bf Hurwitz boundary comet} shows the function $f(n)$ which gives
the number of solutions $p+q=(2,2,2,n)$ with ordered Hurwitz primes $p,q$.  
\label{hurwitzcomet}
}
\end{figure}

Our experiments confirm that the sequences produce primes with the frequency given by 
analogous Hardy-Littlewood density constants. 

\section{Goldbach for Octonions}

Besides $\mathbb{R},\mathbb{C}$ and $\mathbb{H}$, there is a fourth division algebra 
$\mathbb{O}$. It is the space of {\bf Cayley numbers}. Because it is eight
dimensional, it has also been called the space of {\bf Octonions}. It is this name which 
has stuck. The members of $\mathbb{O}$ can either be written as a linear combination of a basis 
\begin{wrapfigure}{l}{4.1cm} \begin{center}
\includegraphics[width=4cm]{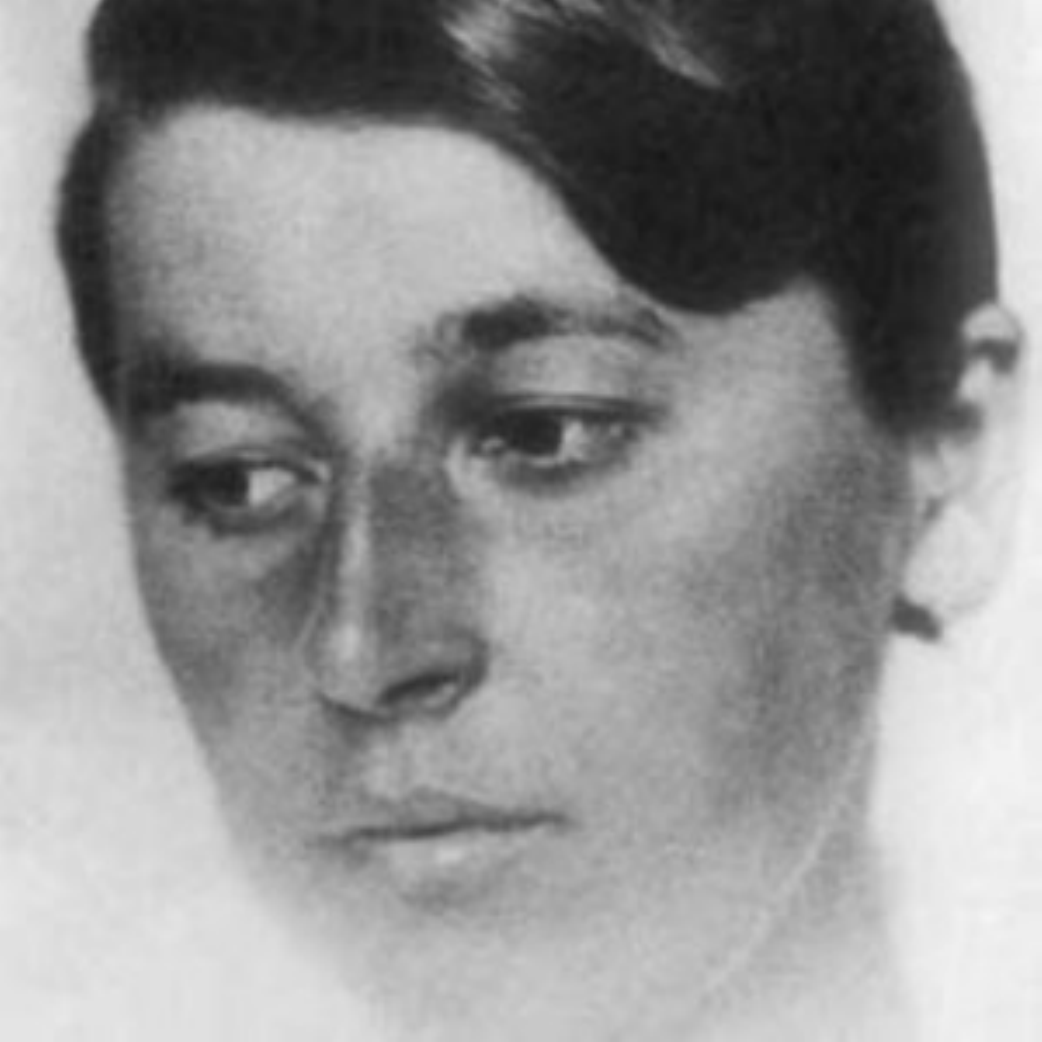}
\end{center} \end{wrapfigure}
$1,i,j,k,l,m,n,o$ or then, according to a suggestion of Cayley-Dixson, as pairs $(z,w)$ of quaternions,
defining $(z,w) \cdot (u,v) = (zu-v^*w,vz+wu^*)$. The algebra is no more associative. \\

In order to do number theory, one has to specify what the {\bf integers} are in $\mathbb{O}$. 
It turns out that there are several classes of integers and even several maximal ones.
First of all there are the {\bf Gravesian integers} $(a,b,c,d,e,f,g)$ which play the role of the
Lipschitz primes in the Hurwitz case. Then there are 
the {\bf Kleinian integers} $(a,b,c,d,e,f,g) + (1,1,1,1,1,1,1,1)/2$ which play the role of the 
Hurwitz primes in $\mathbb{H}$. But then there are more, the
{\bf Kirmse integers} which includes elements for which $4$ of the entries are half integers.
There are $7$ {\bf maximal orders} which Kirmse classified \cite{Kirmse25}. (He counted $8$ which was
later corrected by Coxeter \cite{Coxeter46}). They turn out all to be equivalent 
and produce the ultimate class of integers. They are often called {\bf Octonion} or {\bf Cayley integers}
or then more catchy, the {\bf Octavian integers} \cite{Kirmse25}.
In the mathematical genealogy database, there is no advisor listed, but since he acknowledges Herglotz 
in \cite{Kirmse}, this is our most likely guess. 
Johannes Kirmse wrote his pioneering Octonion paper in 1925.  We could not get hold of
the 1925 paper \cite{Kirmse25} but it is mentioned and discussed in \cite{Coxeter46} who's critics
appears maybe a bit too harsh today when considering that Kirmse explored new ground where 
nobody has been before. \\

The condition $N(z w) = N(z) N(w)$ which assures that
the algebra is a normed {\bf division algebra}, is also called the {\bf Degen eight square identity}. 
Ferdinand Degen discovered it first around 1818 and John Thomas Graves in 1843.
Also Arthur Cayley in 1845 rediscovered the identity.  \\

The primes in the Octonions are either called {\bf Gravesian primes} containing the Lipschitz primes, 
or {\bf Kleinian primes} containing the Hurwitz primes or then the {\bf Kirmse primes}.
After choosing a maximal order, lets just call them {\bf Octavian primes}.
We have followed mostly the nomenclature of \cite{ConwaySmith} but split the three type of 
primes up so that the Octavian primes are made up of three distinct type of primes, similarly 
as the Quaternion primes were made up of two type of primes, the Lipschitz and Hurwitz primes. 
Already the Kirmse integers or the Octavian integers form a lattice called the {\bf $E_8$ lattice}.
It is important as it produces the densest {\bf sphere packing} in $\mathbb{R}^8$. \\

The units form what is now called a {\bf Moufang loop} named after German mathematician
{\bf Ruth Moufang} (1905-1977) \cite{Moufang33}, who was a student of Max Dehn (1878-1952)
and is considered one of the first PhD mathematicians working also in the industry.  Dehn
is also known for the Dehn-Sommerville relations. There is a smaller loop of 16 unit
octonions containing Gravesian integers like $(\pm 1,0,0,0,0,0,0,0,0)$.
The units placed in the unit sphere of $R^8$ form the {\bf Gosset polytope} $4_{21}$ which
was discovered by Thorold Gosset (1869-1962) who as a lawyer without much clients
amused himself as an amateur mathematician. The vertices of $4_{21}$ are the
roots of the exceptional Lie algebra $E_8$ belonging to the 248 dimensional Lie group
$E_8$. As the dimension of the maximal torus is $8$, the root system lives in $R^8$.
One can write points on the sphere of radius $2$ taking vertices $(a,b,0,0,0,0,0,0)$
with $a,b \in \{-1,1\}$ or $(a,b,c,d,e,f,g,h)$ with with entries in $\{-1/2,1/2\}$
summing up to an even number. This lattice $E_8$ has just recently been verified by
Maryna Viazovska to be the densest sphere packing in $R^8$ \cite{Viazovska}. \\

Let us fix now one of the 7 classes of Octonian integers. For the Goldbach statement, we are going 
to look up, it does not matter which. 
Proof. We use the Cayley-Dickson notation which writes an Octonion as a pair of 
quaternions. Assume we have integers of the form $(w,z),(w/2,z/2),(w/2,z),(w,z/2)$ where
$w,z$ are Lipschitz quaternions and $w/2,z/2$ are Hurwitz quaternions. 
Their units of the octonions do not form a multiplicative {\bf group} any more, as multiplication is not associative. 
But it is a {\bf loop} an algebraic structure more primitive than a group in which one does
not insist on associativity. 

When looking for {\bf Octonion Goldbach conjectures} one has to check any of them. 
The most obvious case fails for $K=2$:

One certainly has to distinguish cases. In the simplest case, within Gravesian integers
if $(z,w)=(z_1,z_2,z_3,z_4$,$w_1,w_2,w_3,w_4)$, we want to find 
$(x,y)= (x_1,x_2,x_3,x_4$,$y_1,y_2,y_3,y_4)$
such that both $p=\sum_i x_i^2 +y_i^2$ and $q=\sum_i (z_i-x_i)^2 + (w_i-y_i)^2$ are prime.
The integer $(2,2,2,2,2,2,2,2)$ can not be written as a sum of two Gravesian
primes. Also Kleinian primes do not work:
there are Gravesian integers with entries $\geq 2$ which are not the sum of two Kleinian primes:
Assume $(2,2,2,2,2,2,2n)$ is the sum of two Kirmse primes. Here are the four cases: 
$$(1,1,1,1,1,1,1,2k+1)/2, (3,3,3,3,3,3,3,2n-2k-1)/2 $$
$$(1,1,1,1,1,1,3,2k+1)/2, (3,3,3,3,3,3,1,2n-2k-1)/2 $$
$$(1,1,1,1,1,3,3,2k+1)/2, (3,3,3,3,3,1,1,2n-2k-1)/2 $$
$$(1,1,1,1,3,3,3,2k+1)/2, (3,3,3,3,1,1,1,2n-2k-1)/2 \; . $$ 
Define $m$ by $2m+1=2n-2k-1$. Now this either one of the four pairs are both prime:
$(7+(2k+1)^2)/4=2+k+k^2$ and $(63+(2m+1)^2)/4 = 16+m+m^2$ \\
$(6+9+(2k+1)^2)/4=4+k+k^2$ and $(1+54+(2m+1)^2)/4 = 14+m+m^2$ \\
$(5+18+(2k+1)^2)/4=6+k+k^2$ and $(2+45+(2m+1)^2)/4 = 12+m+m^2$ \\
$(4+27+(2k+1)^2)/4=8+k+k^2$ and $(3+36+(2m+1)^2)/4 = 10+m+m^2$ \\
with $2m+1=2n-2k-1$. But this is not true as in any case all these numbers are divisible by 2.
The failure is related to the fact that the Kirmse integers are not yet a maximal order. 
But also a brute force search over Kirmse does not work for $(2,2,2,2,2,2,2,2,2)$. \\

It appears that it is not obvious how to come up with a conjecture which
is both convincing and also justifiably difficult. Brute force searches are difficult as the volume
of a box of size $r$ grows like $r^8$. Here is a first attempt of get to a conjecture.
We formulate as a question since our experiments did not get far yet, nor do we have an
idea how difficult the statement could be. Anyway, lets denote by $Q$ again the set of all
integer octonions, for which all coordinates are positive.

\question{
There exists $K$ such that every Octavian integer $Z$ with 
coordinates $\geq K$ is a sum of two Octavian primes $P,Q$.}

{\bf Remarks.} \\
{\bf 1)} In physics, one is interested in these structures as they have used them to build a
$10=8+1+1$-dimensional space-time, where time and string parametrization add two more dimensions.
Some theoreticians take the structures of division algebras serious in describing
matter \cite{Gursey,DixonDivisionalgebras}. Some history on their use in quantum mechanics is
\cite{Baez2012,Baez2002}. 

\section{Landau and Bunyakovsky problem}

One of the four problems presented by {\bf Edmund Landau} at the 1912 International
congress of mathematicians claims that there are infinitely many primes
of the form $n^2+1$. The problem is also known as the {\bf fifth Hardy-Littlewood 
conjecture}, but the problem has been raised already by Euler \cite{Pintz}.
It is now part of a more general {\bf Bunyakovsky conjecture} or of the {\bf Schinzel
hypothesis $H$} or more generally the {\bf Bateman-Horn conjecture} \cite{BatemanHorn}
Early laboratory experiments have been done by Euler
who computed all primes of the form $n^2+1$ for $n$ up to $1500$ and saw 
that there are many primes in this list. \\

Despite the extreme simplicity of question: {\bf "are there infinitely
many primes $p$ for which $p-1$ is a square?"}, it is considered a 
\begin{wrapfigure}{l}{4.1cm} \begin{center}
\includegraphics[width=4cm]{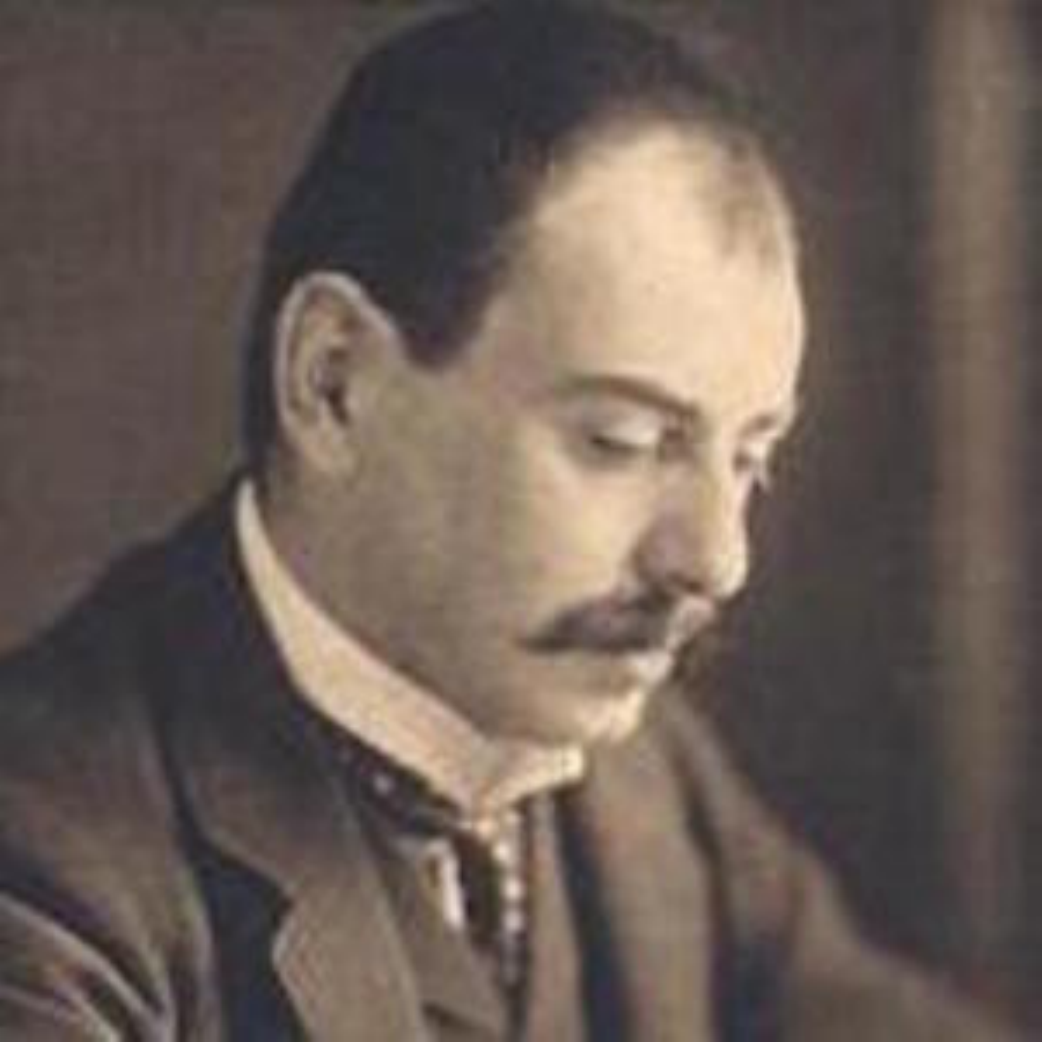}
\end{center} \end{wrapfigure}
``hard" problem. One indication is that despite of having been in the spotlight for 
more than over a hundred years, there is no solution in sight. The tools of analytic number theory 
look helpless. Just to illustrate some attacks, in \cite{FouvryIwaniec} it is shown that there are 
infinitely many Gaussian primes $a+ib$, for which $a$ is a given prime. While this existence result does not 
prove the Landau problem, it gives hope even so the result giving existence is far from 
proving that there are infinitely many. \\

A much more general conjecture was already formulated by Victor Bunyakovsky in 1857. Bunyakovsky was
a student of Cauchy and a gifted educator. His conjecture is that if $f(x)$ is an irreducible
polynomial for which some obvious non-triviality conditions are satisfied, then the range of $f$
has infinitely many primes. \\

\begin{wrapfigure}{l}{4.1cm} \begin{center}
\includegraphics[width=4cm]{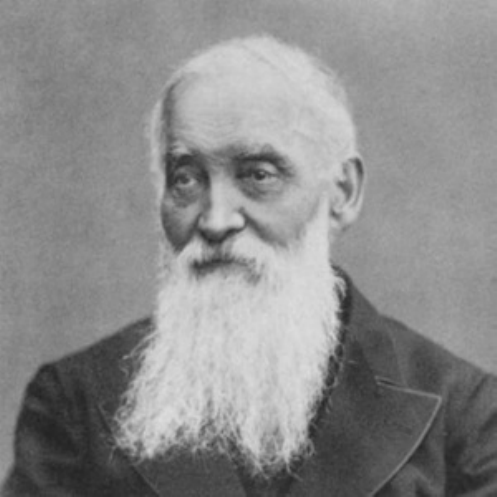}
\end{center} \end{wrapfigure}

It is therefore amazing that Hardy and Littlewood quantized the Landau problem in their landmark paper 
of 1923. They predicted a explicit limiting density $C$ for the number of such primes and more
generally predicted a precise density ration for primes of the form 
$k+ia$ and $k+ib$ as $C_a/C_b$ with
the product $C_a = \prod_{p} (1-\left( \frac{-a}{p} \right)/(p-1))$ over all odd primes.
These are the {\bf Hardy-Littlewood density conjectures}.  \\

The constant $C=C_1$ agrees with the ratio between real Gaussian primes and 
Gaussian primes of the form $n+i$ in the first row. We focussed our own experiments primarily 
on this number. 

The first {\bf Hardy-Littlewood ratio} $C_1$ can be rewritten as
$$ C = \prod_{p \in P_1} [1-\frac{1}{p-1}] \prod_{p \in P_3} [1+\frac{1}{p-1}] = 1.37279 ... \; , $$
where $P_k$ is the set of rational primes congruent to $k$ modulo $4$. 
Intuitively, one can understand this probabilistically. 
The probability to be in a multiplicative remainder
class of $p \in P_3$ is $p/(p-1)$ and the probability to be not in a specific remainder class modulo 
$p \in P_1$ is $1-1/(p-1)$.  \\

More general {\bf Hardy-Littlewood density estimates} predict a density relation between primes of the
form $f(x)$ or $g(x)$ given two irreducible polynomials $f,g$ of the same degree to be given as $C_f/C_g$ with
$C_f = \prod_p \frac{(1-\omega_f(p)/p)}{(1-1/p)}$. Here, $\omega_f(p)$ is the number of solutions
$f(x)=0$ modulo $p$. In order for the estimate to work, one needs that the polynomials a 
positive leading coefficients and $\omega_f(p) \neq p$ for all $p$. This 
especially implies $\omega_f(2)=0$ and even so not irreducible, one can
include $f(x)=x^n$, where $C_f=1$. (See \cite{BatemanHorn,KConrad2003}, where one can even look at the product of 
irreducible cases, combine several conjectures of Hardy and Littlewood). 
The constants $C_k$ are then shorts cuts for $C_{n^2+k}$. 
In the prototype case $f(x)=x^2+1$ one has $\omega_f(p)=1$ for $p \in P_3$ and
$\omega_f(p)=2$ for $p \in P_1$ by quadratic reciprocity so that
$(1-\omega_f(p)/p)/(1-1/p)) = 1-\frac{1}{p-1}$ for $p \in P_1$ and $1+\frac{1}{p-1})$ for $p \in P_3$. \\

The intuition which led to the conjectures of Hardy and Littlewood are of a probabilistic nature:
the key assumption is that solving equations like $a^2=-1$ modulo $p$ or then solving modulo $q$ is pretty
much independent, if $p,q$ are different odd primes. Of course this is not justified. But it is
part of the magic of primes that it seems to work.  \\

Lets try to explain the intuition behind the constant: 
every time a new prime is added, the size of the space changes by $(p-1)/p=(1-1/p)$ because we can only 
take numbers which are not multiples of $p$. 
The product of these size changes is $1/\zeta(1)=0$ and reflects the infinitude of primes. But if we look at the {\bf ratio} 
of solution sets for two polynomials, we don't have to do this {\bf re-normalization} because it happens on both sides.
Now, when looking at solutions of the form $x^2+1$, then whenever $-1$ is a quadratic residue, the
probability decreases by $(1-1/(p-1)) (1-1/p) =(p-2)/p$ but if $-1$ is not a quadratic residue, then 
the probability increases stays the same. Including back the volume change gives $(p-2)/p /(1-1/p) = 1-1/(p-1)$
in the residue case and $1/(1-1/p) = p/(p-1) = 1+1/(p-1)$ in the non-residue case. This explains the formula
for $C$.  \cite{BatemanHorn} explain this skillfully.  Using this frame work,
many of the formulas of \cite{HardyLittlewood1923} make sense, like density formulas for the estimated number of
prime twins. 

\resultremark{
The Goldbach conjecture for Gaussian primes implies the existence of infinitely many Gaussian primes 
of the form $a+i$ or $a+2i$ and especially implies Landau's first problem. 
}

Landau's problem asks whether infinitely many primes exist on the first row of the complex plane. 
One can also ask about existence of primes on rows. This appears much easier but is also open.
We wanted to call it the "highway crossing frog problem"
but the popular US comedian Will Ferrell once said this is a "lame name" so that we call it the 
{\bf frogger problem}. It could be a problem which will be solved first in the arena of Gaussian integers: \\

Assume the rows in the complex plane are the lanes of the "highway" and the primes are the
"gaps between cars". A frog can walk freely horizontally between the highway lanes and hop through gaps.
The question is whether it can hop arbitrarily far in the vertical direction.  \\

Its obviously possible if and only if there exists at least one prime on each highway lane:

\conjecture{
For any integer $a>0$ there exists a rational prime of the form $x^2+a^2$ with integer $x$. 
}

The frogger problem would follow from Landau type
conjectures the existence of infinitely many primes of the form $x^2+a^2$.
Work of Hecke shows that the frog can jump through lines for which $a$ is prime. In general, it appears open. \\

The {\bf Hurwitz frogger problem} asks whether for every integer $a$,
there are integer vectors $(x,y,z)$ for which $a^2+x^2+y^2+z^2$ is a rational prime and
whether for every half integer $a+1/2$, there exists a half integer
$(x+1/2,y+1/2,z+1/2)$ for which $(x+1/2)^2+(y+1/2)^2+(z+1/2)^2+(a+1/2)^2$ is a rational prime.
The first part can be solved: 

\resultremark{
For every integer $a$ there are integer triples $(x,y,z)$ such that $a^2+x^2+y^2+z^2$ is a rational prime.
}

Proof: Fix $a$ and an arithmetic progression $p=u+a^2+k v$ such that $u+a^2,v$ are coprime.
By the {\bf Dirichlet's theorem on arithmetic progressions},
there are infinitely primes $p$ like that. Now $a^2+x^2+y^2+z^2=p=u+a^2+kv$
is $x^2+y^2+z^2 = u+k v$. The {\bf Legendre three square theorem} tells that there is a solution if 
$u+k v$ is not of the form $4^r(8s+7)$. Now just chose $u$ odd with remainder different from $-1$
modulo $8$ and $v$ a multiple of $8$. \\

We don't know yet about the second part which is a Bunyakovsky type problem but for a function 
of several variables: 

\conjecture{
For any integer $a$ there are $(x,y,z)$ such that $a^2+x^2+y^2+z^2+x+y+z+a+1$ is a rational prime. 
}

We have measured numerically that the number of solutions of this problem 
grows like $C_a (n \log(n))^2$ with Hardy-Littlewood type constants $C_a$.

\conjecture{
For any hyper-plane $z=a$ in the space of Hurwitz integers, the number of Hurwitz primes in a ball of radius $r$
grows like $C_a (n \log(n))^2$.
}

Hardy-Littlewood type statements about the asymptotic density of the number of solutions of such 
Diophantine problems appear more approachable since we are in higher dimensions, where
some Landau type problems are already answered. 

{\bf Remarks.}
{\bf 1)} Randomness considerations also add intuition to the believe that factoring large composite
integers is hard because the holy grail of integer factorization is since Fermat the ability to solve
quadratic equations modulo the product $n=pq$. Say, if $x,y$ were both solutions to $x^2+1=0$ modulo $pq$, then
$x^2=y^2$ modulo $pq$ and ${\rm gcd}(x-y,n)$ is either $p$ or $q$. In some sense, the quadratic map on finite
fields shows similar randomness features like the quadratic maps on the field of complex numbers, where one has
the Julia-Fatou story. So, once accepting this intuition that a quadratic maps shuffles the residue system pretty
well, rendering them independent, the formulas of Hardy and Littlewood become transparent. 
The connection with chaotic maps is not an accident, there are basic algorithms like the 
Pollard rho method which use the {\bf iteration of a quadratic} map to get to factors. \\
{\bf 2)} Probability argument using independent can work, Here is an example where one has a product space
what is the probability that two numbers have no common denominator?
The probability to have a common denominator $p$ is then $1-1/p^2$ in the product space $[1,..,n]^2$.
The product of all these probabilities gives $\prod_p (1-1/p^2) = 1/\zeta(2) = 6/\pi^2$.

\section{The Hardy-Littlewood constants}   

Hardy and Littlewood \cite{HardyLittlewood1923} conjectured that
the density of Gaussian primes of the form $a+i$ divided by the density
of Gaussian primes of the form $a+0 i$ approaches a constant.
This statement has attracted the attention of early pioneers in computer experiments like \cite{Shanks60} who
factored numbers of the type $n^2+1$ and was therefore was interested in the density of cases, where
$n^2+1$ is prime. If $\pi_0(n)$ is the number of Gaussian real Gaussian primes in $[0,n]$ and
$\pi_1(n)$ the number of Gaussian primes $a+i$ with $a \in [0,n]$, then $\pi_0 \sim (1/2) {\rm Li}(x)$
and Hardy-Littlewood predicted $\pi_1(x) \sim (C/2) {\rm Li}(x)$.
It is interesting to see the early pioneers go through relatively heavy mathematical gymnastics to
compute the constant $C$ efficiently. Daniel Shanks (1917-1996) \cite{Shanks60} reports in 1959
that A.E. Western rewrote the constant as
$$ C = \frac{3}{4} \frac{\zeta(6)}{\beta(2) \zeta(3)} 
       \prod_{p \in P_1} (1+\frac{2}{p^3-1}) (1-\frac{2}{p(p-1)^2}) \; , $$
where $G$ is the {\bf Catalan constant}. This formula 
gives 5 decimal places already when summing over three primes $p=5,13,17$. 
A bit earlier, in 1922, even before the Hardy-Littlewood article appeared (and probably while helping
to work on the numerical verifications with Hardy and Littlewood), 
A.E. Western \cite{Western22} computed the constant $C$ using further sophisticated identities
involving various zeta values so that one can use two primes $5,13$ only to get $C$
to 5 decimal places! Western was a giant in computation \cite{WesternMiller}.
What a culture had been developed only about a seemingly tangential constant! \\

But one has to remember that these mathematicians had no access to computers.
The Zuse Z3 was built only in 1941, machines like Colossus and Mark I appeared only 
in 1944 but all of them were slow: Mark I for example needed 6 seconds to multiply 
two numbers \cite{NahinNumberCrunching}. By the way, Shanks in 1960 used an 
IBM 704 with a 32K high-speed memory. In contrary, we have today access to machines 
which give each user 500 GBybes of RAM (I have used such a machine, Odyssee, for some computations in this
paper). Shanks needed 10 minutes to factor all $n^2+1$ from $n=1$ to $n=180000$. Today, the command
$"Timing[Table[FactorInteger[n*n+1],{n,180000}]][[1]]"$ on a tiny laptop reports it done in 4 seconds.
When looking at the Hardy-Littlewood paper, the expression
$C \prod_{p>2} (1- \frac{1}{p-1} \left( \frac{-1}{p} \right) )$ involves the
{\bf Legendre symbol} $\left( \frac{-1}{p} \right)$. In this special case,
it is $\left( \frac{-1}{p} \right) = (-1)^{(p-1)/2}$ which is $1$ for primes of the form $4k+1$
and $-1$ for the others. Shanks article of 1960 reveals through how much pain
Mathematicians have gone to compute things effectively before without computers. \\

\begin{figure}[!htpb]
\scalebox{0.3}{\includegraphics{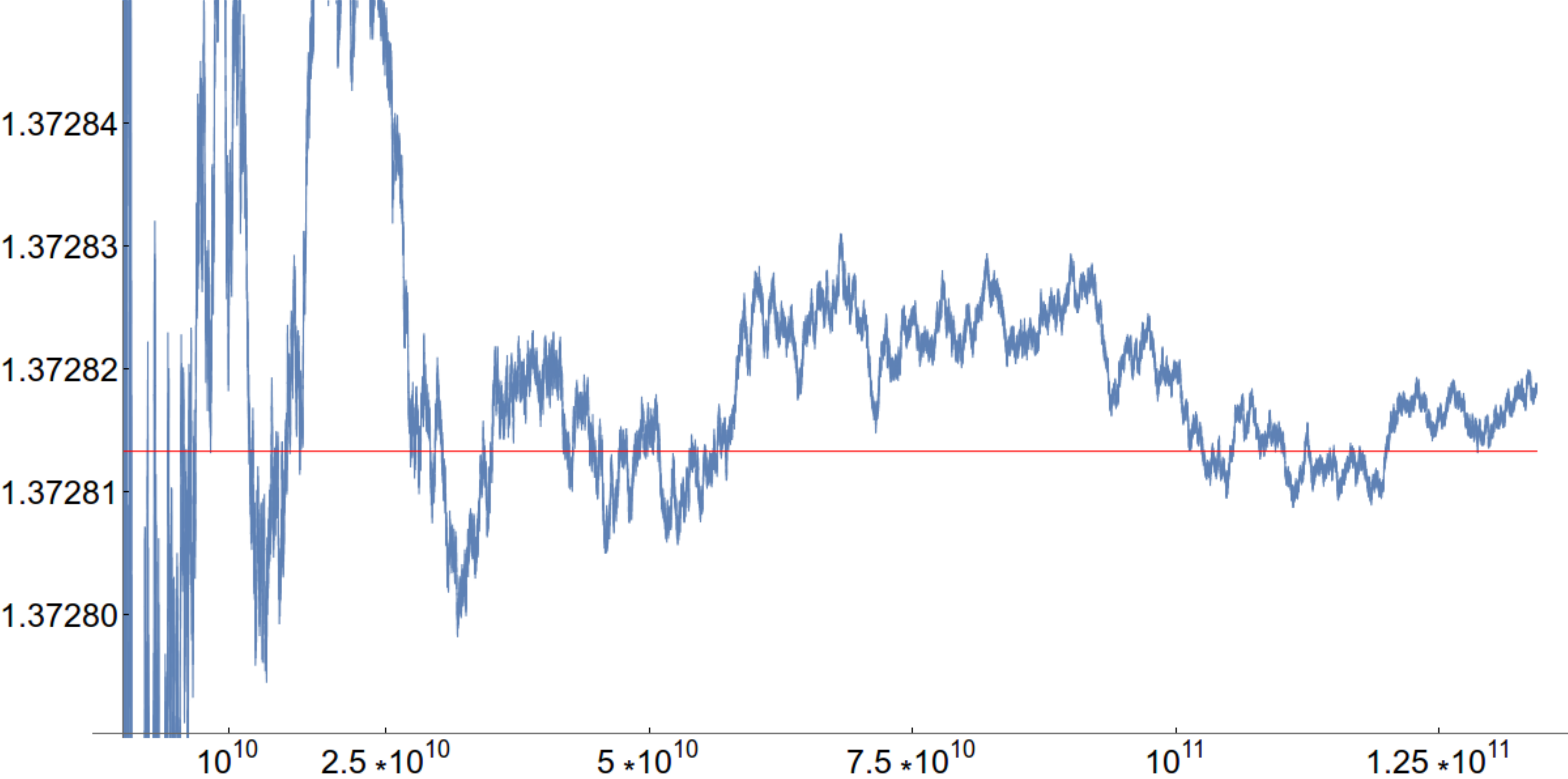}}
\caption{
The convergence to the Hardy-Littlewood constant $C$ giving the
ratio of Gaussian prime density on ${\rm Im}(z)=1$ and ${\rm Im}(z)=0$.
Shanks \cite{Shanks60} gave $C=1.37281346$. This is conjecture E in \cite{HardyLittlewood1923}.
The constant $C$ is almost prophetic as by an open Landau's problem (which currently
appears theoretically beyond reach), one does not even know whether $C$ is positive!
\label{hardylittlewood}
}
\end{figure}

\begin{figure}[!htpb]
\scalebox{0.3}{\includegraphics{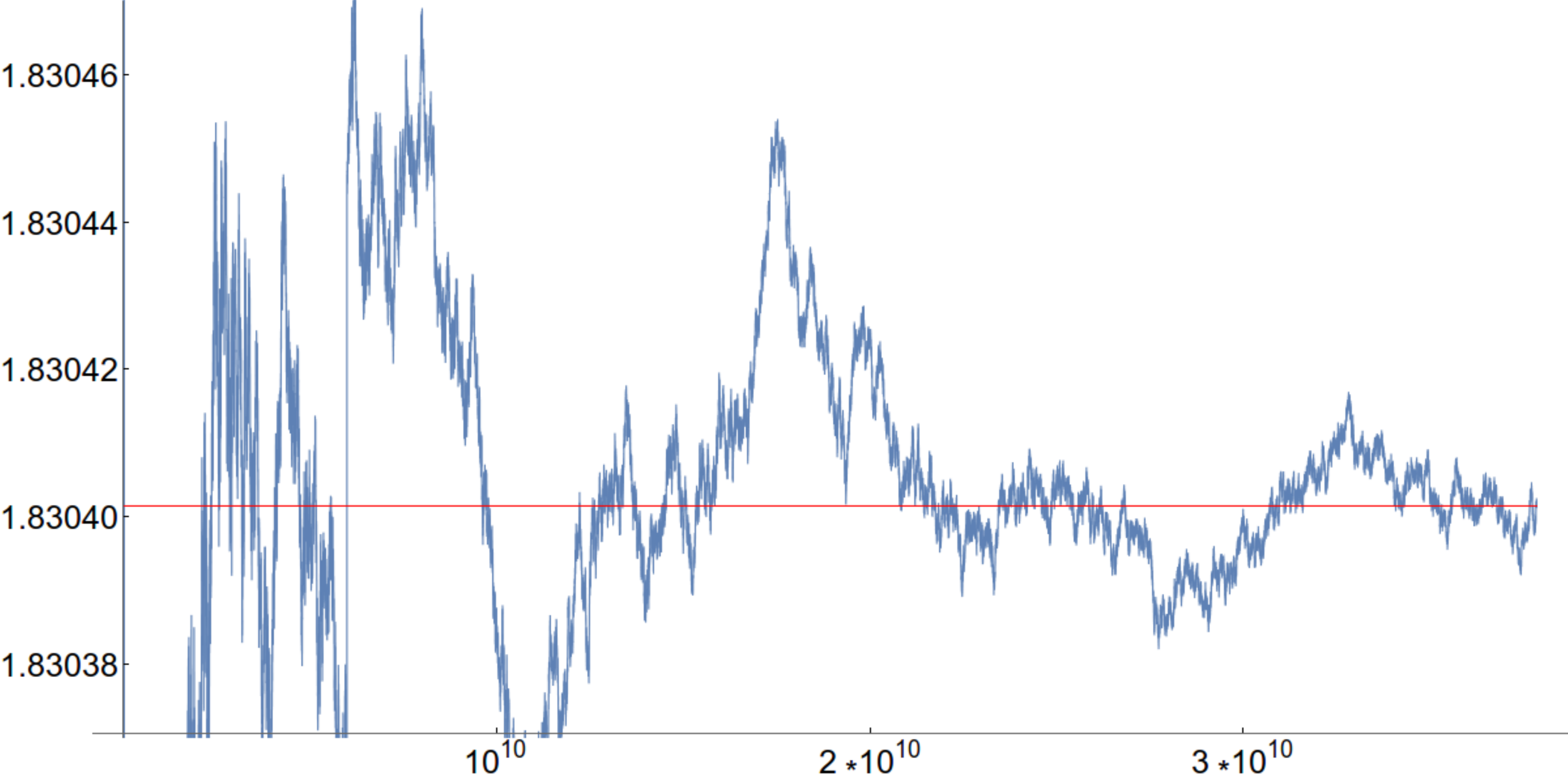}}
\caption{
The convergence to the Hardy-Littlewood constant $C_5$ giving the
ratio of Gaussian prime density on ${\rm Im}(z)=5$ and ${\rm Im}(z)=0$.
In the 5th row of the Gaussian plane, there is a high density of primes. 
}
\end{figure}

Hardy and Littlewood were not without assistance:
On page 62 of their article, they mention that "some of their conjectures have been tested
numerically by Mrs. Streatfield, Dr. A.E. Western and Mr. O. Western". Much has been written about the
influence of computers in mathematical research \cite{Williams82}, the story of the constant $C$ illustrates
that already early in the 20th century, when humans were doing the computations by hand,
the experimental part has been important. \\

We report on some experiments on our own for computing the Hardy-Littlewood ratio $C$.
Wunderlich \cite{Wunderlich} in the 70ies had computed up to $n=14'000'000$. 
Our own runs go up to $140'000'000'000$ which is 10'000 times more. We were surprised to see that 
while increasing $n$ by several orders of magnitudes, we do not get closer to the actually 
predicted constant. Its not that we have 
\begin{wrapfigure}{l}{4.1cm} \begin{center}
\includegraphics[width=4cm]{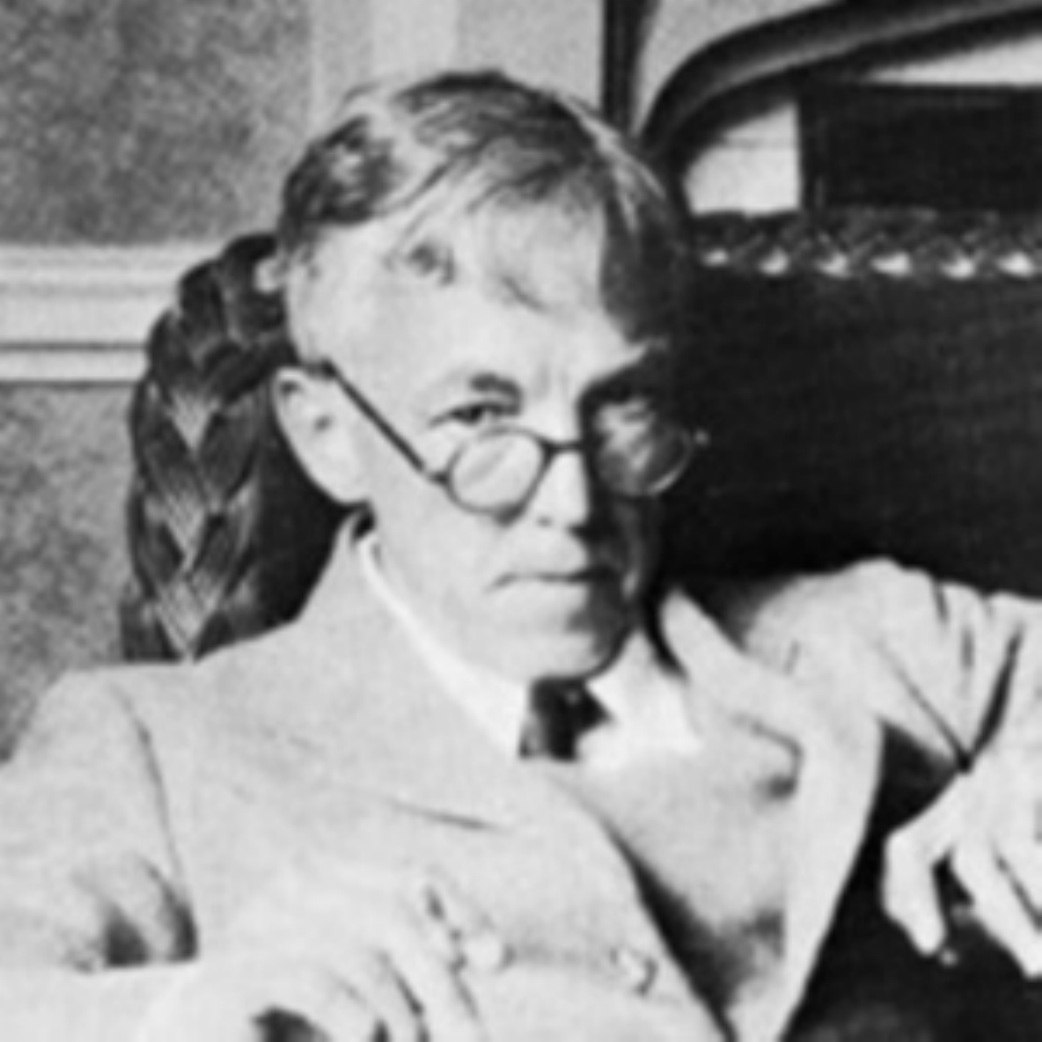}
\end{center} \end{wrapfigure}
to compute more accurately the constant $C$ which was obtained by Hardy and Littlewood 
from sieving considerations, but that the small fluctuations actually do only die out 
{\bf very slowly}, the reason being that we look at the {\bf absolute error} and not the 
{\bf relative error}. As is evident from the papers
of Shanks is that the computing assistants of Hardy and Littlewood got the computations refined
so much to values of zeta functions that they had only to consider $2$ primes in their version of
the product to compute the constant and already got to an accuracy of 5 digits. When checking how close we are
to the constant for $n=2^{32}$, it does not look much better. Note that the cryptologists Dan Shanks (1917-1996) 
and Marvin Wunderlich (1937-2013) primarily focused his computations on factorization and did 
measurements at a time, when the first magnetic card programmable computer HP-65 has 
hit the streets, and when Apple I did not even exist as a concept, and programming 
\begin{wrapfigure}{l}{4.1cm} \begin{center}
\includegraphics[width=4cm]{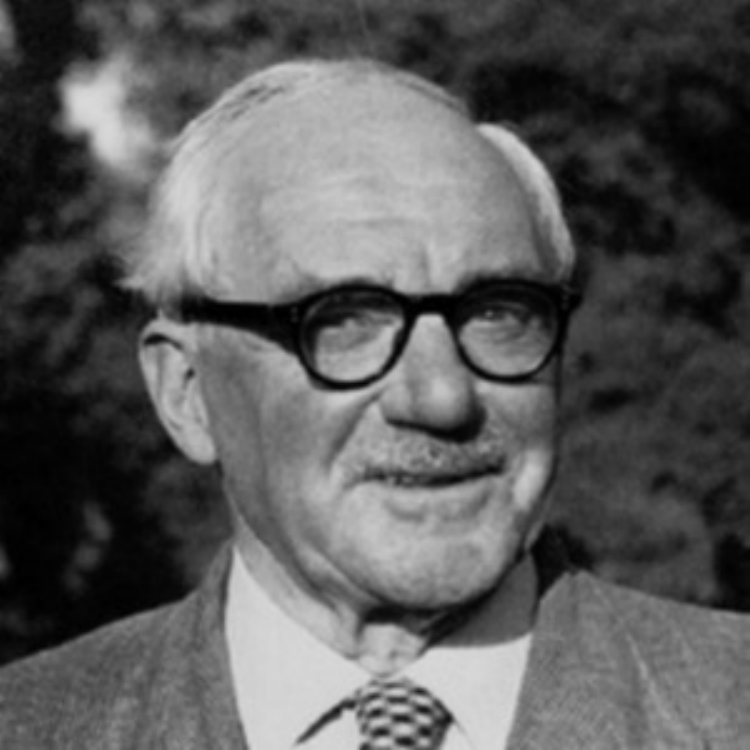}
\end{center} \end{wrapfigure}
languages like Pascal (Niklaus Wirth 1970) and C (Dennis Richie, 1972) had just started to take off. 
Hardy and Littlewood themselves acknowledge the assisted by several human computational collaborators.
Anyway, it is of course still possible that the Hardy-Littlewood claim was too strong and that the
density ratios between different rows of the Gaussian integers remain fluctuating on a small order.
Only the future will show whether Hardy-Littlewood were right.  \\

On $\{1,\dots, n=2^{31}=2'147'483'648\}$, there are $\pi_0(n)=52549599$ Gaussian primes
and on $\{1+i,\dots,1+2^{31}\}$, there are $\pi_1(n)=72139395$. The Hardy-Littlewood
ratio up to this $n$ is $\pi_1(n)/\pi_0(n)=1.37279$.

\begin{figure}[!htpb]
\scalebox{0.8}{\includegraphics{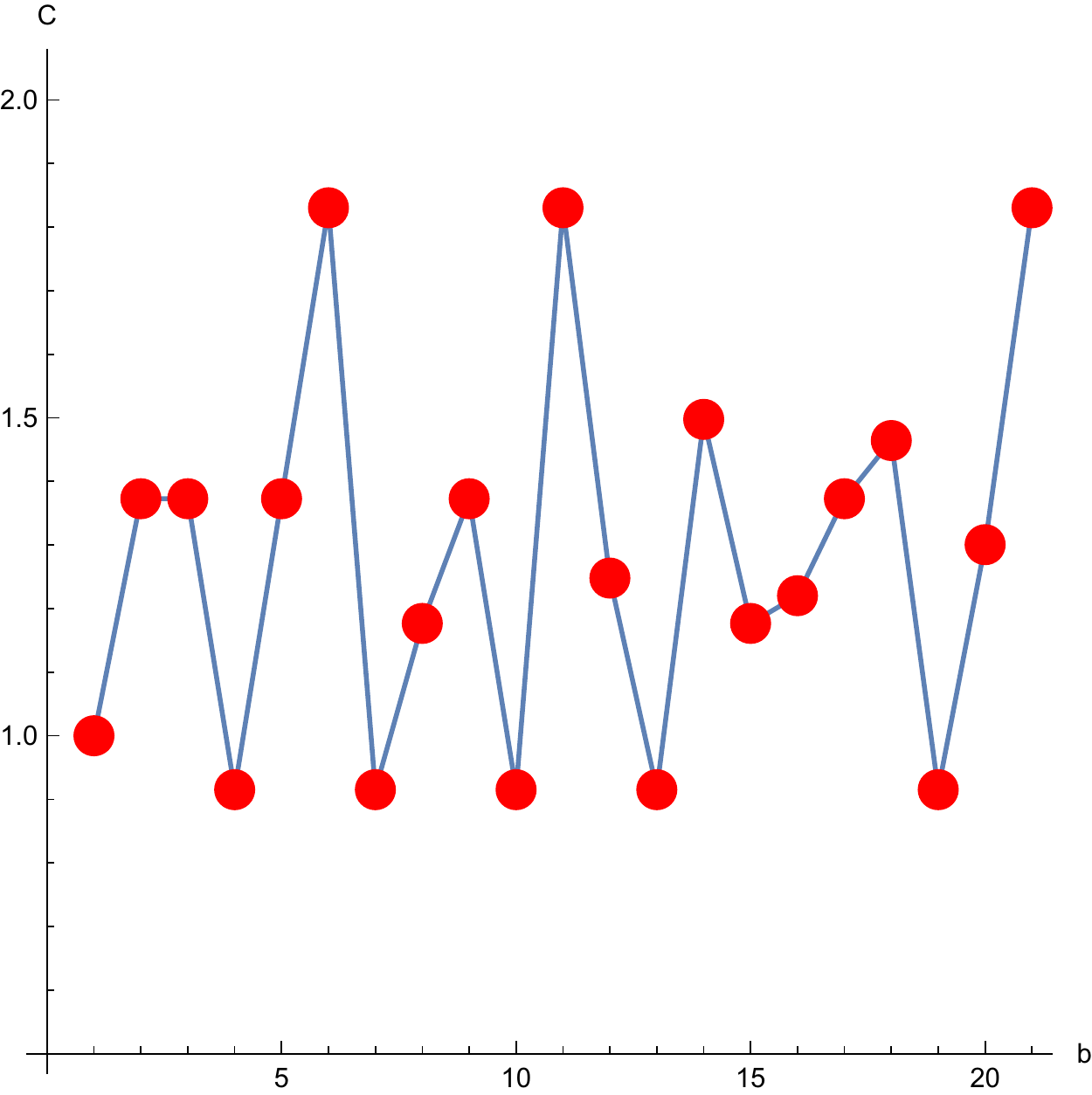}}
\caption{
The Hardy-Littlewood constants $C_k$ giving
the density rations between the first and $k$'th row in the Gaussian primes.
\label{hardyconstants}
}
\end{figure}

To investigate the constant $C_1$, we repeated some measurements of Shanks in 1953 and Wunderlich 1973, 
We could push it further thanks to faster computers. Wunderlich stayed below $10^{24}$ we went to $2^{36}$
which illustrates Moore's law. We see that on the interval $x \in \{1, \dots 2^{36}\}$, there are $1884245341$
Gaussian primes $x+i$ and $1372531868$ real Gaussian primes.  The fraction is $1.3728244749214085$.

\begin{figure}[!htpb]
\scalebox{0.2}{\includegraphics{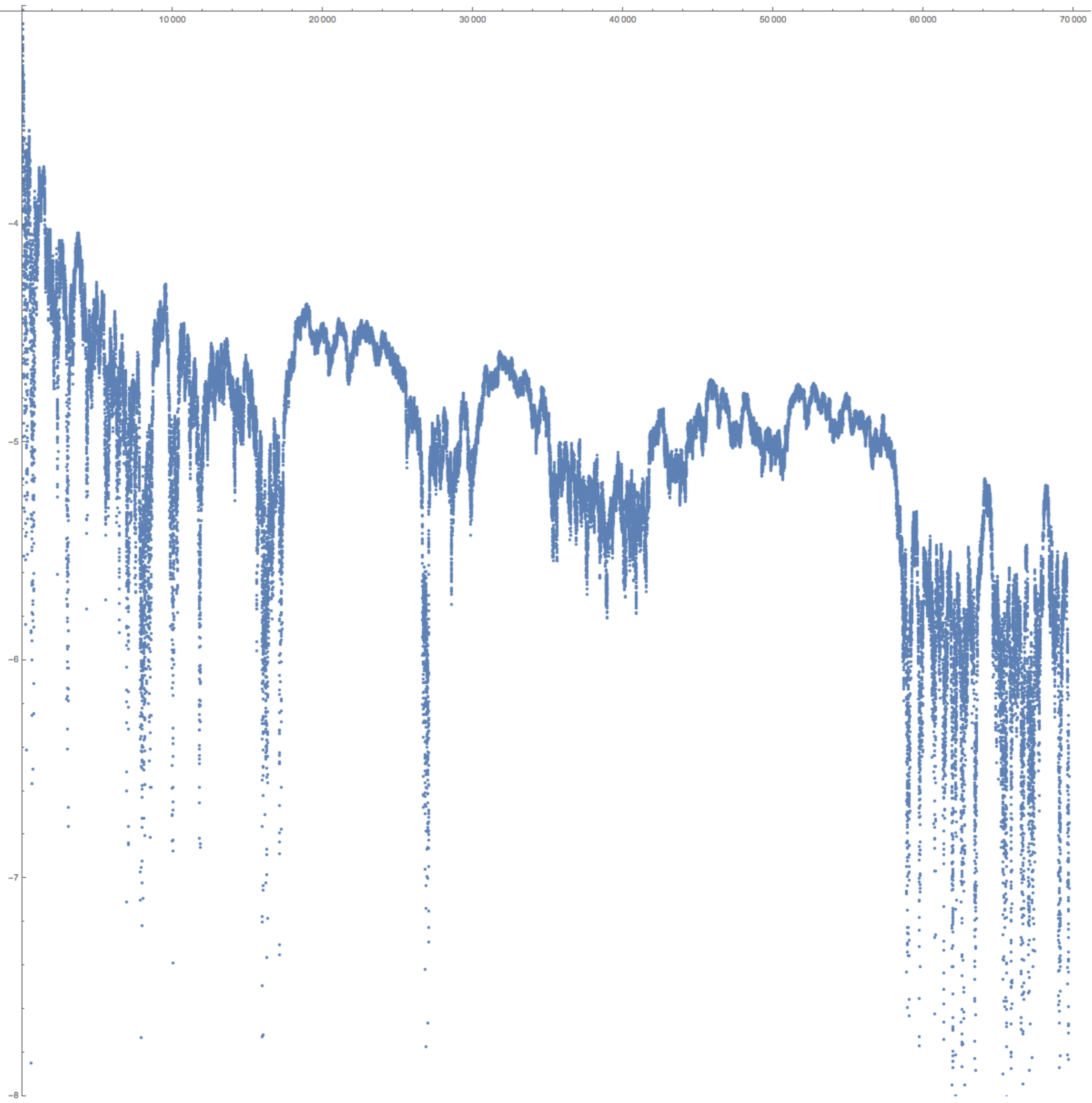}}
\caption{
The deviation from the Hardy constant $C$ when increasing the fraction of primes
of the form $a+i$ and real Gaussian primes. These are up to $10^{11}$ (one unit is $10^6$.).
Convergence is slow. If the Riemann hypothesis holds, we expect an error of the 
order $(\log(x))^2/\sqrt{x}$. For $x=2^{36}$ this is $2 \cdot 10^{-3}$ only.  
\label{hardyconvergence}
}
\end{figure}

As we don't even have 6 digits reliably while Shanks essentially confirmed
5 digits already.  We initially started to doubt the Hardy-Littlewood conjecture until 
realizing that one should not expect a better convergence than in the case of 
the prime counting function itself. Why is the error going to zero so slowly? 
Does it have to do the {\bf Chebyshev bias} which is a
strange phenomenon in the {\bf  prime race}?  No, the reason is much simpler: \\

The Riemann hypothesis is known to be equivalent to 
$|\pi(x)-{\rm Li}(x)| \leq C \sqrt{x} \log(x)$. 
Because ${\rm Li}(x) = x/\log(x)$ means 
$|\pi(x)/{\rm Li}(x)-1| \leq C \log^2(x)/\sqrt(x)$, we expect an error of
the same size also for the convergence to the Hardy-Littlewood constant. 
In order to get $6$ digits reliably, we expect having to go up to $10^{18}$.

\section{Gaussian prime matrices} 

In this section we introduce and investigate some {\bf linear algebra problems} 
for matrices defined by Gaussian primes. 
The idea is to look at a square window in the Gaussian integer lattice and place a $1$
in the corresponding matrix, where we have a Gaussian prime. Otherwise, we place a $0$. 
This construction produces a square {\bf Gauss prime matrix} $A(z,n)$ which can be 
studied experimentally. We
\begin{wrapfigure}{l}{4.1cm} \begin{center}
\includegraphics[width=4cm]{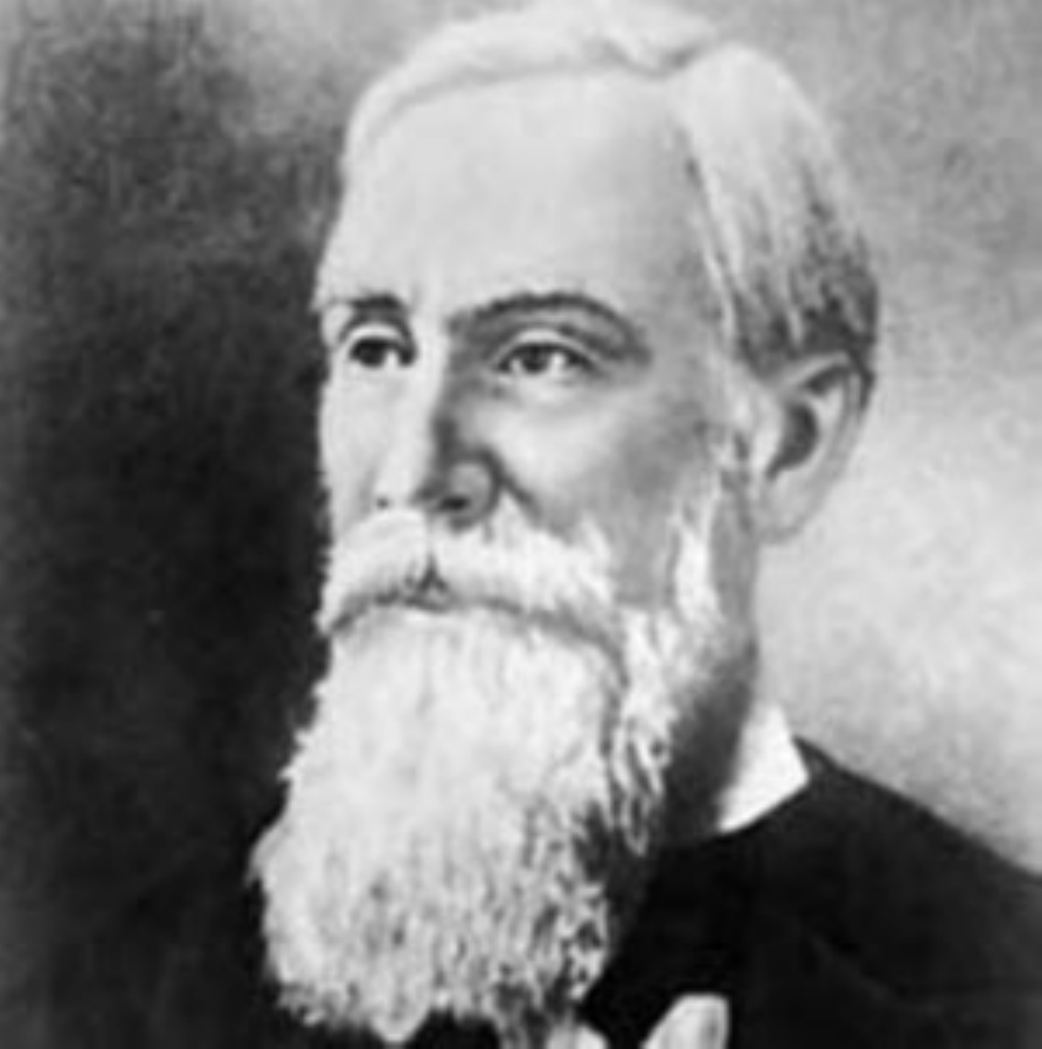}
\end{center} \end{wrapfigure}
noticed empirically for example that two rows with odd distance
are {\bf negatively correlated} and rows of even distance are {\bf positively correlated}. 
We follow here a bit the great Chebyshev, who
has both relations to probability theory and number theory. Chebyshev also noticed a
bias between the prime counting functions $\pi_1$ and $\pi_3$ for $4k+1$ and $4k-1$ primes.
The Bertrand-Chebyshev theorem telling that there is always a prime between $n$ and $2n$ 
finally was a precursor of the prime number theorem.  \\

\begin{figure}[!htpb]
\scalebox{0.6}{\includegraphics{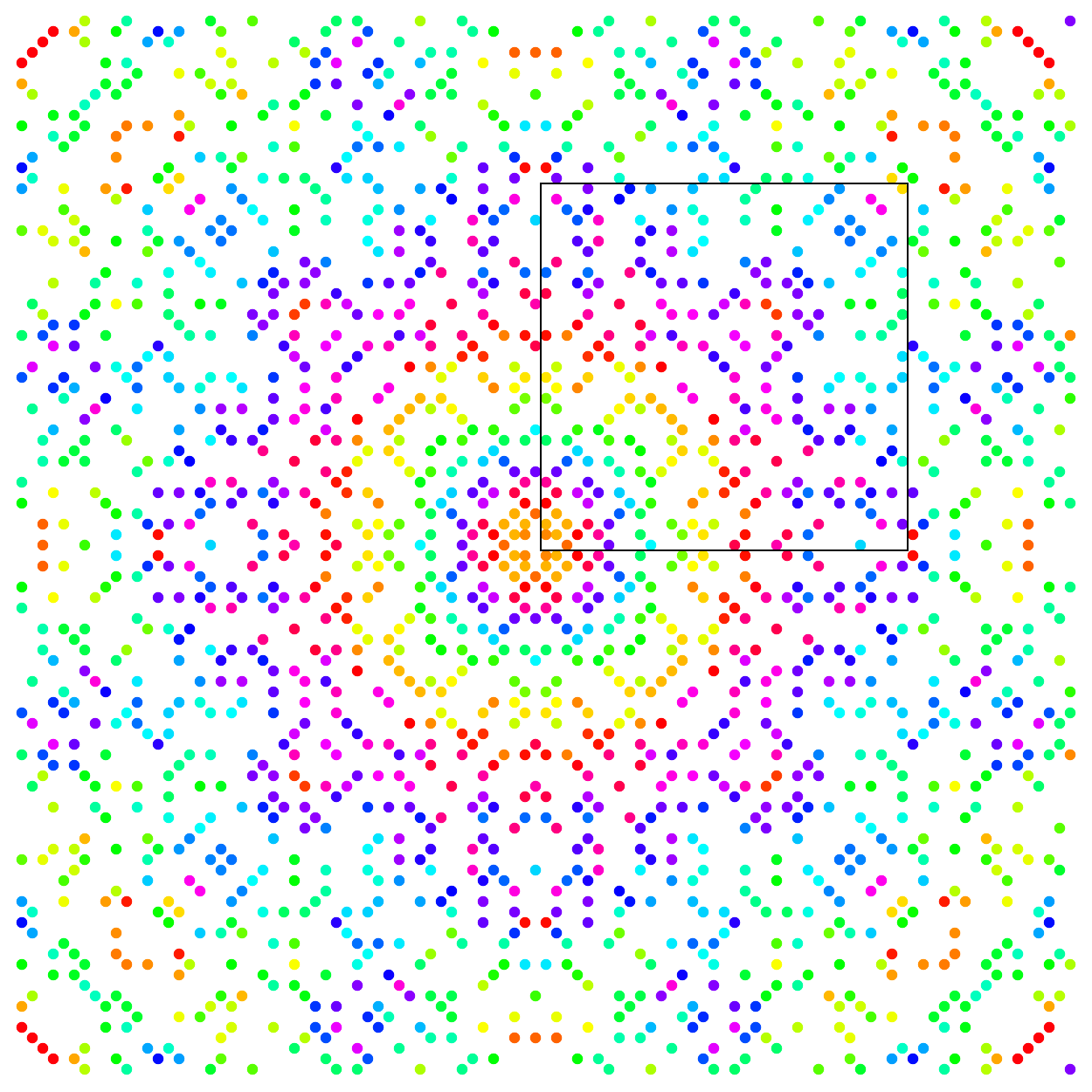}}
\caption{
The picture shows Gaussian primes in $|{\rm Re}(z)|,|{\rm Im}(z)|<50$
color-coded according how many primes there are in a neighborhood of radius $4$.
Outlined is the part produced to generate the matrix $A(35)$.
}
\label{figure2}
\end{figure}

We start with a simple but curious {\bf super symmetry} which has the effect that
half of the determinants disappear. This symmetry is based on the fact that for a 
Gaussian prime $a+ib$ with $a^2+b^2 >2$ exactly one of the $a$ or $b$ is even. 
Since all except $4$ Gaussian primes fit this  "checkerboard" pattern",
it gives some intuition why the correlation claim is reasonable: the claim would obviously 
hold if the Gaussian primes were generated by a random process. 
For matrices not affected by the symmetry, we see that they become invertible 
if they are sufficiently large and that the determinant grows with a definite super exponential 
growth rate. Also here, a disprove of Landau's problem on the infinitude of primes $n^2+1$ 
would lead to arbitrary large matrices singular $L(k+i \cdot 0,n)$ as the bottom row would then be zero.
Since the problem of determinants is linked to a hard problem in number theory, the 
problem of invertibility of the matrices $A(n)$ is least as hard. \\

Let us start to describe the set-up in more detail: 
for any positive integer $n \in \mathbb{N}$ and any {\bf Gaussian integer} $z \in \mathbb{Z}[i]$,
define the $n \times n$ matrix $A_{kl}(z,n)=1$ if $z+k+i l$ is a
{\bf Gaussian prime} and $0$ else. We call it a {\bf Gaussian prime matrix}. For non-zero and even $z$, that is if
$z$ is a non-zero multiple of $1+i$, then $A(z,n)$ is the adjacency matrix of a triangle free graph $G(z,n)$
which for large enough $n$ is a non-planar of chromatic number $2$. When looking at spectra, the non-selfadjoint
case is more interesting. The simplest non-self adjoint case $z=1$ defines the matrix
$A(n)=A(1,n)$. See Figure~(\ref{figure1}). Let $\mu(A) = \frac{1}{n} \sum_\lambda \delta_{\lambda}$
denote the {\bf density of states} of a matrix $A$. It is the normalized uniform
discrete Dirac point measure supported on the discrete spectrum of $A$ in the complex plane
$\mathbb{C}$. In random matrix theory, one calls it the  {\bf empirical measure}.
Denote by $\rightharpoonup$ weak-* convergence of measures meaning
$\mu_n \rightharpoonup \mu$ if $\int_{\mathbb{C}} f(z) \; d\mu_n(z) \to 
\int_{\mathbb{C}} f(z) \; d\mu(z)$ for any continuous function $f$ of compact support.
Let $P$ denote the diagonal matrix with entries $P_{jj}=(-1)^j$. Since for every
$A=A(z,n)$ with $z=a+ib \neq 0,a \geq 0, b \geq 0$ and even $a+b>0$,
the {\bf anti-commutation relation} $\{P,A\}=0$ holds, the spectrum of $A(z,n)$
has then a {\bf reflection symmetry} $\sigma(A) = -\sigma(A)$ like Dirac matrices in
physics. This symmetry assures that for odd $n$, the matrix $A(z,n)$ always has a kernel,
leading in that case to some unexpected linear relations between columns of the Gaussian prime
matrices. In the other cases, for $A(a+ib,n)$ with $a,b \geq 0,a+b$ odd or $(a,b)=(0,0)$, we
detect a {\bf threshold} $n_0(z)$ so that the matrices are invertible for all $n>n_0$.

\resultremark{The matrix $A=A(a+ib,n)$ with $a \geq 0, b \geq 0$ and even $a+b>0$,
satisfies the anti-commutation relation $\{P,A\}=PA-AP=0$. Its spectrum 
satisfies the symmetry $\sigma(A) = -\sigma(A)$.}

Proof.  The statement that $A(z,n)$ anti-commutes with the matrix
$P={\rm Diag}(1,-1,1,\dots)$ if $z=k+il$ with $k+l$ odd and $k,l>0$
follows from the fact that any Gaussian
prime $n+im$ has the property that exactly one of the two numbers $n,m>0$ is even if $n,m$
are not both $1$. The number $n^2+m^2=p$ is then a prime congruent $1$ modulo $4$.
Now $P A P$ changes to its negative. We see that $A$ is then conjugated to its negative
implying in the odd $n$ case that there is a zero eigenvalue. \\

The {\bf characteristic polynomial}
$p(A,x)={\rm det}(A-x)=p_0 (-x)^n + p_1 (-x)^{n-1} + \cdots + p_n$
of $A$ defines a real-valued piecewise constant function
$$  f_n(x) = \frac{\log|p_{[n x]}(A(n))|}{\log|{\rm Det}(A(n))|} \;  $$
on the interval $[0,1]$. We call it
the {\bf characteristic polynomial function} of $A(n)$. Here, ${\rm Det}(A(n))$
is the {\bf pseudo determinant} \cite{cauchybinet}, the product of non-zero eigenvalues, which is
up to a sign the last nonzero coefficient $p_k$ of the characteristic polynomial;
the function $t \to [t]$ is the floor function rendering the largest integer smaller or equal
to $t$. In many statistical settings like for Erd\"os-Renyi graphs or random matrices,
the coefficient functions $f_n(x)$ is observed to converge uniformly to a limit.
Let $Q_n,R_n$ be the matrices in the {\bf QR decomposition}
$A(n)=Q_n R_n$ and let $q_n(x,y) = Q_n([n x]+1,[n y]+1)$ on
$[0,1) \times [0,1)$. The row vectors $X_k = (q_n(1,k), \dots, q_n(n,k))$
have constant $L_2$ norm $1$ and are pairwise perpendicular.
As usual for data, the {\bf expectation} of $X$ is
${\rm E}[X]=\frac{1}{n} \sum_k X(k)$, the {\bf covariance} as
${\rm Cov}[X,Y] = {\rm E}[X Y]-{\rm E}[X] {\rm E}[Y]$, the {\bf standard deviation} as
$\sigma[X] = \sqrt{Cov[X,X]}$ and the {\bf correlation} is
${\rm Cor}[X,Y] = {\rm Cov}[X,Y]/(\sigma(X) \sigma(Y))$. We can apply these notions
especially to the row vectors
$R_k(n) = \{ 1+k i, 2+k i, \dots, n+k i\}$ of the Gaussian prime matrices.
The function ${\rm Li}(x)$ is the {\bf Eulerian logarithmic integral}
${\rm Li}(x)=\int_2^x dt/\log(t)$.
For the graphs $G(n) = (\{2,\dots, n+1\}, \{(a,b) \; | \; a+ib \; {\rm prime} \}$
the Euler characteristic decreases monotonically like $-C n^2/\log(n)$. The graphs are bipartite
and so triangle free. We have $\chi(G(n))=n-\pi(n)/2$, where $\pi(n)$ counts the number
of primes in $\{ (a+ib) \; | \; 2 \leq a \leq n+1, 2 \leq b \leq n+1 \; \}$.  \\

Our measurements indicate that there are real constants
$C_0,C_1,C_2,C_3$, a rotationally symmetric measure
$\mu = \rho(r) dr d\theta$ with smooth $\rho$ in the complex plane and
a measure $\nu = \sigma(x) dx$ with smooth $\sigma$ on the interval $[0,1]$ such that:

\begin{center}
\fbox{
\parbox{12cm}{
\parbox{5.5cm}{
{\bf A) Invertibility} \\
${\rm det}(A(n))>0$ for $n>28=C_0$. 
}
\parbox{6.1cm}{
{\bf B) Determinant} \\
$\frac{1}{n \log(n) \log(\log(n))} \log|\det(A(n)|)  \to C_1$
}
}
}

\fbox{
\parbox{12cm}{
\parbox{5.5cm}{
{\bf C) Trace} \\
$\frac{1}{{\rm Li}(n)} {\rm tr}(A(n))  \to C_2$
}
\parbox{6.1cm}{
{\bf D) Eigenvalues} \\
$\mu(\frac{C_3 \sqrt{\log(n)}}{\sqrt{n}} A(n)) \rightharpoonup \mu$
}
}
}

\fbox{
\parbox{12cm}{
\parbox{5.5cm}{
{\bf E) Characteristic polynomial} \\
$f_n(x)dx \to \sigma(x)$. 
}
\parbox{6.5cm}{
{\bf F) QR factorization} \\
${\rm E}[X_y] \to 0$ for $y \in (0,1)$. 
}
}
}

\fbox{
\parbox{12cm}{
\parbox{5.5cm}{
{\bf G) Row correlation} \\
${\rm Cor}[R_k(n),R_l(n)] \to 0$ for $k \neq l$. 
}
\parbox{6.5cm}{
{\bf H) Covariance sign} \\
${\rm sign}({\rm Cov}[R_k(n),R_l(n)]) \to (-1)^{k-l}$
}
}
}
\end{center}

\vspace{5mm}
Computer code to investigate all these statements is included in the TeX source 
code of this file.\\

\begin{figure}[!htpb]
\scalebox{0.22}{\includegraphics{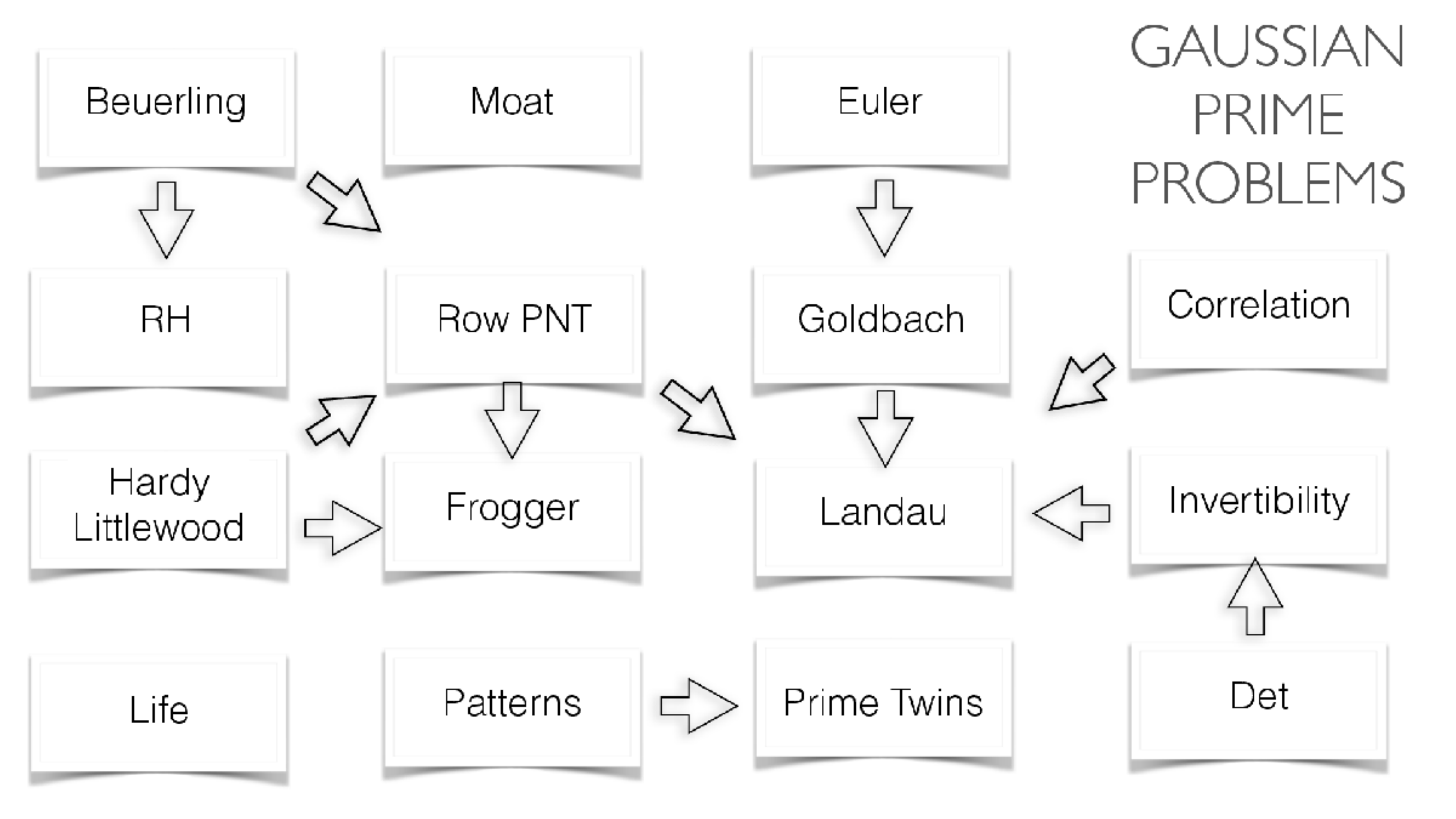}}
\caption{
Some problems about Gaussian integers. 
}
\label{problems1}
\end{figure}

\begin{figure}[!htpb]
\scalebox{0.5}{\includegraphics{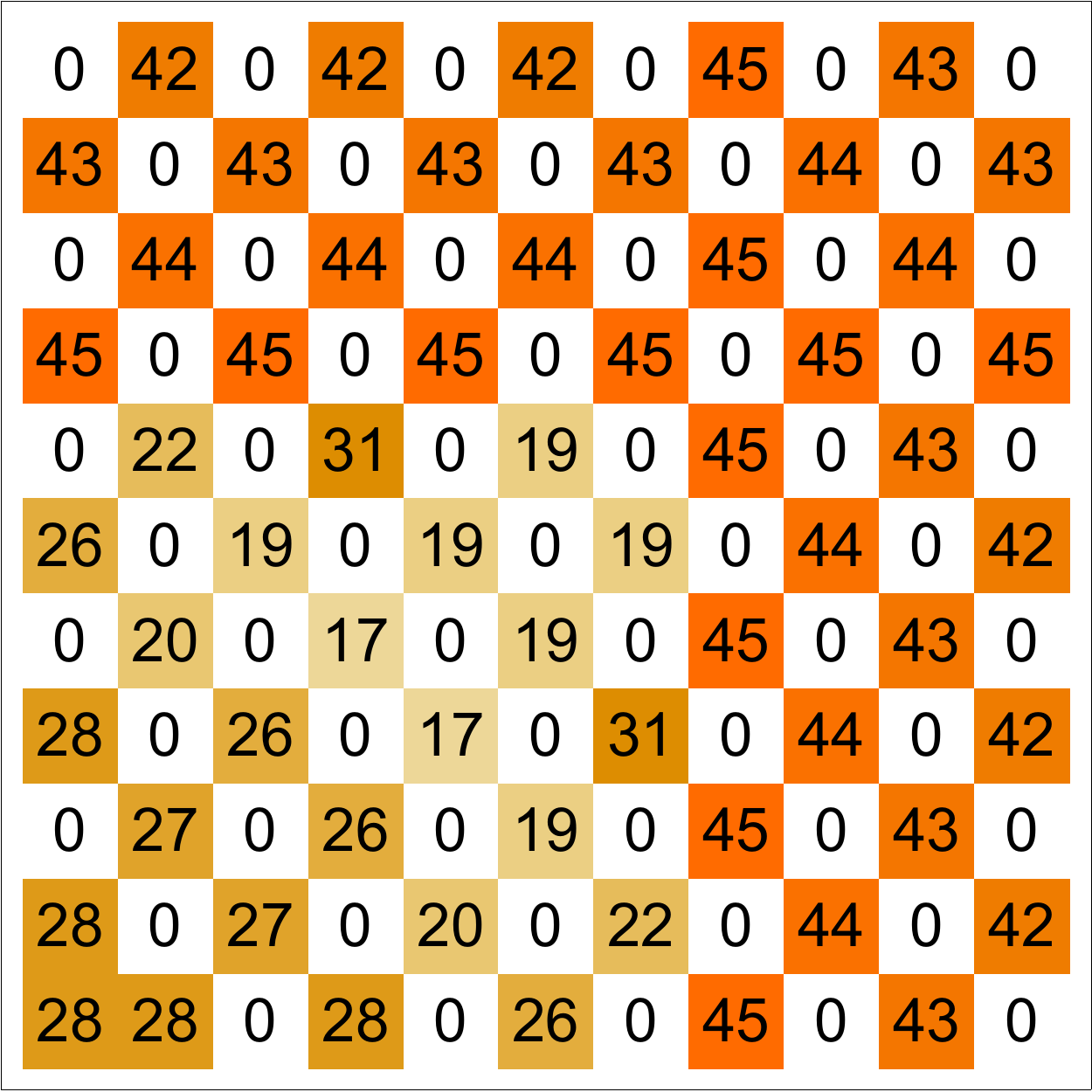}}
\caption{
As an illustration to A) we compute the apparent
non-invertibility thresholds of a matrix $A(z,n)$ with
$z=a+ib$ and $0 \leq a \leq 30, 0 \leq b \leq 30$. An entry $0$ means
that the matrix $A(a+ib,n)$ is singular for arbitrary large $n$. We know from
the anti-commutative relation $\{P,A\}=0$ which gives a spectral symmetry
$\sigma(A)=-\sigma(A)$ in the complex plane forces in the odd $n$ case to have a
zero eigenvalue of $A(z,n)$.
}
\end{figure}

\begin{figure}[!htpb]
\scalebox{0.12}{\includegraphics{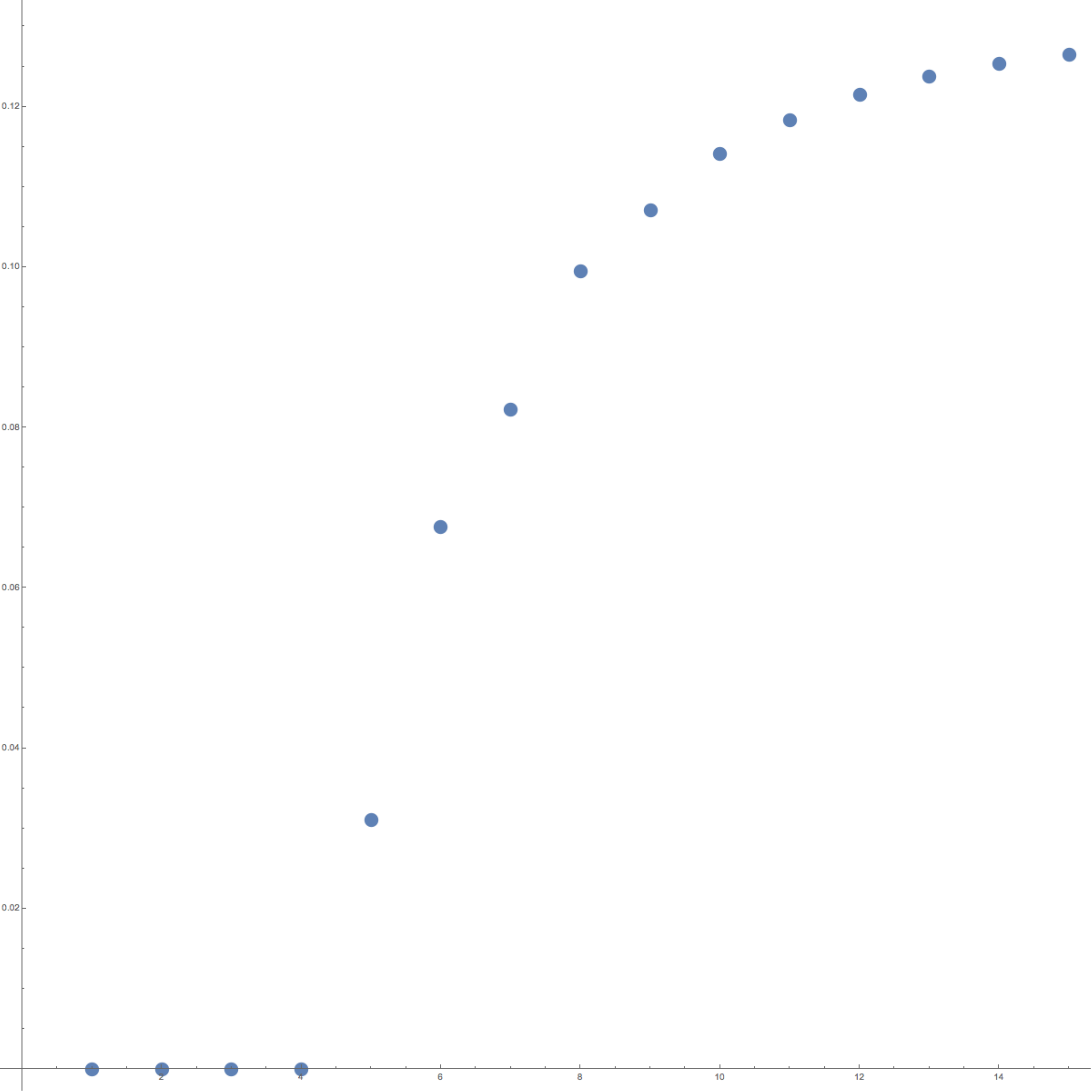}}
\scalebox{0.12}{\includegraphics{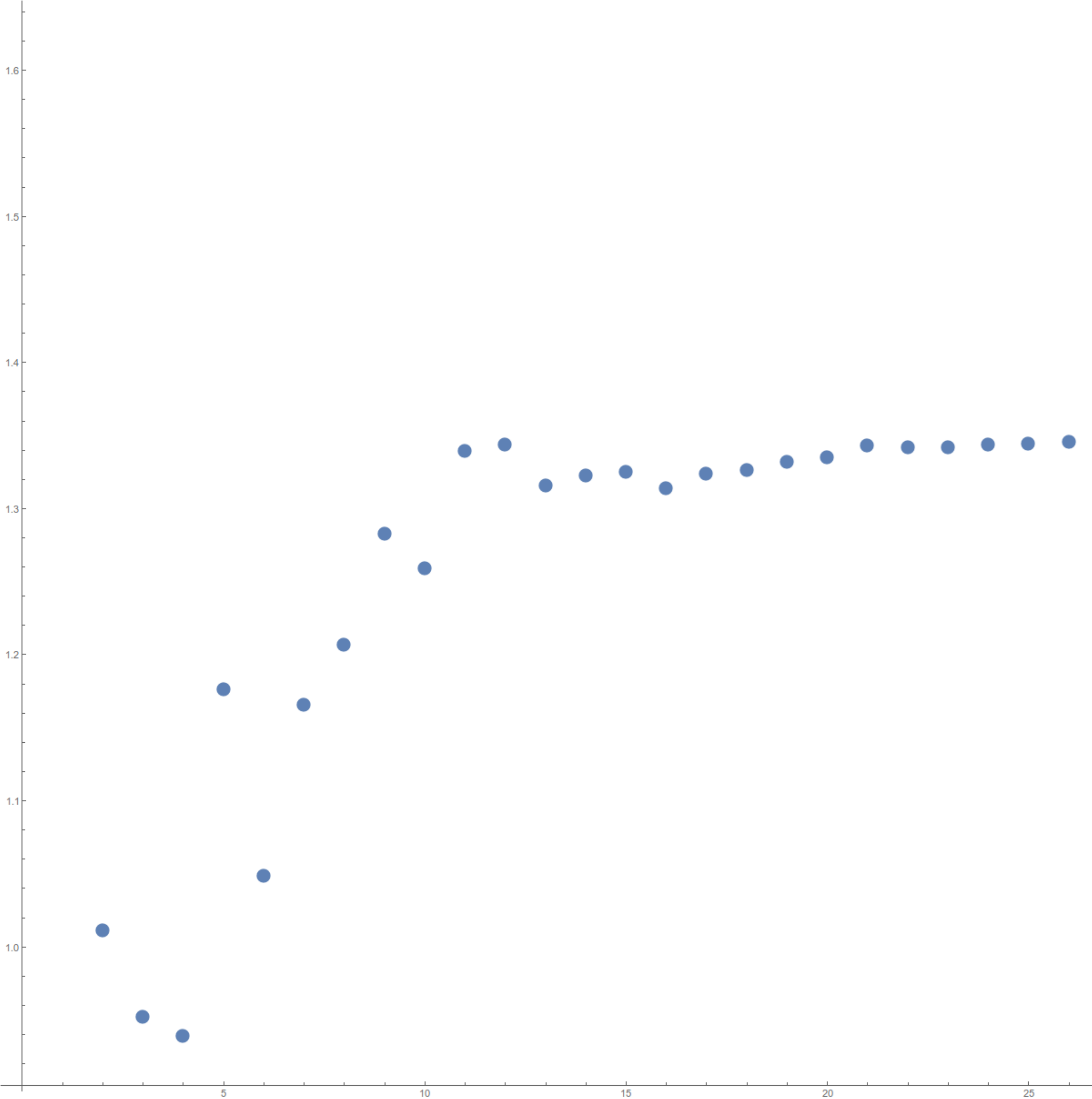}}
\caption{
To the left we see the values of the determinant function in B) for
$n=2^k$ with $4 \leq k \leq 16$.
To the right we see the value of the trace function in C) for up
to $n=2^{29}$.
}
\end{figure}

\begin{figure}[!htpb]
\scalebox{0.12}{\includegraphics{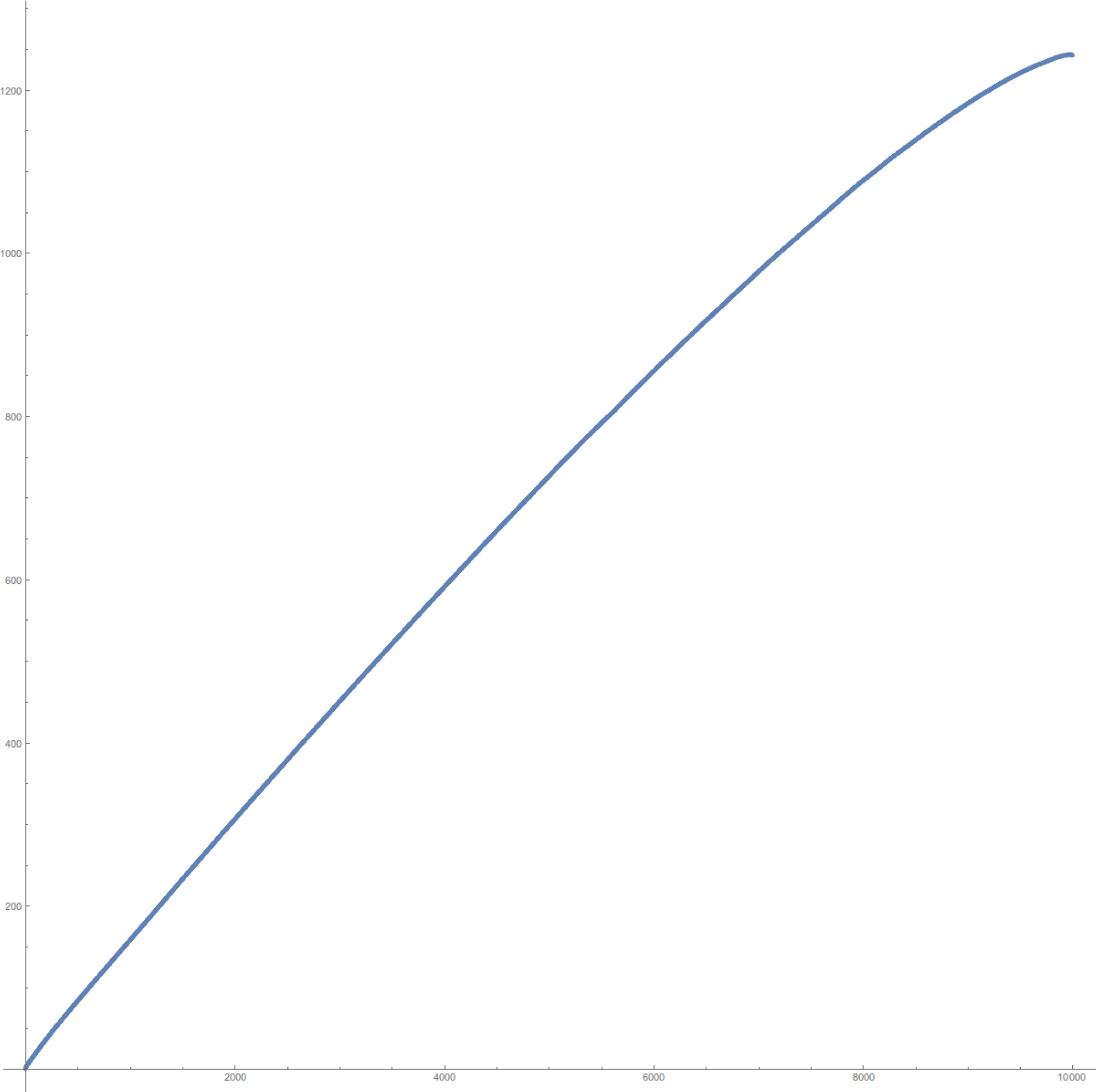}}
\scalebox{0.12}{\includegraphics{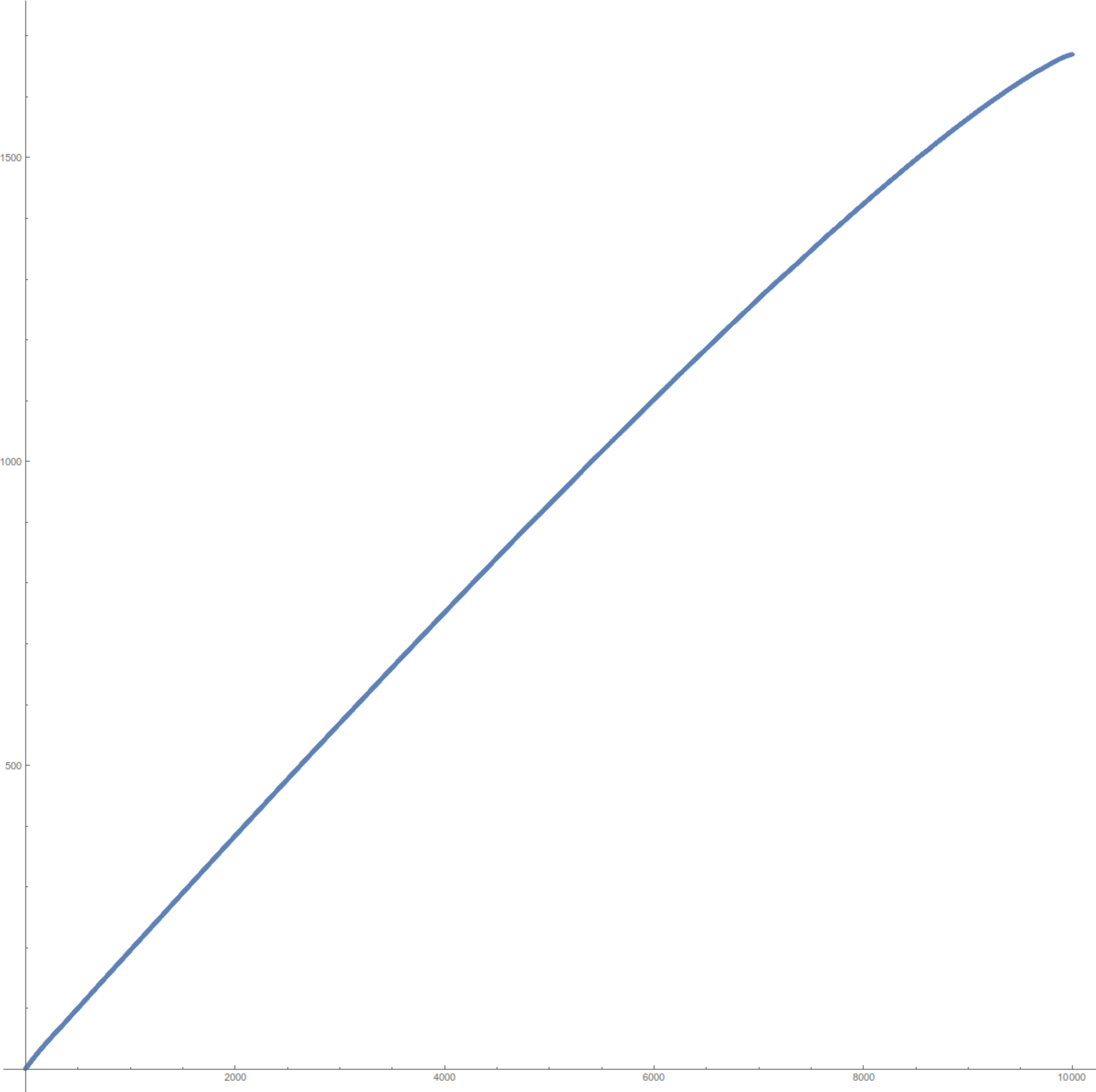}}
\caption{
To the left, we see the characteristic polynomial function $f_n(x)$ for $A(n)$ with $n=10000$
By B), the value $p(0)={\rm det}(A(n))$ to the left converges to a constant.
Statement F) claims that all scaled coefficients converge.
To the right we see the characteristic polynomial function for a random
0-1 matrix of the same size. As the coefficients of the characteristic
polynomial of a matrix is is $p_k(A)={\rm tr}(\Lambda^k A) = \sum_{|P|=k} \det(A_P)$,
here the sum is over all $k \times k$ sub matrices $A_P$ of $A$, one knows
in the random case that $f_n(x) \to x$ for all $x$.
}
\end{figure}

{\bf 1)} The numerical evidence has been accumulated only with relatively small 
$n$, namely $n \leq 2^{15}$, where we compute with $32768 \times 32768$ matrices.
We are confident about A) and predict it is the most reachable of
all these statements. But it is certainly not easy as the generalized problem with 
matrices $L(k+0 i,n)$ would imply the unresolved {\bf Landau problem} about the existence 
of infinitely many rational primes of the form $n^2+1$. We contemplated about these problems 
while writing a linear algebra exam for the Math 21b service course at Harvard.
The actual exam only featured the problem to compute ${\rm det}(A(5))$ and
mentioned  A) as the "21b conjecture". The threshold A) will grow for matrices $A(z,n)$ 
with larger $|z|$, where primes are more sparse as its easy to see that there 
are arbitrary large disks in $\mathbb{Z}[i]$ without primes.  \\

{\bf 2)} In B)-D), smaller order growth rates like $\log(\log(\log(n))$ might not be 
easy to detect in experiments. Anyhow, proving bounds like B) or C) appear theoretically 
out of sight. We do not even know whether the {\bf logarithmic potential}
$\frac{1}{n} \log|\det(A(n)|)$ diverges. This later divergence statement would definitely 
deserve the name conjecture (especially, when taking $\limsup$) 
because we see it diverge logarithmically. \\

{\bf 3)} For E) the result holds for random matrices, as one can express $p_k$ 
through minors. We also expect part F) to be true for random matrices. 
The characteristic polynomial function $f_A(x)$ is interesting in more general setups. 
We see clear convergence in classes of matrices like if $A(n)$ is the Kirchhoff matrix
of complete graphs $K_n$, circular graphs $C_n$, wheel graphs $W_n$ or random graphs. 
To make the statement reasonable also in more graph theoretical setups, we 
divide by the Pseudo determinant. In our case, we could take the determinant. \\

{\bf 4)} Matrices like $A(k(1+i),n)$ are symmetric.
The symmetric case $A(0,n)$ is the only where one still has a finite 
threshold $28$ for zero determinant, all cases $A(a+ai,n)$ have a kernel for odd $n$
as the super-symmetry relation $\{ A(n),P \}=0$ implies then a spectral symmetry $\sigma(A(n))=-\sigma(A(n))$
which for odd cases implies the existence of a zero eigenvalue.  We see that $A \to A P$ 
switches the sign of half of the eigenvalues and that $A \to P A$ switches the other half. 
In the symmetric case $B(n)=A(-n-ni,2n)$, we have ${\rm rank}(B(n)) = n$ for $n>28$.
As the spectrum is real, D) becomes a semi-circle law which should be compared in the random 
case with Wigner's theorem. \\

{\bf 5)} The divergence rate of ${\rm tr}(A(n)))$ is a statement about arithmetic progressions 
of Gaussian primes on lines parallel to the coordinate axes. We measure what is expected
from having the primes equally distributed on diagonal lines. 
Questions about arithmetic progressions is attributed to John Leech
\cite{Guy}. But density questions seem difficult as one does not even 
know whether for fixed $a,b$, infinitely many primes of the form 
$k (a+i b)$ exist. For any prime $a$, work of Hecke shows that there are infinitely many 
primes of the form $a+ib$. The linear algebra questions posed here are related as if 
there was an $a_0$ for which on the line $a_0+ib$ only finitely many primes existed, then 
there exists $b_0$ such that $A(a+ib,n)$ is singular for $a_0-n<a<a_0,b>b_0$. 
A statement like that for any $a+ib$ in the first quadrant with $a+b$ odd, there exists
$n_0$ such that $A(a+ib,n)$ is nonsingular for all $n>n_0$ would prove the infinitude 
of Gaussian primes on any axes parallel to to the coordinate axes. \\

{\bf 6)} Observation D) is intuitively explained by {\bf Girko's law} which holds for random 
$\{0,1\}$ matrices, where each entry takes value $1,0$ with some probability $p$ going to
zero for $n \to \infty$. Probabilistic thinking in that context are not new.
For random percolation considerations, see \cite{Vardi}. More generally, one has
looked at generalized primes and corresponding Beuerling integers \cite{Beuerling,HilberdinkLapidus}. 
The random case in some sense investigates the 
statistics obtained from generalized Gaussian Beuerling primes. \\

{\bf 7)} Case E) is also observed for random matrices. One can write the
coefficients of the characteristic polynomial in terms of minors \cite{cauchybinet} and
get $f_n(x) \to x$. 
Special coefficients in the characteristic polynomial are the trace and the determinant, covered in
more detail in statements B),C). One can get similar convergence statements for
each of the coefficients $p_k$ of the characteristic polynomial. If the picture in E) is 
correct, each of the $p_k$ grows in the same way than the determinant, with the exception
of the trace which grows at a slower pace. \\

{\bf 8)} Situation F) is very much expected
as the QR decomposition ``decorrelated" the columns of a matrix. It could be
interesting but as the initial remark about spectral symmetry shows, there
are unexpected linear relations between columns. The question about QR decomposition
might be the least interesting, but we looked at QR decompositions of matrices $A(n)$
when looking for problem sources for exam problems. \\

{\bf 9)} In the random matrix case, the expectation of having an invertible
matrix goes to zero exponentially fast. Gaussian prime window matrices
$A(n)$ behave similarly as random matrices. The rate C) is what one expects
from the prime number theorem, as there is no reason why on the diagonal
of $A(n)$, the distribution of the Gaussian primes should be different than on
any other side diagonal. \\

{\bf 10)} Since $\frac{1}{n} \log(\det(A(n))) = \frac{1}{n} \sum_\lambda \log(\lambda)$ is the
{\bf Riesz potential} $\Delta \mu$ of the measure $\mu(A(n))$, problems B) and D)
are related and B) is the logarithmic potential at 0. Also the minimal absolute
value of eigenvalues appears to be what one would expect from random
matrices, where $\log(\det(A(n))/(n \log(n))$ converges. The claim that for large enough 
$n$ the Riesz potential is bounded away from $0$ by a definite constant is much 
weaker than observation B) but appears also out of reach. \\

{\bf 11)} The mathematics of estimating determinants of large matrices appears frequently 
in Hamiltonian dynamics or solid state physics. Here is an example of an open problem
in ergodic theory which I had worked myself for a long time without success: 
if $x(k+1) - 2 x(k) + x(k-1) = c f(x(k))$ is the recursion of the symplectic
map $T(x,y)=(2x-y+c f(x),y)$ on the torus $\mathbb{T}^2$ and 
$L(n)u_k= (c f'(x_k)+2) u_k + u_{k-1} + u_{k+1}$ 
linearizes a finite piece of orbit of length $n$, then the Green-McKay-Weiss formula 
${\rm tr}(dT^n)-2={\rm det}(L(n))$ for the Jacobian matrix $dT^n$ of the iterate $T^n$ leads 
to the Thouless formula assuring that 
$\lambda(x_0,x_1) = \lim_{n \to \infty} (1/n) \log({\rm det}(L(n)))$ is the 
{\bf Lyapunov exponent} of the orbit starting at $(x_0,x_1)$, measuring sensitivity with respect to
changes of initial conditions. By Oseledecs theorem, the limit exists for almost all $(x_0,x_1)$. 
Integrating $\lambda(x_0,x_1)$ over the normalized Lebesgue measure on 
$\mathbb{T}^2=\mathbb{R}^2/(2\pi \mathbb{Z})^2$ is the Kolmogorov-Sinai entropy
of $T$. Pesin theory shows that if it is positive, then 
$T$ to a measure theoretical factor of a Bernoulli system on a set of positive Lebesgue measure. 
In the Gaussian prime case, where one deals not with band matrices, the scaling factor of the logarithmic  
potential is $n \log(n)$ not $n$, but the difficulty of establishing a limit is similar than for
Lyapunov exponents: the situation is deterministic and so of pseudo random nature. It is 
neither random nor integrable and in both cases, the difficulty is in the complement of 
integrable situations. Some progress about the Lyapunov exponent of the Standard map has been done
in \cite{Bourgain2012}.\\

{\bf 12)} In random cases (probability theory) or integrable cases 
(harmonic analysis), there are usually attacks. 
For Jacobi matrices (like for the Anderson model) has been covered by 50 year old 
analysis F\"urstenberg-Kesten. Also integrable almost periodic situations
like almost Matthieu allow Fourier theory (Aubry duality) or subharmonic methods  
(Herman, Sorets-Spencer) are understood. There are more parallels: while almost 
periodic Jacobi matrices have in general Cantor spectrum, the eigenvalues of almost 
periodic matrices like $A_n(i,j) = \cos(i \alpha  + j \beta +\theta)$ with rationally 
independent $\alpha,\beta,\theta$ appear to have a fractal limiting law in the complex 
plane.  The Gaussian integers behave a bit like
a mostly non-uniformly hyperbolic dynamical system, which both show order and randomness at 
the same time. 
Both, in the ergodic and number theoretical setups, 
the difficulty is in the analysis is the complement of integrable structures: the complement 
of KAM islands playing the role of the Erastostenes sieve picture. 
This theme \cite{TaoStructureRandomness} also appears within dynamical system theory 
\cite{MoserStableRandom}. The interplay between ergodic theory and number theory 
is hardly new \cite{FurstenbergRecurrence} and has got through a renaissance 
in the last decade.  

\section{Circular laws}  

We also look at the eigenvalues of the Gaussian prime matrices and observe a 
{\bf Wigner-Ginibre-Girko circular law} emerge, even so the density appears different than for 
{\bf random matrices} with independent coefficients as random variables,
where the law has been confirmed and a uniform density on the disc emerges. 
The idea is to take a large matrix defined by Gaussian primes and draw the 
eigenvalues in the complex planes. We see a uniform angular distribution to emerge but the radial
distribution of the spectrum is different: the density inside is slightly larger than for 
random matrices. Furthermore, we look at the coefficients of the characteristic polynomial 
and observe that, like in many other random or pseudo random situations, that they converge 
to a concave limiting function if properly rescaled. 
The Gaussian primes appear to be random enough to smooth out the angular
distribution. We see such phenomena in non hyperbolic situations which are located between
random-Anosov-Anderson type or integrable almost periodic-KAM-Mathieu type situations. 
In any case, establishing a Girko type law for Gaussian prime matrices looks not easy. 
Even the random case is technically hard and it was only recently achieved by Tao and Vu 
\cite{TaoVu} for matrices for which the entries can have discrete distribution.
The circular law has been found by {\bf Jean Ginibre} in 1965 and 
{\bf Vyacheslav Girko} in 1984. 

\begin{figure}[!htpb]
\scalebox{1.0}{\includegraphics{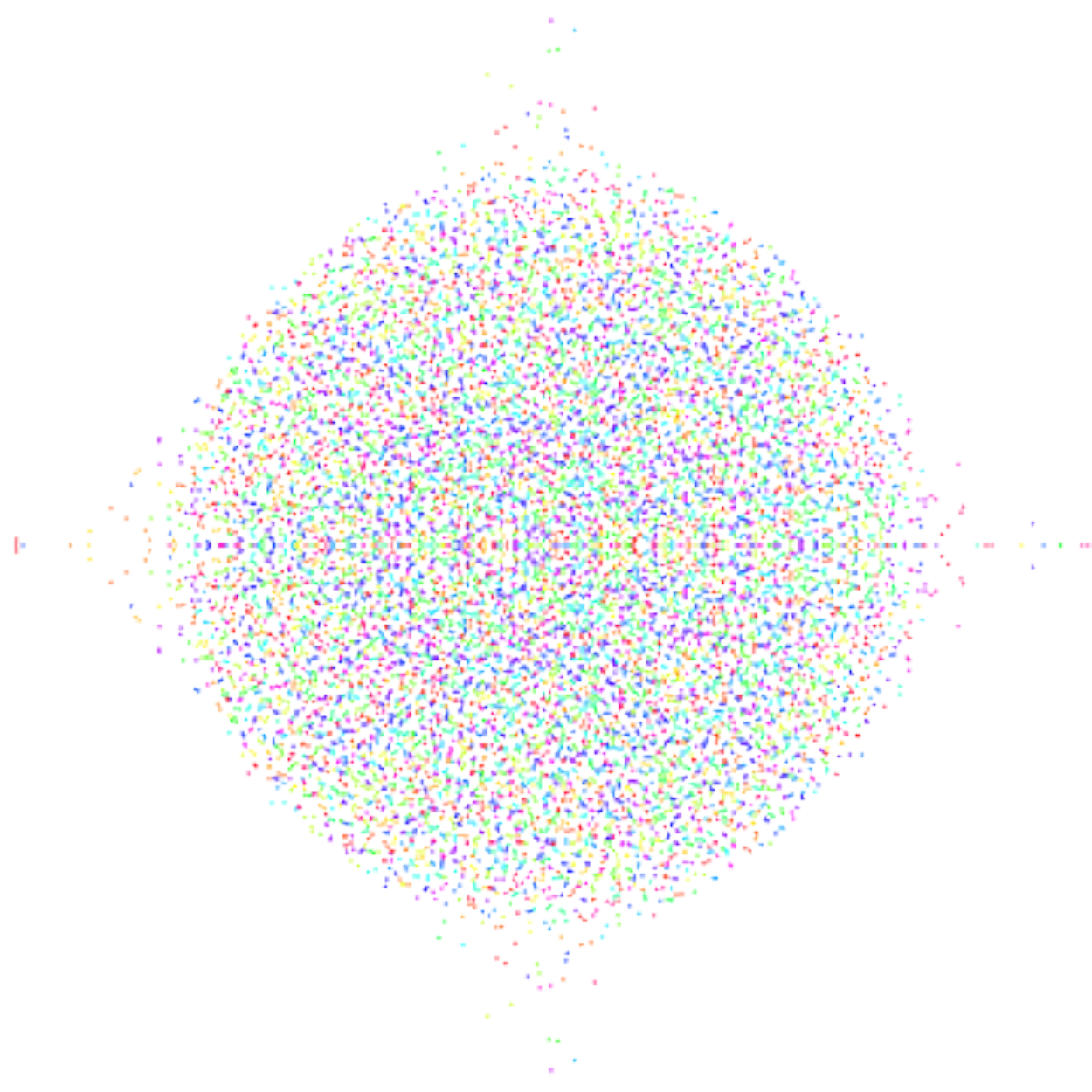}}
\caption{
We see the eigenvalues of the matrix
$A(10000)$ in the complex plane. The colors encode the
distance to the nearest neighbor, the idea being to detect any
{\bf clustering} patterns. Some {\bf spectral leaking}
near the coordinate axes can be observed but it appears to be
too small to be statistically relevant leading to a violation of
statement D). Still, there appears to be a larger inner density. }
\end{figure}

\begin{center} 
\includegraphics[width=4cm]{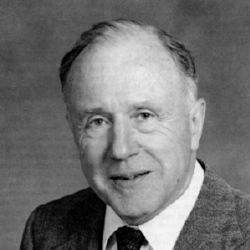}
\includegraphics[width=4cm]{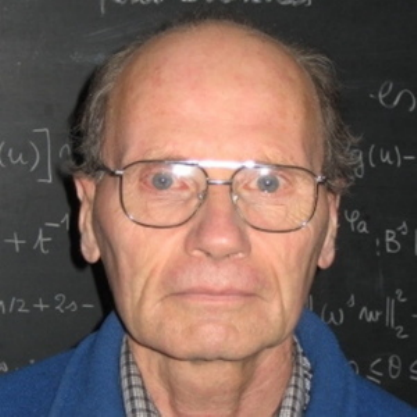}
\includegraphics[width=4cm]{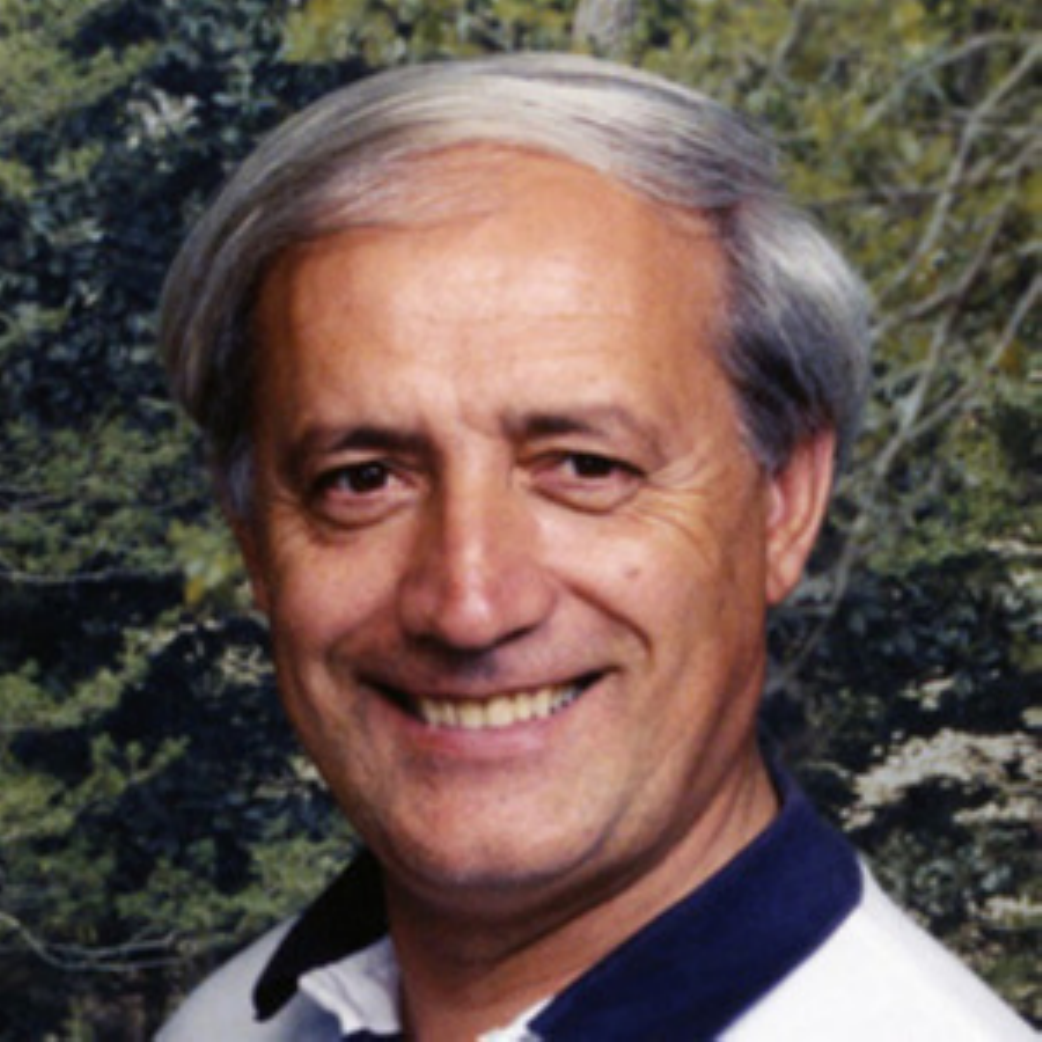}
\end{center} 

In the self-adjoint case, one has a semi-circle law.

\begin{figure}[!htpb]
\scalebox{1.2}{\includegraphics{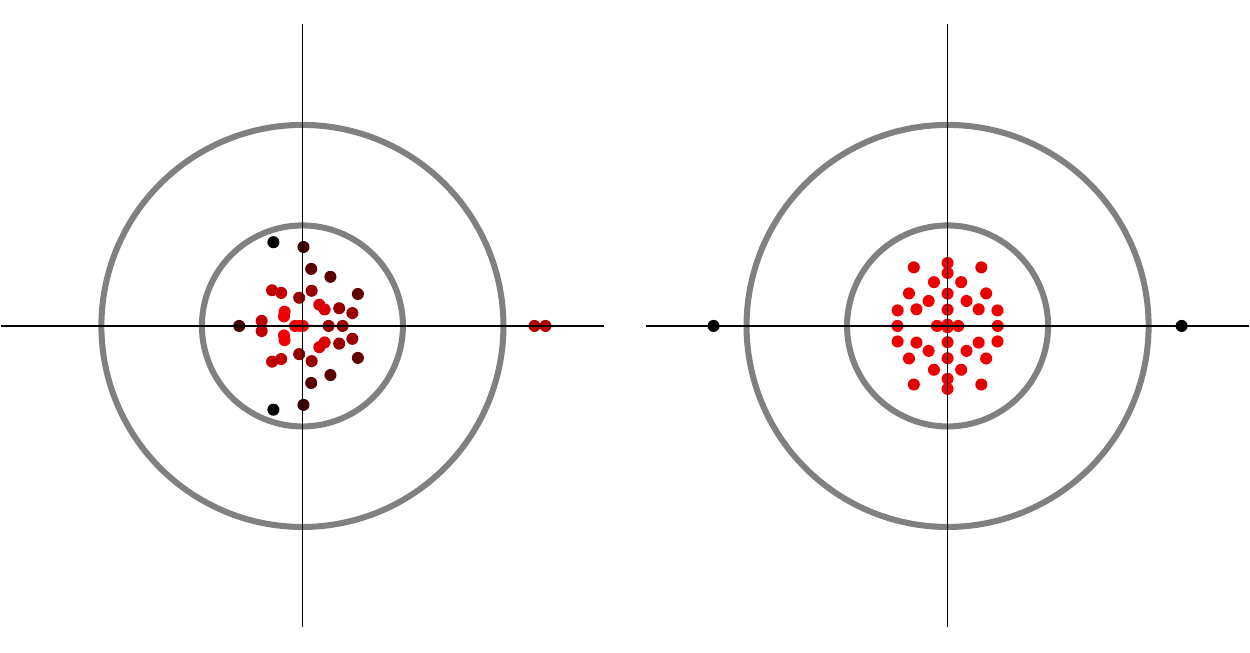}}
\caption{
We see the eigenvalues of the matrix $A(1,201)$ in the complex plane to the left
and the eigenvalues of the matrix $A(2,201)$ to the right. In the second case,
there is a super symmetry present which is reflected in an additional symmetry
$\sigma(A) \to -\sigma(A)$. The symmetry group of the left picture is $Z_2$
for the right picture $Z_2 \times Z_2$. 
}
\end{figure}

\begin{figure}[!htpb]
\scalebox{0.26}{\includegraphics{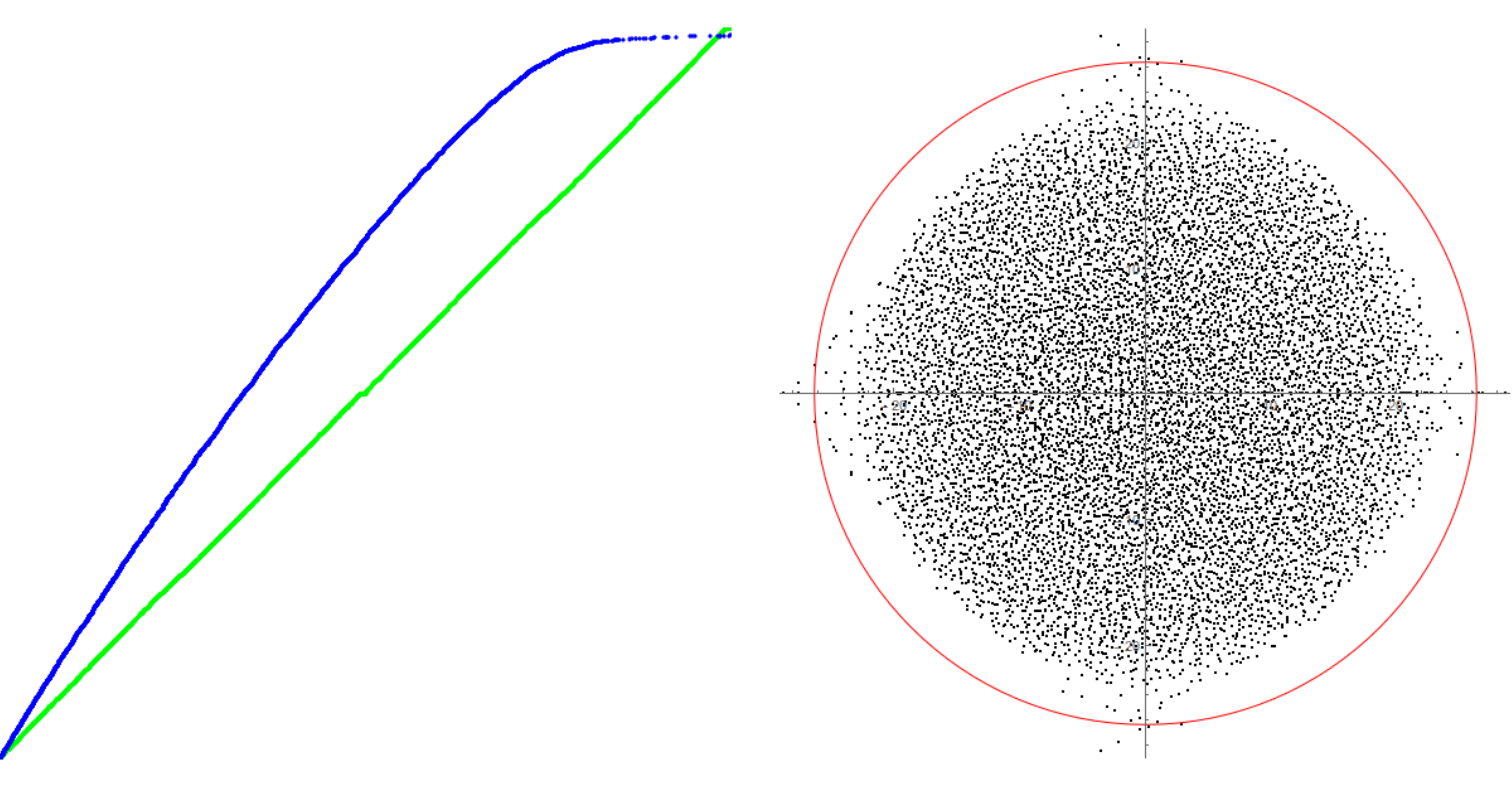}}
\caption{
This figure shows the cumulative distribution functions of the argument and
the absolute value of the eigenvalues.
We see that the density decreases towards the boundary and that it is a bit larger
inside that in the random uniform case.
The actual eigenvalues are shown to the right together with a circle
of radius $\sqrt{n}/\log(n)$. The {\bf spectral radius}
of $A(n)$ grows like $n/\log(n)$, if the Hardy-Littlewood conjectures hold.
In comparison, the spectral radius of random $\{0,1\}$ matrices grows 
linearly due to the presence of a single maximal real eigenvalue of size $n/2$.
}
\end{figure}

\begin{figure}[!htpb]
\scalebox{0.26}{\includegraphics{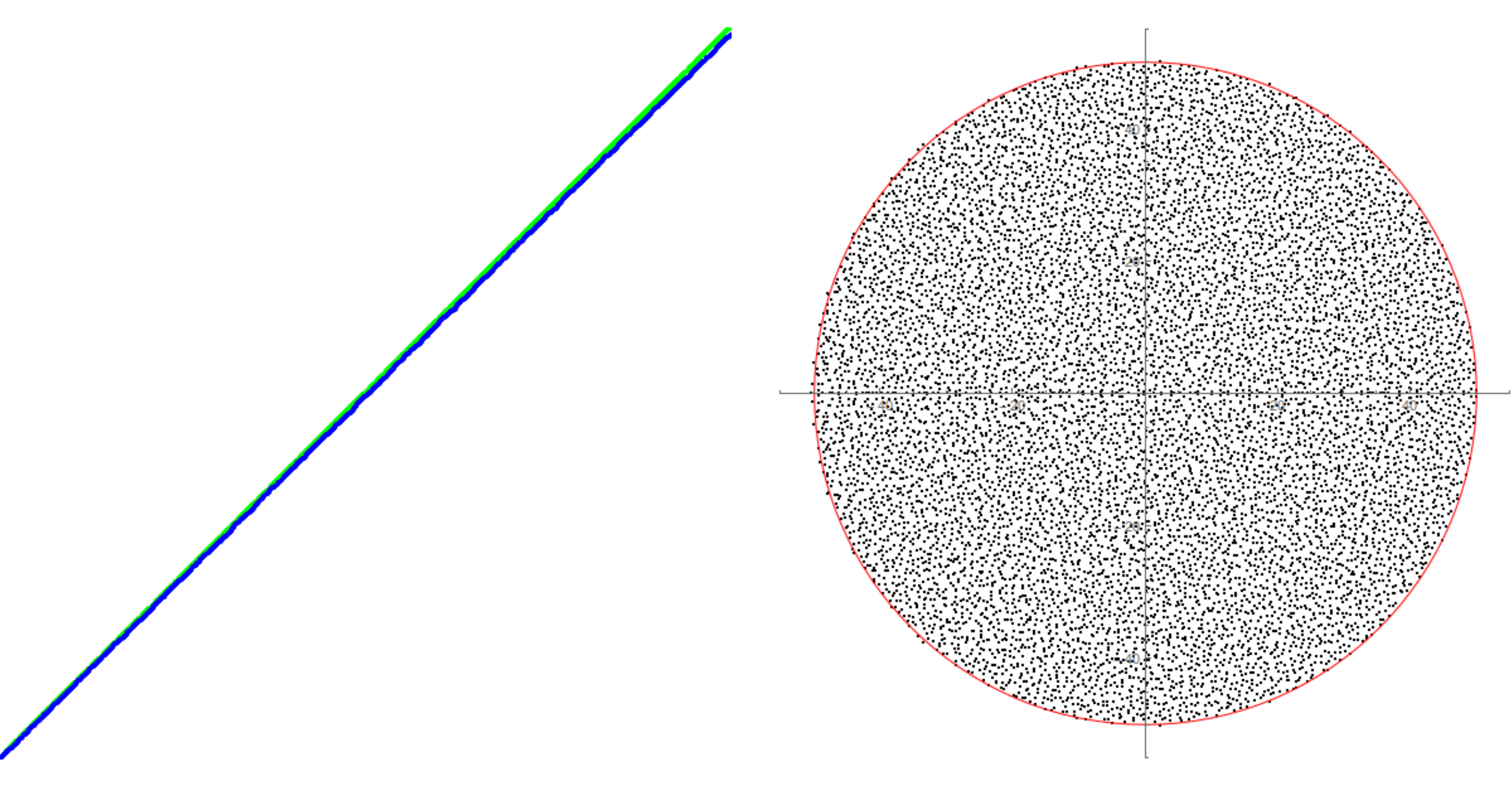}}
\caption{
As a comparison, we see a random case, where the matrices have
been chosen to take random $\{0,1\}$-values with probability $1/2$.
Also drawn is the circle of radius $r=\sqrt{n}/2$.
Girko's law for such random matrices is known to hold \cite{TaoVu}.
}
\end{figure}

{\bf Remark.} \\
{\bf 1)} The laws for random matrices were first found by 
Wigner \cite{Wigner58} in the selfadjoint case 1958, the {\bf circular law} in the non-selfadjoint
random was discovered by Ginibre in 1965 for the Gaussian case, by Girko in 1984 for more general
case. Under rather minimal assumptions like matrices taking finitely many values is covered in
\cite{TaoVu}

\section{Gaussian prime graphs}   

There is a topological connection between Gaussian primes and graphs. We will look 
in this section at graphs whose Euler characteristic is related to the 
{\bf prime number theorem} for Gaussian primes. \\
 
We like to see the adjacency matrices of a class of bipartite graphs defined by 
Gaussian primes. Define $G(n)$ as the finite simple graph 
which has as the vertex set the integers $\{2, \dots, n+1 \; \}$ and for which 
two vertices $a,b$ are connected if $a+ib$ is a Gaussian prime.  \\
Since Gaussian primes $k+il$  have the property that exactly one of the $k,l$ 
are odd, the graphs $G(n)$ are all {\bf bipartite}. They consequently have no triangles and
their Euler characteristic is $|V|-|E|$ and by Euler-Poincar\'e given by $b_0-b_1$, 
where $b_0$ is the number of connectivity components and $b_1$ is the {\bf genus}. 
Since $|V(n)|=n$ and $|E(n)|=\pi(n)$ counts the number of Gaussian primes in 
$[2 \times n+1] \times [2 \times n+1]$ we have a class of graphs for which we 
know the behavior of the Euler characteristic well. \\

What is the exact relation between the prime counting function for rational primes and Gaussian
primes? Since every rational prime $p$ in $\{2\} \cup P_3$ corresponds to exactly $4$ primes on the      
circle $|z|=p$ and every rational prime in $P_1$ corresponds to exactly $8$ primes on 
the circle $|z|=\sqrt{p}$ the growth rates are related.
Let $\pi_G(r)$ is the {\bf Gaussian prime counting function} giving the number of Gaussian primes
in $\{ |z|^2 \leq r \}$. We have now

\resultremark{
The Gaussian prime counting function satisfies
$\pi_G(x)=4+8 \pi_1(x)+4 \pi_3(\sqrt{x})$. 
}

This implies that it behaves like $8 \pi_1(x) 4 \pi(x)$, where $\pi(x)$ is the 
prime counting function of the rational integers. 

\begin{figure}[!htpb]
\scalebox{0.9}{\includegraphics{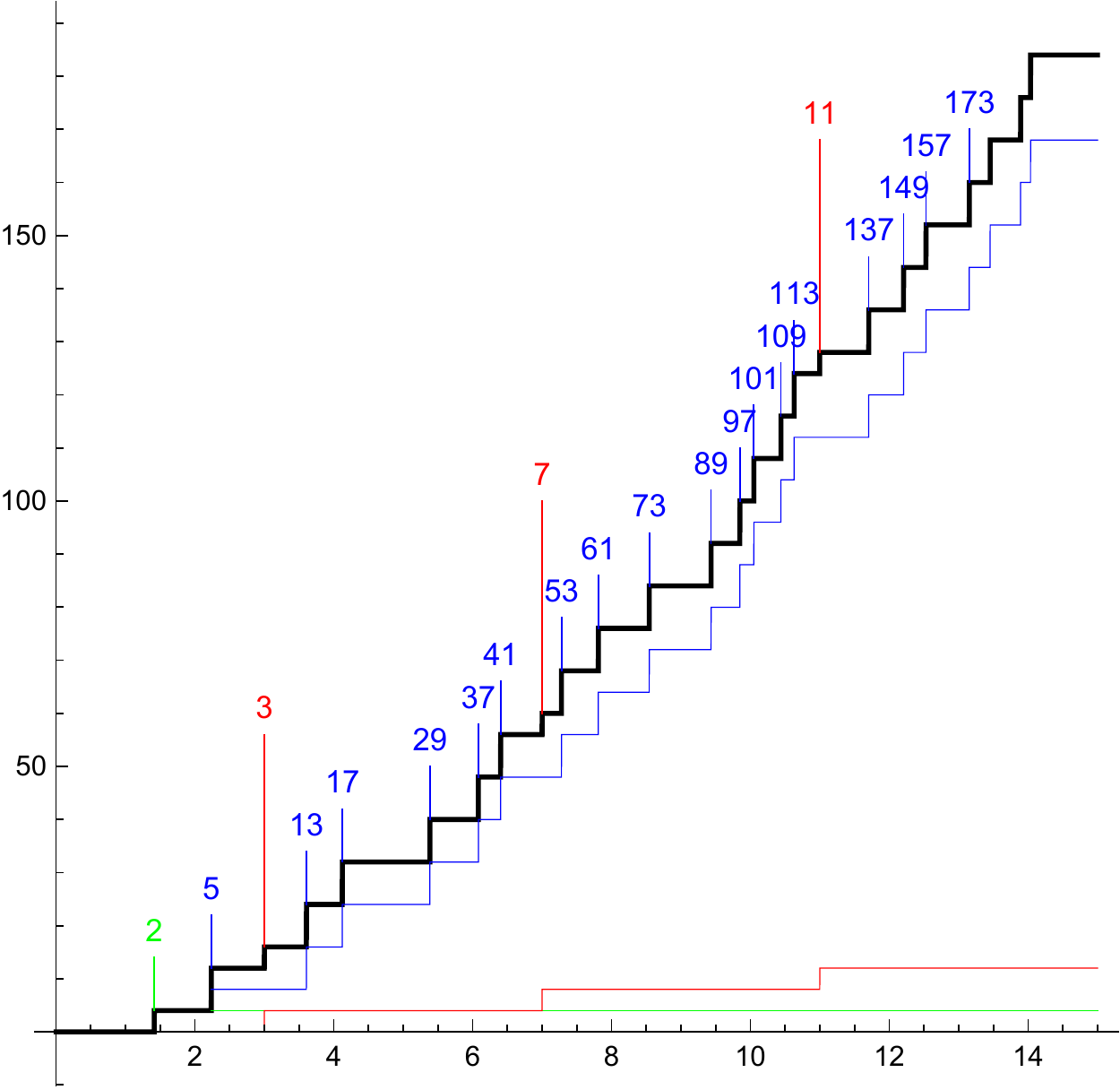}}
\caption{
The prime counting function $\pi_G(r)$ for Gaussian primes is the sum
of three functions $\pi_2(\sqrt{r})$ counting the even primes,
$8 \pi_1(r^2)$ and $4 \pi_3(r)$, where $\pi_k(r)$ counts the
number of primes giving remainder $k$ when divided by $4$.
Since $\pi_2$ is bounded by $1$ and $\pi_1(r)$ and $\pi_3(r)$
grow in the same way, but $\pi_1$ evaluates at $r^2$,
the Gaussian prime counting function $\pi_G$
is essentially governed by the growth of $\pi_1(r^2)$. In other words,
one can pretty much neglect the even primes and the primes on the axes.
The error $\pi_G(r)-{\rm Li}(r)/2$ is described by the Riemann
hypothesis. Also the problem of estimating the lattice point error
$\pi_0(r)-\pi r^2$ is open but of a different kind: there are no primes
involved and the number of solutions $x^2+y^2=n$  for non-prime $n$
can be much larger.
}
\label{howmanyprimes}
\end{figure}

We know from the prime number theorem $\pi(x) \sim {\rm Li}(x)$ and 
$\pi_1(x) \sim \pi(x)/2$ that $\pi_G(x) \sim 4 {\rm Li}(x)$. 

\resultremark{
The Euler characteristic of $G(n)$ grows like $\pi_G(x)$.
}

Proof: 
The graphs are bipartite. This prevents any triangles to appear
in the graph so that the Euler characteristic is $\chi(G(n))=|V|-|E|=n-\pi(n)$, where $\pi(n)$
is the number of Gaussian primes in the rectangle $[2,n+1] \times [2,n+1]$ of the 
complex plane.  The genus of these graphs is therefore $b_1=1-|V|+|E| = 1-n+\pi_G(n)$. \\

It follows that the Riemann hypothesis for Gaussian primes has a topological
interpretation as the genus is directly linked to the prime counting function $\pi$.

\begin{figure}[!htpb]
\scalebox{1.0}{\includegraphics{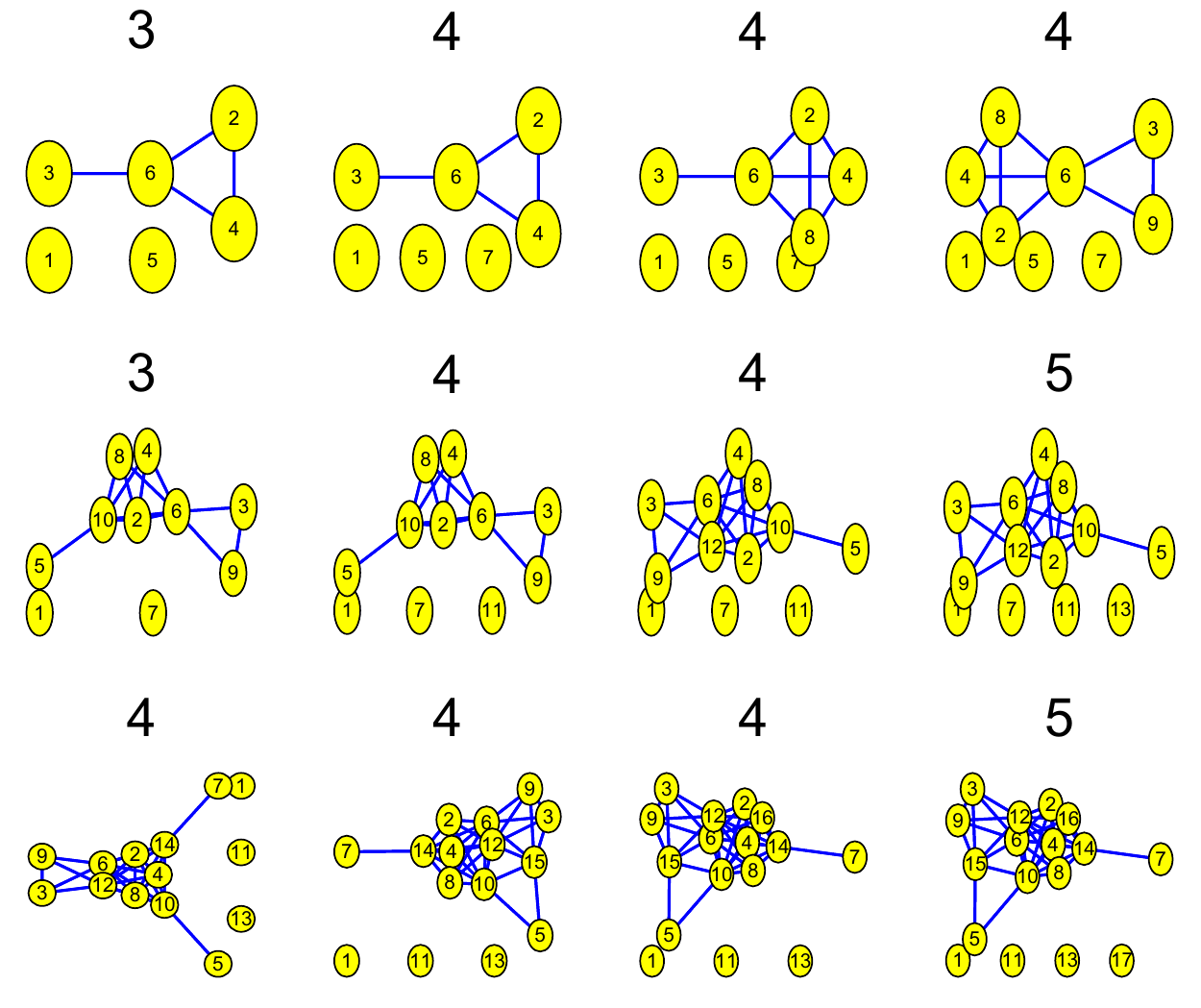}}
\caption{
We see the bipartite graphs $G(n)$ for $n=4, \dots 19$. The vertex set of $G(n)$
is $\{ 2,\dots, n+1 \}$ and two vertices $k,l$ in $G(n)$ are connected if $k+il$ is a
Gaussian prime. The label above a graph gives its Euler characteristic $b_0-b_1=|V|-|E|$.
The Euler characteristic $\chi$ of $G(n)$ is directly related to the Gaussian
prime counting number: $\chi(G(n))=n-\pi_G(n)$, where $\pi_G(n)$ is the number of Gaussian
primes in $[2,n+1] \times [2,n+1]$. Similarly as in a picture drawn by Knauf \cite{Knauf98}
but much more elementary, the Riemann hypothesis is here related to topological
properties of graphs. }
\label{figure1}
\end{figure}

\begin{center}
\includegraphics[width=4cm]{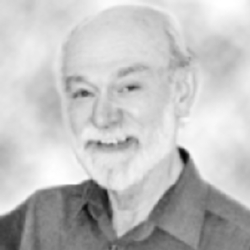}
\includegraphics[width=4cm]{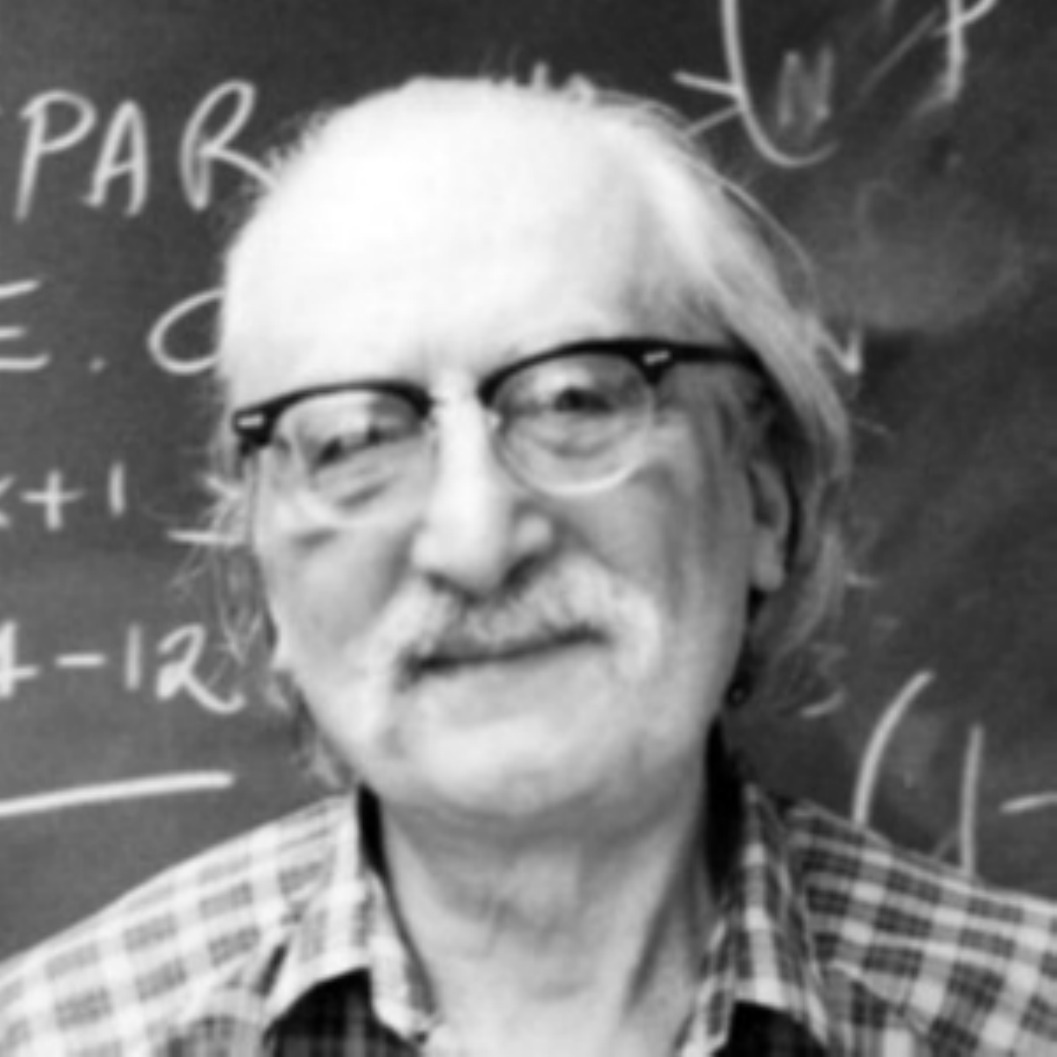}
\end{center} 

Similarly, we can look at graphs defined by Lipschitz or Hurwitz integers. 
Let $H(n)$ denote the graph with vertex set $a+ib$ with $1 \leq a,b \leq n$  
for which two vertices $a+ib, c+id$ are connected, if $a+ib+jc+kd$ is a
Lipschitz prime,  which means that  $a^2+b^2+c^2+d^2$ is prime. 
Similarly, define $L(n)$ denote the graph with the same vertex set for which 
two vertices $a+ib,c+id$ are connected, if $(a+ib+jc+kd)/2$ are Hurwitz
primes (meaning $(a^2+b^2+c^2+d^2)/4$ being prime. Since even and odd can 
not be connected, the graph $H(n)$ has at least two component. But we can conjecture
to have enough Hurwitz primes so that

\conjecture{ 
The graphs $H(n)$ have exactly two components for $n \geq 4$. 
The graphs $L(n)$ has exactly one component for $n \geq 2$. 
}

Similarly as before for the Gaussian graphs, the Lipschitz graphs
are bipartite so that the Euler characteristic is is equal to 
$|V|-|E|$, where $|V|$ are the number of vertices and $|E|$ the 
number of edges. 

\begin{figure}[!htpb]
\scalebox{0.4}{\includegraphics{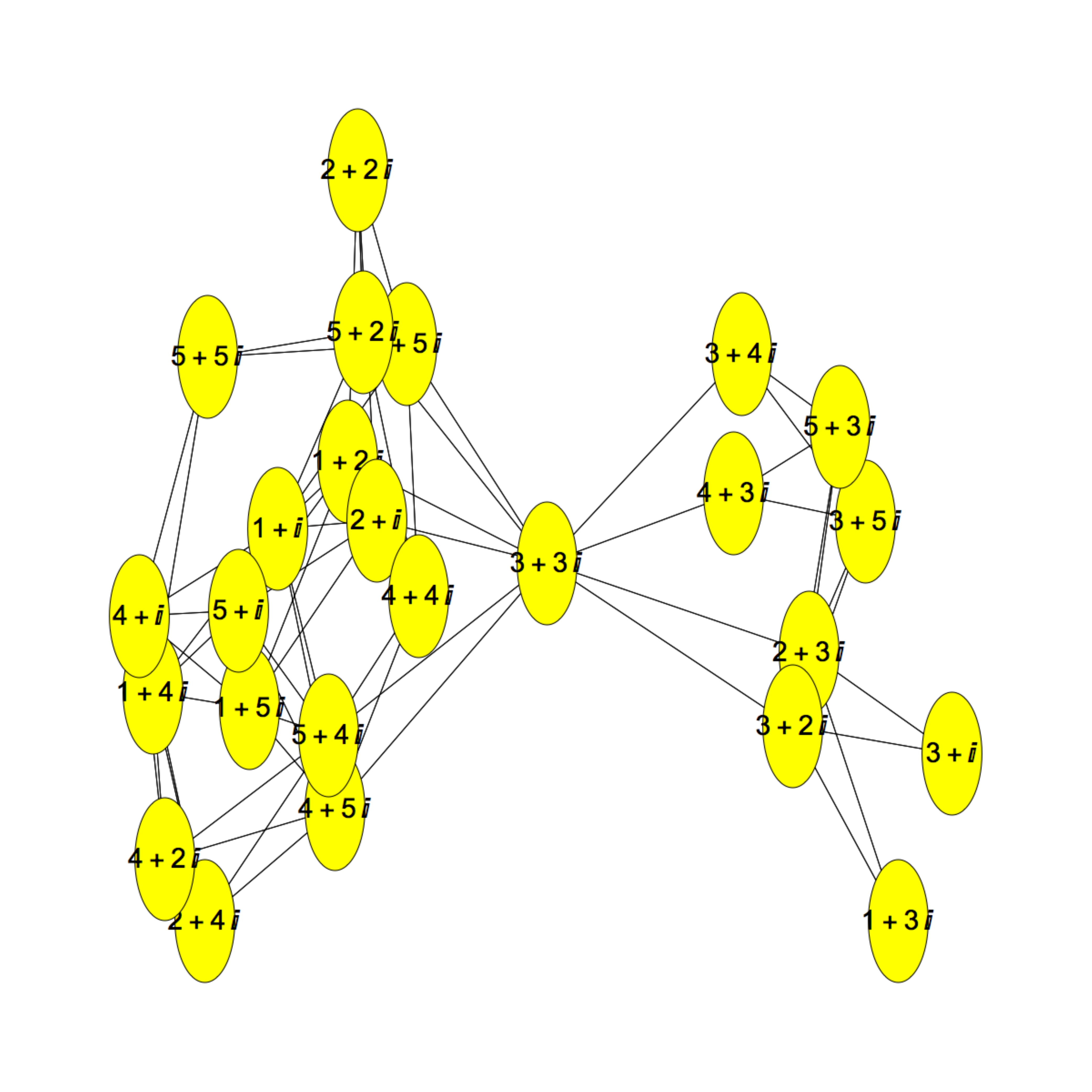}}
\caption{
In this graph connects Gaussian integers $a+ib, c+id$ with $a>0,b>0$
for which $a+ib+jc+kd$ is a Lipschitz prime. 
\label{hurwitzgraph}
}
\end{figure}

\begin{figure}[!htpb]
\scalebox{0.4}{\includegraphics{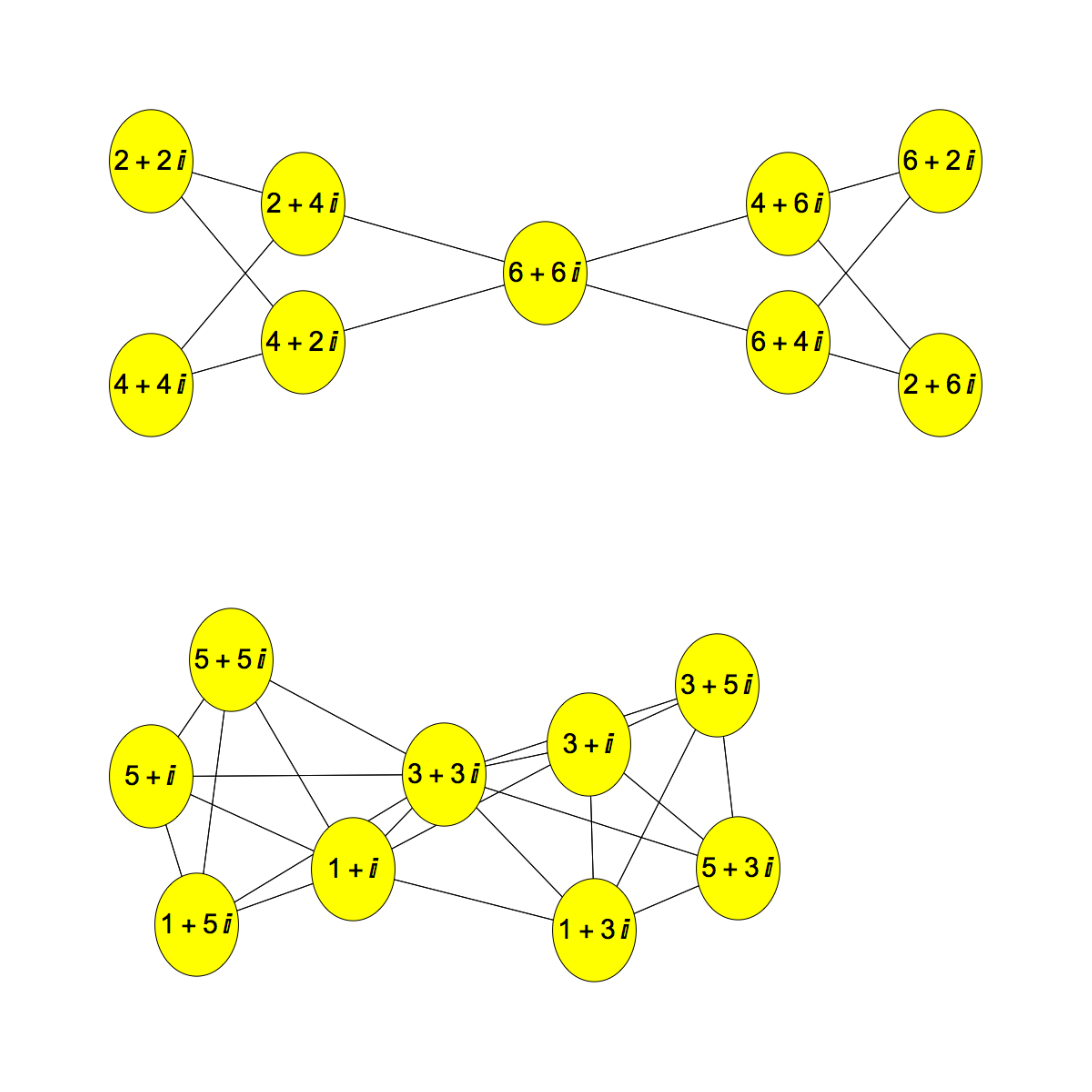}}
\caption{
Here, we connect Gaussian integers $a+ib,c+id$ with $a>0,b>0$ if $(a+ib+jc+kd)/2$
is a Hurwitz prime.
\label{hurwitzgraph}
}
\end{figure}

{\bf Remarks.}
{\bf 1)} Graphs can be introduced differently. 
For every integer $n$ one can look at the countable simple graph $G$
for which the Gaussian primes form the vertex set and where two
primes are connected if their distance is $\leq n$.
The {\bf Gaussian moat problem} of Motzkin and Gordon \cite{Guy}
asks whether the connected component containing the prime $1+i$
can become infinite for some $n$, \cite{GethnerWagonWick,Conrad}.
The moat problem looks independent of the linear algebra problem posed here.
But when seen like this, it is also an eigenvalue problem.  \\

{\bf 2)} The prime counting function should be compared with the
{\bf Gaussian circle problem} which deals with the number of Gaussian integers within a disc
of radius $x$. If $\pi(x)$ counts the number of primes there,
a result of Koch implies with the Riemann hypothesis that 
$|\pi(x) -{\rm Li}(x)/4|  \leq C \sqrt{x} \log(x)$
for a constant $C$ and sufficiently large $x$. The reason is that half of the primes remain on the axes and
the other half gets their norm squeezed by the square root and multiply by $8$.
For the Gauss lattice problem on the other hand,
one believes that for every $\epsilon>0$, one has $|N(r)-\pi r^2| \leq r^{1/2+\epsilon}$
for large enough $r$.  \\

{\bf 3)} The monotonicity question for the Euler characteristic of the graphs $G(n)$
would follow from that fact that for any sufficiently large $a$, there exists $2 \leq b \leq a$ for
which $a+i b$ is a Gaussian prime. A counter example
would also lead to counter examples in A). Of course, if Landau's 4th problem about the existence
of infinitely many primes $p$ for which $p-1$ is a square is false, then there was a counter
example. But having the density of primes on rows quantitatively so well nailed down by the
Hardy-Littlewood constants, there is hardly anybody who would doubt that there are
Gaussian primes on every row $\{ n + b i \; | \; n \in \mathbb{Z} \}$ of the complex integers. 

\section{Gaussian Zeta function}   

In the context of prime numbers, it is impossible to avoid mentioning zeta functions and especially
the {\bf Riemann hypothesis}, which is widely considered the most important open problem in mathematics.
For any algebraic structure containing primes and multiplicative norm $N$, one has a zeta function
$\sum_p N(p)^{-s}$. A bit more than a decade ago, the literature for the Riemann hypothesis 
aiming at the larger public has started to grow: 
\begin{wrapfigure}{l}{4.1cm} \begin{center}
\includegraphics[width=4cm]{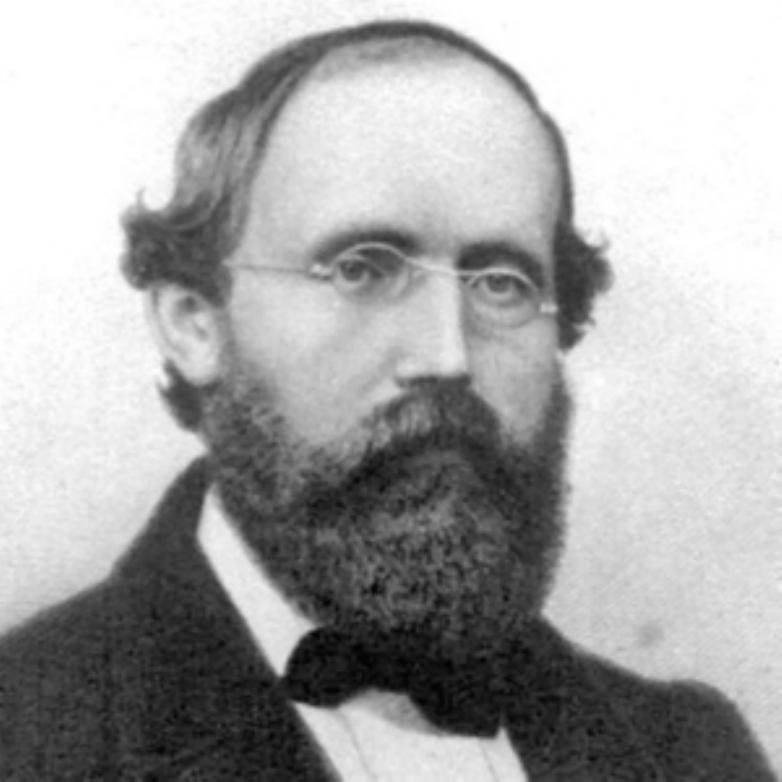}
\end{center} \end{wrapfigure}
\cite{Derbyshire,MusicPrimes,Sabbagh,Rockmore,Watkins2010,VeenCraats,MazurStein}.
Popular expositions have appeared even before the Millennium problems were offered.
Remarkable are the gorgeous paper \cite{Bombieri92} and \cite{Conrey}, which is a paper 
which won a prize for expository writing. Some of the books mention also Gaussian primes. 
What about the zeta function for Gaussian integers or Hurwitz integers? 
As the zeta function is an example of a Dedekind zeta function for function fields,
{\bf Riemann hypothesis for Gaussian integers} is part of the 
{\bf generalized Riemann hypothesis}. The story is quite similar than in the case of the
standard zeta function. \\

The {\bf Dedekind zeta function} of the Gaussian integers is
$$ \zeta_G(s)=\sum_n N(n)^{-s} \; , $$ 
where $n$ ranges over all nonzero Gaussian integers $n$ in the Euclidean 
domain $\mathbb{Z}[i]$ and where $N(a+ib)=a^2+b^2$
is the usual arithmetic norm. As $\mathbb{Z}[i]$ is a unique factorization domain, 
the {\bf Euler product formula} 
$\zeta_G(s)=4 \prod_{p \in Q} (1-N(p)^{-s})^{-1}$ still holds, where $Q$ is set
of Gaussian primes in a half open quadrant. It is the {\bf Euler golden rule} 
for $\zeta_G(s)$. Note that the product is over all Gaussian primes in a quadrant
only, similarly as for rational primes, where one only takes the product over positive
primes. Let $P_k$ denotes the set of rational primes which give reminder $k$ 
modulo $4$. The Gaussian zeta function relates with the {\bf Dirichlet beta function}
$\beta(s) = \prod_{P_1}(1-p^{-s}) \prod_{P_3} (1+p^{-s})$
which is also called the {\bf Dirichlet $L$-function} $L(s,\chi_{4,3})$ for the 
character $\chi_{4,3}$. \\
While the Riemann hypothesis for Gaussian primes is part of the 
{\bf generalized Riemann hypothesis}, the two statements for 
rational and Gaussian primes are known to be equivalent \cite{McPhedran2013}. 
Already Dirichlet knew:

\resultremark{(Dirichlet)
The Gauss Zeta function relates with the usual Zeta function 
by $\zeta_G(s) = 4 \zeta(s) \beta(s)$.}

Dirichlet knew such factorizations for quadratic number fields. We 
will look at it below also in the case of the Eisenstein zeta function which belongs
to an other ring of integers.  \\

The Riemann hypothesis appears first a bit stronger than the Riemann hypothesis because
$\zeta_G(s)$ involves also the {\bf Dirichlet beta function}
As the $\zeta$ function, also the $\beta$ function has both trivial as well as nontrivial
roots. But it is known that the non-trivial roots of $\beta(s)$ are on the critical line 
if the zeta function has the roots of the critical line \cite{McPhedran2013}. 
The Gaussian Riemann hypothesis is therefore 
equivalent to the standard Riemann hypothesis.  \\

The proof of the above formula for $\zeta$ and $\zeta_G$ appears in many places 
like \cite{Schlackow} page 16, exercise 4.12.h in \cite{VinceLittle} or 
\cite{GaussianZoo}). 

There are many parallels between $\zeta$ and $\zeta_G$: 
the {\bf Euler golden keys} are
$$  \beta(s) = \prod_{p \in P_1} (1-p^{-s})^{-1} 
               \prod_{p \in P_3} (1+p^{-s})^{-1}, 
    \zeta(s) = \prod_P (1-p^{-1})^{-1} \; , $$
where $P$ is the set of rational primes. The functional equations of the
two Dirichlet L-functions are very similar: 
\begin{eqnarray*}
\beta(1-s) &=&      2^s      \pi^{-s}   \sin(\pi s/2) \Gamma(s) \beta(s)  \\
\zeta(1-s) &=& 2\pi 2^{-s}   \pi^{-s}   \cos(\pi s/2) \Gamma(s) \zeta(s)   \; ,
\end{eqnarray*}
the {\bf functional equation for the Gaussian zeta function} is therefore even simpler
and given by 
$$ \zeta_G(1-s) = \sin(\pi s) \Gamma(s)^2 \pi^{-2s} \zeta_G(s) \; . $$
The {\bf reduced Gaussian zeta function}
$\xi_G(s) = \zeta_G(s) \Gamma(s)/\pi^s$ is now invariant under the evolution $s \to 1-s$
and $s(1-s) \xi_G(s)$ has an {\bf analytic continuation} to the entire plane.  \\

The equivalence of the hypothesis for $\mathbb{R}$ and $\mathbf{C}$ implies
statements about the growth rates of the corresponding {\bf Mertens function}. 
This function is defined in the same way as for rational integers:  \\

First of all, because $\mathbb{C}$ is a division algebra implying $N(n)$ to be multiplicative,
the formula $1/\zeta(s) = \sum_{n>0} \mu(n)/n^s$ for the usual Zeta function becomes now
$$ \frac{1}{\zeta_G(s)} = \sum_{z, N(z)>0} \frac{\mu(N(z))}{N(z)^s}  \; . $$
with the {\bf Gaussian M\"obius function} $\mu_G(s)$ which is 
$0$ if it contains a square prime factor and $1$ if $n$ has an even number of 
different prime factors and $(-1)$ if there is an odd number of different prime factors. 
But now, for real Gaussian integers of the form $z=4k+1$, we have $\mu(z)=1$. For
$z=2$, we have $\mu(2)=0$ as $2$ is composed of two primes which are conjugated. 
In the same way than the {\bf rational Mertens function} $M(n)=\sum_{k=1}^n \mu(k)$,
define the {\bf Gaussian Mertens function} $M_G(n) = \sum_{N(k) \leq n} \mu_G(k)$. 

\begin{figure}[!htpb]
\scalebox{0.2}{\includegraphics{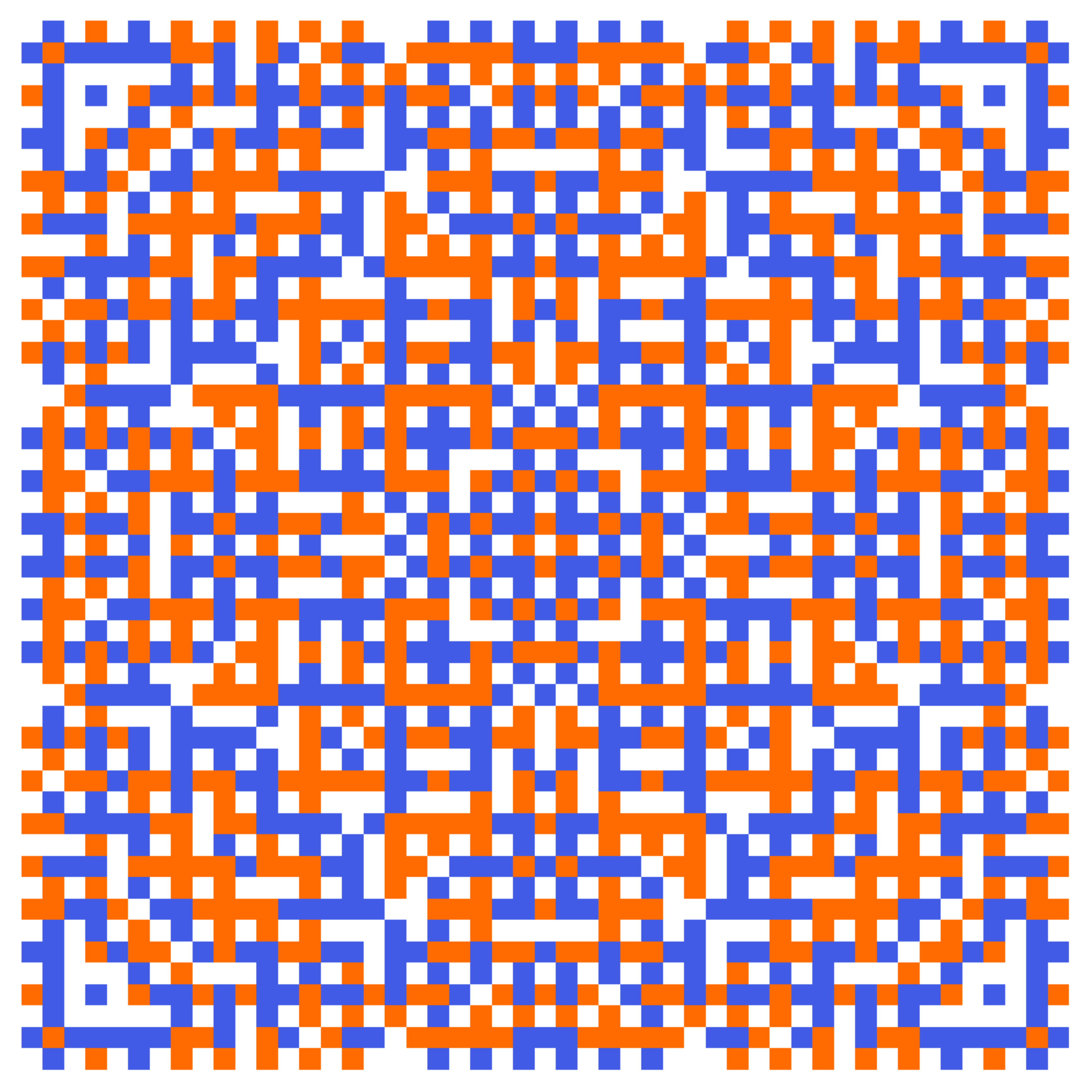}}
\caption{
The Moebius function for Gaussian integers is defined in the same way
as for rational integers.
As it is zero on even Gaussian integers $2a+2bi$ because $2 = (1+i)(1-i)$
and $2a+2bi= (1+i)^2 (1-i) (a+ib)$ contains a square. 
}
\end{figure}

The equivalence of the Riemann hypothesis means that also 
than $M_G(n) \leq n^{1/2+\epsilon}$ for every $\epsilon>0$ and large
enough $n$. 

\begin{figure}[!htpb]
\scalebox{0.12}{\includegraphics{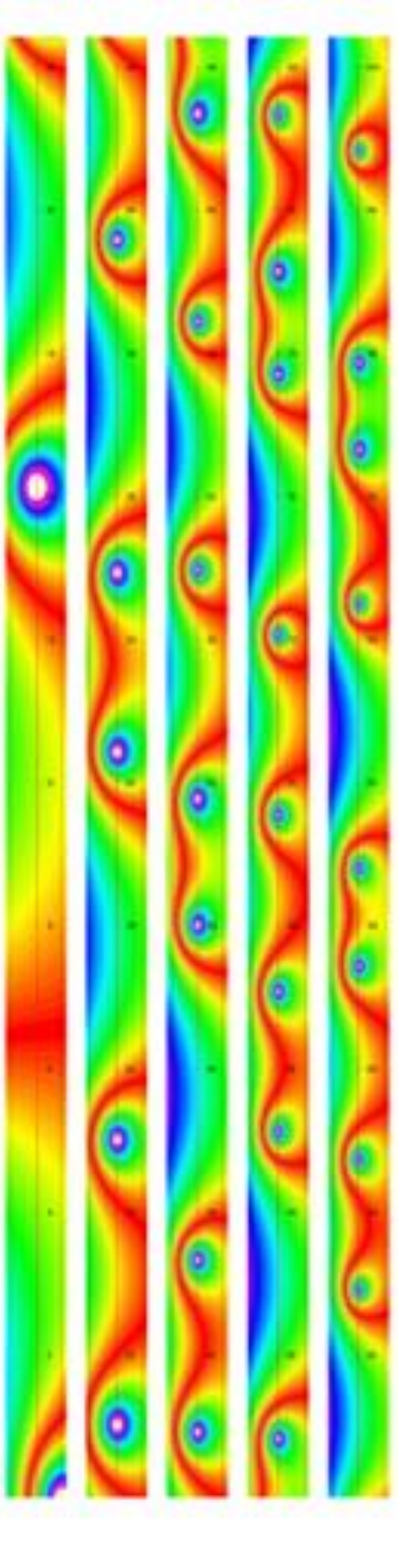}}
\caption{
A contour plot of the Riemann zeta function in the
{\bf critical strip} $z=a+ib$ with $0 \leq a \leq 1, 0 \leq b \leq 100$.
The nontrivial roots are all located on the line ${\rm Re}(z)=1/2$.
}
\end{figure}

\begin{figure}[!htpb]
\scalebox{0.12}{\includegraphics{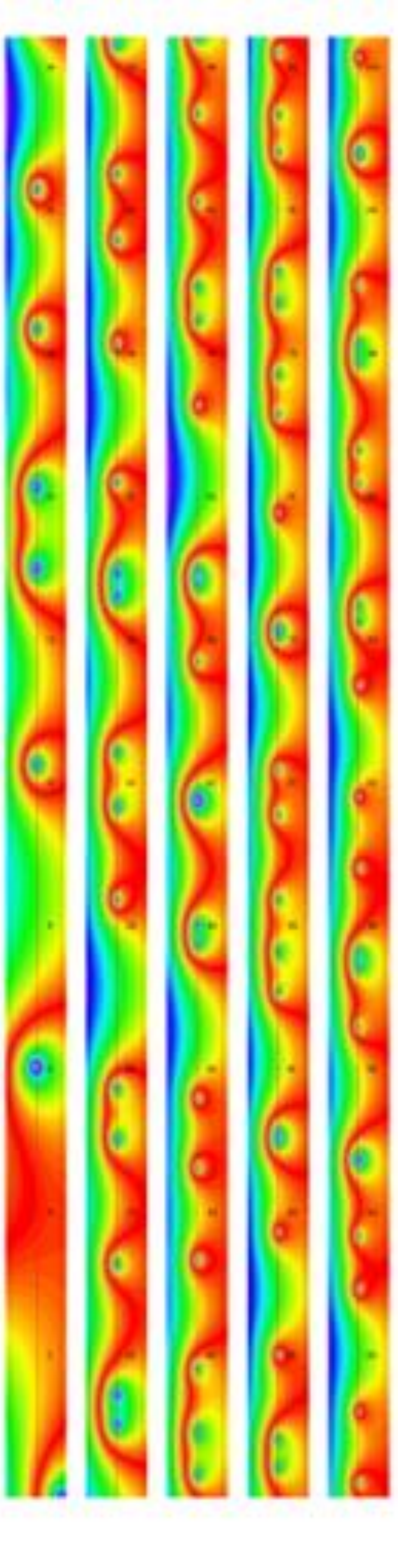}}
\caption{
A contour plot of the {\bf Dedekind Zeta function} $\zeta_C(z)$  of the
Gaussian primes. Its known that the all its roots are on the critical
line if the Riemann zeta function has this property. 
The zeta function is a product of the ordinary zeta function and the
Dirichlet Beta function. }
\end{figure}

\begin{figure}[!htpb]
\scalebox{0.2}{\includegraphics{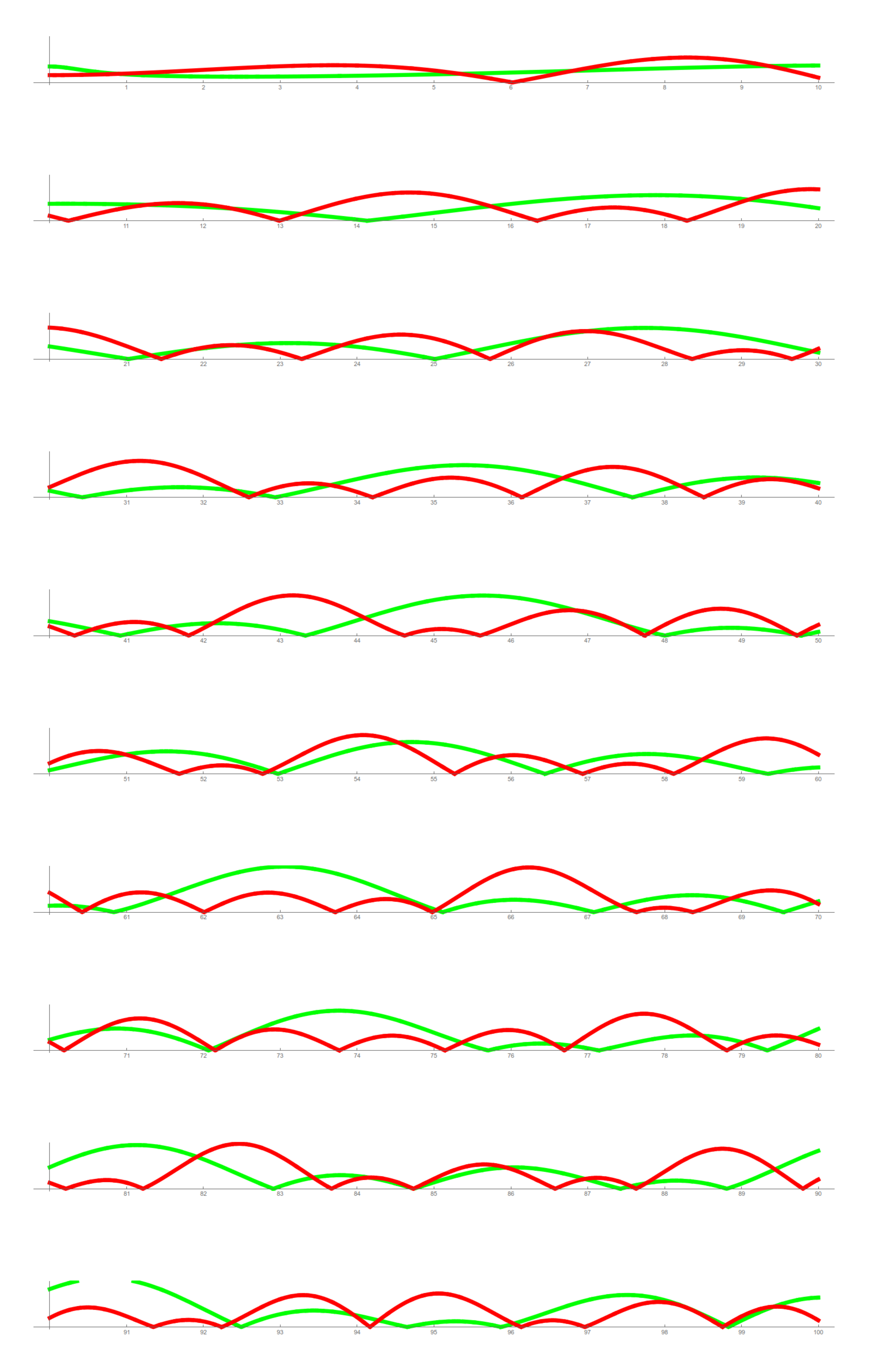}}
\caption{
The absolute value of the $\zeta$ and Dirichlet $\beta$ function
on the critical line ${\rm Re}(s)=1/2$ for $0 \leq {\rm Im}(s) \leq 100$.
The roots have been computed for example in \cite{Ossicini}.
}
\end{figure}

A highlight of Riemann's theory is the fact that the {\bf Chebyshev function} 
$$\psi(x) = \sum_{p^k<x} \log(p) = \sum_{n \leq x} \Lambda(n) $$
which logarithmically counts primes, satisfies the {\bf Riemann-Mangoldt formula}
$$  \psi(x) = x - \sum_w \frac{x^w}{w} - \log(2\pi) - \frac{1}{2} \log(1-x^{-2})  \; , $$
where $w$ runs over the non-trivial roots of $\zeta$, 
and where $\log(2 \pi)=\zeta'(0)/\zeta(0)$ comes from
the simple pole at $1$ and $\log(1-x^{-2})/2$ is the contribution of the trivial zeros
$-2,-4,-6,\dots$ and $x^w/w=e^{\log(x) a + \log(x) i b}/w$ is the contribution from the nontrivial zeros.
Pairing complex conjugated roots $w_j=a_k+ib_j=|w_j|e^{i \alpha_j},\overline{w}$ 
gives a sum of functions
$f_j(x) = e^{\log(x) a_j} 2 \cos( \log(x) b_j-\alpha_j)/|a_j+i b_j|$. This is the
major theme in books about the Riemann hypothesis for the more general audience like
\cite{VeenCraats}. The functions $f_j$ are the tunes of the {\bf music of the primes},
as they guide the distribution of primes \cite{Bombieri92,Sautoy,BorweinChoiRooneyWeirathmueller}.
One can hardly be too much excited about this formula: \cite{GranvilleSoundarajan} state
"What an unexpected and delightful identity". Lets look at the picture of these functions: 

\begin{figure}[!htpb]
\scalebox{0.2}{\includegraphics{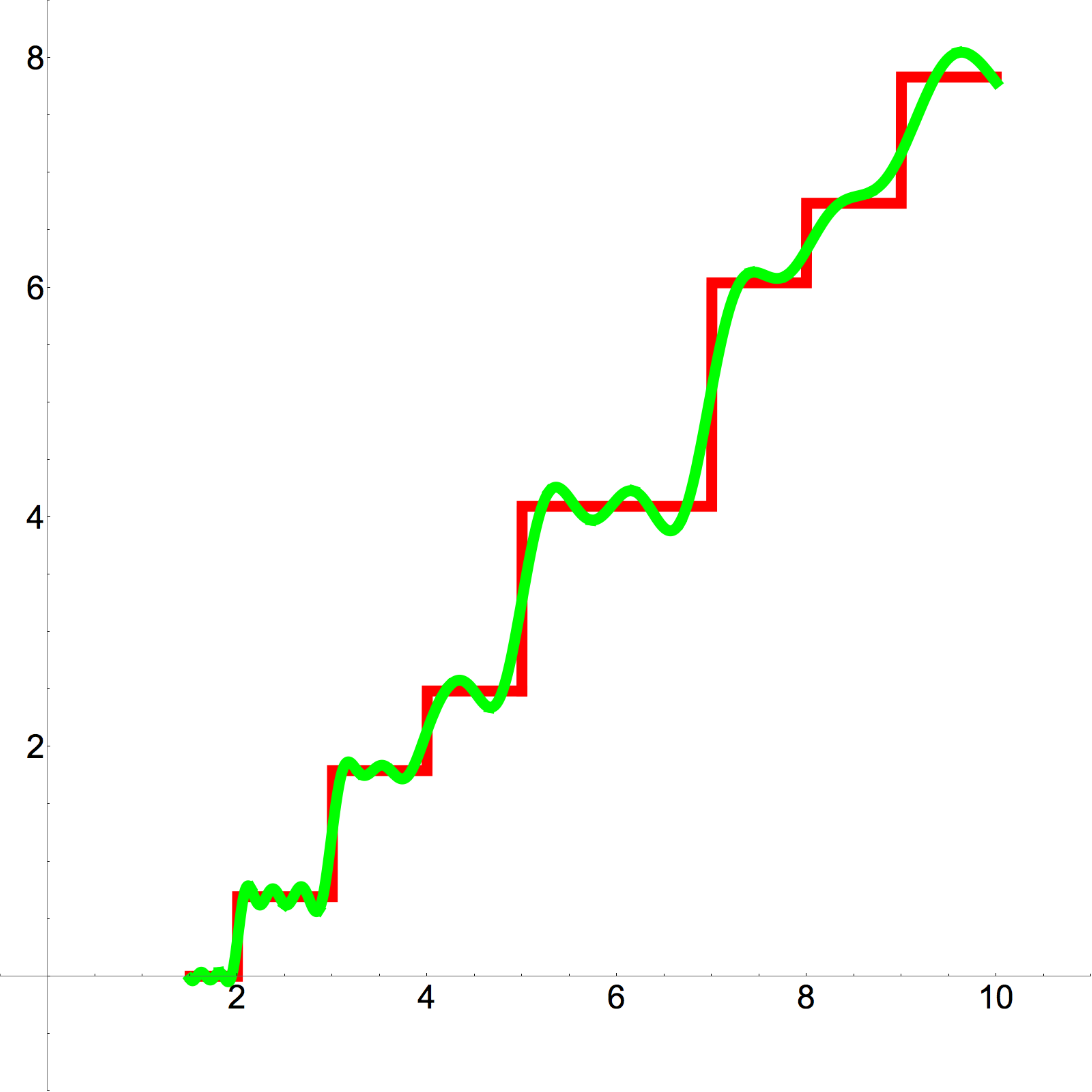}}
\caption{
The graph of the Chebyshev function $\psi$ with a finite Fourier approximation
given by the Riemann-Mangoldt formula. 
}
\label{figure1}
\end{figure}

For quaternions, for which unique prime factorization only holds modulo
meta-commutation, permutation and recombination,
one still has the golden key formula $\zeta(s) = \prod_p (1-1/N(p)^s)^{-1}$ 
where the product is over all Hurwitz primes in the fundamental region of $\mathbb{H}$. 
This leads to relations between primes and the roots of zeta. Both in the Gaussian
as well as the Hurwitz case, there are direct relations to the story in the rational prime case.
In the Hurwitz case, the zeta function is just shifted by $1$.
The upshot is that for rational, Gaussian or Hurwitz arithmetic, the Riemann hypothesis 
are equivalent. 

\section{Greatest common divisor matrices}

A rather unexpected relation between number theory and matrices appears 
for the {\bf greatest common divisor matrices} $A_{ij}(n,s)={\rm gcd}(i,j)^s$ 
introduced in 1876 by Henry John Smith  \cite{SmithGCD}
\cite{SmithGCD}. Smith was interested in various fields of 
mathematics. He also discovered the Cantor set well before Cantor \cite{Stewart89}. \\

The Smith matrices are finite matrices have entries $A(n,s)_{ij} = {\rm gcd}(i,j)^s$.  
For example, 
$$S(7,s) = \left[
                 \begin{array}{ccccccc}
                  1 & 1 & 1 & 1 & 1 & 1 & 1 \\
                  1 & 2^s & 1 & 2^s & 1 & 2^s & 1 \\
                  1 & 1 & 3^s & 1 & 1 & 3^s & 1 \\
                  1 & 2^s & 1 & 4^s & 1 & 2^s & 1 \\
                  1 & 1 & 1 & 1 & 5^s & 1 & 1 \\
                  1 & 2^s & 3^s & 2^s & 1 & 6^s & 1 \\
                  1 & 1 & 1 & 1 & 1 & 1 & 7^s \\
                 \end{array}
                 \right] \; . $$
The determinants is explicitly known in terms of Jordan functions:
$$ {\rm det}(A(n,s) = \prod_{m=1}^n  m^s \prod_{p|m} (1-p^{-s}) =\prod_{k=1}^n \phi^s(k) \; . $$
For $s=1$ already Smith has obtained 
${\rm det}(A(n,1) = \prod_{k=1}^n \phi(k)$, where $\phi$ is the Euler function
$\phi(n) = n \prod_{p|n} (1-1/p)$. \\

\begin{wrapfigure}{l}{4.1cm} \begin{center}
\includegraphics[width=4cm]{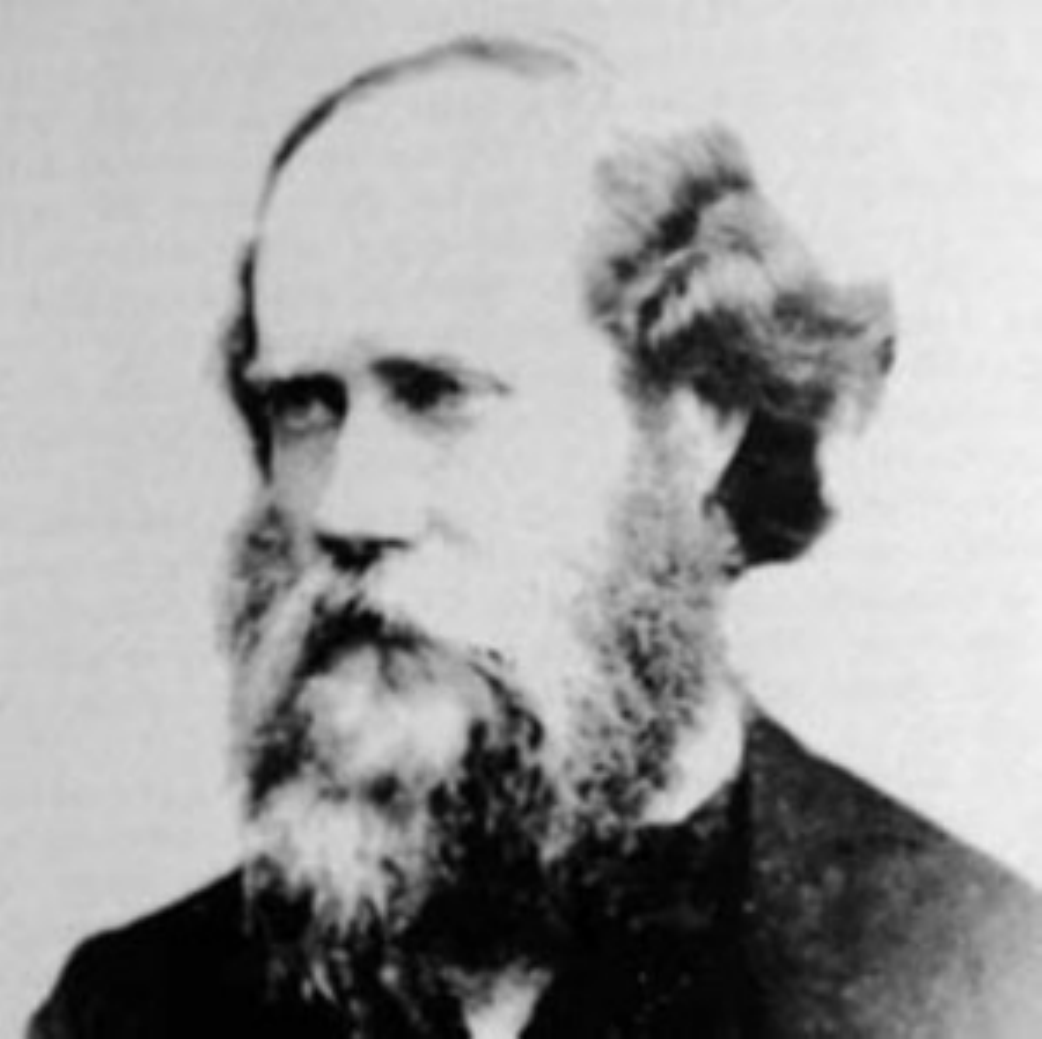}
\end{center} \end{wrapfigure}
More generally, Smith has looked at the matrices $A_{ij} = {\rm gcd}(x_i,x_j)^s$ 
defined for any finite set $S=\{x_1, \dots, x_n\}$ which is {\bf factor closed}, 
meaning that every factor of an element in $S$ must be in $S$. In that case, the determinant is 
$\prod_{k=1}^n \phi^s(x_k)$.  \\

Note that the roots of the function $h(s) = \det(A(n)^{-s})$ are all on 
the line ${\rm Re}(s)=0$ because it is a product of terms
$(1-1/p^s) = (1-\exp(-s \log(p)))$ which are zero exactly for ${\rm Re}(s)=0$ and
${\rm Im}(s) = 2\pi/\log(p)$. \\

\begin{figure}[!htpb]
\parbox{14.8cm}{
\scalebox{0.16}{\includegraphics{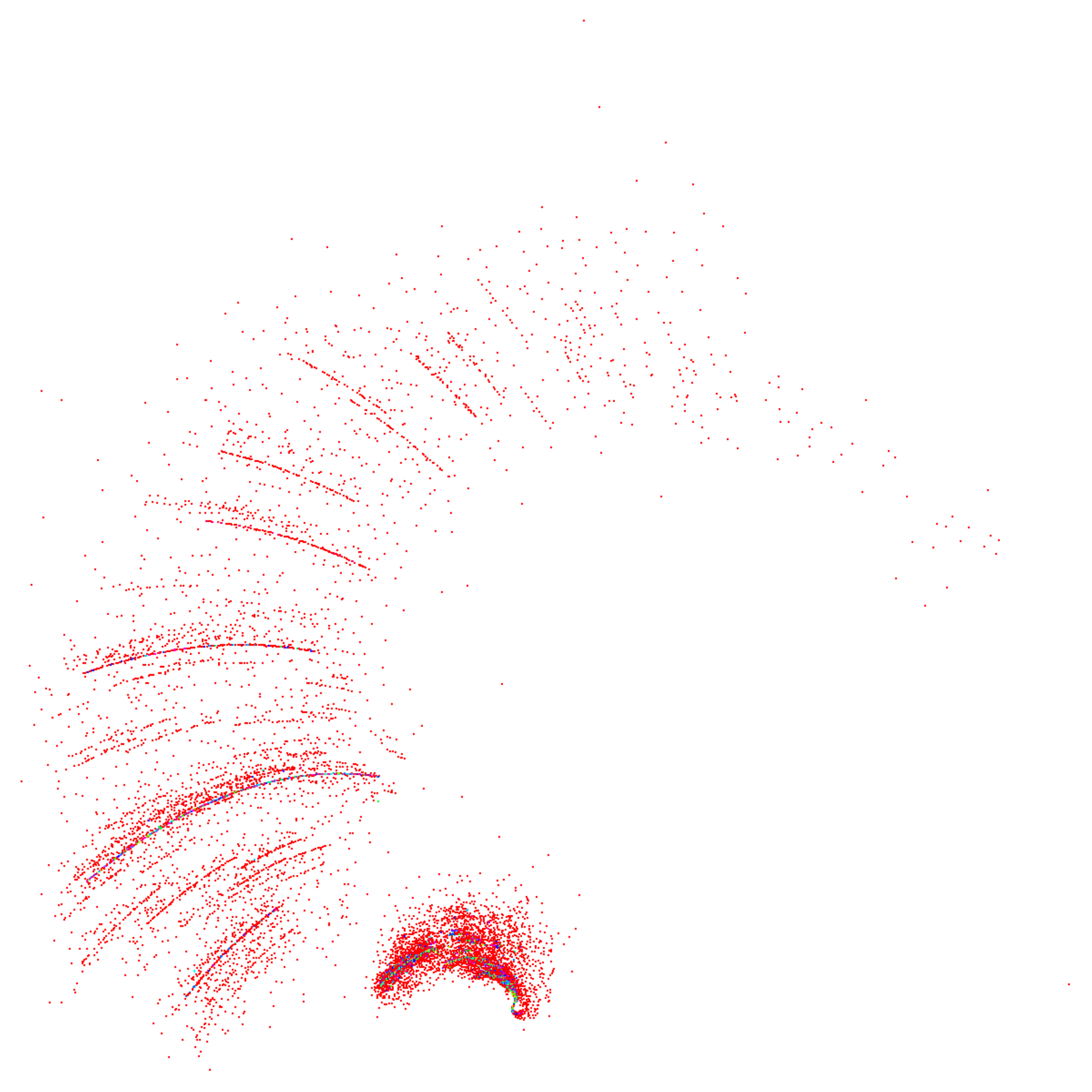}}
\scalebox{0.16}{\includegraphics{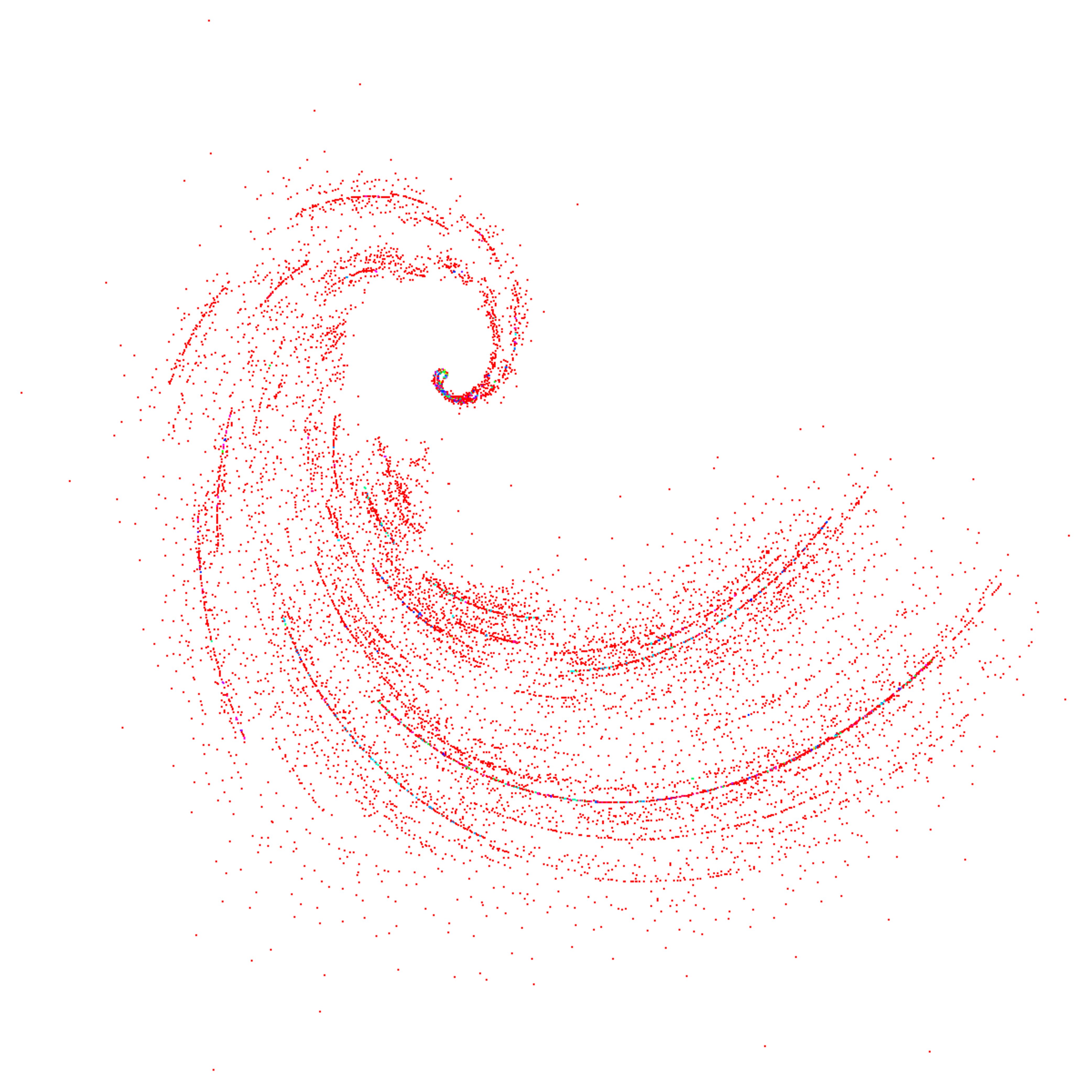}}
}
\parbox{14.8cm}{
\scalebox{0.16}{\includegraphics{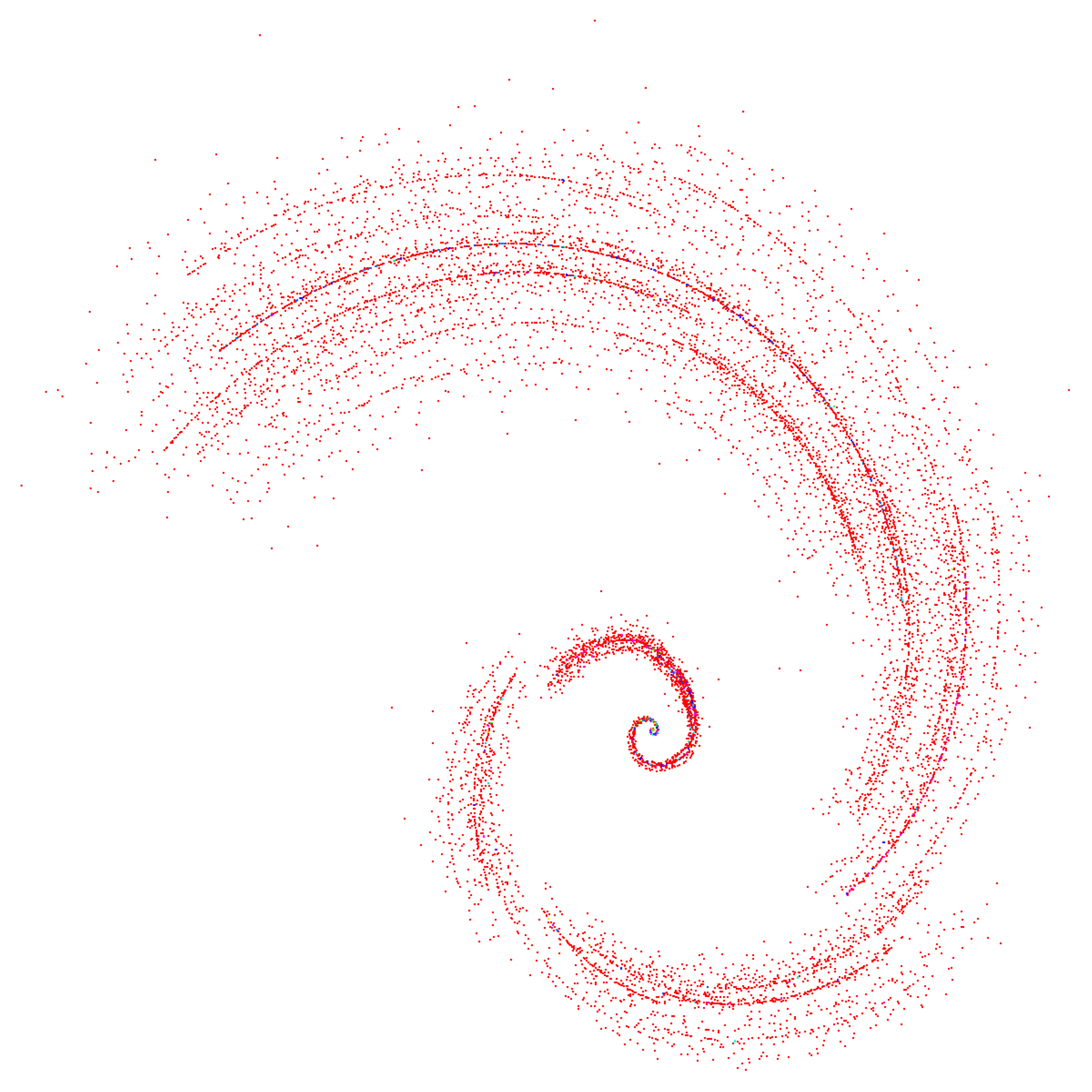}}
\scalebox{0.16}{\includegraphics{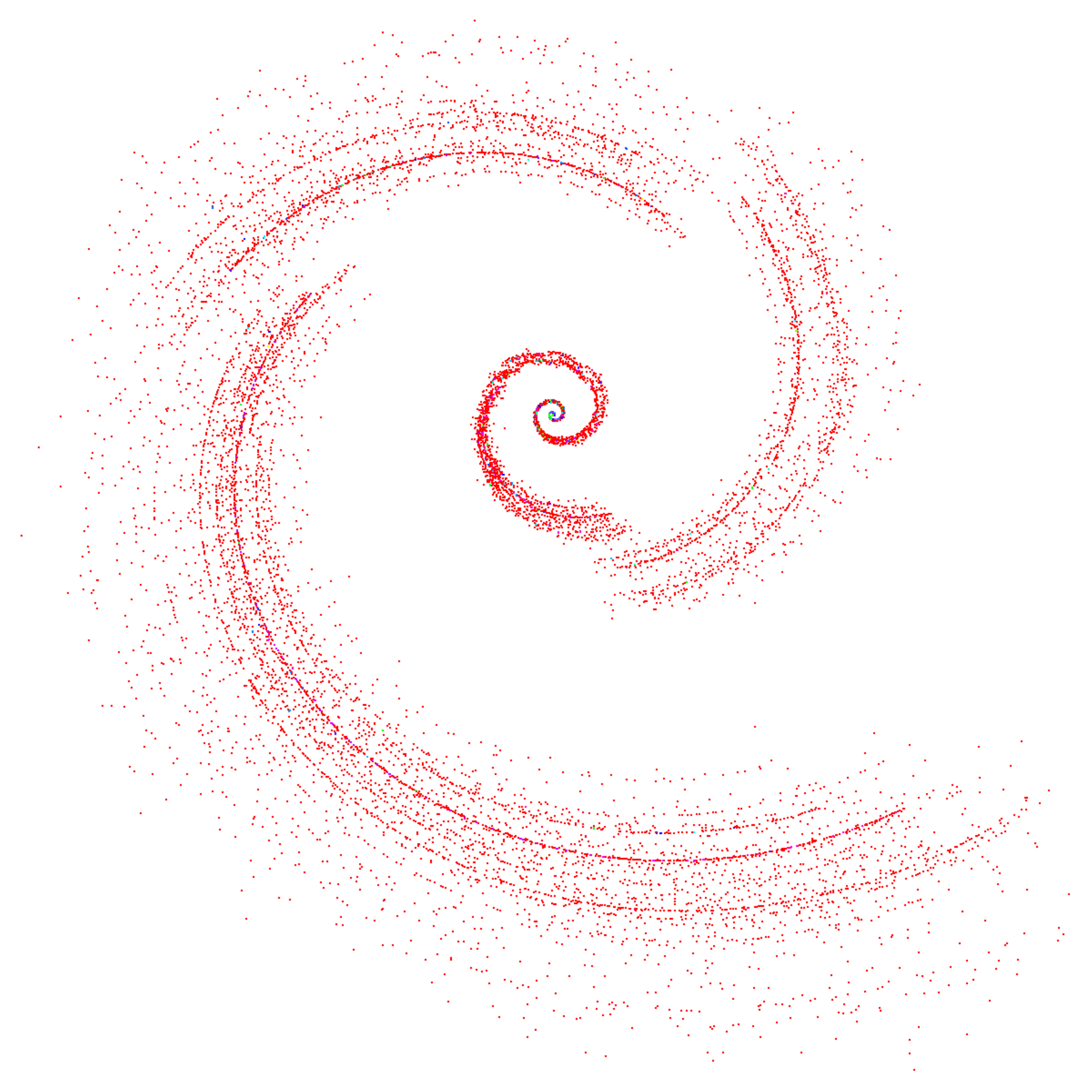}}
}
\caption{
The spectrum of the GCD matrix $A(n)^s$ with $A_{ij}(n)= {\rm gcd}(i,j)$
for the values $n=10'000$ for the parameter $s=1+ki$ for $k=1,2,3,4$. 
}
\end{figure}

We can also look at adjacency matrices of graphs $G(n)$ with vertex set $\{1, \dots, n\}$
for which two integers $k,l$ are connected if ${\rm gcd}(k,l)>1$.  Let us call them 
{\bf GCD graphs}. Since the rule making the connection is transitive, the graph
is topologically quite trivial and its Euler characteristic is the number of 
connected components. Because these connected components must consist of primes
larger than $n/2$ and smaller or equal to $n$, we have: 

\resultremark{
The Euler characteristic of $G(n)$ is $2$ plus the number of primes in 
the half open interval $(n/2,n]$.}

We look now at the vertex degrees. 
Look at a prime $p \leq n$. It is connected to $[n/p]-1$ other integers, where
$[t]$ is the floor function giving the largest integer smaller or equal to $t$. 
A product of two primes $pq$ is connected to $[n/p] + [n/q] - [n/(pq)]-1$ 
other integers. For example, in $G(30)$ the vertex $2$ has degree $[30/2]-1=14$,
the vertex $3$ has degree $[30/3]-1=9$ and the vertex $6$ has degree $19=14+9-4$. 
Now, also the vertex degree of a prime power $p^k$ is $[k/p]$ as it is connected
to the same integers than $2$. And the vertex degree of two prime powers $p^k q^l$
are all equal to $[n/p] + [n/q] - [n/(pq)]-1$. 
The vertex $p^n q^m r^k$ with three primes $p,q,r$ has degree
$[n/p]+[n/q]+[n/r]-f[n/(pq)]-f[n/(pr)]-f[n/(qr)]+f[n/(pqr)]-1$. 
By the {\bf Euler handshake lemma} the sum of the vertex degrees is twice to the
number of edges. This sequence is the sequence A185670 
and explicitly given as $n(n-1)/2 -\sum_{i=1}^n \phi(i)+1$. 

\resultremark{
The edge degree of $G(n)$ is $n(n-1)/2 - \sum_{k=2}^n \phi(k)$. 
}

This formula is attributed to Reinhard Zumkeller in \cite{A185670} and
as a formula counting the number of "non connections" then take this away from 
$n(n-1)/2$. It is amusing that for the determinants of the GCD matrices, 
we had the product of Euler totients $\det(A(n)) = \prod_{k=1}^n \phi(k)$. 
Now, the number of edges of of the GCD graphs $G(n)=(V(n),E(n))$ 
is related to the sum of Euler totients $\sum_{k=1}^n \phi(k)$. This function 
is called the {\bf totient summatory function} $\Phi(n)$, which grows
like $(3/\pi^2) n^2$. Coming back to the GCD matrices, the Euler summatory totient 
function is just the sum of the entries of the last column of the matrix $A(n)$.  \\

For the Kirchhoff Laplacian of the graph $G(n)$, the trace is the sum of the 
vertex degrees. We see that the Euler summatory totient function has a spectral 
interpretation of the Laplacians of a sequence of graphs. 

\begin{figure}[!htpb]
\scalebox{1.0}{\includegraphics{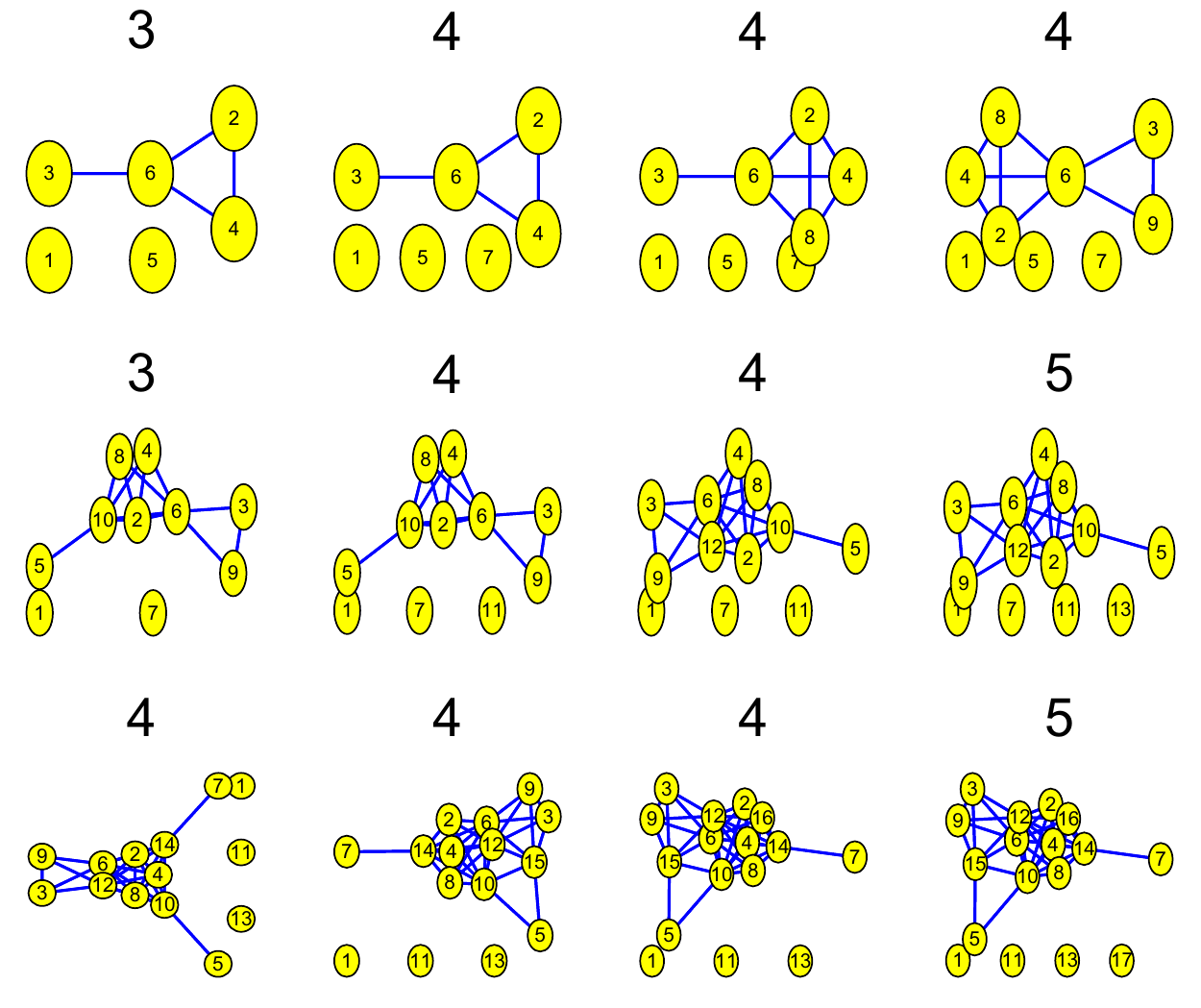}}
\caption{
We see the $G(n)$ for $n=4, \dots 19$. The vertex set of $G(n)$
is $\{ 1,\dots, n \}$ and two vertices $k,l$ in $G(n)$ are connected if $gcd(k,l)>1$.
The label above a graph gives its Euler characteristic $\sum_k v_k(G)$ where $v_i(G)$
counts the number of complete subgraphs $K_k$ of $G$. 
}
\label{figure1}
\end{figure}

{\bf Remarks}: \\
{\bf 1)} It is fitting that similar formulas appear in the context of quaternions. 
\cite{Hurwitz1919} proves that if $n$ is an integer, then the number of 
quaternions $z$ congruent to $1$ modulo $n$ is 
$$ n^3 \prod_{p|n} (1-\frac{1}{p^2}) = n \phi^2(n) \; . $$

{\bf 2)} Much more is now known about GCD matrices. One has formulas for the inverse
\cite{BourqueLigh} which uses an explicit formula $E L(n,s) E^T$, with the unimodular
lower triangular matrix $E_{ij}$ which is $1$ if $i$ divides $j$ and $0$ otherwise. 
Since $E^{-1}_{ij} = \mu(i/j)$ if $i$ divides $j$ and $0$ else, the inverse is
explicitly given as $(E^{-1})^T L(n,s)^{-1} E^{-1}$. Note that since $E$ is not orthogonal,
the explicit formula does not provide a diagonalization. The structure of the eigenvalues
of $S(n,s)$ is therefore not so clear. \\

{\bf 3)} Some properties of these determinants were rediscovered by 
Juan Jose Alba Gonzalez who communicated
it to Omar Antolin, who showed it to me. The matrices became a fixture in our 21b 
linear algebra course of spring 2015, as they appeared in homework, project and exams.
We have asked students to look at the eigenvalue structure of these matrices. 
While working on that Mathematica project, two of our students, Isabelle Steinhaus and
Jerry Nelluvelili discovered the spiral for the parameter $s = 1+4i$.
Similarly as Julia sets in complex dynamics, the spectral pictures are
parameterized by a complex parameter. We would like to know for example whether for some $s$
and $n \to \infty$, the spectrum converges to some fractal set.
We know from the Smith formula that $\det(A(n)^s/n^s$ is a product of factors $(1-1/p^s)$,
where $p$ runs over all primes dividing $n$.
In some sense, the determinants of the GCD matrices lurch up to the zeta function if $n$ is
the product of the first $k$ primes and $\zeta(s)$ is a limit of determinants.

\section{Cellular automata}    

An other famous problem in additive number theory is the {\bf twin
prime conjecture}. It was first formulated by Polignac in 1849 \cite{Pintz} and mentioned
by Kronecker in 1901 and Maillet in 1905.  \\

Cellular automata were first introduced by Gustav Hedlund, who not only made early contributions
to the calculus of variations and ergodic theory but also worked in
symbolic and topological dynamics. The Hedlund-Curtis-Lyndon theorem assures that any continuous
map $T$ on a product space $A^L$ of a compact topological space $A$, 
where $L$ is a lattice such that $T$ commutes with all translations, must be given by a 
{\bf local rule}. In other words, such a topological dynamical system must be a cellular automaton. \\

One can see the set of Gaussian integers as an element in the linear
space $X$ of all functions from the two dimensional lattice $Z^2$ to the
field $Z_2$. When seen like this, the set $X$ is a configuration over the alphabet
\begin{wrapfigure}{l}{4.1cm} \begin{center}
\includegraphics[width=4cm]{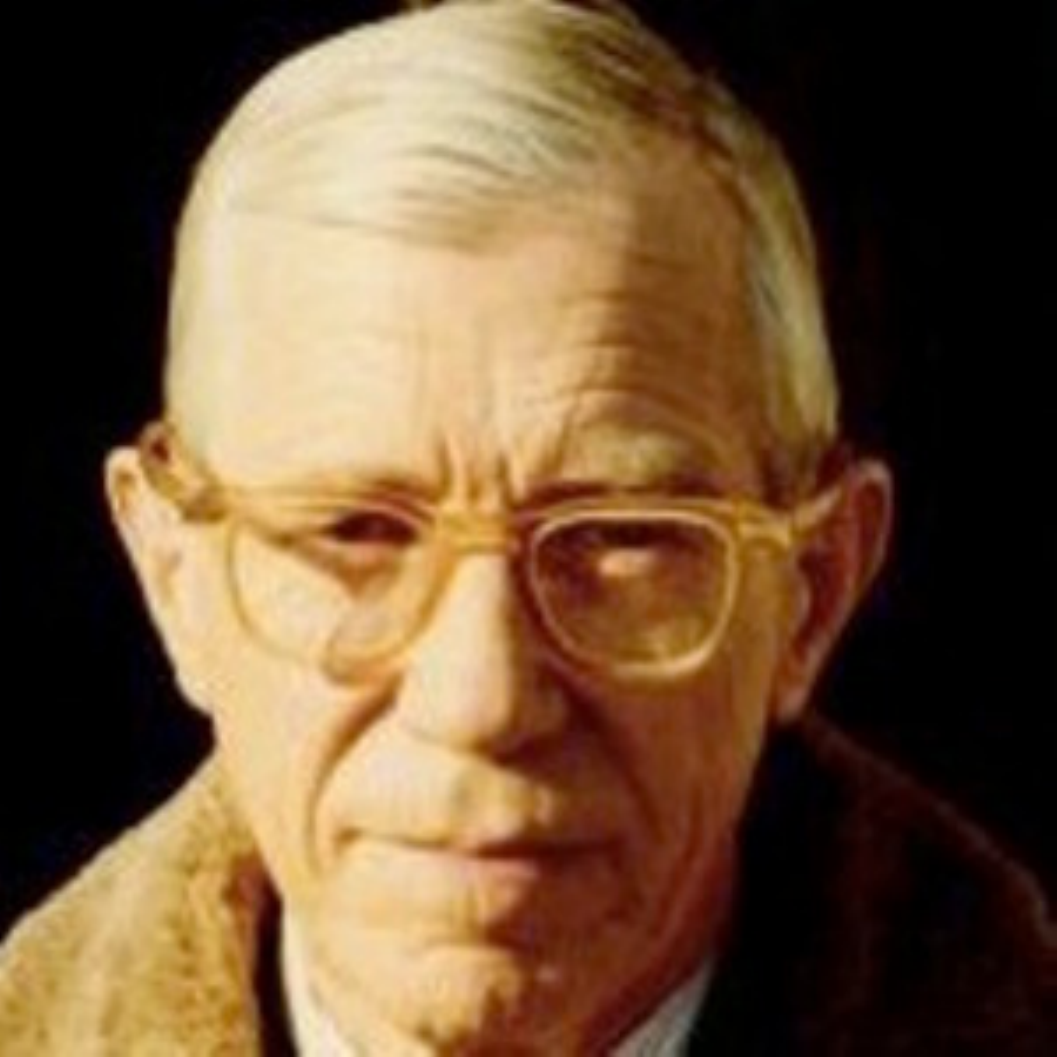}
\end{center} \end{wrapfigure}
$A=\{0,1\} = Z_2$ or as a compact metric space $X=A^{Z^2}$ with the product topology. 
Gustav Hedlund was the first to consider continuous maps from $X$ to $X$ which commute
with all translations. They are now called {\bf cellular automata}. The most famous by far
is the {\bf game of life}. What happens if we apply this to the Gaussian prime initial condition $x$? 
Let us call a region inside a configuration $x$ "alive" for a cellular automaton $T$, if it 
is not a fixed point of $T$. It is alive, if it moves. One can conjecture that "the region which
is alive is infinite". The prime twin conjecture would imply that.  \\

Hardy and Littlewood have precise predictions about the number of classical prime twins.
What about the complex cases?  The {\bf twin prime conjecture} for Gaussian primes
asks for the existence of infinitely many {\bf Gaussian primes twins}, pairs of primes for
which the Euclidean distance is $\sqrt{2}$. \cite{HolbenJordan}.
While one does not know whether infinitely many Gaussian prime twins exist, one
can estimate that there are asymptotically $C r/\log^2(r)$ of them in a ball of
radius $r$ \cite{GaussianZoo}. Building the graph with the set $P$ of Gaussian primes as 
vertex set, where two are connected if their distance is $\sqrt{2}$, 
the prime twin conjecture whether infinitely many components
of length $2$ exist.  \cite{HolbenJordan} also conjectured infinitely many connectivity components
 of size $3$ and $4$ and point out that there are only finitely many quintuplets. 
The {\bf prime gap problem} is formulated analogue in the complex:
define $g_n$ as the radius of the largest punctured disc without primes in $|z| \leq n$. For
rational primes, finite bounds for $\liminf g_n$ are known and since Zhang's proof the bound
has decreased to 246. This means that there are infinitely many prime pairs of distance 246.
\cite{Polymath2014}. The twin gap problem could also be asked for Gaussian primes and smaller and
smaller bounds searched until the twin prime problem is solved.  \\

\begin{figure}[!htpb]
\scalebox{0.5}{\includegraphics{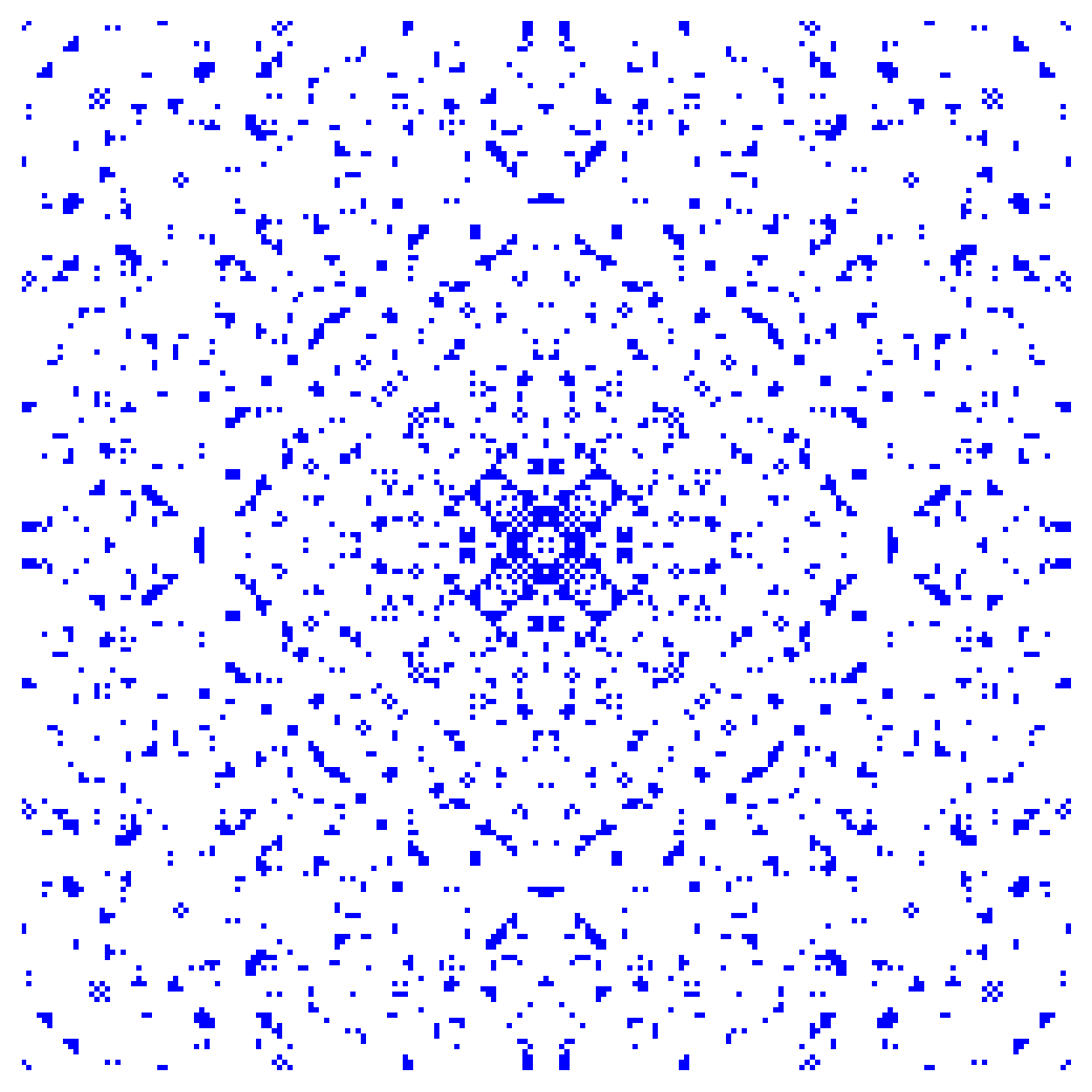}}
\scalebox{0.5}{\includegraphics{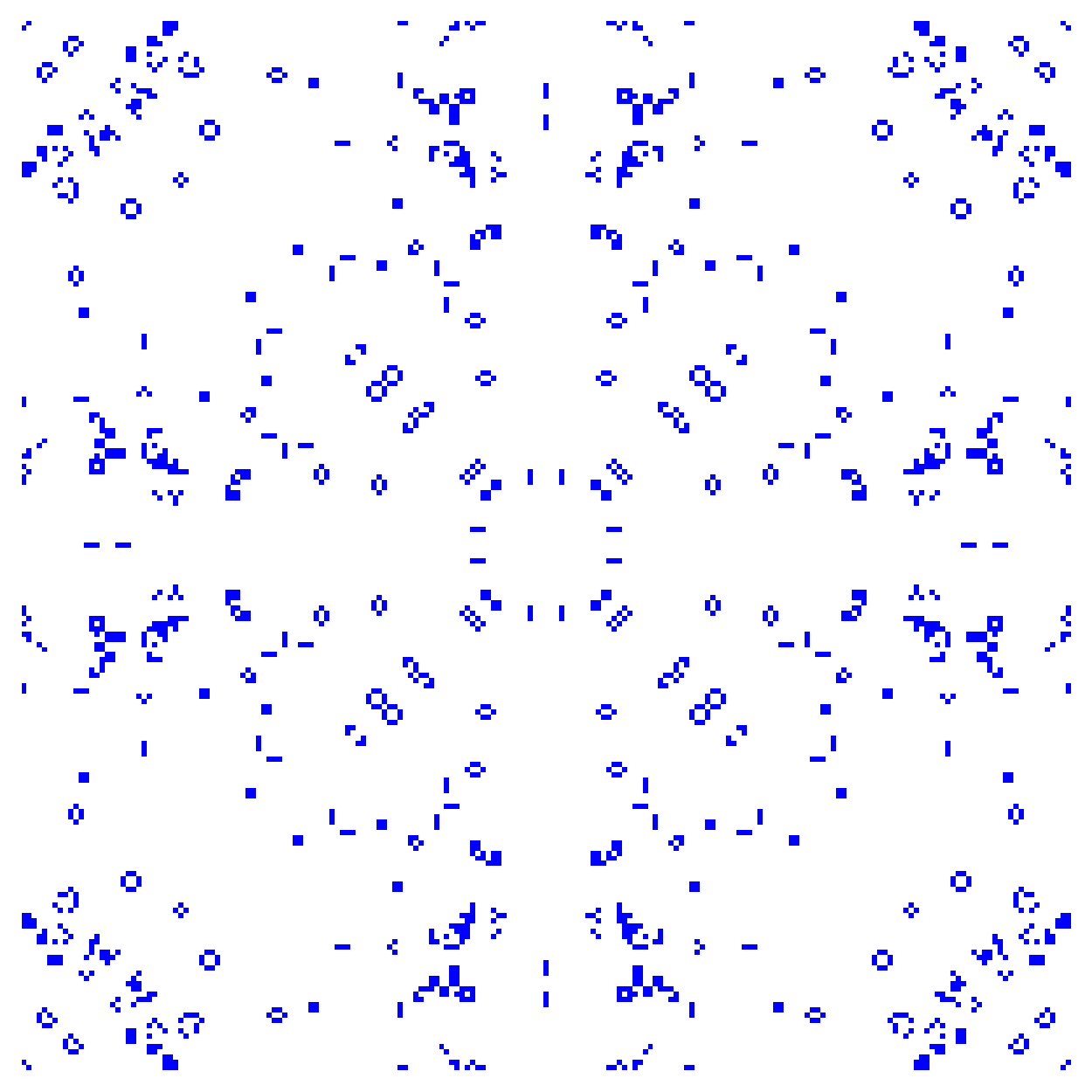}}
\caption{
When seeing the Gaussian primes as a configuration in $A^{Z^2}$,
where $A=\{0,1\}$ is the alphabet, we can apply cellular automata
on it.
The first picture shows the prime configuration after applying the Game of life
once, the second after three times. If the prime twin conjecture holds,
then there is life arbitrary far away from the origin. ``Life" in a
region is a configuration which moves when applying the time map.
\label{life}
}
\end{figure}

One can also use cellular automata to illustrate the moat problem. Just apply a map which gives
$1$ if there is a neighboring cell alive and 0 else. After applying this a couple of times,
we can look for connected components. The computed moats of course depend on the cellular
automaton used but its clear that if the moat conjecture is true, then also after applying
such a CA a couple of times, the connected components are still bounded. Indeed,
for any cellular automaton, and every configuration $x$: if arbitrarily large moats exist
for $x$, then for any finite $m$, arbitrarily large moats exist for $T^m x$. 
Now, like with additions of complex numbers, also the application of cellular automata
maps is not compatible with the multiplicative structure of the complex numbers.

\begin{figure}[!htpb]
\scalebox{0.2}{\includegraphics{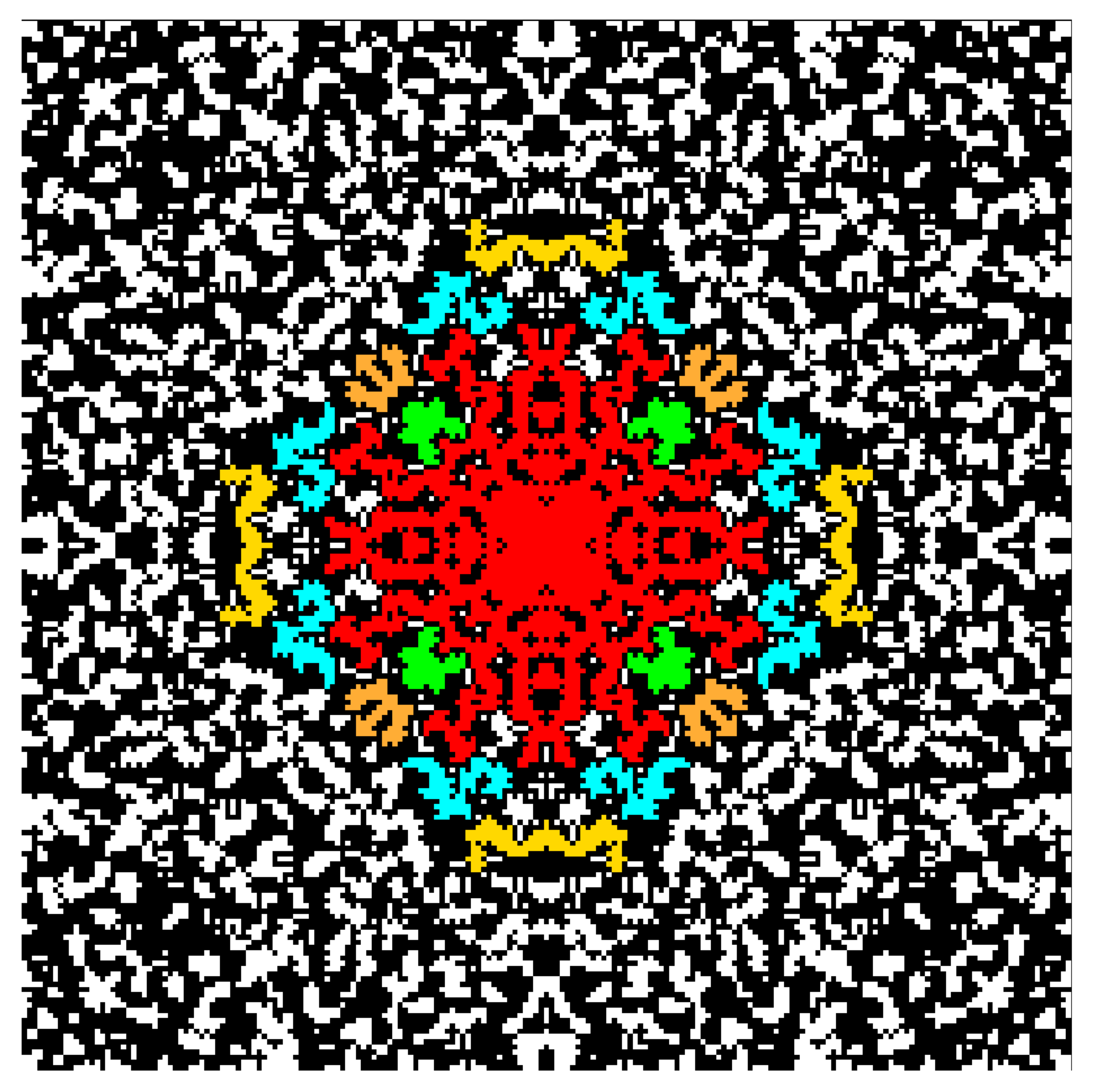}}
\scalebox{0.2}{\includegraphics{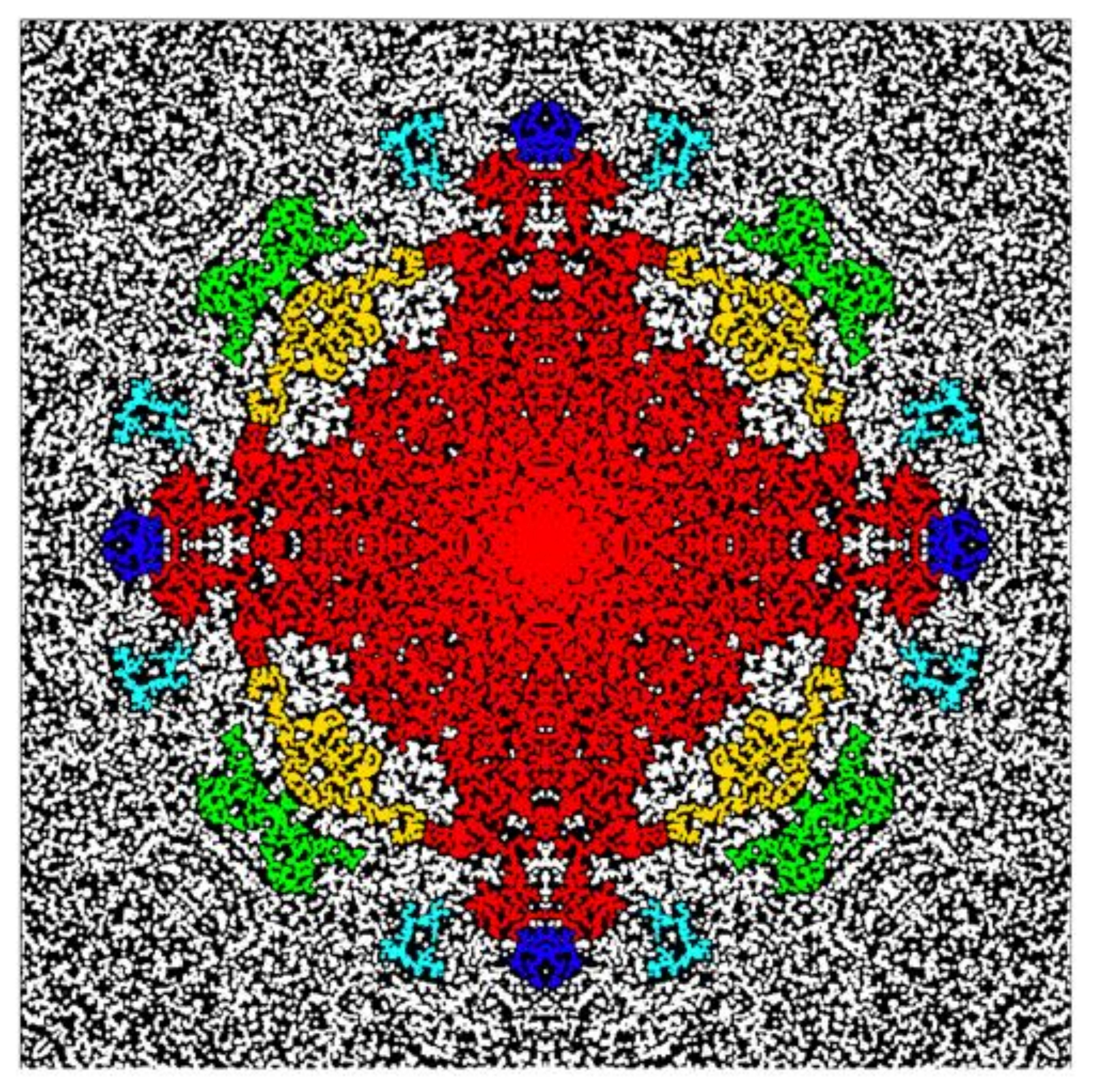}}
\caption{
Illustration of the moat problem: we applied a cellular automaton
map and looked at the central connected component. The first
picture is after applying the map once, the second after applying
the map two times. In principle, one can compute like that any
moat. The problem is that it is inefficient. The best known moats
have also been computed using computer assistance but are much
more sophisticated.
\label{moat}
}
\end{figure}

So, here is the conjecture about ``Prime life arbitrary far, far away": 

\conjecture{When applying the game of life cellular automaton to the
Gaussian integers, there is motion arbitrary far away from the origin.  }

If the conjecture is false, there are only finitely many living creatures 
in the Gaussian prime  and that would be rather sad. 
The problem looks hard. It would follow from the prime twin conjecture for
Gaussian primes as it would produce {\bf blinkers} which are time periodic
configurations.  \\

{\bf Remarks.} \\
{\bf 1)} The {\bf Gaussian moat problem} is covered in
\cite{JordanRabung,GethnerStark,GethnerWagonWick,Wagon}.
Possible Gaussian moat paths to infinity have some argument restrictions
\cite{Loh}. \\

{\bf 2)} Percolation problems for Gaussian primes have been compared before with the case of random matrices
\cite{Vardi}, for which the distribution $\pi(x) \sim {\rm Li}(x)$ matches the Gaussian prime distribution.
The spectral situation for the Gaussian pseudo random case appears similar as the random case.
The random case can be seen as generalized primes leading to {\bf Beuerling zeta functions}. \\

\section{Almost periodic matrices} 

The last topic is related to an other passion of Hardy and Littlewood, as they  studied 
almost periodic and especially quasi-periodic structures related to number theory. \\

These structures have been studied in solid state physics in the form of {\bf quasi crystals},
which are special almost periodic structures. An example of an extensively studied model is the 
\begin{wrapfigure}{l}{4.1cm} \begin{center}
\includegraphics[width=4cm]{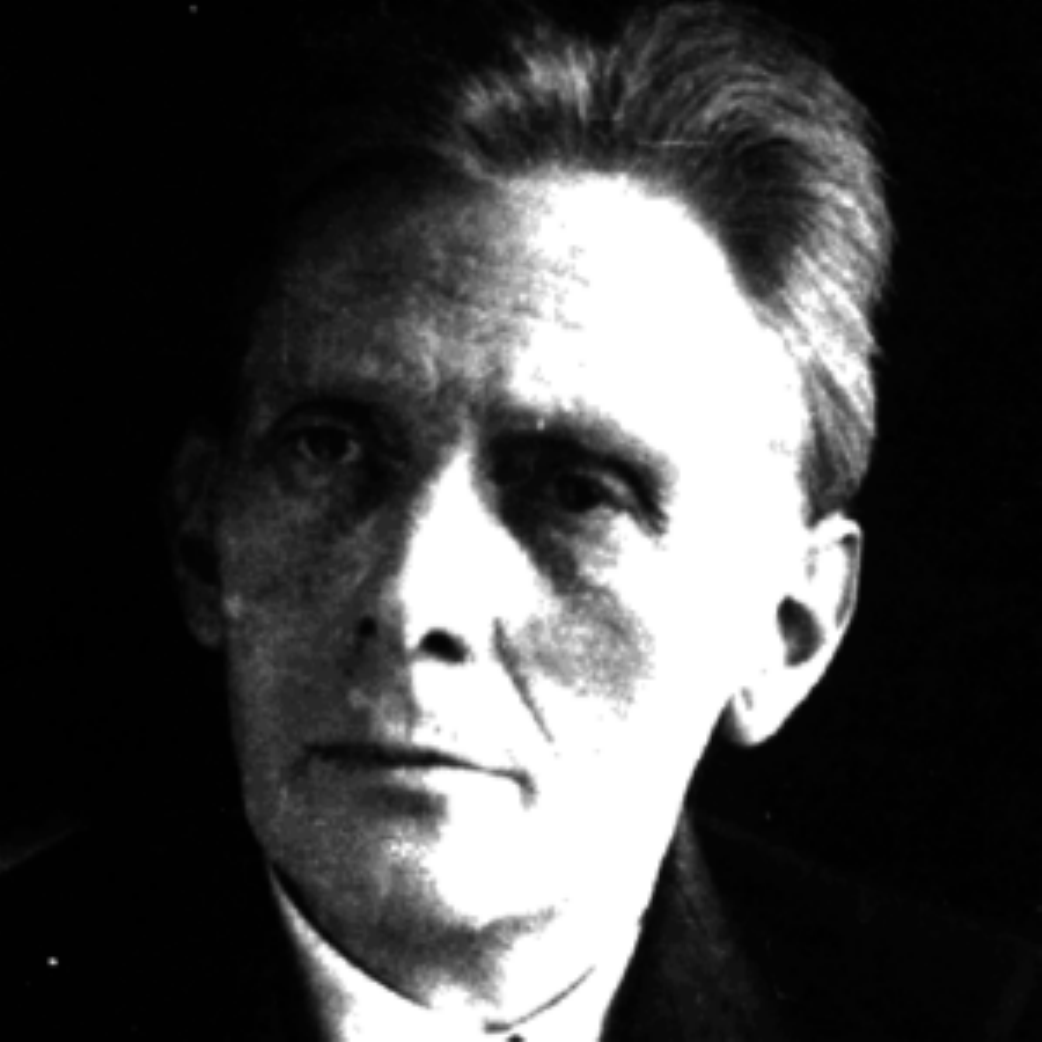}
\end{center} \end{wrapfigure}
{\bf almost Mathieu operator} $(Lu)_n=u_{n+1}+u_{n-1}+V_n u_n$ with $V_n=c \cos(n \alpha)$,
where the almost periodicity is quasi-periodic. 
In its analysis number theoretical Diophantine properties relate to spectral
properties of these Jacobi matrices. When plotting the picture of spectra for all $\alpha$ 
in the interval $0$ to $\pi$ one sees the Hofstadter butterfly. \\

We look at almost periodic matrices defined by irrational rotations defined by 
$$   A_{km} = \cos(k m \alpha  + m \beta) \; , $$
where $\alpha,\beta$ are irrational numbers. We took $\alpha=\sqrt{p}, \beta=\sqrt{q}$, where
$p,q$ are primes. The spectrum of these non-selfadjoint $n \times n$ matrices often feature
an unusual spectral structure. We call them {\bf snowflake spectra}. \\

The non-selfadjoint {\bf van der Monde matrices}
$$   B_{km}=\exp(i(k m \alpha + m \beta)) = (z^k w)^m  $$
with $z=e^{i \alpha}, w=e^{i \beta}$ also produce interesting spectra in the complex plane
but they look more like unions of curves. We have $A = {\rm Re}(B)$ by Euler's formula. \\

Since the matrices $A$ are real, the spectrum of $A$ is symmetric with
respect to complex conjugation. The spectrum of $A(n)$ is however {\bf not} symmetric with respect to 
$\lambda \to -\lambda$ or $\lambda \to i \lambda$ but in the large, this symmetry emerges.
One can always see these model in a probabilistic setup and study the matrices
$$   A_{km}(\theta) = \cos(\theta + k m \alpha  + m \beta) \; , $$
experimentally. Here is a trivial upper bound: 

\resultremark{
The spectral norm of $A$ is bounded above by $n$.
}
Proof; The spectral norm $||A||_2$ is the square root of the maximal eigenvalue 
of $A^T A$. The norm $||A||_1$ is the maximal $l_1$ norm of its columns
and $||A||_{\infty}$ is the maximal $l_1$ norm of the rows. 
One has $||A||_2^2 \leq ||A||_1 ||A||_{\infty} \leq n^2$ so that 
$||A||_2 \leq n$. \\
Asymptotically, we expect the spectral radius of $A$ to be bounded by $\sqrt{n} +R(n)$
where $R(n)$ goes to zero for $n \to \infty$ at least if $\alpha,\beta$ are Diophantine
like square roots of primes. 
Since every column $w_k$ has Euclidean length $\leq \sqrt{n}$,
for every vector $v$ of length $1$, we have $Av = \sum_k v_k \vec{w}_k$
with length $\leq (\sum_k v_k^2 \sqrt{n}) + \sum_{l,k} v_k v_l \vec{w}_k \vec{w}_l
\leq \sqrt{n} + R$ where $R$ is expected to zero since the different columns 
become more and more orthogonal. 

\resultremark{
The matrix $A$ is the real part of a Van der Monde matrix for the numbers $a_k = z^k w$.
}

We know therefore the determinant 
$\prod_{k<l} (a_l-a_k) = \prod_{k<l} (w (z^l-z^k))$
which is $w^{n(n-1)/2} \prod_{r=1}^{n-1} f_r(z)$ with $f_r(z)=\prod_{m=1}^{r-1} (1-z^m)$. 
If we look at the absolute value and take the logarithm, this is a sum of 
Birkhoff sums which can be estimated from above by $\log(n!)$, as long as
$\alpha$ is strongly Diophantine like the golden mean. 

In the case $z=e^{i \alpha}$ and $w=e^{i \beta}=1$, the complex matrices are
$$  B_{km} = e^{i k m \alpha} \;  $$
for which the determinant is $\prod_{k<l} (z^l-z^k)$. Taking
absolute values $\prod_{k,l,k<l} |1-z^{k-l}|$ which is $f_1(z) \dots f_{n-1}(z)$
if $f_n(z) = \prod_{k=1}^n (1-z^k)$. We know that if $\alpha$ is the 
{\bf golden ratio}, then $|f_n(z)| \leq \log(n)$. The proof of the later uses that
$S_k(\alpha) = \sum_{j=1}^k g(j \alpha)$ 
$G(x)=\log(2-2 \cos(2\pi x)) =2 \log|2 \sin(\pi x)|=2\log|1-e^{2 \pi i x}|$. 
Now, the Birkhoff sum of its derivative $G'(x) = \cot(\pi x)$ has a universal 
property.  \\

So, we know at least how to handle the determinant in the complex case for special  
irrational rotations. 


\resultremark{
For $\alpha$ the golden ratio and $\beta=0$, the absolute value of the 
determinant of $B$ is bounded above by $\log(2^n n!)$. 
}

We see numerically a very similar behavior for the real case $A$

\conjecture{
For $\alpha$ the golden ratio and $\beta=0$, the absolute value of the
determinant of $A$ is bounded above by $\log(2^n n!)$.
}

If we animate the spectra of $A$ by changing the rotation numbers $\alpha,\beta$
we see the eigenvalues move around. One can observe {\bf eigenvalue repulsion}. 
Let us illustrate this in the simplest case $n=2$, where we have the $2 \times 2$ matrix $A$
given as
$$   A= \left[
                 \begin{array}{cc}
                  1 & 1 \\
                  \cos (\beta ) & \cos (\alpha +\beta ) \\
                 \end{array} \right] 
      =  \left[
                  \begin{array}{cc}
                   1 & 1 \\
                   a & b \\
                  \end{array}
                  \right] \; . $$
It has the characteristic polynomial $x^2 - (1+b) x + (b-a)$. By changing
$\alpha$ and $\beta$ we can achieve that the eigenvalues approach each other
on the real axes, then bounce off and scatter away in the complex plane. \\

\begin{figure}
\scalebox{0.20}{\includegraphics{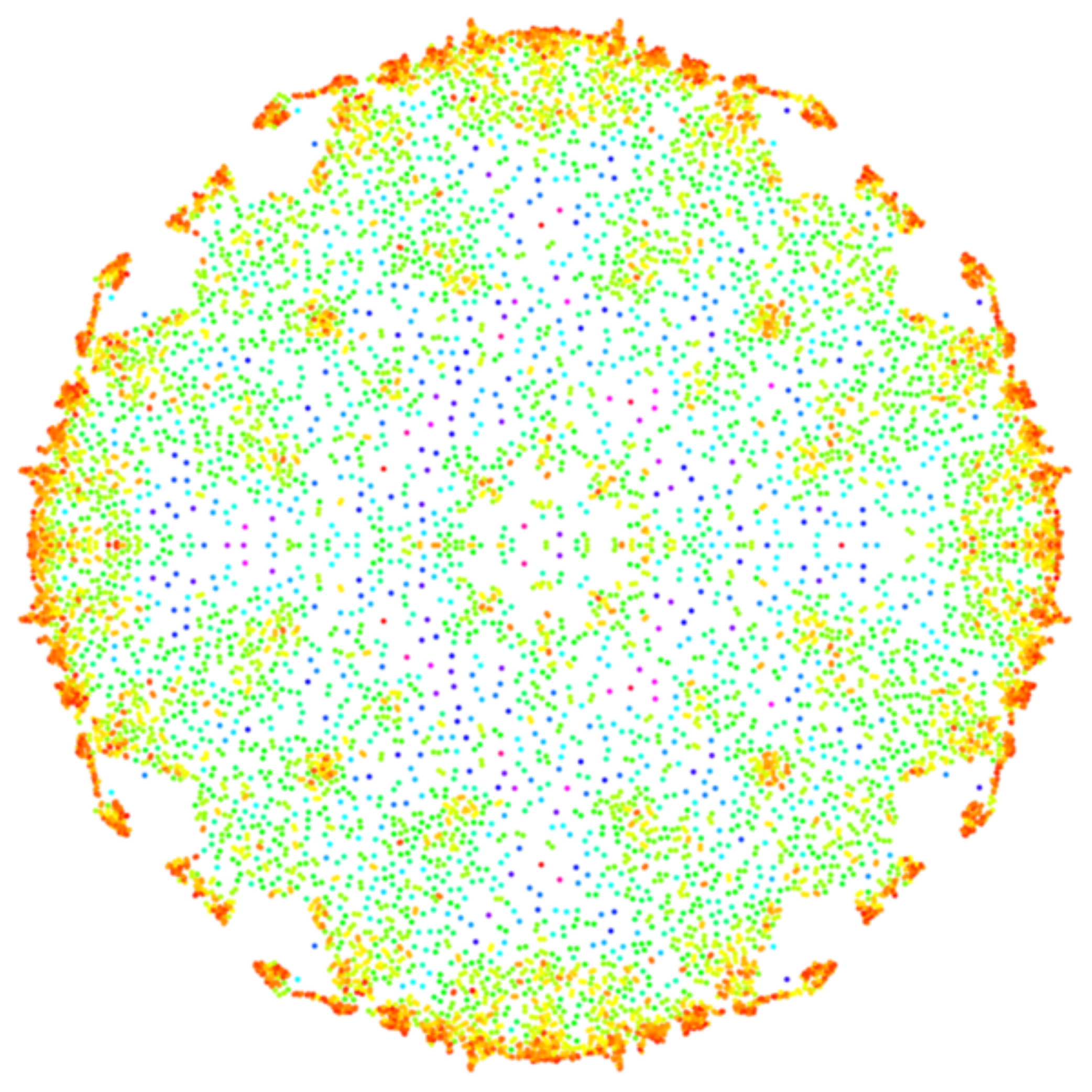}}
\scalebox{0.20}{\includegraphics{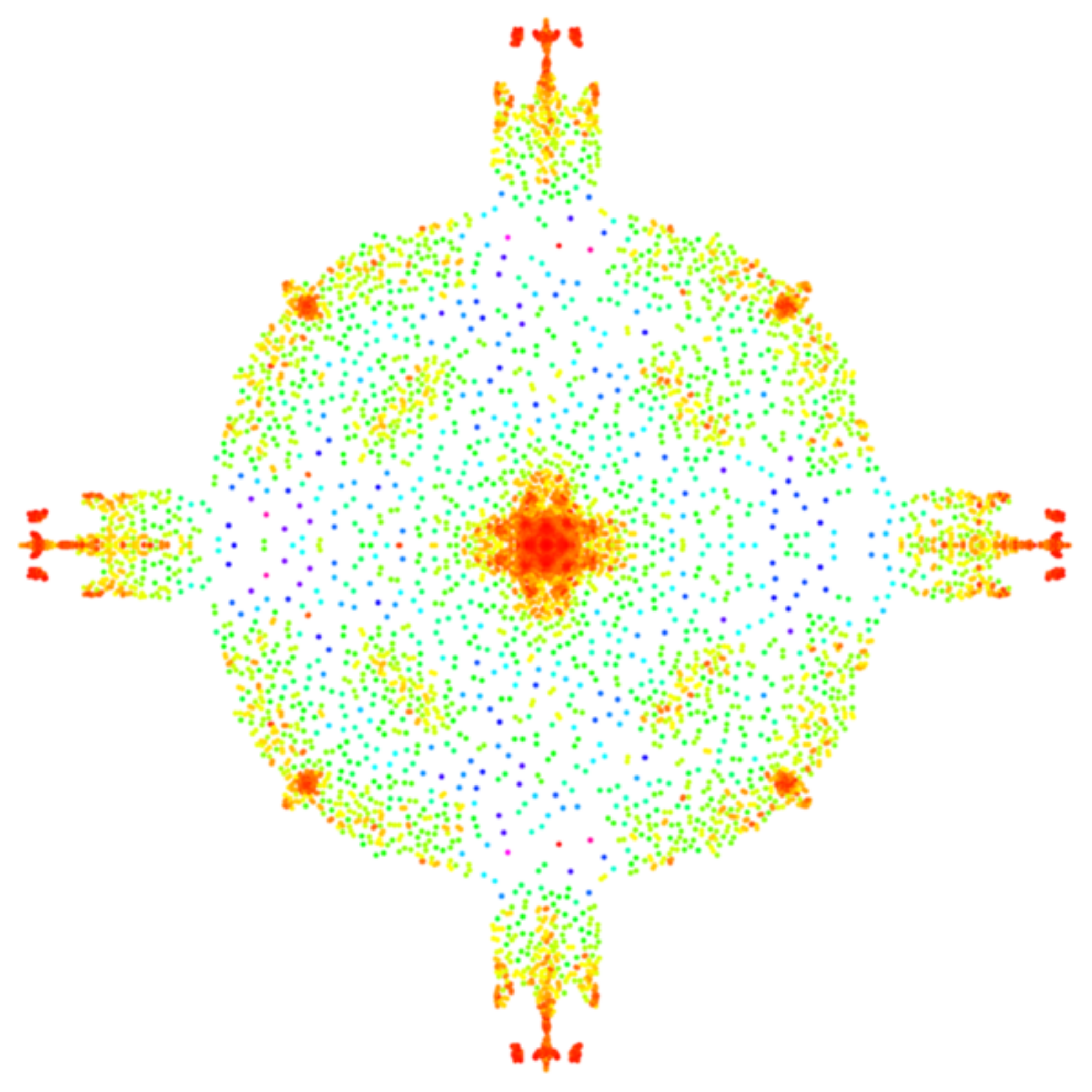}}
\scalebox{0.20}{\includegraphics{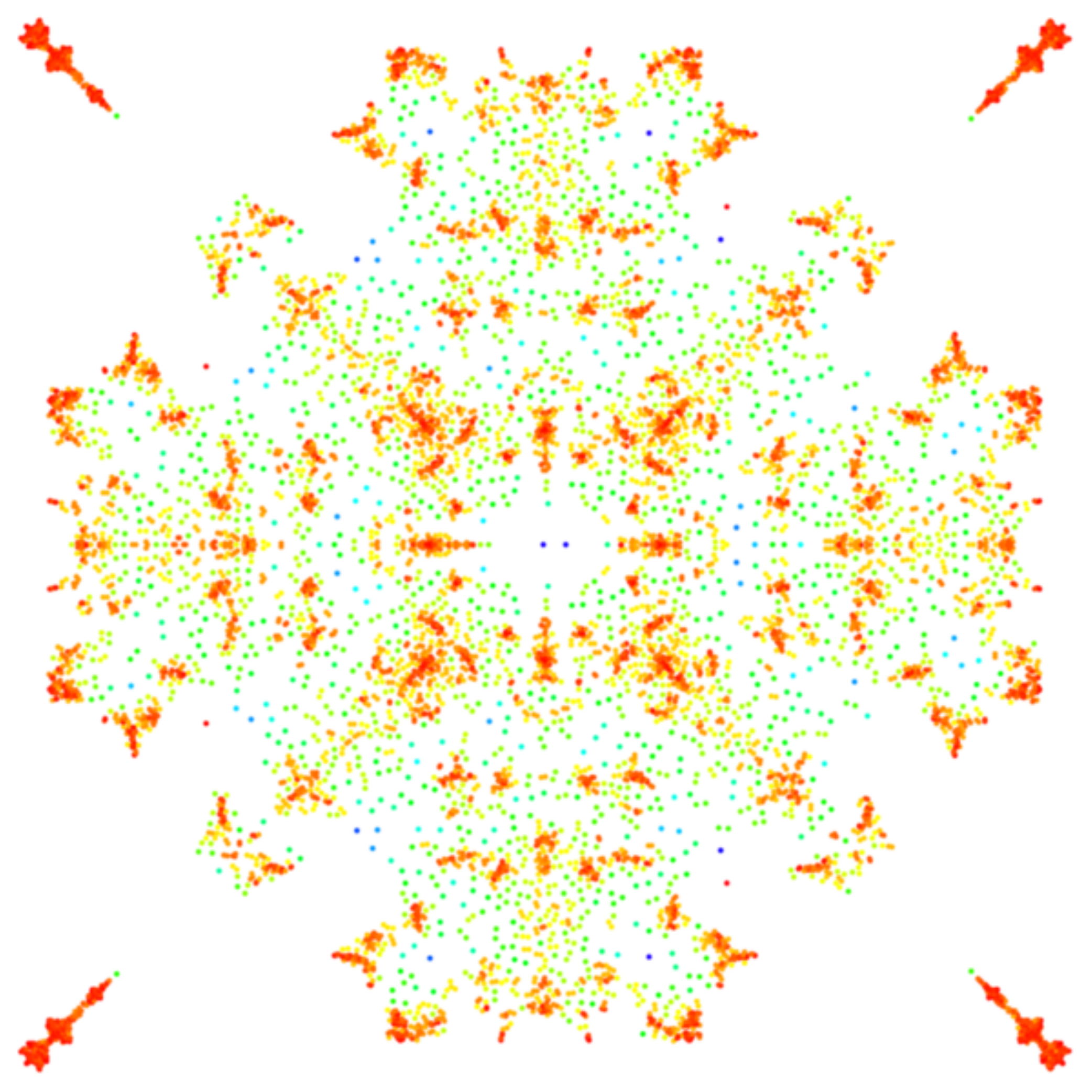}}
\scalebox{0.20}{\includegraphics{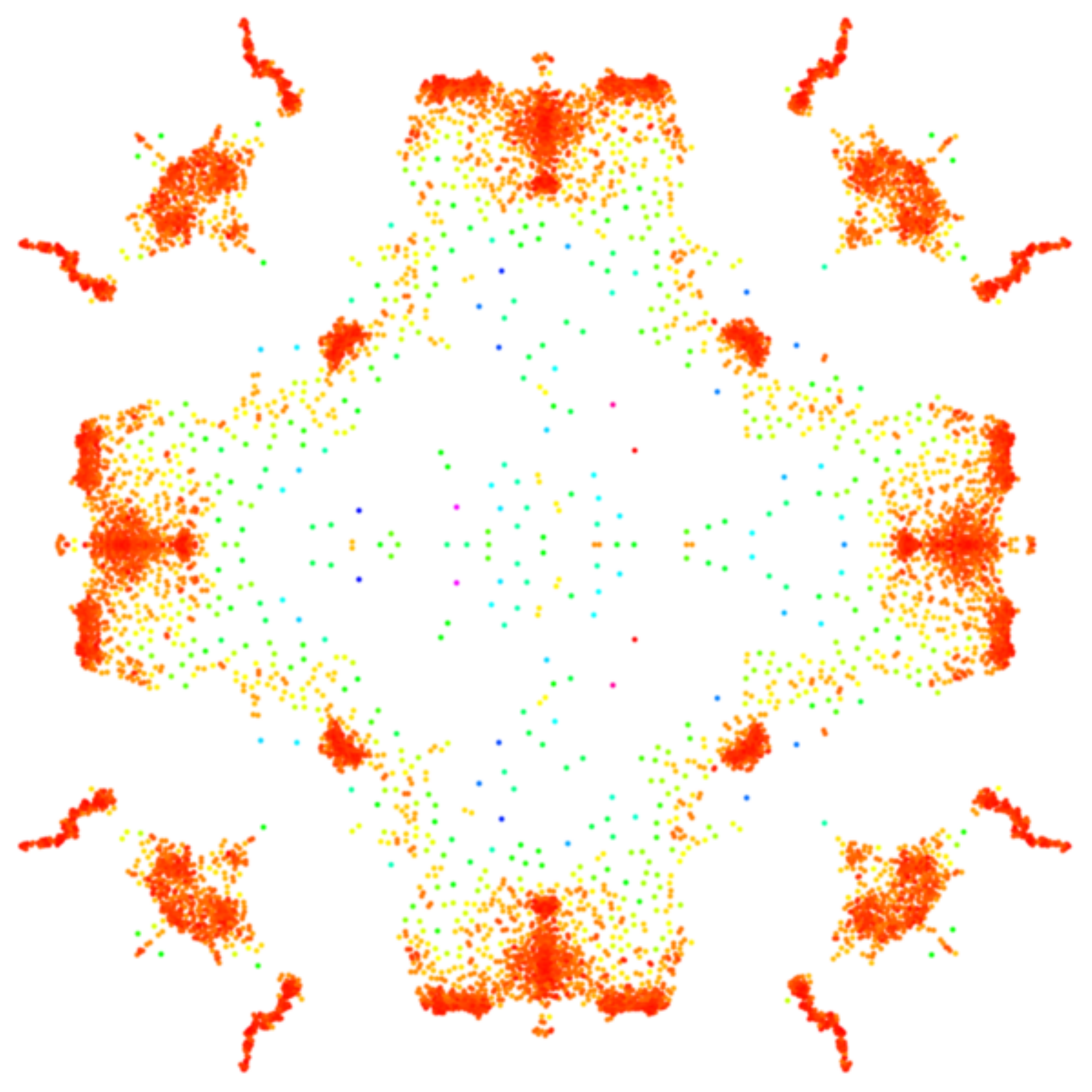}}
\scalebox{0.20}{\includegraphics{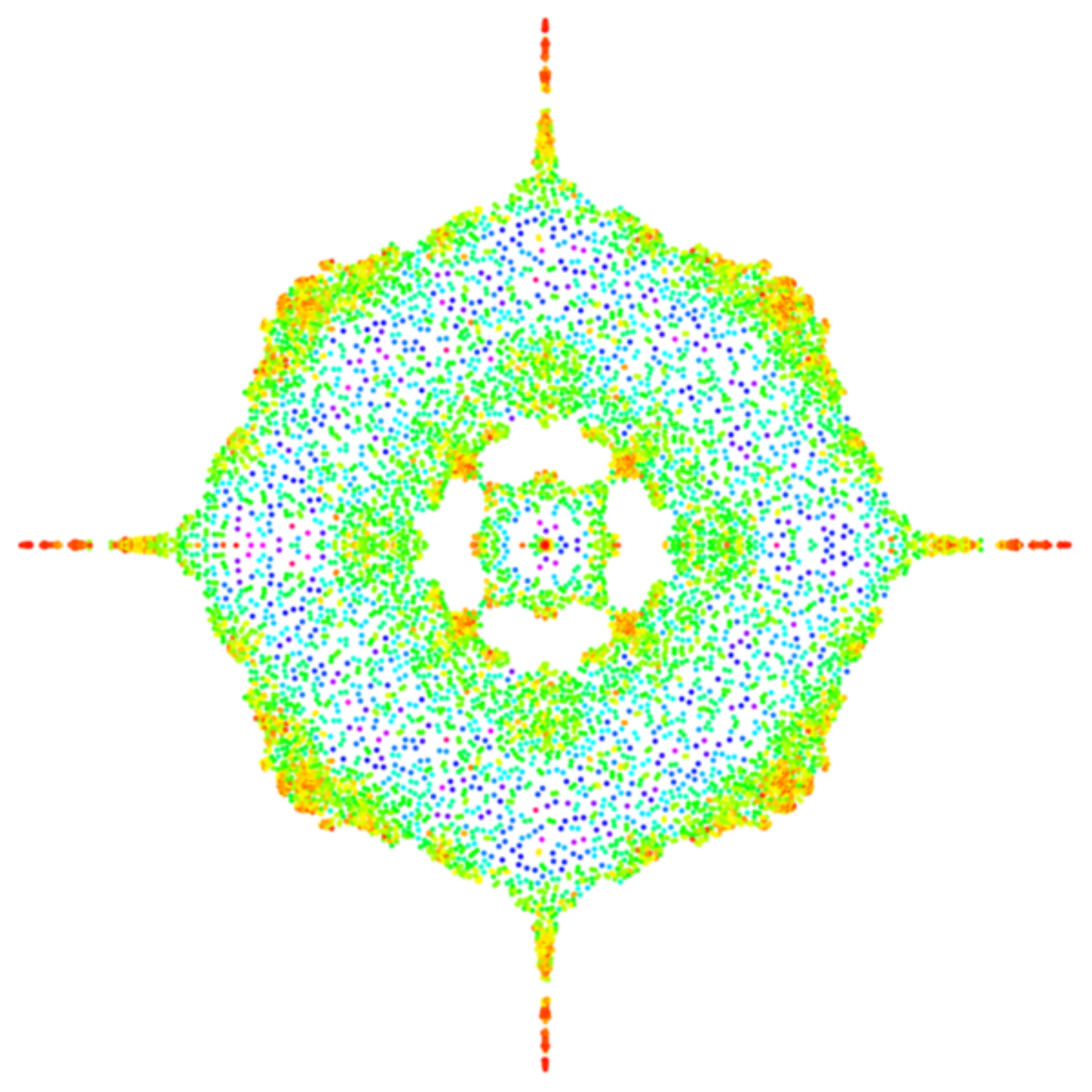}}
\scalebox{0.20}{\includegraphics{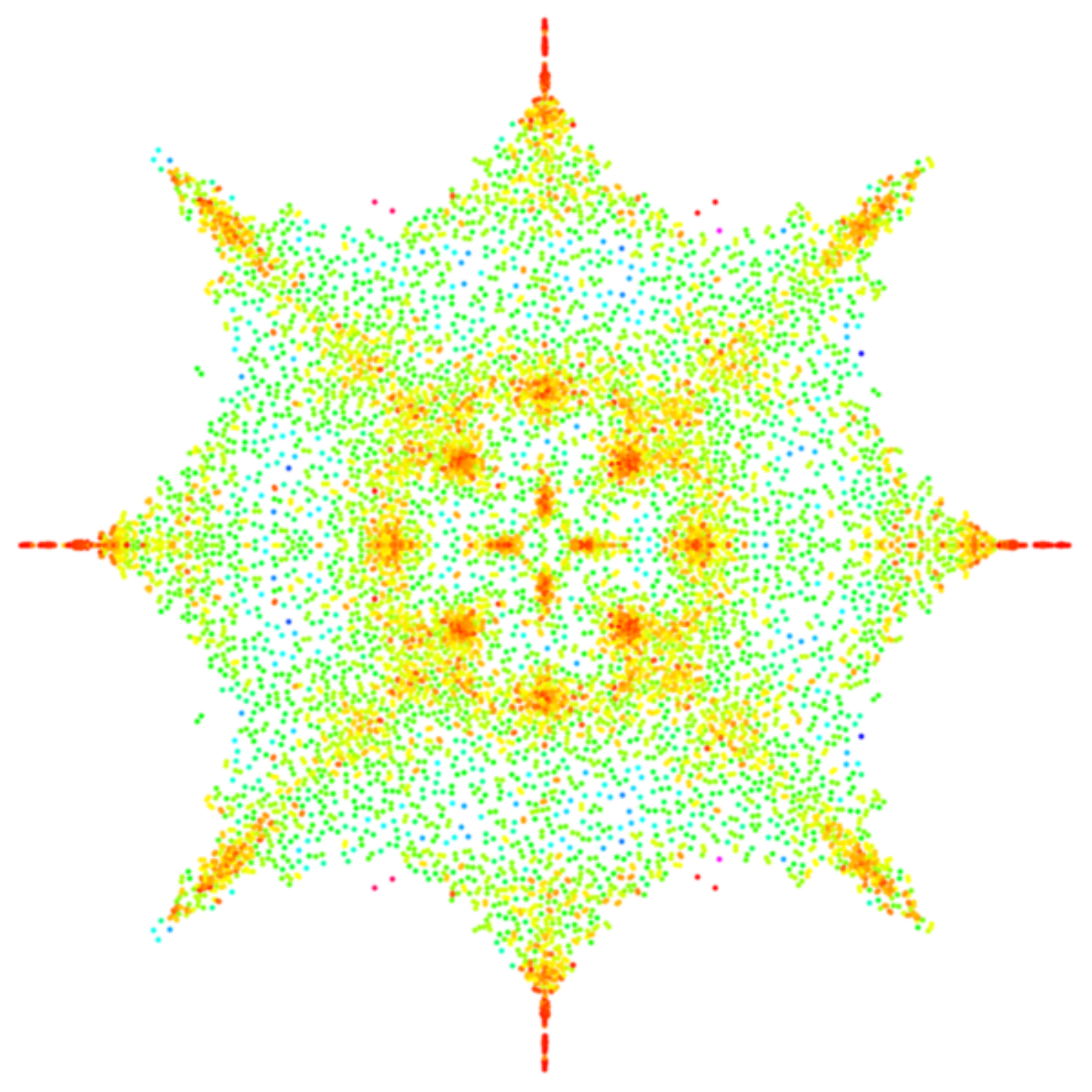}}
\scalebox{0.20}{\includegraphics{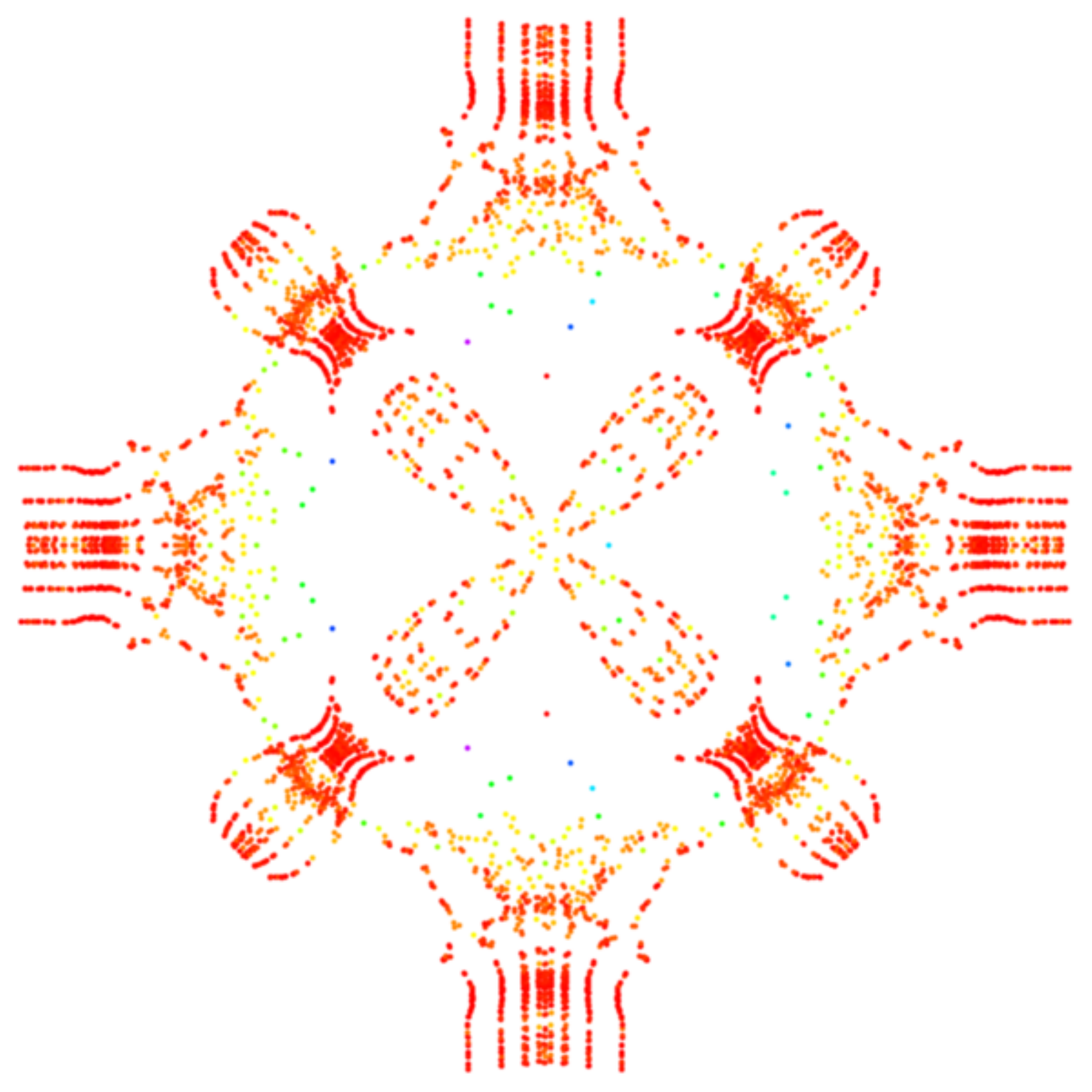}}
\scalebox{0.20}{\includegraphics{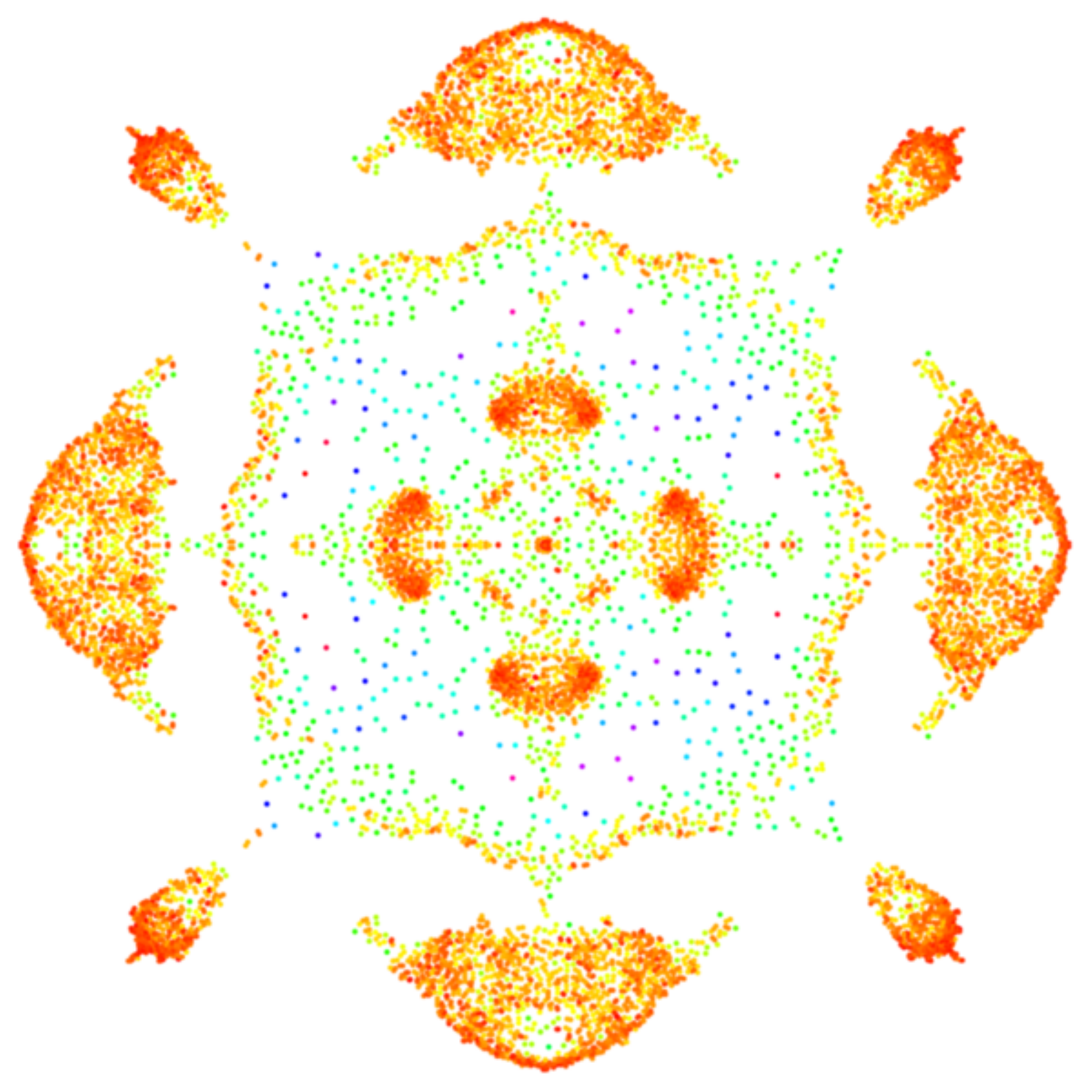}}
\caption{
\label{snowflake}
The spectrum of $10000 \times 10000$ matrices $A_{ij} = \cos(i^2 \alpha  + j \beta)$
with algebraic $2\pi \alpha,2\pi \beta$. The color of an eigenvalue depends
on the minimal distance of the eigenvalue to the rest of the spectrum.
}
\end{figure}

\begin{figure}
\scalebox{0.20}{\includegraphics{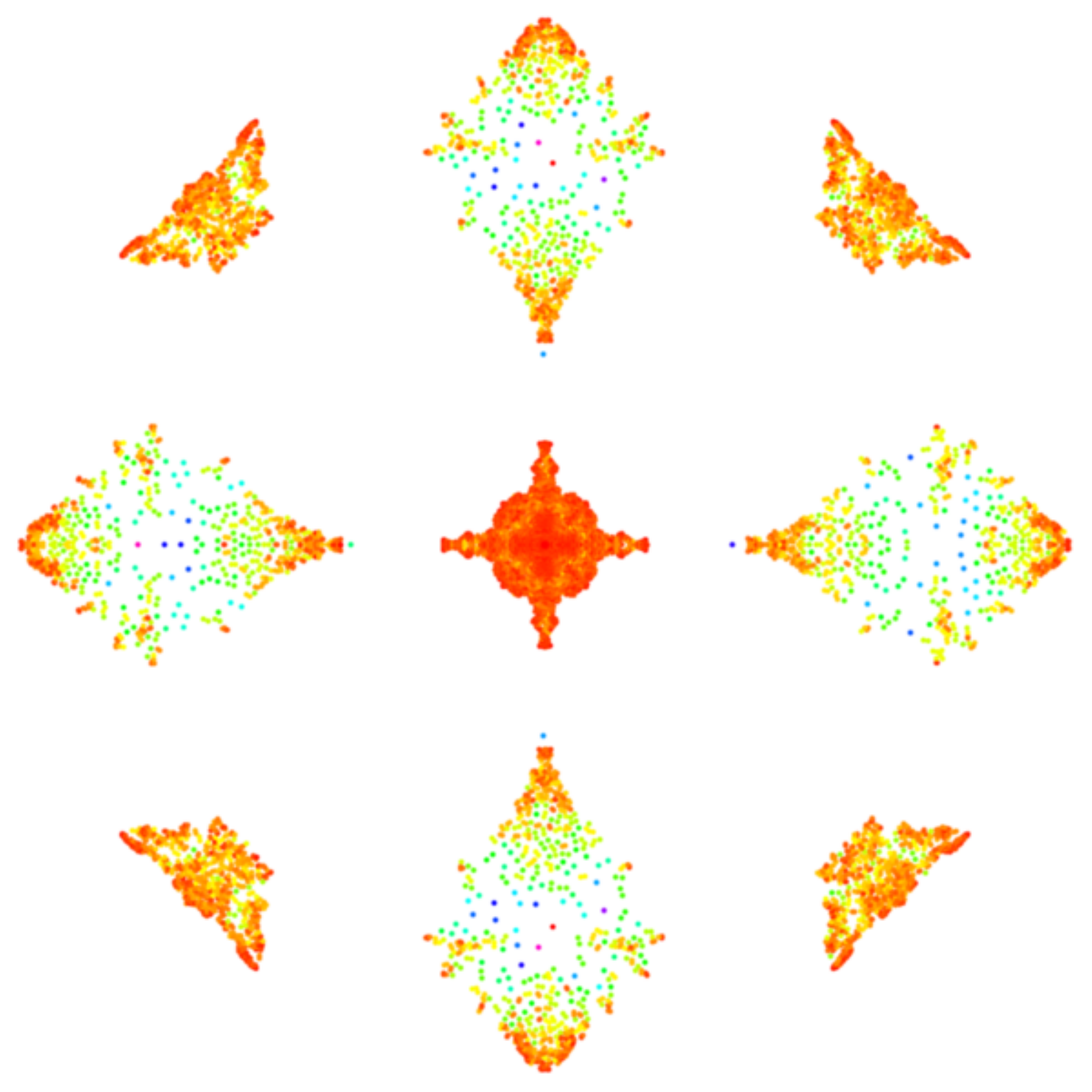}}
\scalebox{0.20}{\includegraphics{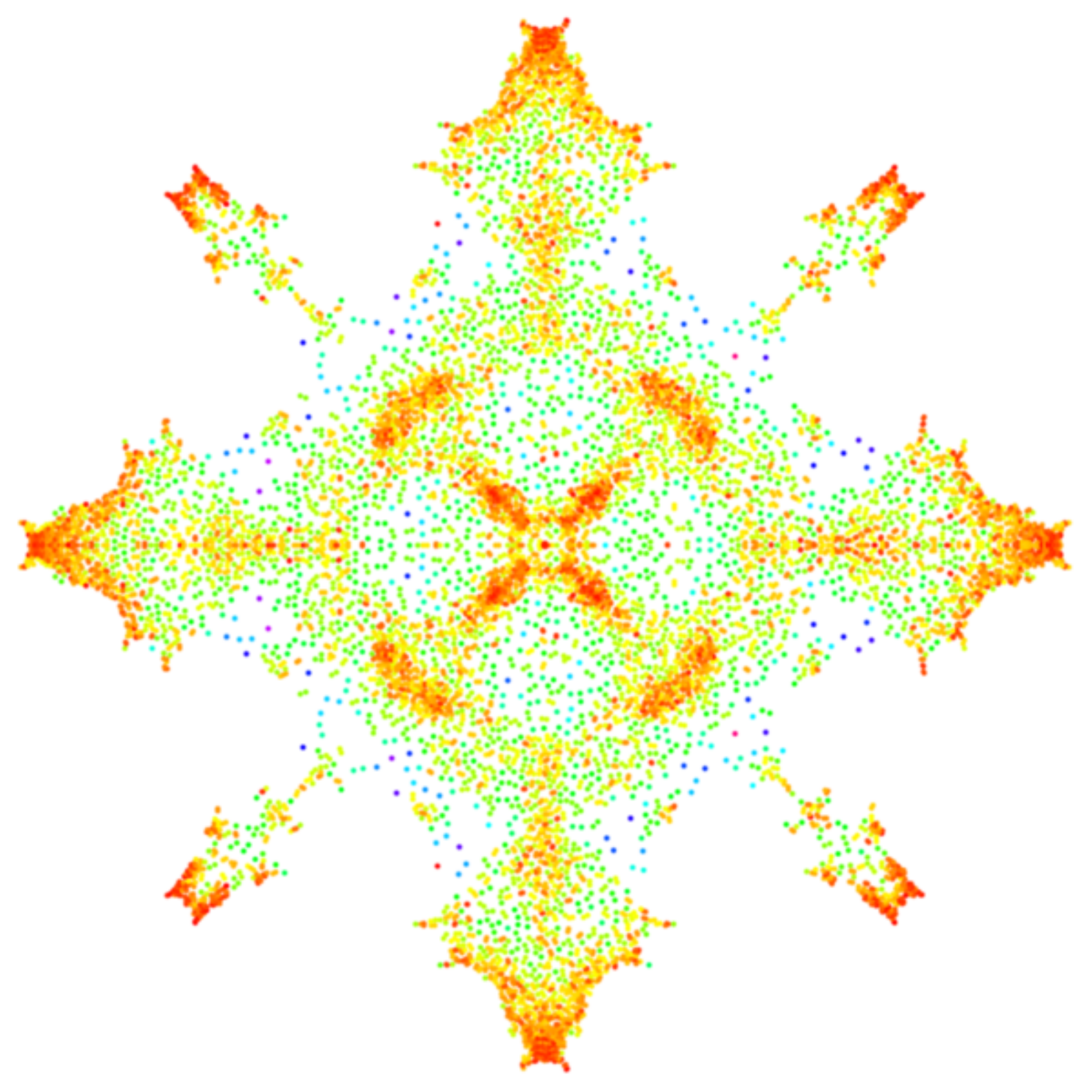}}
\scalebox{0.20}{\includegraphics{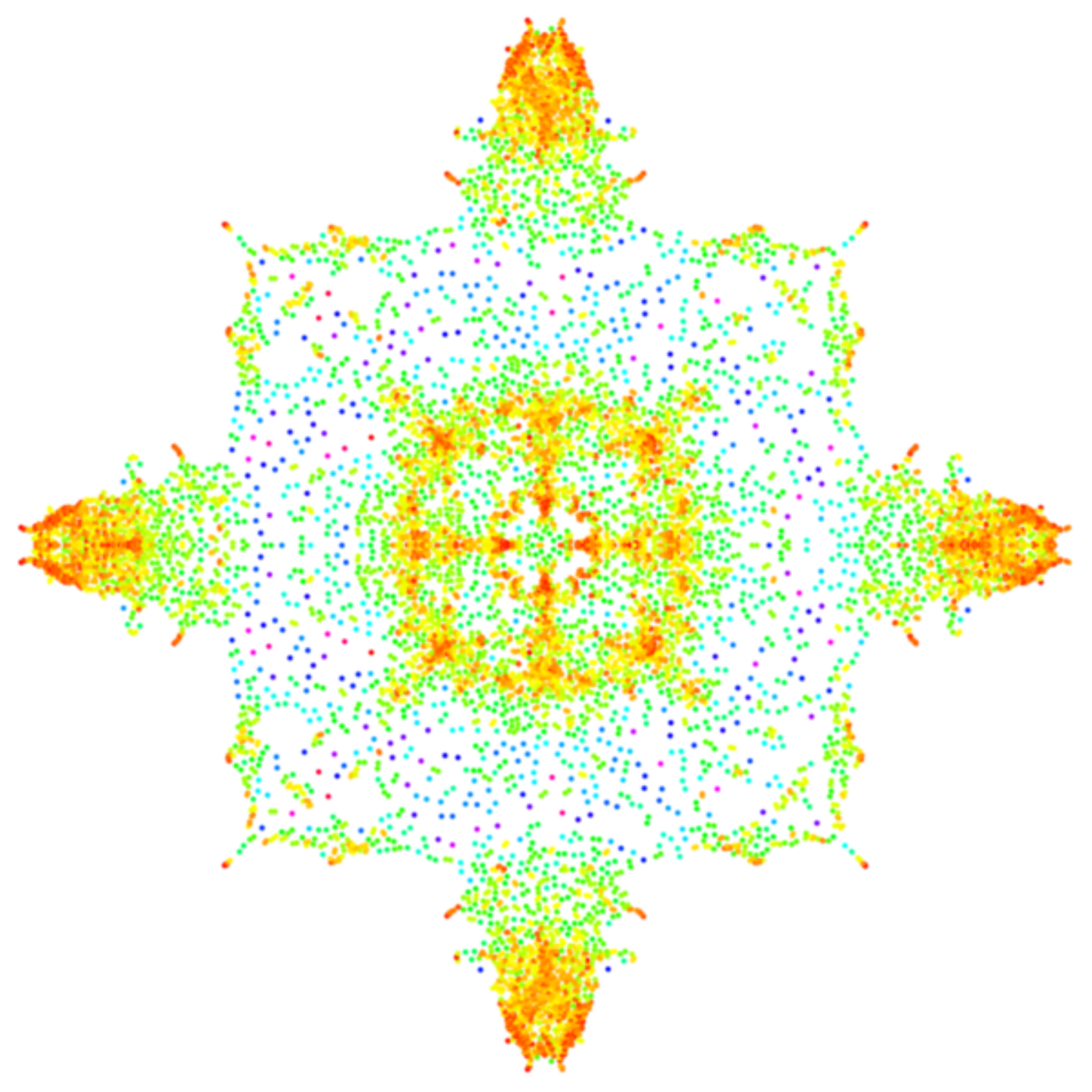}}
\scalebox{0.20}{\includegraphics{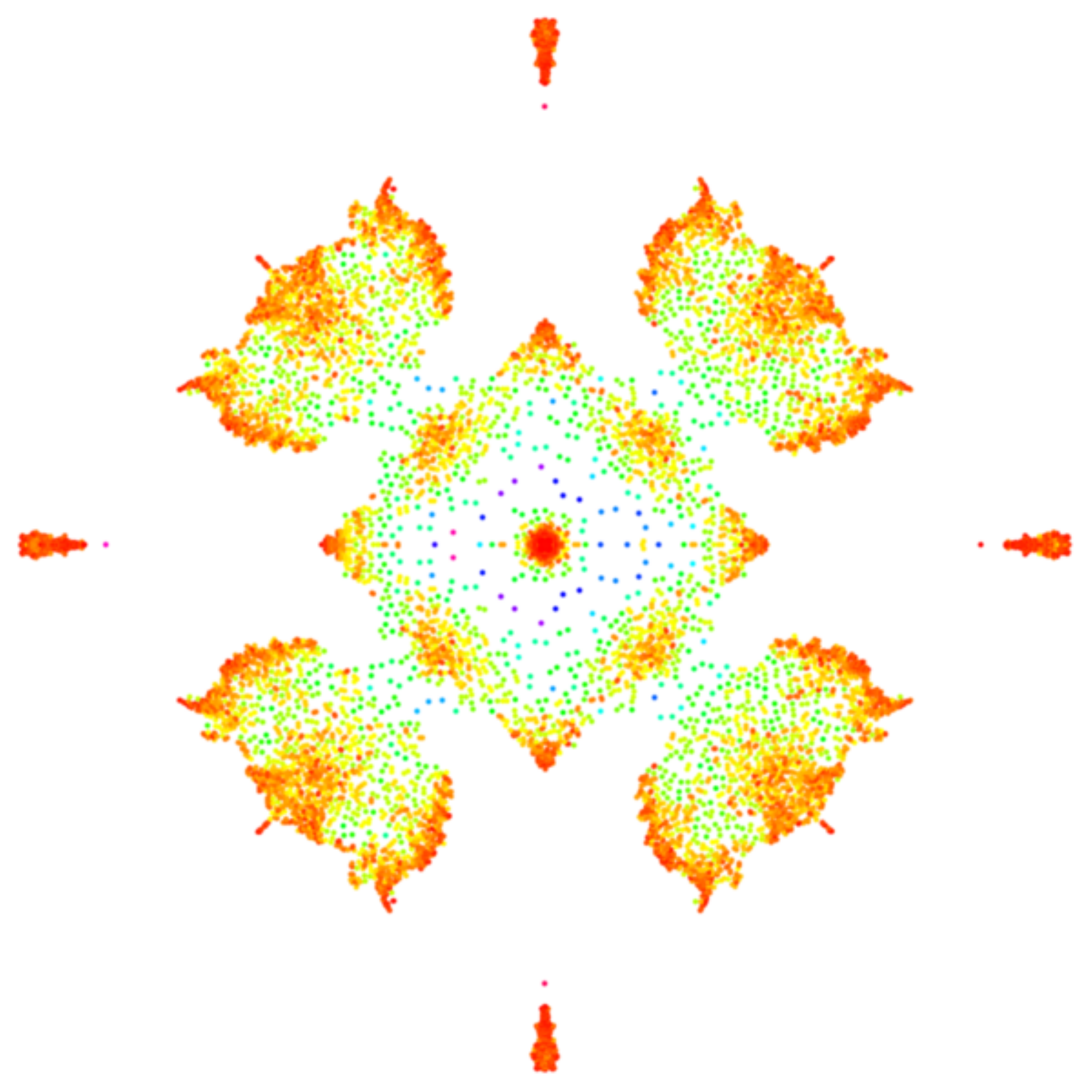}}
\scalebox{0.20}{\includegraphics{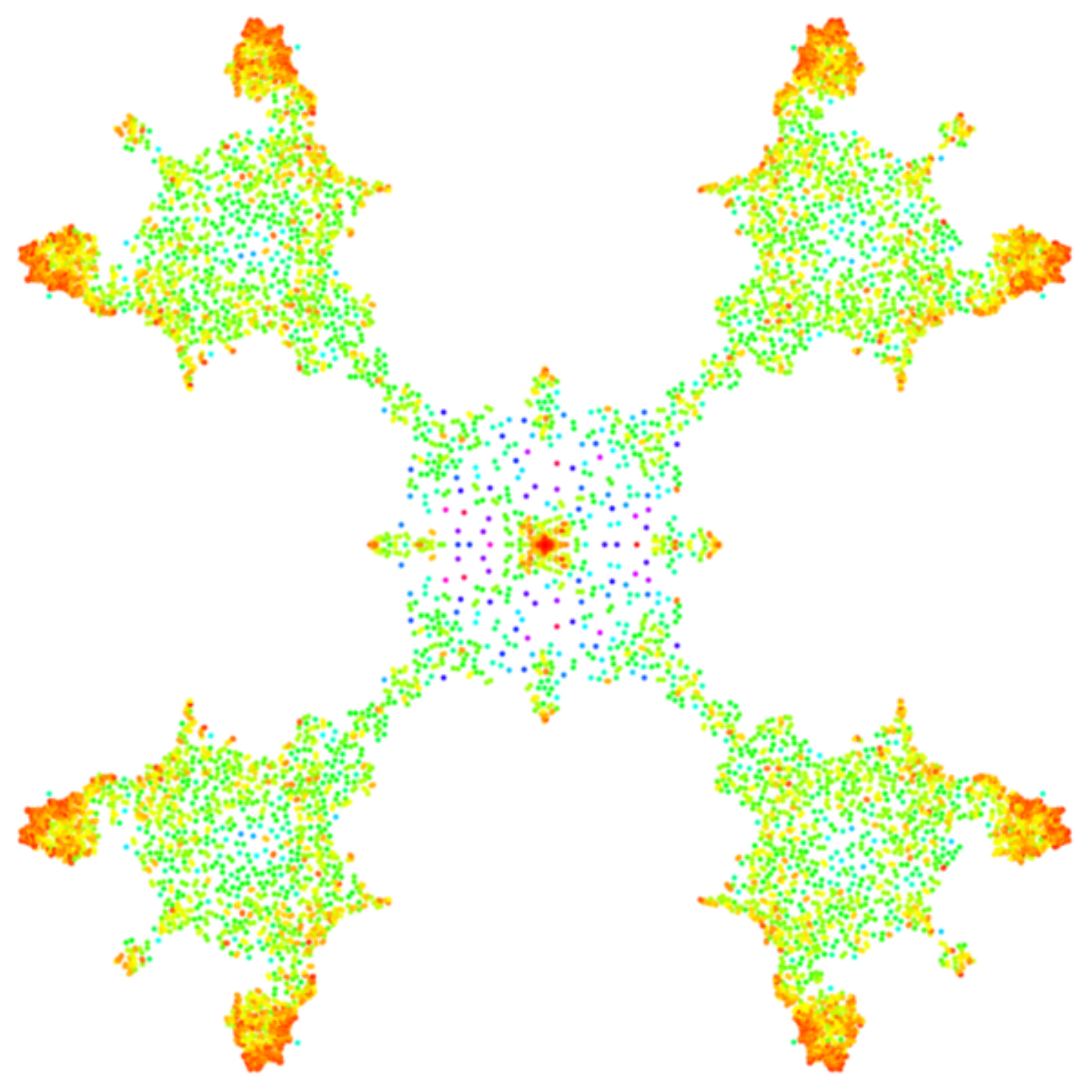}}
\scalebox{0.20}{\includegraphics{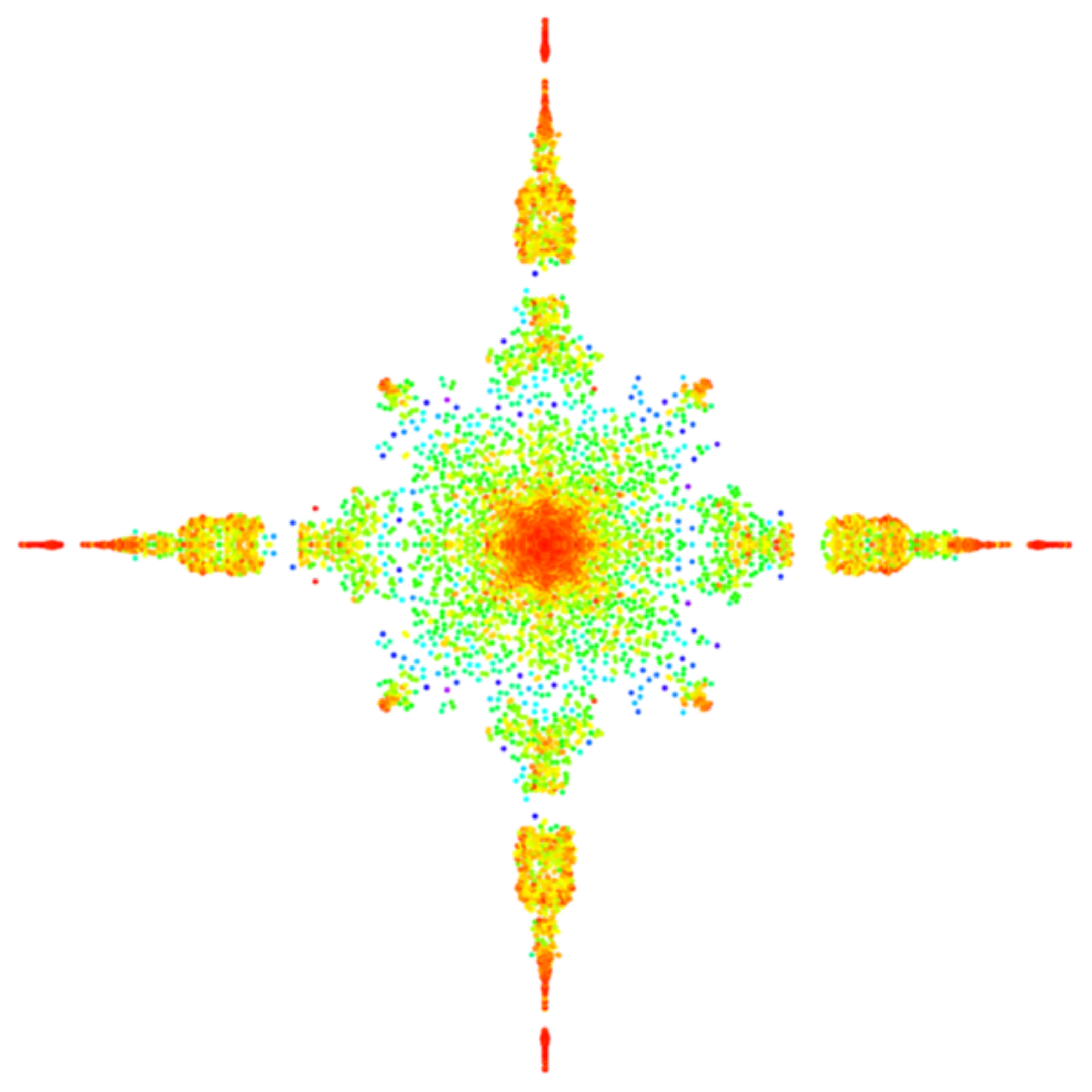}}
\scalebox{0.20}{\includegraphics{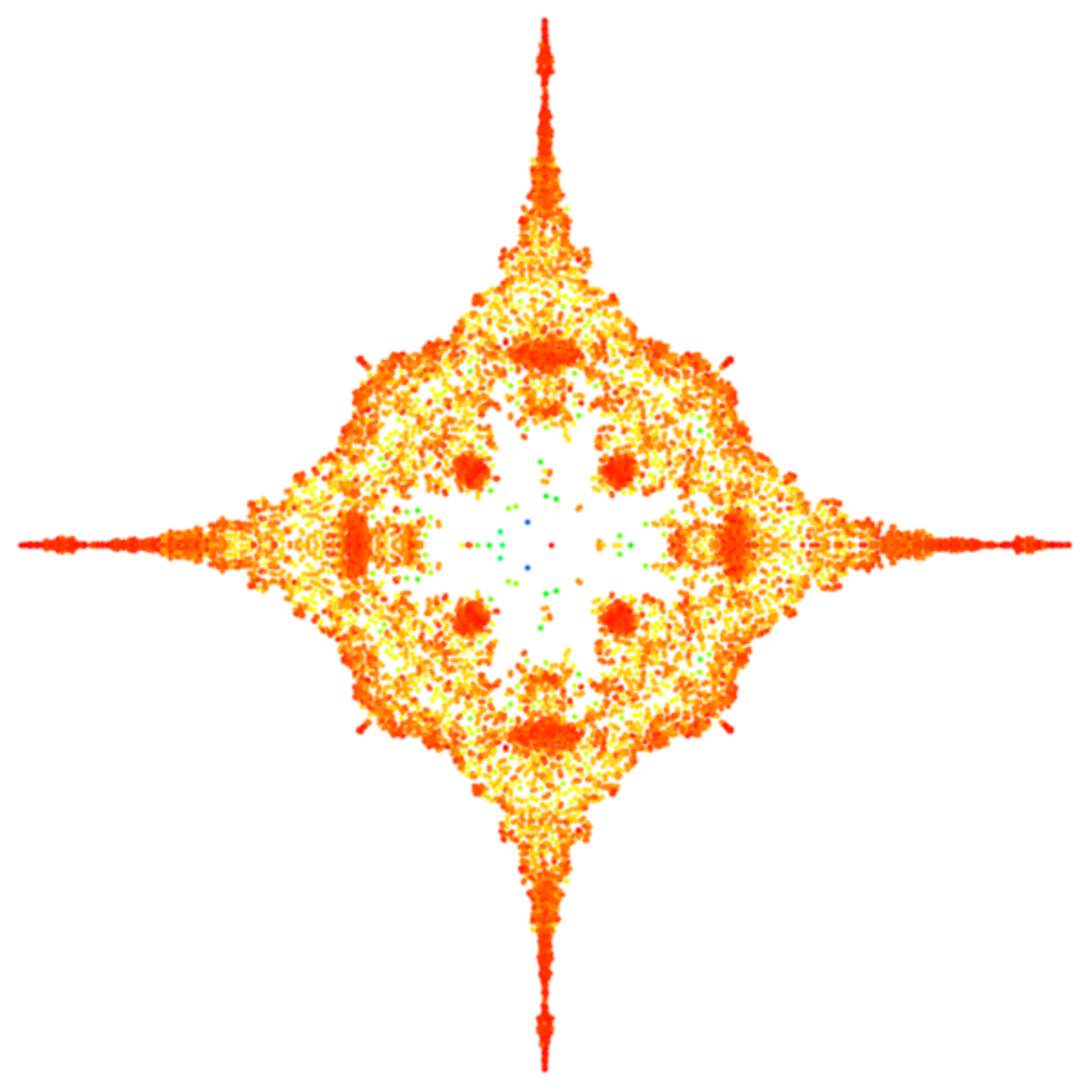}}
\scalebox{0.20}{\includegraphics{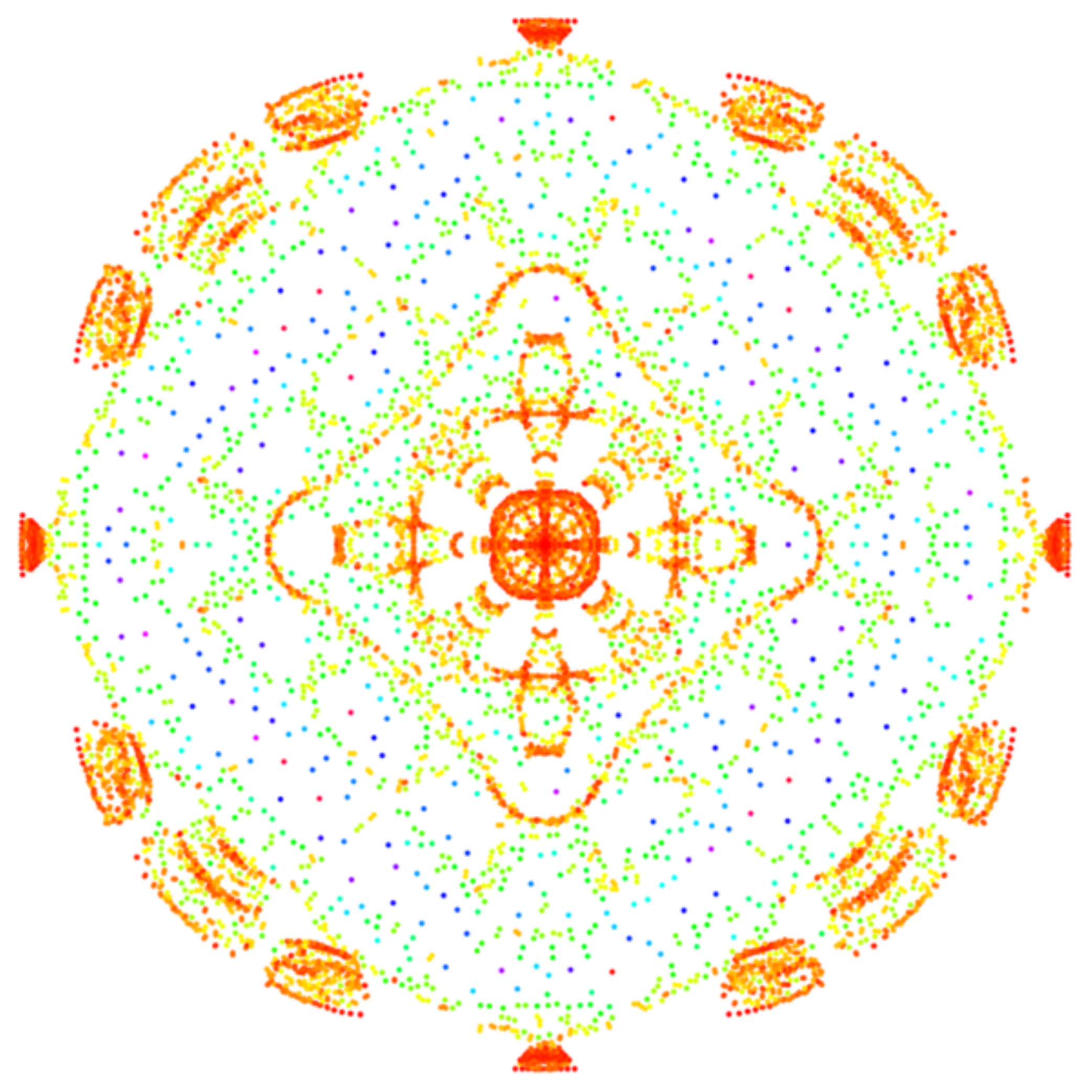}}
\caption{
\label{snowflake}
The snow flake spectra of these matrices were featured in a movie made 
for a talk in 2008 and were used in the spring 2010 for a Mathematica project. 
}
\end{figure}

How come, almost periodic notions appear in number theory? One reason is the approximation
of real numbers by fractions like the continued fraction expansion. How fast this can be done
depends on Diophantine properties. When looking at Diophantine equations or problems in the
geometry of numbers, ergodic processes related to irrational translation appear naturally too.
Here is where Aleksandr Yakovlevich Khinchin comes in who made profound contributions in
this area \cite{Chinchin92} and also was a master of probability. Actually, his law of iterated
logarithms has been put into connection with the Riemann hypothesis as popularized first in
\cite{Feller} and since been used in popular texts about the Zeta function like
\cite{Derbyshire}. \\

Unexplored is:

\conjecture{
The empirical measure of the spectra of $A(n)$ converges weakly for $n \to \infty$.
}

In any case, even a rough explanation for the strange spectral phenomena is missing. 

\section{Eisenstein Primes}   

In the ring $R=\mathbb{Q}$ of rational numbers, the {\bf integers} where the
subring containing the multiples of $1$. {\bf Primes} were
the integers different from $1$ which can not be written as a product of 
smaller primes. With this definition, examples of primes are $2,3,101,-5$ or $-11$. 
When factoring out the {\bf units}, the multiplicative subgroup $U=\mathbb{Z}_2=\{-1,1\}$ of units
one has a list of primes is $2,3,5,7, \dots$ as $-5$ and $5$ are now identified. An 
different way to get rid of the ``negative primes" is to look at {\bf prime ideals} 
in the ring $Z$. The ``primes" are now no more elements in the ring but they are
{\bf ideals}, additive subgroups $I$ of the ring which absorb other integers in the sense that
for all $n$ in $I$ and every integer $m$, also the product $nm$ is in $I$. Each nonzero 
integer in $Z$ generates the ideal of all multiples of the integer. An ideal $P$ different
from the ring is called a {\bf prime ideal} if for any pair $n,m$ in the ring, the property
that $nm$ is in $P$ implies that either $n$ or $m$ is in $P$. In the ring $Z$ of ordinary
integers the ideals are of the form $(n) = \{ \dots, -2n,-n,0,n,2n \dots \} = n Z$, 
where $n$ is an integer and the prime ideals 
are of the form $(p)$ where $p$ is a prime. The more abstract ideal setup has been 
developed by Ernst Kummer in 1844, while attempting to prove Fermat's last theorem.
It was then generalized by Dedekind in 1876 to the notion we use today.
The setup might be artificial at first, but the setup is unavoidable when generalizing 
primes to number fields. 

\begin{figure}[!htpb]
\scalebox{0.80}{\includegraphics{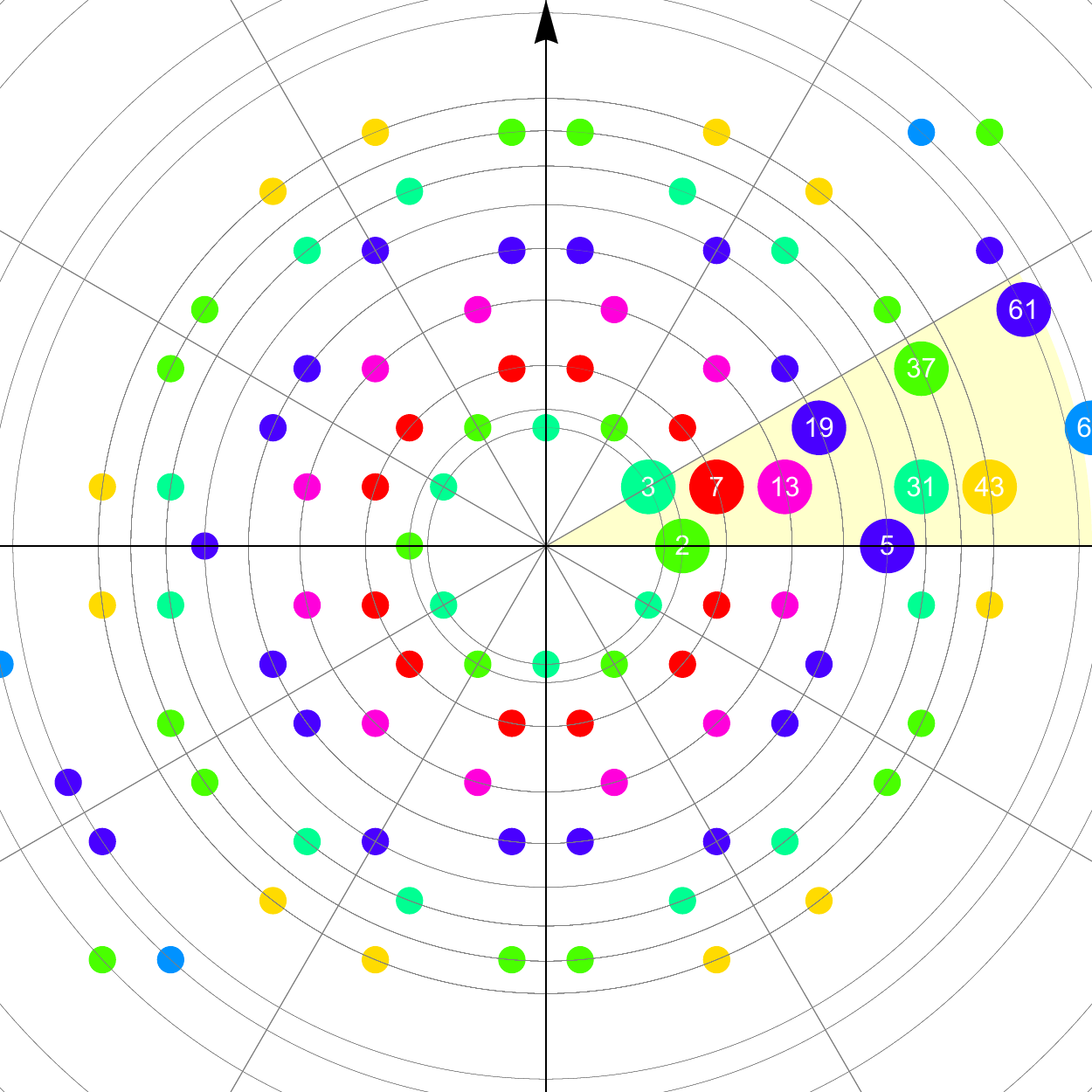}}
\caption{
The Eisenstein primes in the fundamental region cover the 
rational primes. The ordering is different than for the Gaussian primes. 
}
\label{circles}
\end{figure}

\begin{figure}[!htpb]
\scalebox{0.8}{\includegraphics{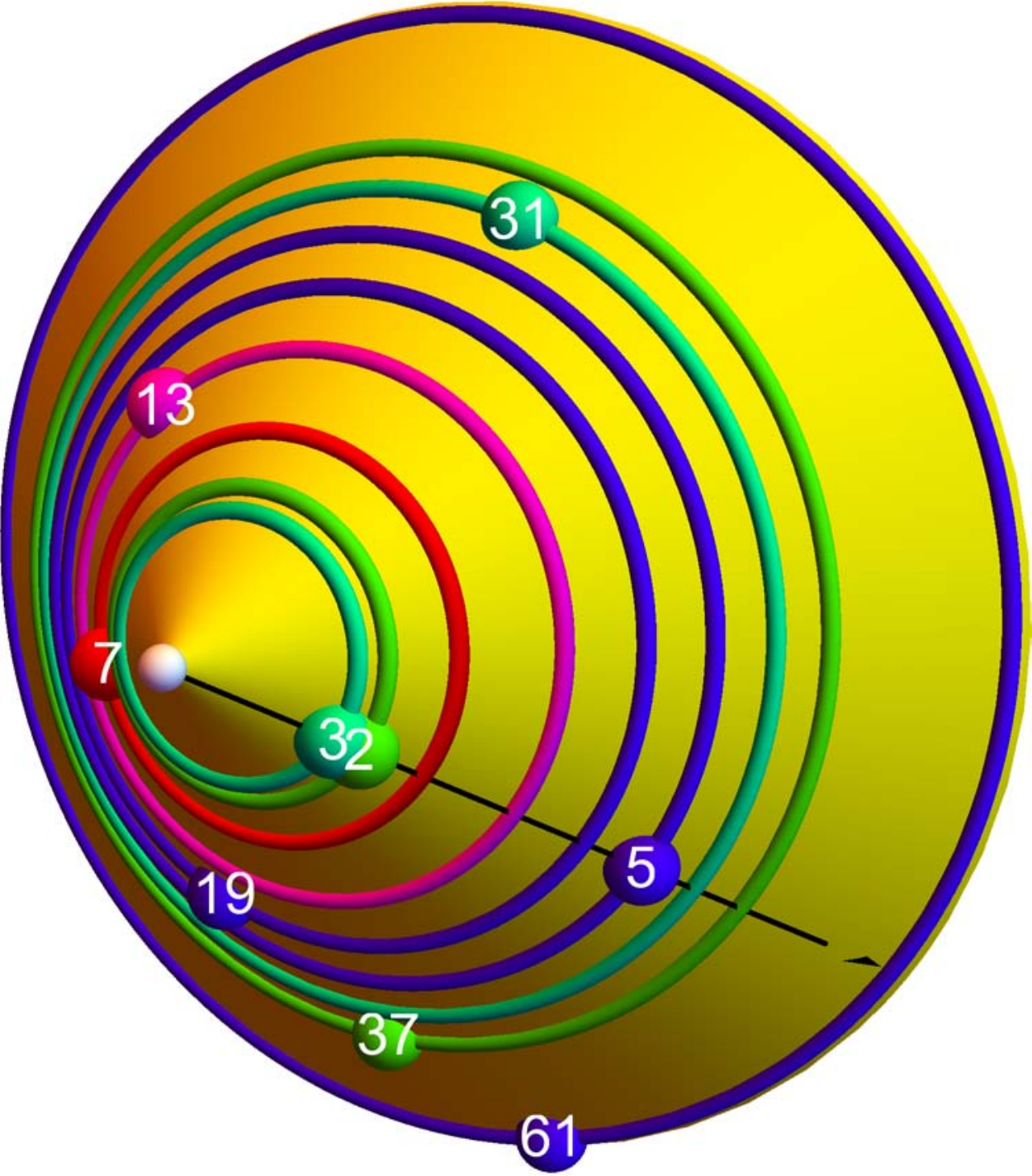}}
\caption{
Eisenstein primes cover the rational primes.
}
\end{figure}

If $R$ is an {\bf algebraic number field}, an algebraic extension of $\mathbb{Q}$, 
the {\bf ring of integers} $O$ in $R$ is the set of roots in $R$ of a 
monic polynomial with coefficients in $\mathbb{Z}$.
The ring of integers $O$ is known to always be a {\bf Dedekind ring}: the product of any two nonzero elements 
is always nonzero and every ideal different from $R$ can be written as a product of prime ideals. 
While unique prime factorization can fail, like 
in the ring of integers $R=Z[\sqrt{-5}]$ in the number field $\mathbb{Q}[i \sqrt{5}]$ 
obtained from $\mathbb{Q}$ by adding the root of $x^2+5$. The integer $6$ has there 
two different factorizations: $6=2 \cdot 3 = (1+i\sqrt{5})(1-i\sqrt{5}$. But when 
looking at ideals, the unique prime factorization is restored. If we write $(a_1 \dots, a_k)$ for the
ideal generated by $a_1, \dots, a_k$, then  Now $(2),(3),(1+i \sqrt{5}),(1-i \sqrt{5})$ are 
all ideals but they are not all prime ideals. There are "smaller primes" in the form of
"larger ideals" like $p = (2,1+i \sqrt{5})$ or 
$q = (3,1+i \sqrt{5})$. With these, the factorization $6=pq$ has become unique. The {\bf principal ideals}
generated by one element had to be enlarged and include {\bf fractional ideals}. 
The {\bf order} of a ring $R$ is a sub-ring $O$ which is an algebra over the field $\mathbb{Q}$ 
which is the free Abelian group generated by a basis of $O$. The notion of ``order" has been important
already in non-commutative cases like the ring of quaternions, where the Lipschitz integers 
do produce an order, but not a maximal one forcing the Hurwitz integers. 
The example of the Octonion ring shows that the maximal orders are not necessarily unique.
In an example like $\mathbb{Q}[i]$, the order $\mathbb{Z}[2i]$ is not maximal neither. 
To summarize: one can define ``primes", if there is a norm which is multiplicative $N(x) N(y) = N(x y)$
implying that one has to be in a division algebra. This works for number fields in the complex plane.
The {\bf integers} of a number field are the maximal orders in that field.
Sometimes these rings of integers have a unique factorization by itself,
sometimes, it is necessary to use {\bf fractional ideals} which are not necessarily principal.
The {\bf ideal class group} is the quotient of the set of fractional ideals divided
by the set of principal ideals. Its cardinality is the {\bf class number} of $O$. \\

But what about $O=Z[\sqrt{-3}]$? One has the ideal $p=(2,1+i \sqrt{3})$ but $p^2 = (2) p$ 
and the fact that $p$ is not equal to $(2)$ shows that there is no unique factorization
into ideals. The ring $O$ appears not to be a Dedekind ring. Did we not just state that any
ring of integers $O$ of an algebraic number field like $\mathbb{Q}[\sqrt{-3}]$ is a Dedekind ring?
Yes, but the ring $Z[\sqrt{-3}]$ is {\bf not} a ring of integers. Indeed, there is a larger ring which 
needs to be considered. In other words, $Z[\sqrt{-3}]$ is not a {\bf maximal order}. 
There is a larger ring $Z[(1+\sqrt{-3})/2]$ which now contains the roots of the
polynomial $z^3=1$ and is called the ring of {\bf Eisenstein integers}. This ring has already 
been constructed by Euler. The Eisenstein integers are now a {\bf ring of integers} and unique
factorization holds. The pictures of the primes in this ring are very beautiful as they feature
a hexagonal symmetry. Why are the Eisenstein integers a ring of integers?
We have to show that if $\alpha = (1+\sqrt{-3})/2$, then every element 
$a+b \alpha$ is a root of a monic polynomial. One can give it explicitly: 
as $x^2-(2a+b)x + (a^2+ab+m^2)$.  (Note that we use $\alpha$ for which both real and imaginary part
are positive, one usually looks at $\overline{\alpha}$ as a basis but for Goldbach, the positive
sign is better). 
We have seen above that the class number of $Z[i \sqrt{5}]$ is $1$, but
that the class number of $Z[(1+i\sqrt{3})/2]$ was $2$. The class number quantifies,
how far away the ring is from a unique factorization domain. 
While in general not understood yet, one knows the class numbers 
for quadratic number fields $Z[\sqrt{d}]$ with negative $d$. \\

An Eisenstein integer $z=a+b \alpha$ is prime if and 
only if $p=N(z)=z \overline{z} =a^2+b^2+ab$ has one of the following properties: either
$p$ is prime and congruent to $0$ or $1$ to $3$ or then $\sqrt{p}$ is prime and congruent to 
$-1$ modulo $3$. The Eisenstein integers show a similar {\bf dichotomy} like the Gaussian 
integers but with the $4$ replaced by $3$. Obviously, also $3$ plays a special role now.

\resultremark{
There is a natural bijection between Eisenstein primes and rational primes.
}

The picture is the same: identify ${\rm arg}(z)=0$ with ${\rm arg}(z)=2\pi/6$ to get a cone. Now the 
prime $2$ represented by $2+0 w$ and $3$ represented by $1-w$ as well as any prime of the form
$3k-1$ are on the gluing line, while the primes of the form $3k+1$ are in the interior. 
Again, we get for every of the later a unique angle $0<\alpha<\pi/6$ which is the argument
of ${\rm arg}(a+b w)$ with $p=a^2+b^2+ab$. Again we don't know of a fast way to compute that angle
for a given prime of the form $3k+1$. \\

Let us define $Q$ to be the open {\bf hexant} $Q=\{ x + y w \; | \; x>0,y>0 \}$, where 
$\omega = (1+\sqrt{-3})/2$. Note that we use the cube root of $-1$ and not the cube root 
$(-1+\sqrt{-3})/2$ as usual. The reason is Goldbach, were it is more convenient to work with 
$w$.  When experimenting with these numbers,
we notice that every Eisenstein integer in the open quadrant
in distance 2 or larger from the boundary is a sum of two Eisenstein primes. 
Almost! There are two exceptions. Not on the boundary of integers of the form $2+k w$ 
which seems most vulnerable but on the row $3+k w$!

\begin{figure}[!htpb]
\scalebox{0.6}{\includegraphics{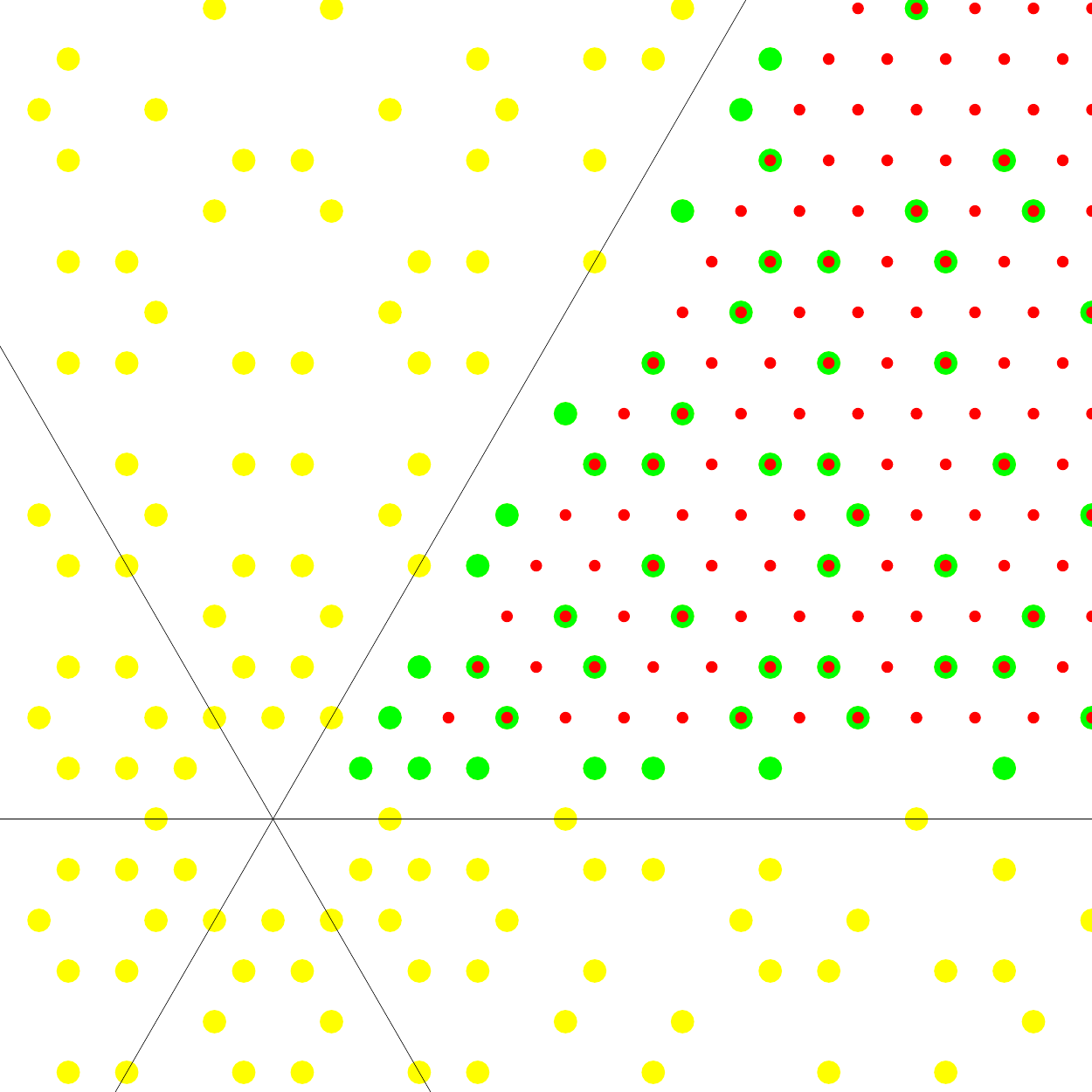}}
\caption{
The Eisenstein Goldbach conjecture looks good at first for coefficients $>1$. 
Here we see the initial part. But there are
integers $(109 +3w)$ and $(121 +3w)$ (as well as their mirrors) which can not be written as the sum
of two positive Eisenstein primes $a+bw + (c+dw)$ with positive $a,b,c,d$. 
These are the bad {\bf Eisenstein ghost twins}. They are not visible on this picture
}
\label{circles}
\end{figure}

\begin{figure}[!htpb]
\scalebox{0.3}{\includegraphics{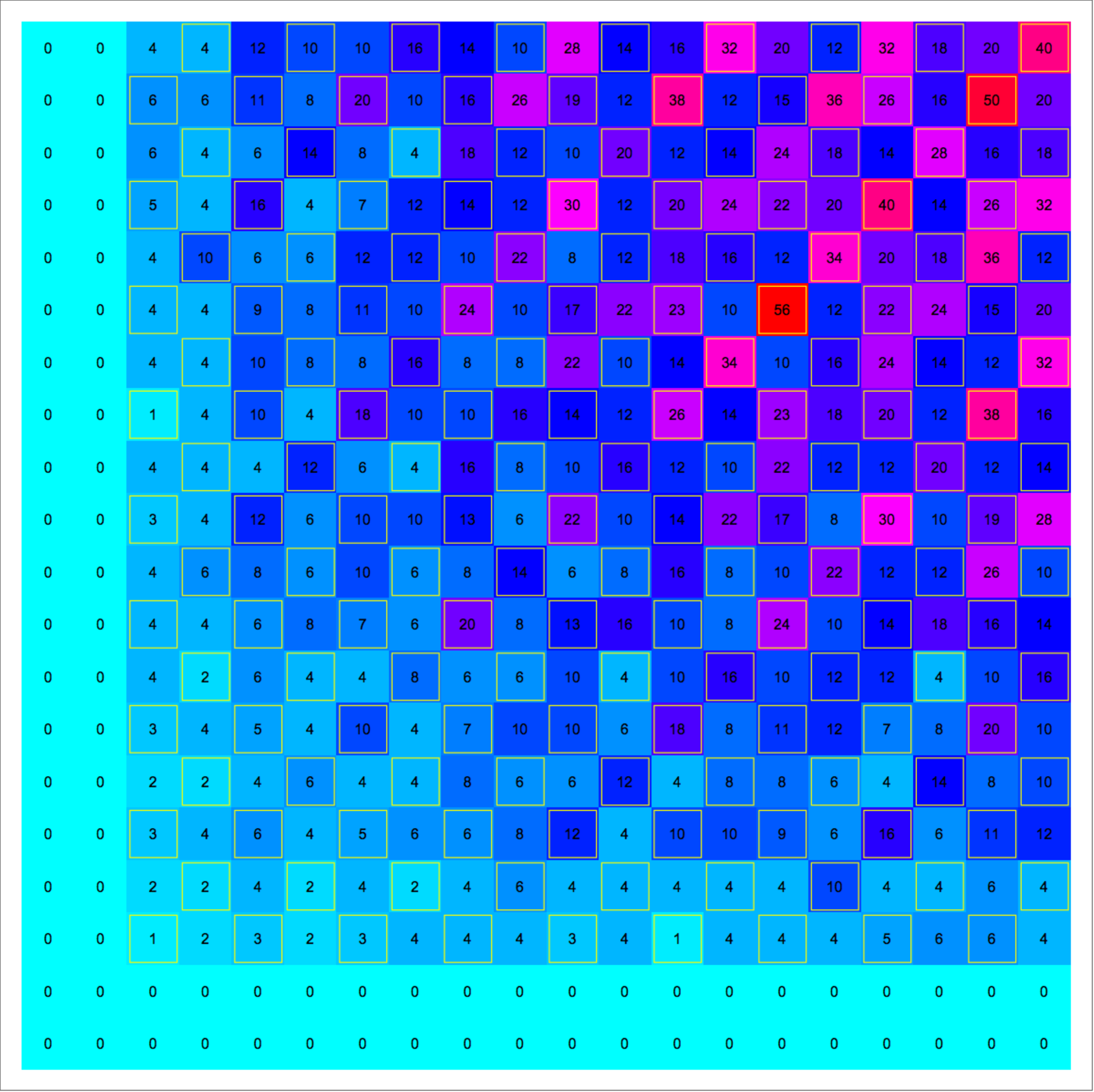}}
\caption{
Here we see how many times an Eisenstein integer $a+w b$ can be written 
as a sum of two Eisenstein primes.
}
\label{circles}
\end{figure}

\begin{figure}[!htpb]
\scalebox{0.3}{\includegraphics{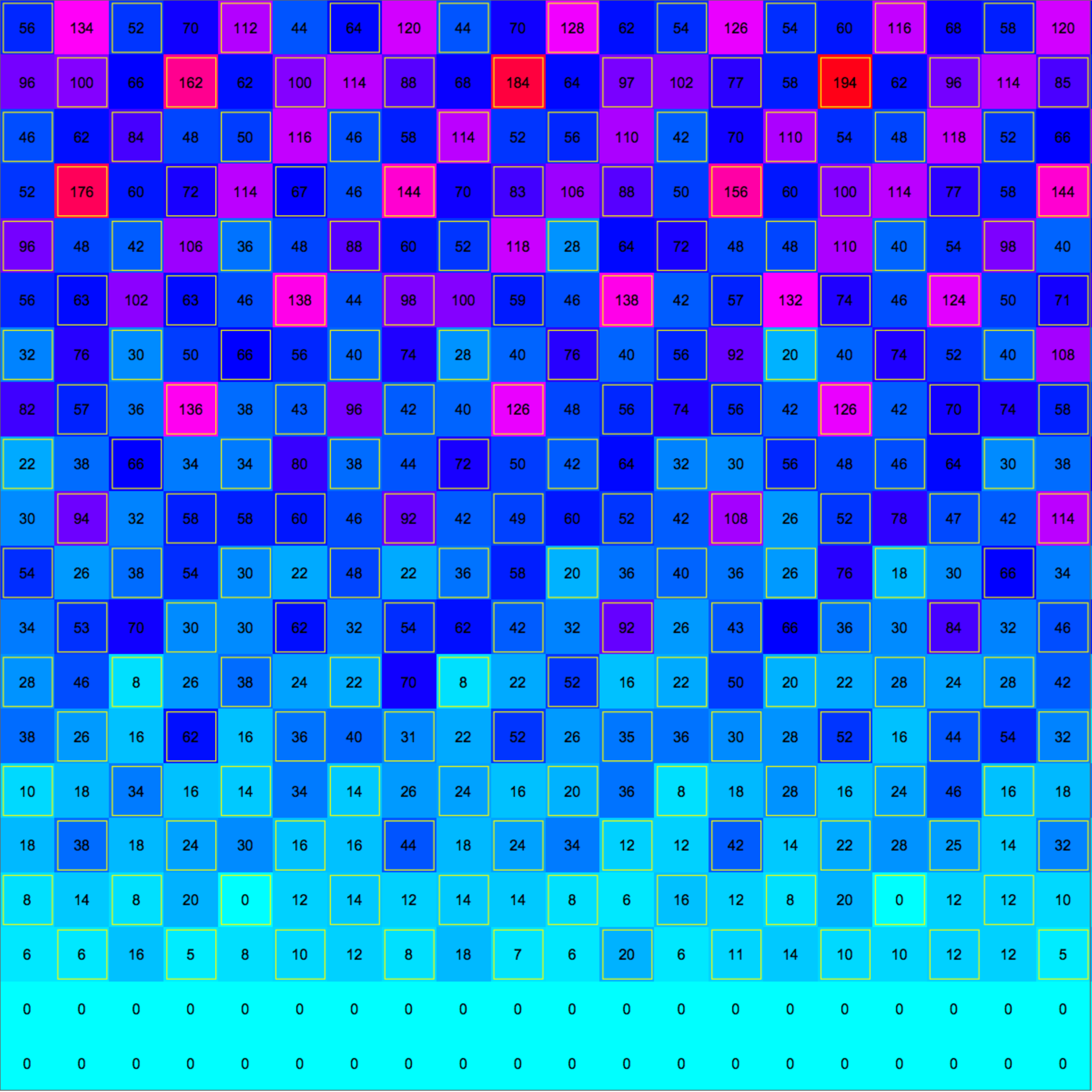}}
\caption{
The Eisenstein ghosts $(109 +3w)$ and $(121 +3w)$ (as well as their mirrors) which can not be written as the sum
of two positive Eisenstein primes $a+bw + (c+dw)$ with positive $a,b,c,d$.
These are the bad {\bf Eisenstein ghost twins}. Can you find the ghosts in the matrix? They are 
caught in the interior of the matrix and not at the boundary. 
}
\label{circles}
\end{figure}

We don't even need the condition that the integers are {\bf even} except for the
two ghost examples (we identify $a+b w$ with $b+a w$ when counting). 
The Eisenstein primes are so dense that their sum seems to 
cover the sector: 

\conjecture{
Every Eisenstein integer $z=a+b w$ with $a>2,b>2$ is a sum of two
Eisenstein primes in $Q$.
}

and

\conjecture{
Every Eisenstein integer $z=a+b w$ with $a>1,b>1$ is a sum of two 
Eisenstein primes in $Q$ except for the two Eisenstein twin ghost examples 
$3+109 w, 3+121 w$.
}

\begin{figure}[!htpb]
\scalebox{0.80}{\includegraphics{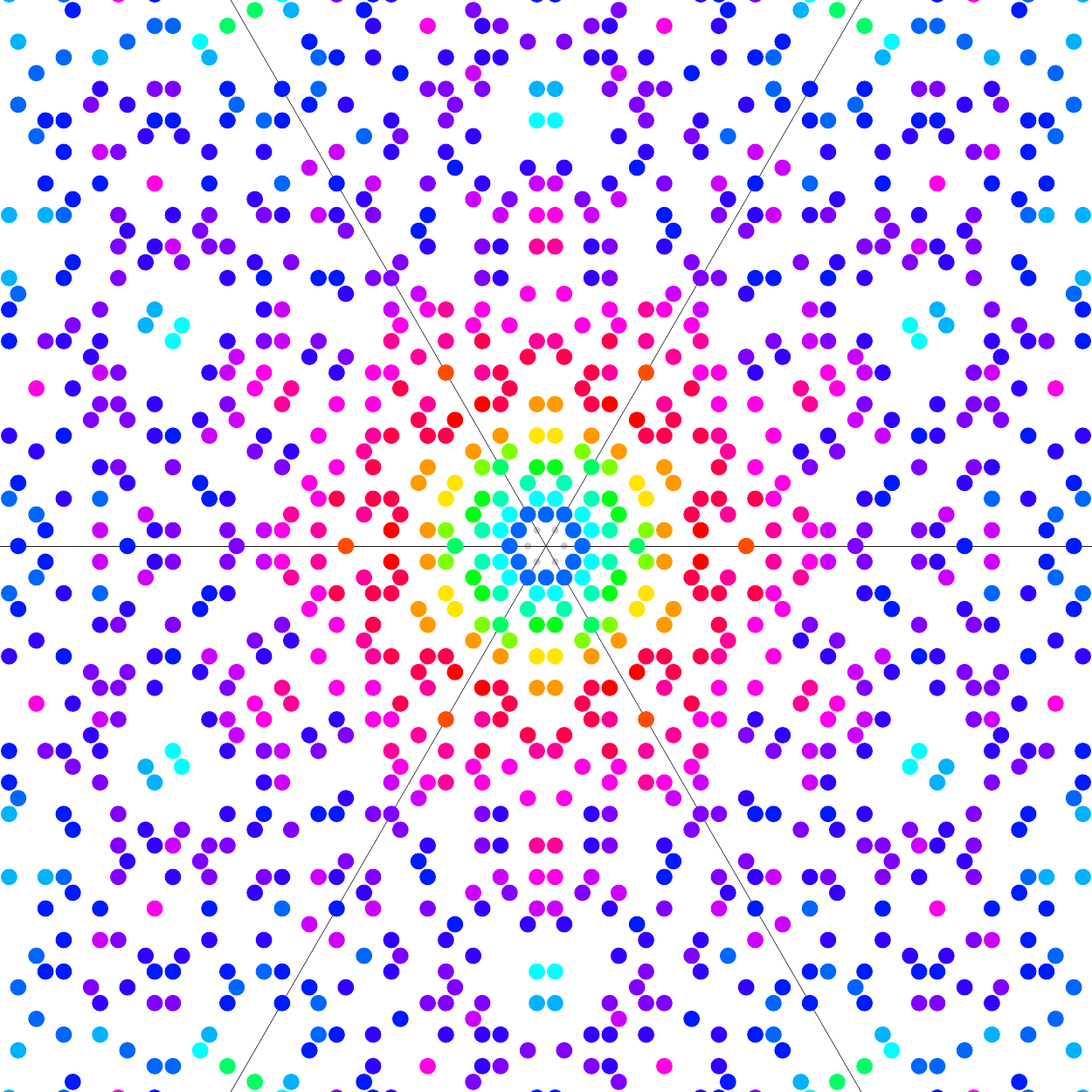}}
\caption{
The Eisenstein primes.
}
\label{circles}
\end{figure}

\begin{figure}[!htpb]
\scalebox{0.80}{\includegraphics{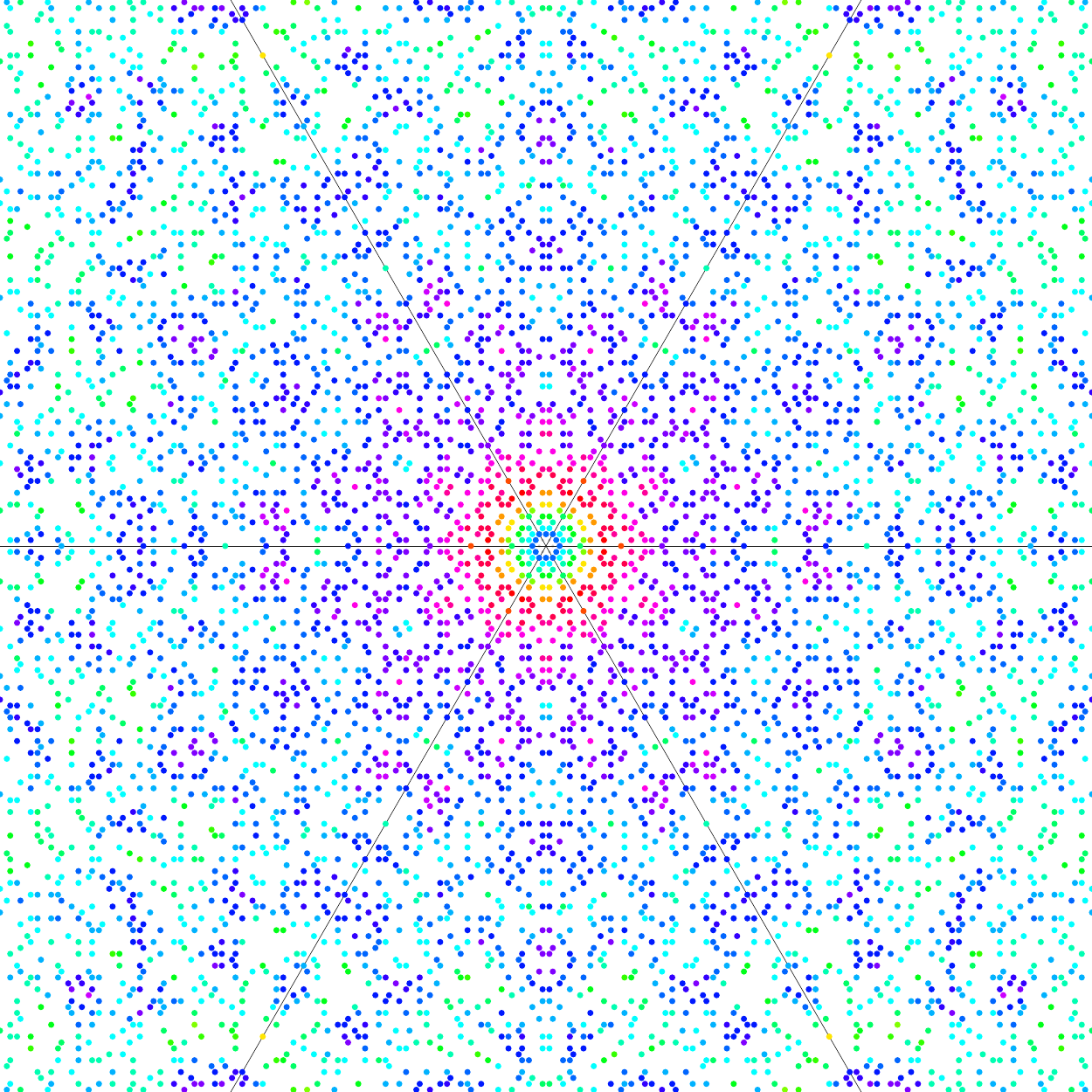}}
\caption{
A larger view on the Eisenstein primes.
}
\label{circles}
\end{figure}

\conjecture{
Except for two counter examples, every rational integer $n$ larger than $1$ is a sum $x+y$ where $p=x^2+x+1$
and $q=y^2+2y+4$ are both rational primes. The two counter examples are $n=109$ and $n=121$. 
}

Again, we do not have to worry that a positive proof of the Goldbach question for Eisenstein integers is
``easy" (except of course if there would exists an obvious counter example)
To see this, look at the boundary case and consider the Eisenstein integers
of the form $n+2w$. If the Goldbach conjecture is true, then every $n+4w$ with odd $n$
is a sum of two Eisenstein primes $a+w, c+3w$. 
But since $x^2+x=1$ is a cyclotomic polynomial, the Bunyakovsky condition
is satisfied. In other words, we must have infinitely many primes of the form $x^2+x+1$,
a problem which is considered similarly hard than the Landau problem: 

\resultremark{
If Goldbach conjecture holds for Eisenstein integers, then the Bunyakovsky conjecture holds for $x^2+x+1$. 
}

Goldbach is even stronger at the boundary than Bunyakovsky: 
the density of the primes of the form $n^2+n+1$ is so large that one can 
reach every rational integer larger than $2$:
one could call this the {\bf boundary Eisenstein Goldbach conjecture}: 

\conjecture{
Every rational integer $n>1$ is a sum $n=x+y$, where $x^2+x+1$ and 
$y^2+y+1$ are both rational primes. 
}

There is a surprise on the next column as there are two ghost counter examples. We call them 
the {\bf Eisenstein ghost twins} (not to confuse with Eisenstein prime twins, 
which are neighboring Eisenstein primes). 

\begin{figure}[!htpb]
\scalebox{0.4}{\includegraphics{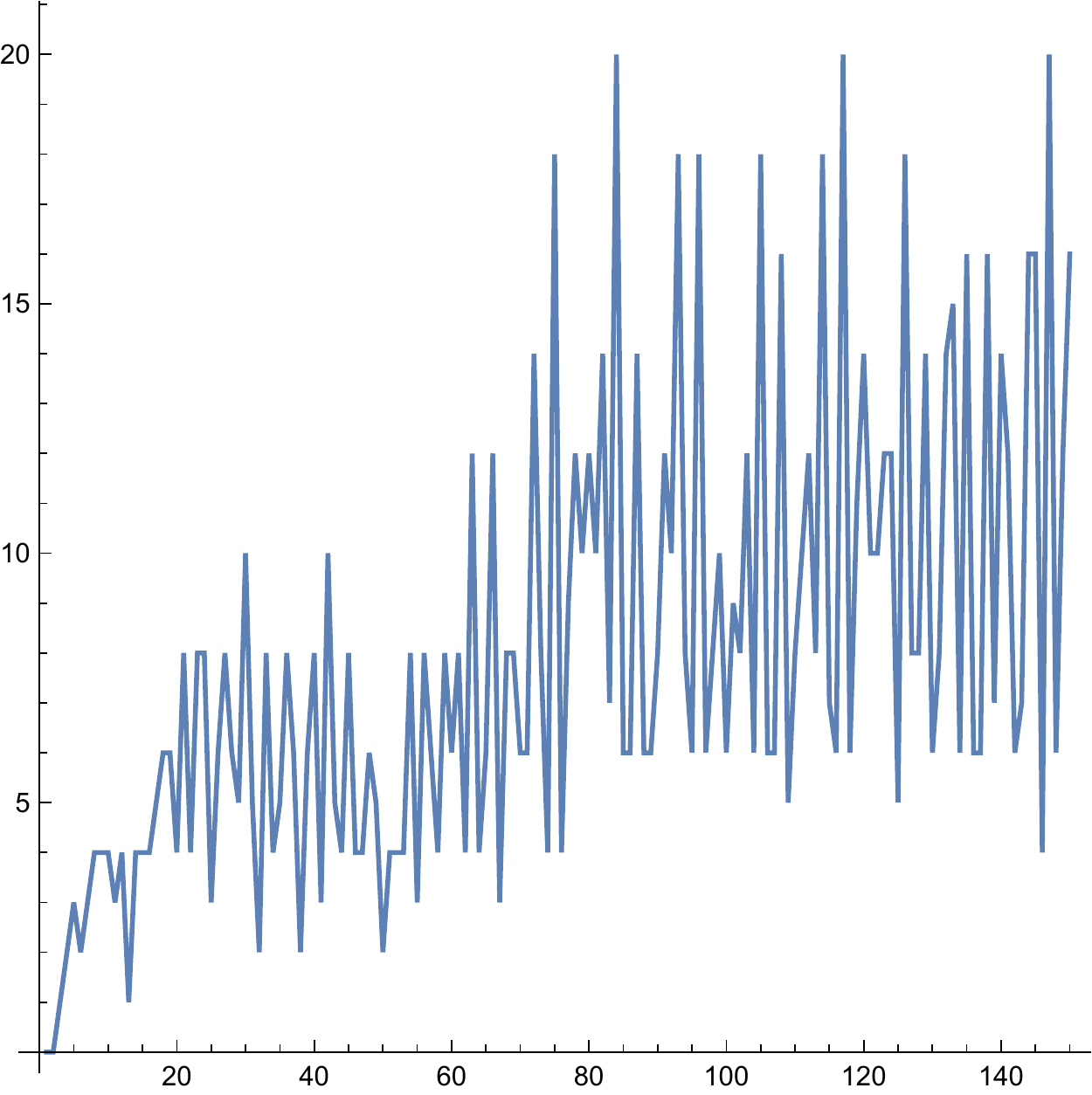}}
\scalebox{0.4}{\includegraphics{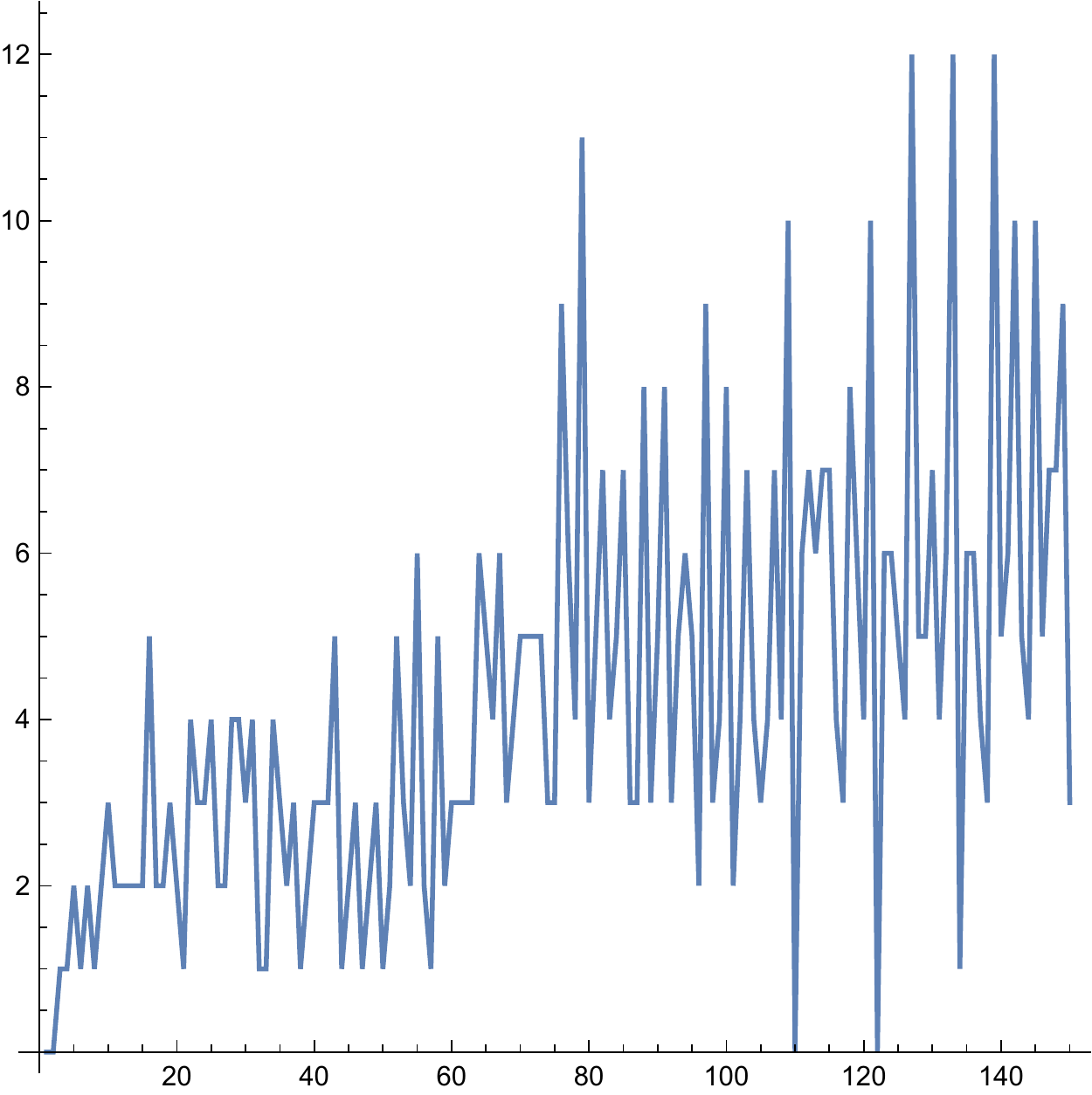}}
\scalebox{0.4}{\includegraphics{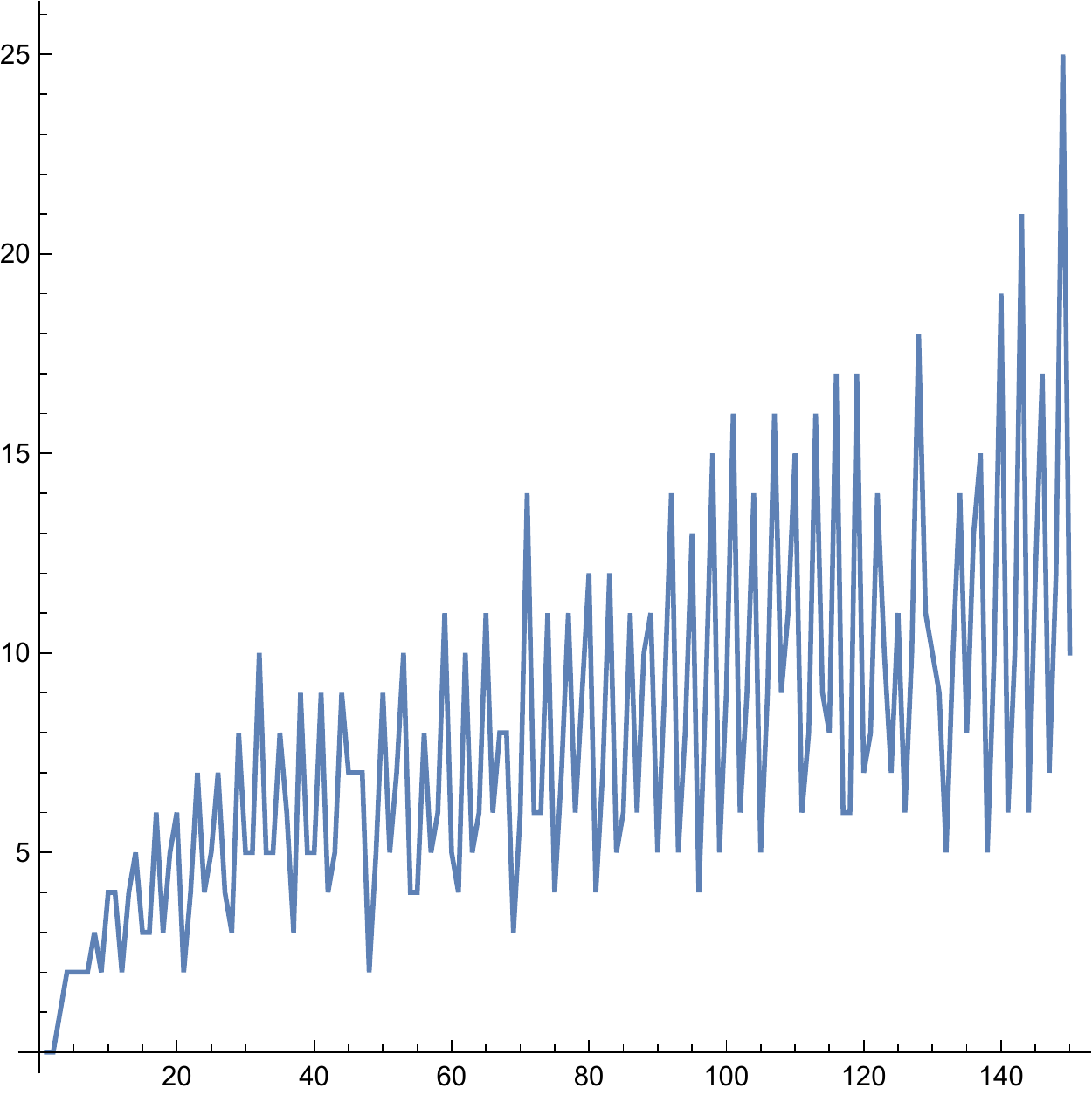}}
\scalebox{0.4}{\includegraphics{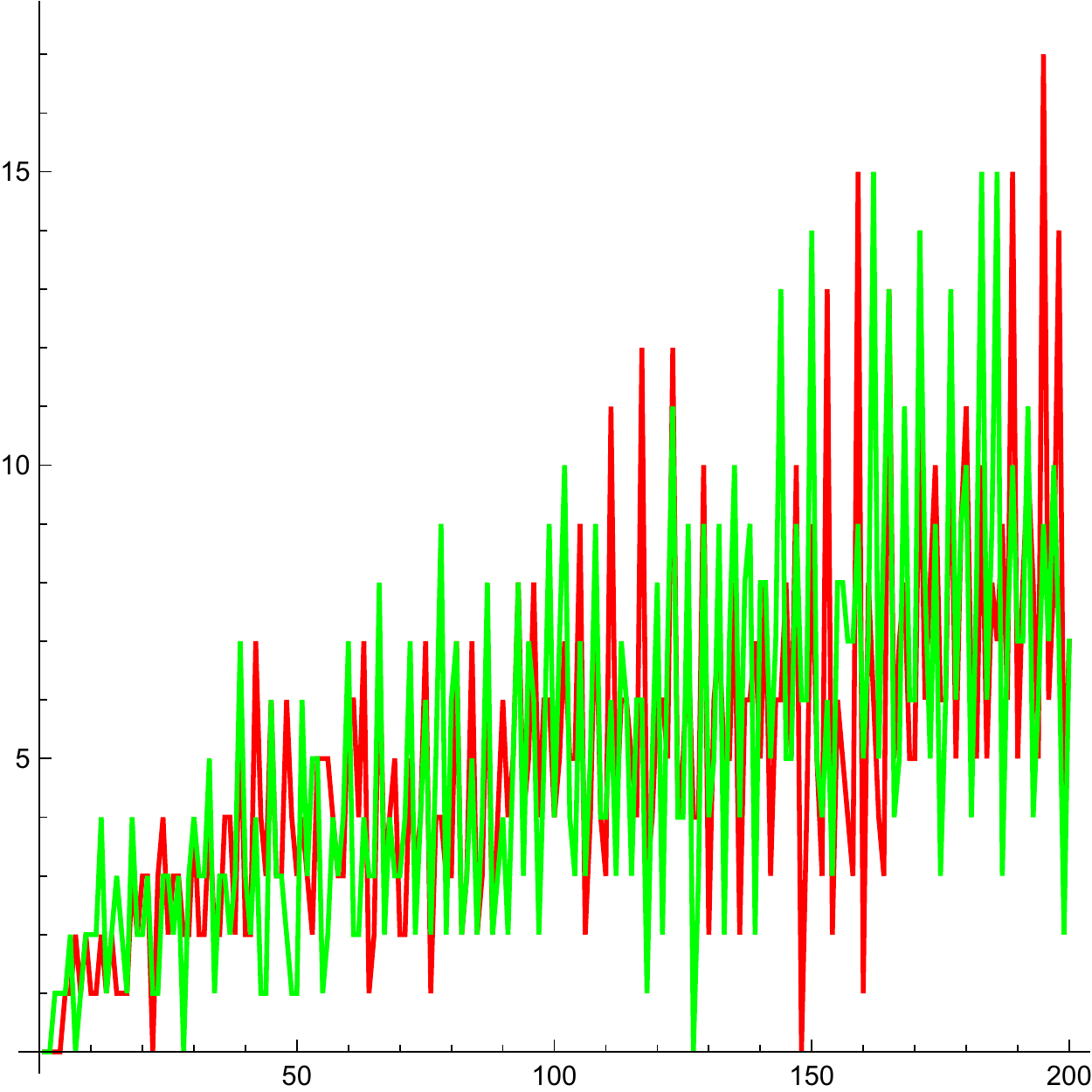}}
\scalebox{0.4}{\includegraphics{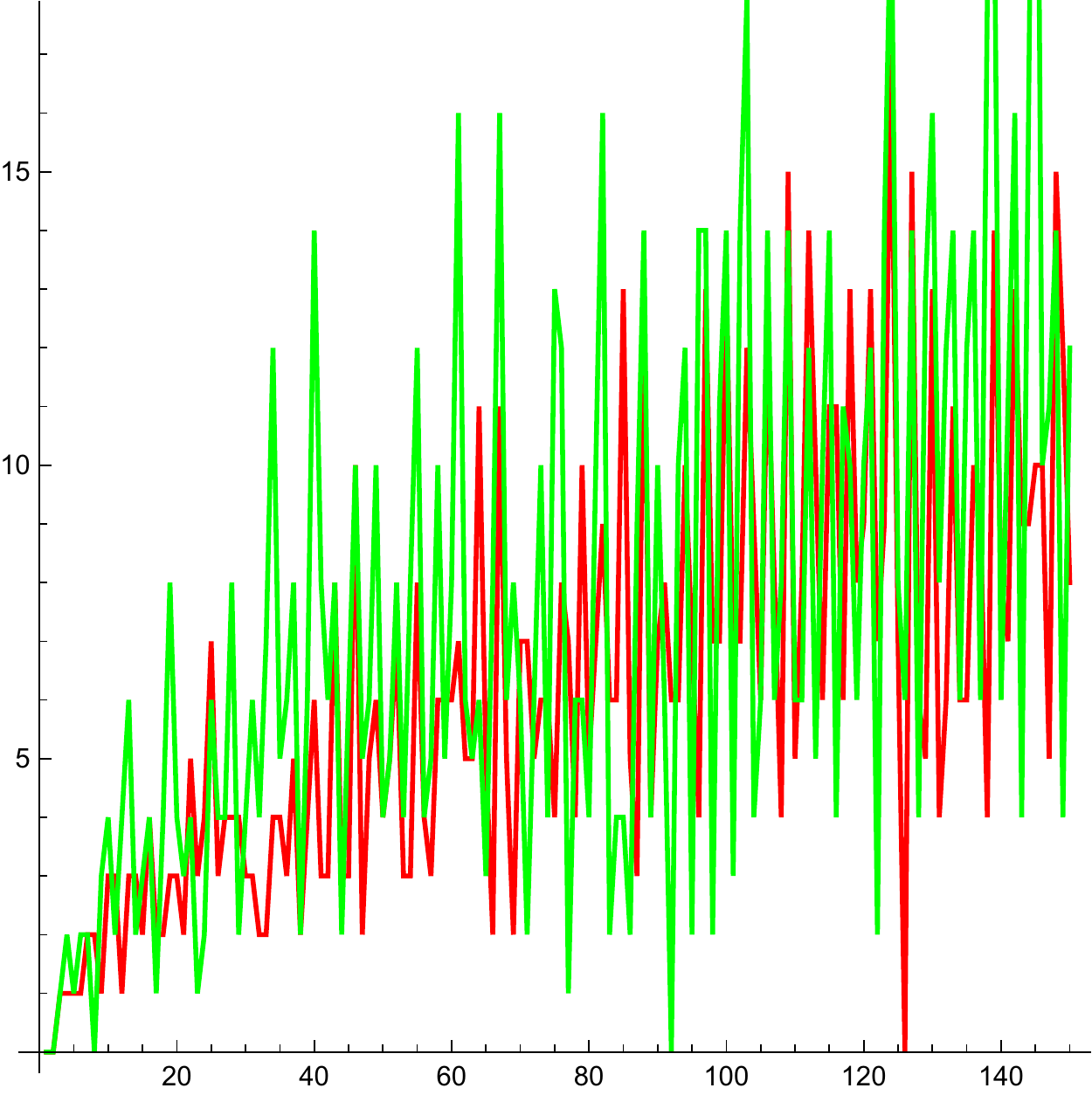}}
\scalebox{0.4}{\includegraphics{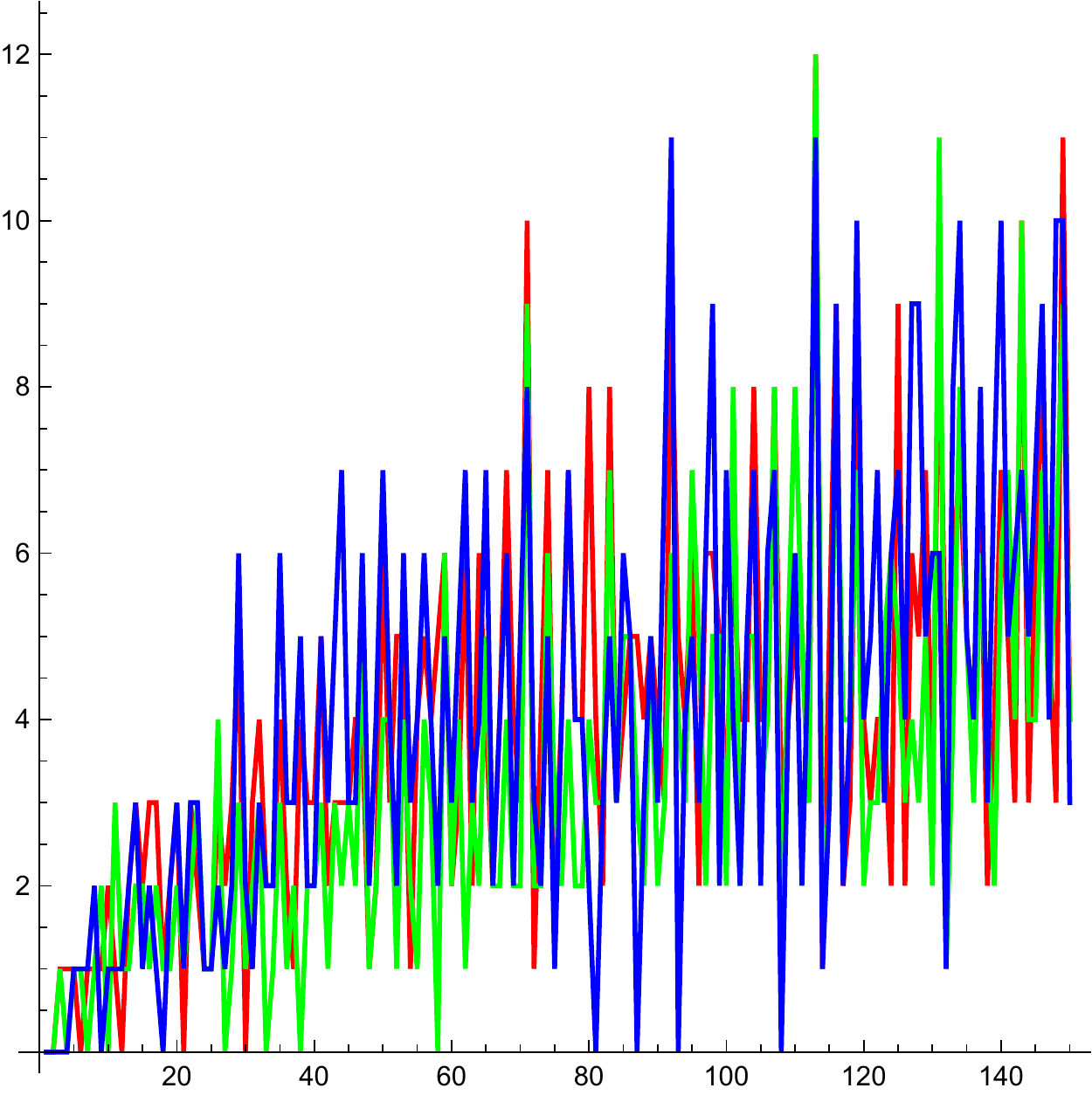}}
\caption{
The number of times an Eisenstein integer $a+ b w$ can 
be written as a sum of Eisenstein primes. For $b=2$, we
measure that this is always possible: there are rational primes of the
form $p=x+w, q=y+w$ which add up to $n=a+2w$. For $b=3$
we see two bad counter examples, the {\bf Eisenstein twins}
which are $3+109w$ and $3+121w$. For $b=4$ we appear already
fine by writing as a sum $p=x+w$ and $q=y+3w$. For
$b=5$ there are gaps again for pairs $p=x+w, q=y+4w$ as
well as pairs $p=x+2w, q=y+3w$, but they don't overlap.
For $b=5$ the decomposition $p=x+w$ and $q=y+4w$ is not
possible for $125+5w$ and $p=x+3w$ and $q=y+3w$ is not
possible for $91+5w$ but again they don't intersect. 
}
\end{figure}

When going further, it appears that every integer larger than $1$ is the sum of 
$n=x+y$ where $p=x^2+x+1$ and $q=y^2+3y+9$ are rational primes. 

On the next row $5+k w$ again, there are exceptional cases, 
of twins when writing $n$ as a sum of primes of the form $p=x^2+x+1$ and $y^2+4y+16$
where $n=21,147$ are not possible or in the form of $p=x^2+2x+4$ and $q=y^2+3y+9$
where $n=6,27,126$ are not possible. But now the two sets don't intersect and we 
appear always to be able to write the Eisenstein integer $n=5+kw$ as a sum of two 
Eisenstein primes, either in the form $n=(1+x w) + (4+ y w)$ or then in the form 
$n=(2+x w) + (3 + y w)$. As further away we are from the boundary, as more possibilities
we have and the chance gets smaller and smaller. 

{\bf Remark}:
{\bf 1)} Having seen so many Goldbach versions, one can ask for more. Goldbach conjectures
have been formulated in polynomial rings. One could also change the field and look
for example at  p-adic fields or then use an other division algebra. By a theorem of Frobenius, 
there are only the quaternions and octonions left. One could look at number fields and rings 
of integers inside the quaternions and octonions. Also to this there is literature but
one has first to make sense of ideal theory in non-commutative and non-associative setups.  \\

Lets at last look at the Zeta function in this case. 
The {\bf Eisenstein zeta function} is an example of a Dedekind zeta function. 
Since the norm of $z=n+m w$ is $n^2+nm+m^2$, it is 
$$ \zeta_E(s) = \sum_{(n,m) \neq (0,0)} \frac{1}{(n^2+nm + m^2)^s} \; . $$
As Eisenstein integers form a factorization domain, the {\bf Euler product formula}
$$  \zeta_E(s) = \prod_{p \in P} (1-N(p)^{-s})^{-1}   $$
still holds, where $p$ runs over all Eisenstein primes as well as
$$  \zeta_E(s)_{-1} = \sum_{(n,m) \neq (0,0)} \frac{\mu(n+m w)}{(n^2+nm + m^2)^s}   $$
where $\mu(z)$ is the {\bf M\"obius function} on Eisenstein integers, encoding the
parity of the prime factorization. Again, as realized by Dirichlet for quadratic number fields:

\begin{figure}[!htpb]
\scalebox{0.12}{\includegraphics{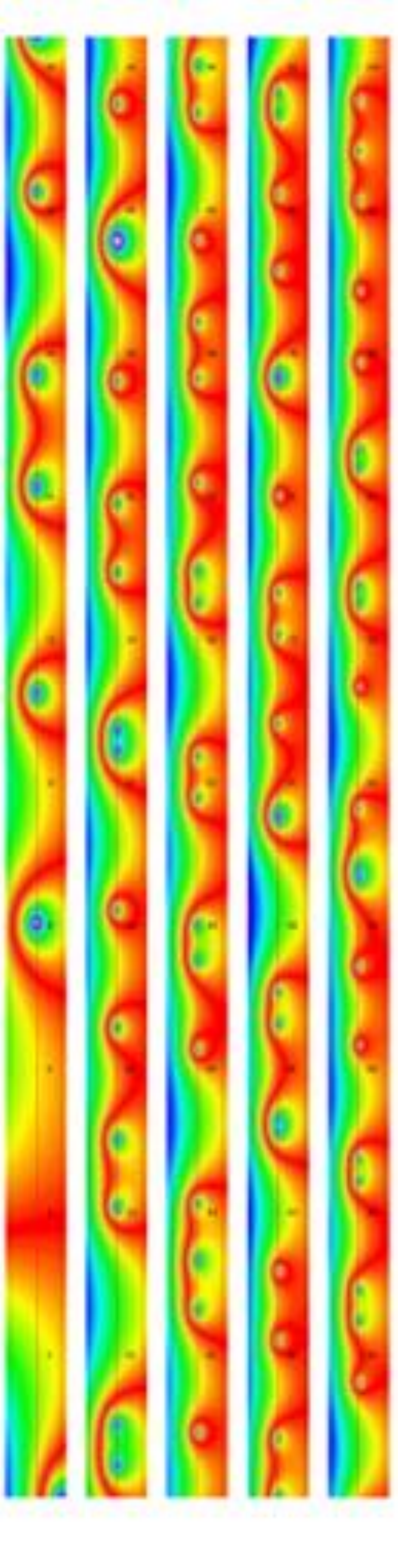}}
\caption{
A contour plot of the {\bf L-function} $\zeta_E(z)$  of the
Eisenstein primes. Part of the generalized Riemann hypothesis 
predicts that all roots are on the critical line. Due to the 
product property, it would imply the classical Riemann hypothesis. 
}
\end{figure}

\begin{figure}[!htpb]
\scalebox{0.2}{\includegraphics{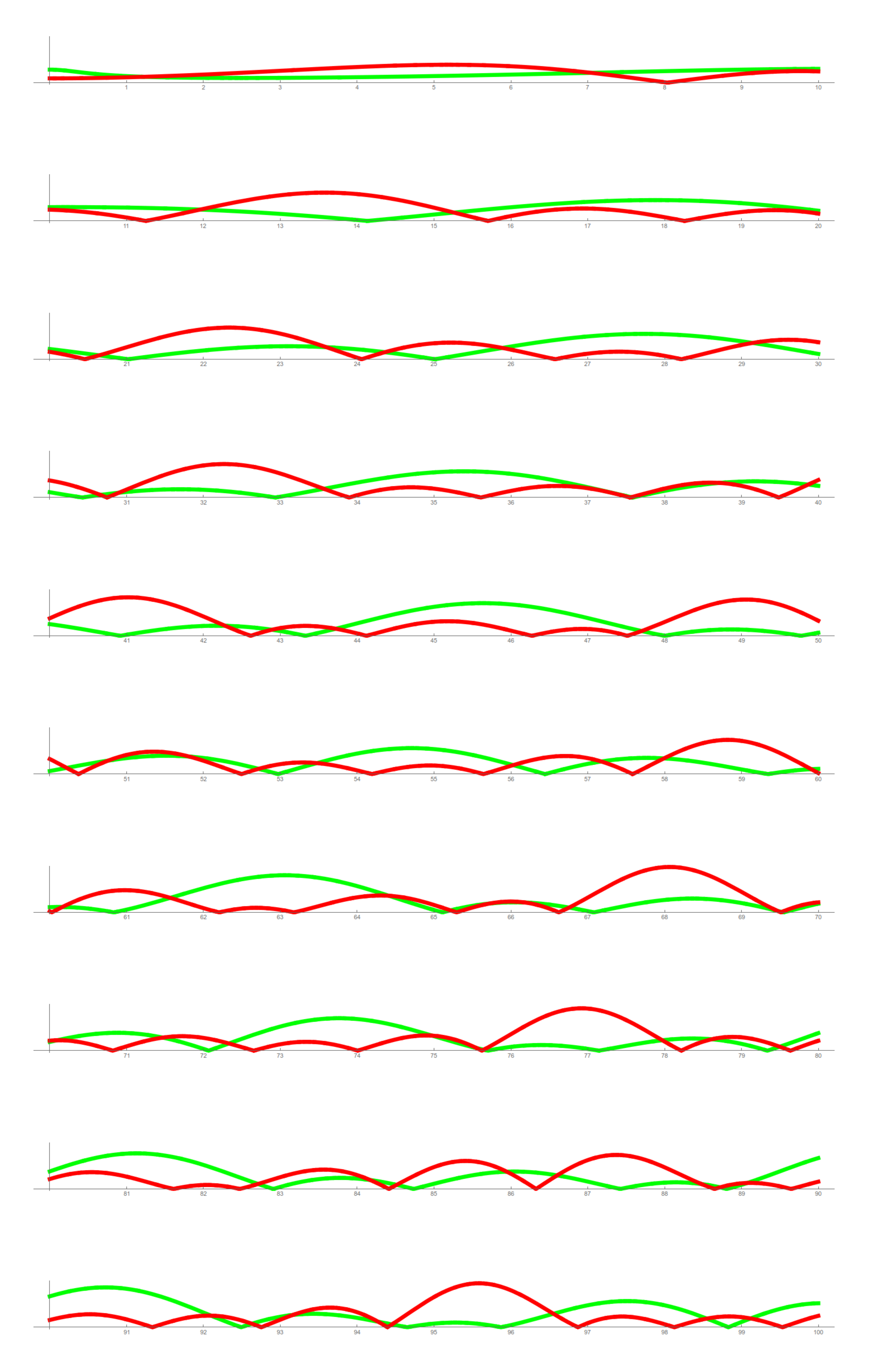}}
\caption{
The absolute value of the $\zeta$ function for the rational integers 
and the zeta function for the Eisenstein integers. 
}
\end{figure}

\resultremark{(Dirichlet)
There is a factorization
$$ \zeta_E(s) = 6 \zeta(s)  \beta_{3,2}(s)  \; . $$
}
where $\beta_3(s)$ is a $L$-function with character $\chi_{3,2}(n)=1$ if $n \in P_1$ 
and $\chi(n)=-1$ for $n \in P_2$ and $\chi(3k)=0$. The verification is the same. 

\section{A physics allegory}

Since the unit sphere in quaternions is $SU(2)$ which relates to 
{\bf Pauli matrices} in various ways, which plays a role in fundamental interactions,
there has been early motivation to apply {\bf particle phenomenology} in division algebras. 
But currently, the standard model is still king and no quaternion nor octonions are
needed. In order not to be misunderstood, we don't think the following will change
that. But we believe the following allegory illustrates that using an accelerator to smash 
{\bf building blocks of matter} onto each other and experimenting 
with the {\bf building blocks of the integers} is closely related.
One should see it more as {\bf entertainment}. We firmly believe that any physical 
theory of value must predict or explain something which can not be done otherwise. 
However the story which follows, renders some of the number theory in the complex and
in the quaternion case more pictorial: \\

{\bf Leptons}: we think of Gaussian integers as
a collection of {\bf leptons}, where the individual
Gaussian primes are indecomposable {\bf Fermions}. A prime $z$ of type $4k+1$ together with an
opposite charged particle $\overline{z}$ form an {\bf electron-positron} pair. We gauge the integers
with units $U=\{1,i,-1,-i\}$ to be in the sector $\pi/2<{\rm arg}(z) \leq \pi/2$. The {\bf charge}
of a lepton $z$ is defined as the sign of the argument of $z$ in the branch $(-\pi,\pi]$. 
Think of the logarithm of the norm $N(z)$ as {\bf mass}. A rational positive prime $4k+3$ is called a {\bf neutrino} 
as it is lighter: its momentum $|z|$ is prime while for $4k+1$ primes which are electrons or
positrons,  the energy $N(z)$ is prime. A neutrino is {\bf neutral} as it is located 
on the real axes. The largest known prime for example is a Mersenne prime and so a neutrino.
An integer $n=p_1 \dots p_k$ is a {\bf Lepton configuration}. The fact that the Gaussian primes 
form a {\bf unique factorization domain} translates into the statement that any lepton configuration 
can be decomposed uniquely into such leptons as well as a neutral mystery particle $2$ which is its
own anti-particle. The uniqueness holds only modulo {\bf gauge transformations} which act here as
multiplications by units. We will see that in the Quaternion case, this fact is no more the case, 
because that is, where the Hadrons will come in explaining why quarks form Baryons and Mesons.
The electron-positron pair is {\bf not bound} together: there is no unit which maps one into the other.
Factoring out the symmetry of units renders the factorization unique. The product $(-3) (-7)$ for
example is gauge equivalent to the product $3 \cdot 7$. Let us now move from primes to rational integers and
call a {\bf rational integer} $n \in \mathbb{Z}$ a {\bf Boson configuration} if it contains an even
number of Fermionic prime factors counted with multiplicity, otherwise it is a Fermion. Mathematically,
$n \in {\bf N}$ is a Fermion if and only if its {\bf Jacobi symbol} 
$(-1|n)=\left( \frac{-1}{n} \right)$ is $-1$.
Otherwise, if $(-1|p)=1$, it is a Boson. The Gauss law of {\bf quadratic reciprocity result} 
tells now that two odd primes $p,q$ satisfy the {\bf commutation relations}
$(p|q) =(q|p)$ if at least one of them is a Boson and that the
{\bf anti-commutation relation} $(p|q)=-(q|p)$ hold exactly
if both primes $p,q$ are Fermions. In other words, if we look at the 
Jacobi symbol as an operator $p \cdot q$, then Bosons commute with everything else, but the sign changes,
if we switch two Fermions. 
The {\bf two square theorem Fermat} telling that an integer $n$ can be represented as $a^2+ b^2$
if and only if $n$ is a Bosonic integer can be interpreted as the fact that a Bosonic rational prime
defines two leptons $a+ib,a-ib$, where $a^2+b^2=p$. The positron and electron are anti-particles
of each other, but they are not equivalent since there is no gauge from one to the other. If we factor out the
{\bf gauge symmetries} given by the {\bf units}, then the factorization aka particle decomposition of the
Lepton set is unique. This is the
{\bf fundamental theorem of arithmetic} for Gaussian integers. It can be proven from the rational
case using the $1-1$ identifications on the orbifold $\mathbb{C}/D_4$ so that we get the rational primes. 
The {\bf Pauli exclusion principle} is 
encoded in the form of the {\bf Moebius function} $\mu_G(n)$ which is equal to $1$ if a Gaussian integer 
$n$ is the product of an even number of different Gaussian primes, and $-1$ if it is the product 
of an odd number of different Gaussian primes and $0$, if it contains to identical particles. 
Again this {\bf particle allegory} is already useful as a {\bf mnemonic} to
remember theorems like the two square theorem, or the quadratic reciprocity theorem:
"Quadratic Reciprocity means that only Fermion primes anti-commute $(p|q)=-(q|p)$.
Fermat's two square theorem assures that Bosonic rational primes $p=a^2+b^2$
are composed of two Gaussian primes $a \pm i b$. The others are all real, light and neutral."  \\

\begin{figure}[!htpb]
\scalebox{0.8}{\includegraphics{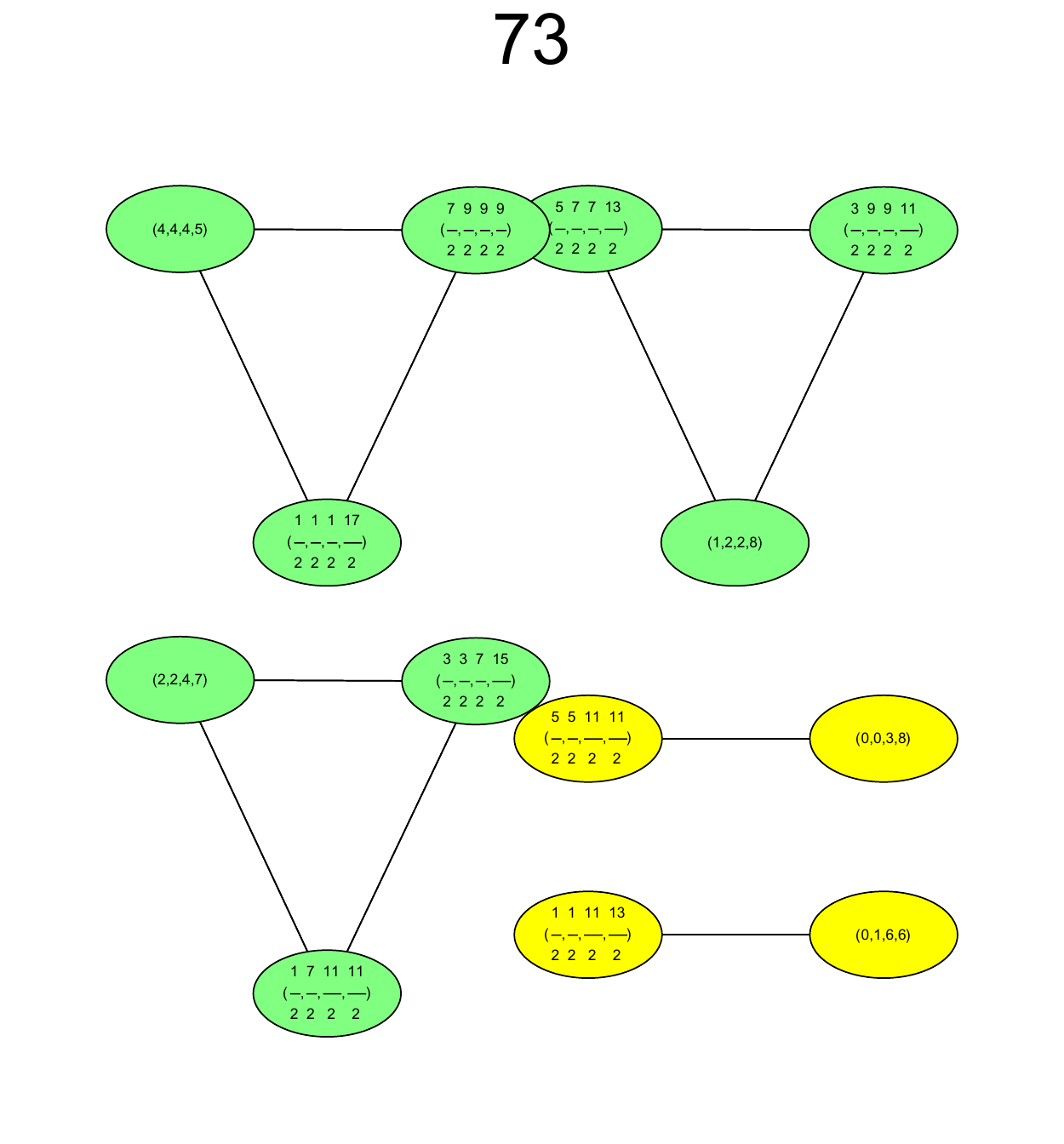}}
\caption{
The hadrons belonging to the prime 73. There are three Baryons and 
two Mesons.  
}
\label{HardyLittlewood}
\end{figure}

\begin{figure}[!htpb]
\scalebox{0.8}{\includegraphics{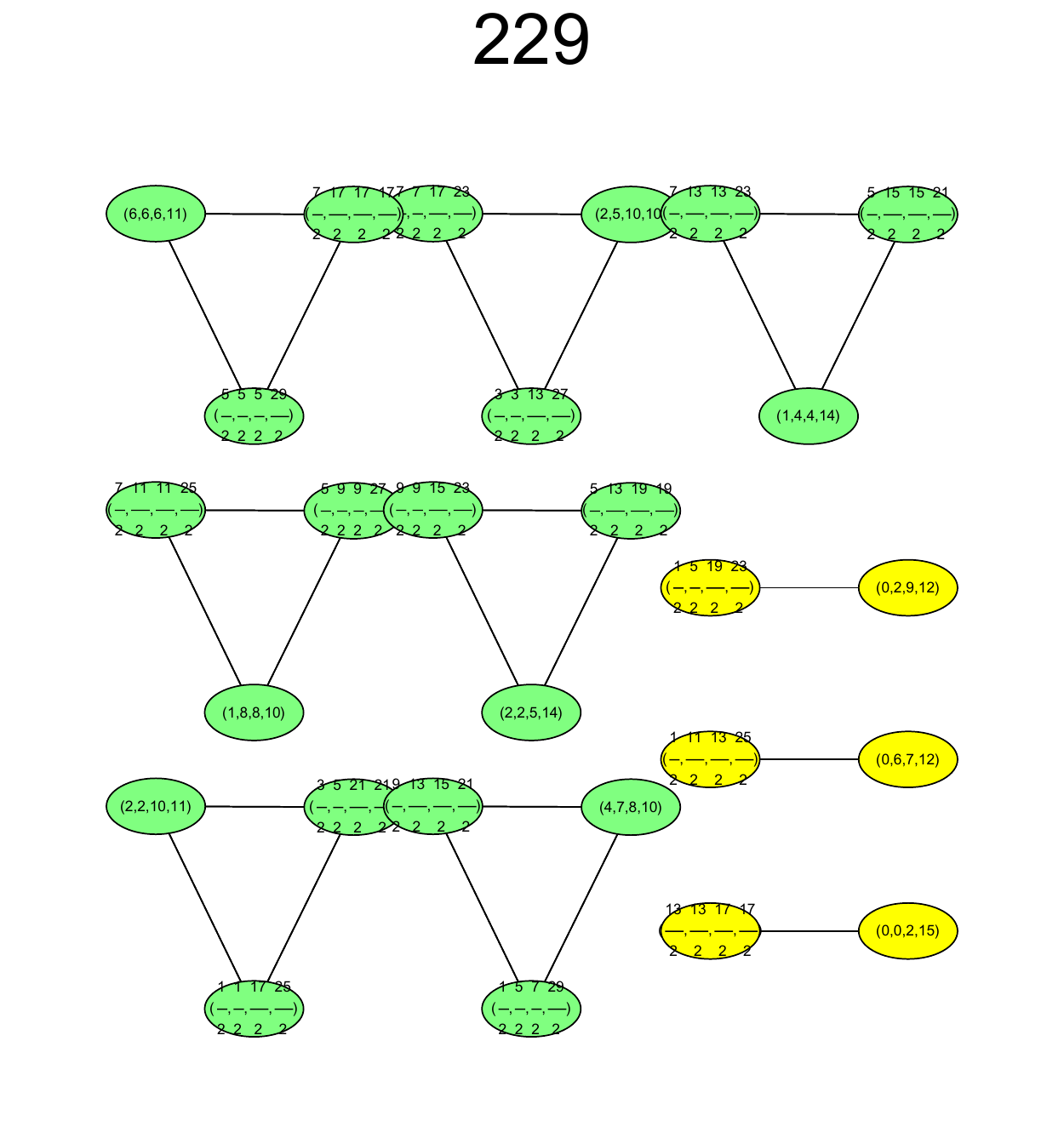}}
\caption{
The hadrons belonging to the prime 229 features 7 Baryons and 3 mesons. 
}
\label{HardyLittlewood}
\end{figure}

{\bf Hadrons:} \\
Hadrons are quaternion primes $z$ with norm $N(z)$ different from $2$. The prime $2$ is special 
also in quaternions. As Hurwitz already showed that despite non-commutativity, one can place them
outside: every integer quaternion $z$ is of the form $z=(1+i)^r b$ where $w$ is an odd 
integer quaternion. There are two symmetry groups acting on hadrons. 
One is the group $U$ of units, the other is the group $V$ generated by coordinate permutations and conjugation. 
The group $U$ has 24, the group $V$ has $48$ elements. The groups are no more
contained into each other like in the complex case, where $U=\mathbb{Z}_4$ was a subgroup of $V=D_4$. If
we look at the $U$ equivalence classes first and then look at the orbits of $V$, we see that some
particles are fixed under $V$ or then that they come in pairs. We will interpret this as {\bf particle
anti-particle pairs}.  We can however also look at the $V$ equivalence classes first and then 
look at the orbits of $U$, then we see 1 or 2 or 3 particles combined. This is remarkable. 
We call a quaternion prime equivalence class with is not a single element
a {\bf hadron}. Each hadron to an odd prime $p$ consists either of two or three 
{\bf quarks}. Quarks can be either Lipschitz or Hurwitz primes. Each hadron contains a 
single Lipschitz prime. 
The {\bf Lagrange four square theorem} assures that there are no neutrini particles among Hadrons.  
So, we can think of Quaternion primes as {\bf quarks} which form equivalence classes in the form of 
Baryons or Mesons. The conjugate quaternion is the {\bf anti particle}. This can be seen from 
the Hurwitz factorization theorem which shows that in such a factorization, one never can have
factors $z \overline{z}$ near each other. The allegory is that they would ``annihilate" into a shower
of smaller particles. Going from one factorization to an other is a rather complex interaction process
changing the nature of some Baryons involving gauge bosons. 
Baryons are Fermions and Mesons are Bosons. Like in the complex case, we have 
a mystery $p=2$ case, which has only one equivalence class. Lets call it the $2$ particle even so we
would like to associate it with something real like Higgs because it is neutral, light, its own 
anti-particle and can give more mass to other particles by multiplying with it. 
The {\bf group of units} contains $8$ particles of the form 
$(\pm 1,0,0,0),(0,\pm 1,0,0), (0,0,\pm 1,0),(0,0,0,\pm 1)$ which modulo $U$ are all equivalent
to the neutral $(1,0,0,0)$, the $Z$-Boson. 
There are $16$ remaining units. Modulo $V$ they are all equivalent to $(1,1,1,1)/2$. This has
a positive charge and is the {\bf $W^+$ boson}. Its conjugate is the {\bf $W^-$ boson}.
Now lets look at a Meson $(pq)$ containing a Lipshitz prime $p$ and Hurwitz prime $q$. 
Since $p$ and $q$ are equivalent in $U$, there exists a permutation, possibly with a
conjugation, such that $\overline{p}$ is gauge equivalent to $q$. If a conjugation is 
involved, then $p$ and $q$ have the same sign of charge, otherwise opposite. Because there
is a Lipschitz prime involved, it is not possible that all three quarks have the same
charge. So, two must have one charge, and one the other. 
Lets postulate that the {\bf charge} of a Lipschitz quark is $\pm 2/3$. 
Since we have identified modulo $V$ we can assume that it is positive. 
In the Meson case, the charge of the other particle is $1/3$ if it has the same charge 
and $-2/3$ if it has opposite charge. 
In the Baryon case, if the two other quarks have the same charge sign, they have charge $-1/3$.
If one has the same charge than the Lipschitz one, then both have charge $2/3$ and the other $-1/3$. 
The charge of an equivalence class is the sum of the charges. We have now assigned a charge in a gauge 
invariant way: a Lipschitz quark $(a,b,c,d)$ has charge $+2/3$ if $a \leq b \leq c \leq d$ and $-2/3$
if it is obtained from that by switching two coordinates. The structure of the equivalence classes 
assures a compatible choice so that the total charge is an integer. 
In the Meson case, we observe that one of the Lipschitz primes is 
located on the three or two dimensional coordinate plane. If a coordinate is zero, then we can
not perform flipping operations and we have to see whether we have to flip it in the {\bf Gaussian 
subplane}. In this picture, all Hadrons have charge $0,-1,1$. There are no Hadrons of charge $2$. 
\footnote{Particle physicists mention {\bf beta uuu-hadrons} but seem not detected in 
experiments. Also strange quark matter consisting of more than 3 quarks have not been observed.}

The non-uniqueness of prime factorization allows us to see going from one to the other factorization
as a {\bf particle process}. It involves the gauge bosons. 
It only becomes only unique modulo {\bf unit migration}.
(see \cite{ConwaySmith}). This means that if $x$ is a Hurwitz integer and $N(x)=p_1 \dots p_n$
then $x= (P_1 u_1) (u_1^{-1} P_2 u_2)  \dots  (u_{n-1}^{-1} P_n)$, where the $u_j$ are units
and $P_i$ are Hurwitz primes. In other words, the factorization becomes unique if
we look at it on the Meson/Baryon scale but it depends on the order. It becomes unique
when including meta-commutation: the prime factorization of a nonzero Hurwitz integer
is unique up to meta-commutation, unit migration and recombination, the process
of replacing $P \overline{P}$ with $Q \overline{Q}$ if $P,Q$ have the same norm.  \\
Let us look at some examples:

\begin{figure}
\scalebox{0.30}{\includegraphics{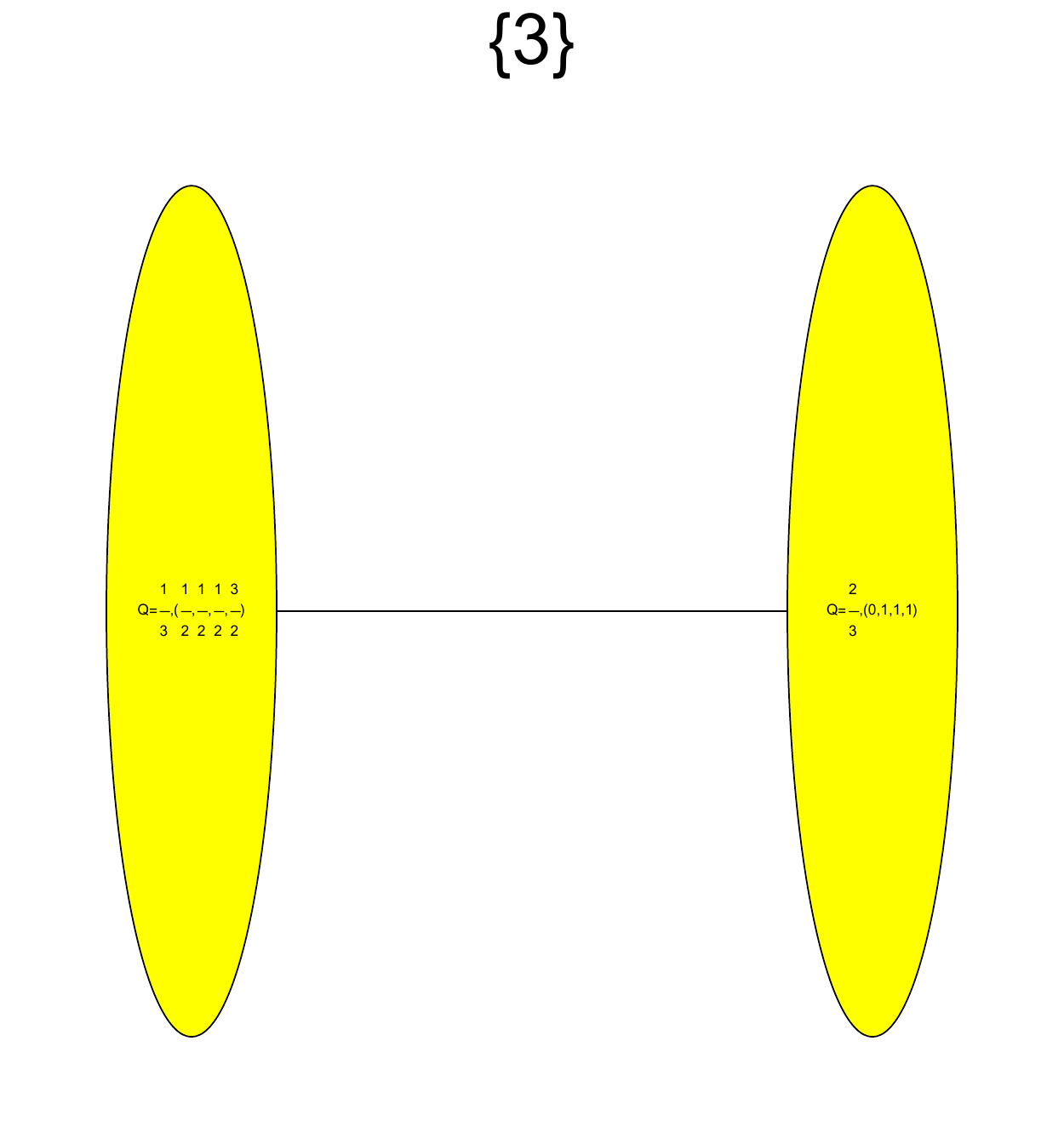}}
\scalebox{0.30}{\includegraphics{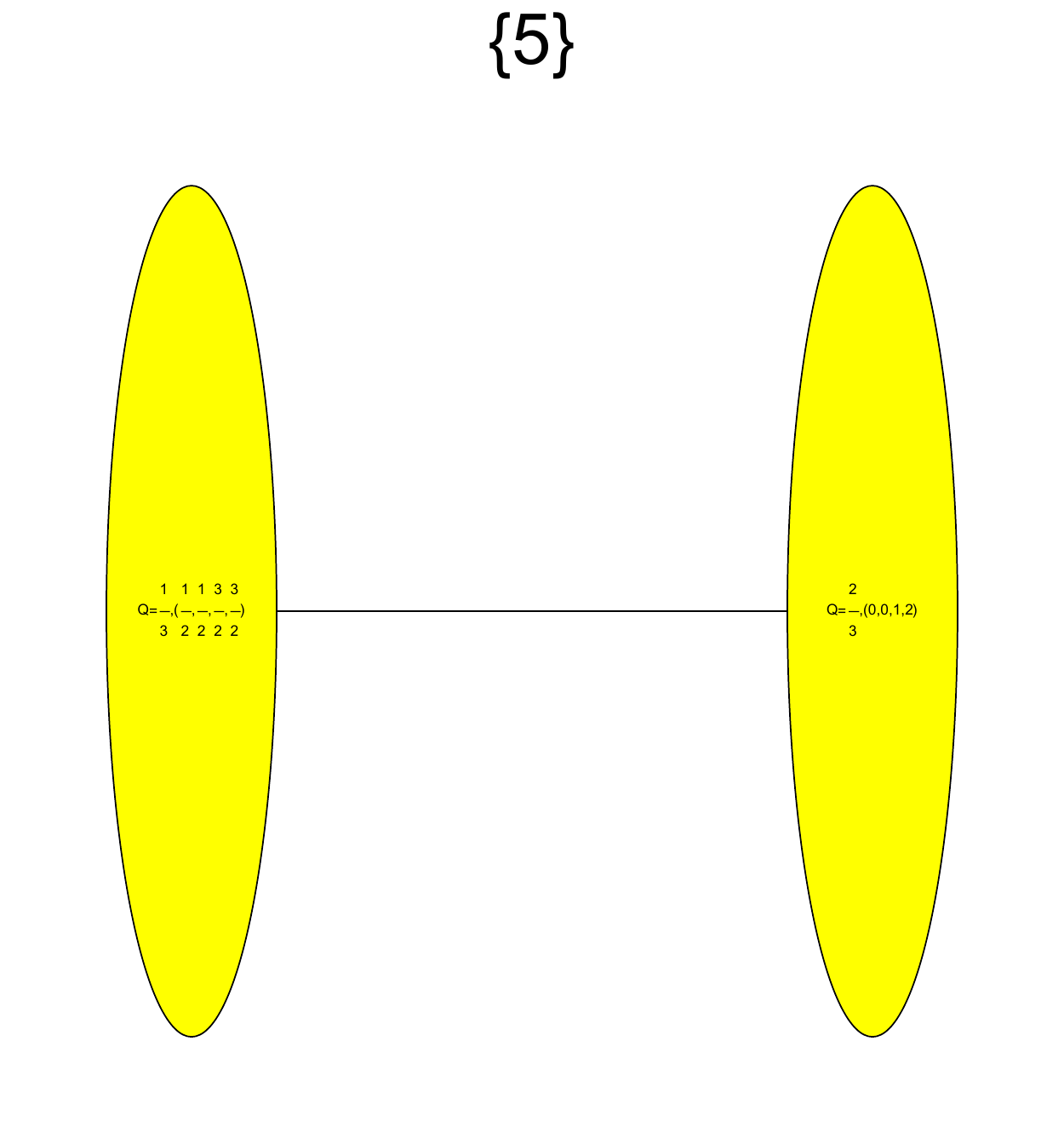}}
\scalebox{0.30}{\includegraphics{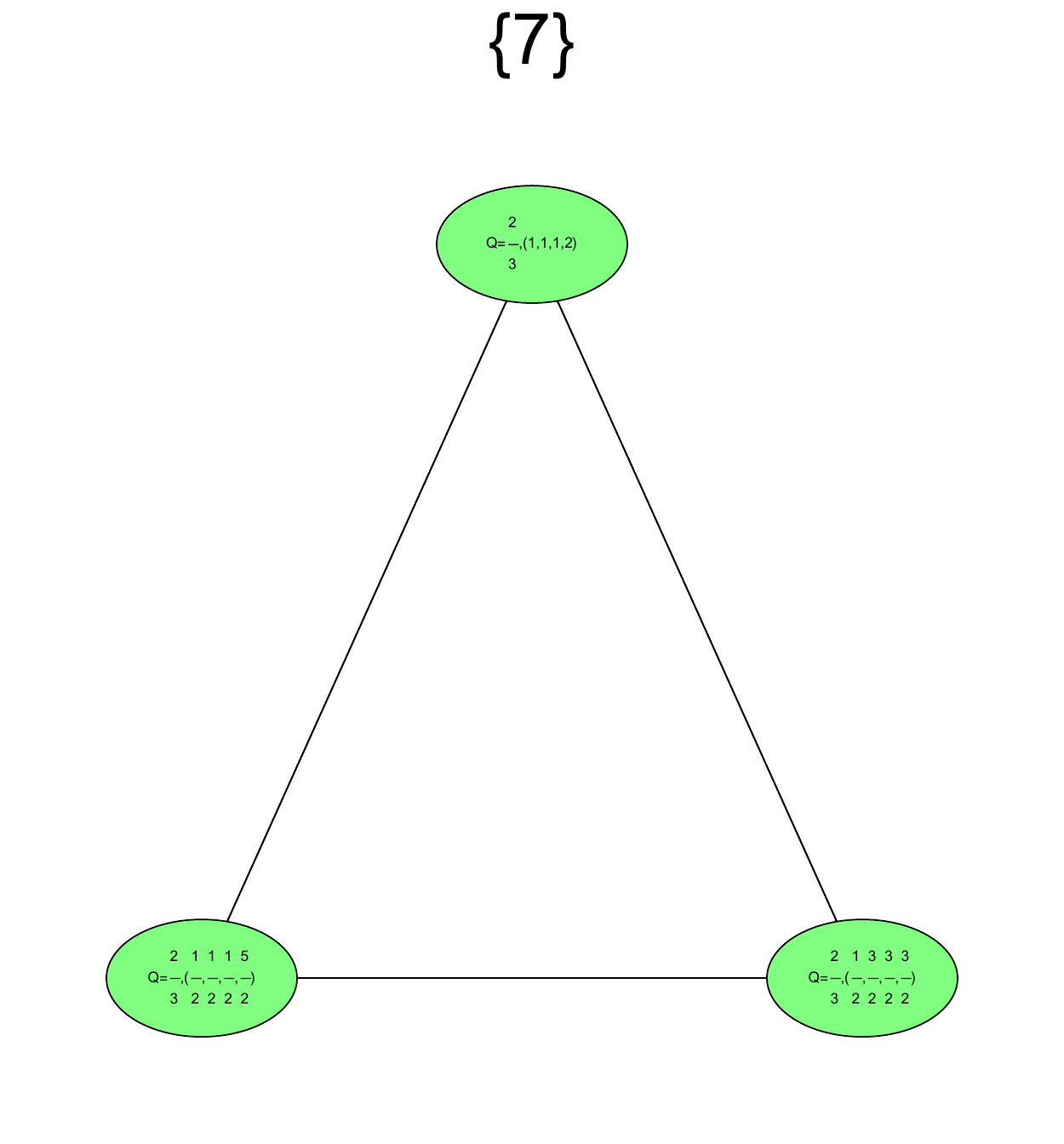}}
\scalebox{0.30}{\includegraphics{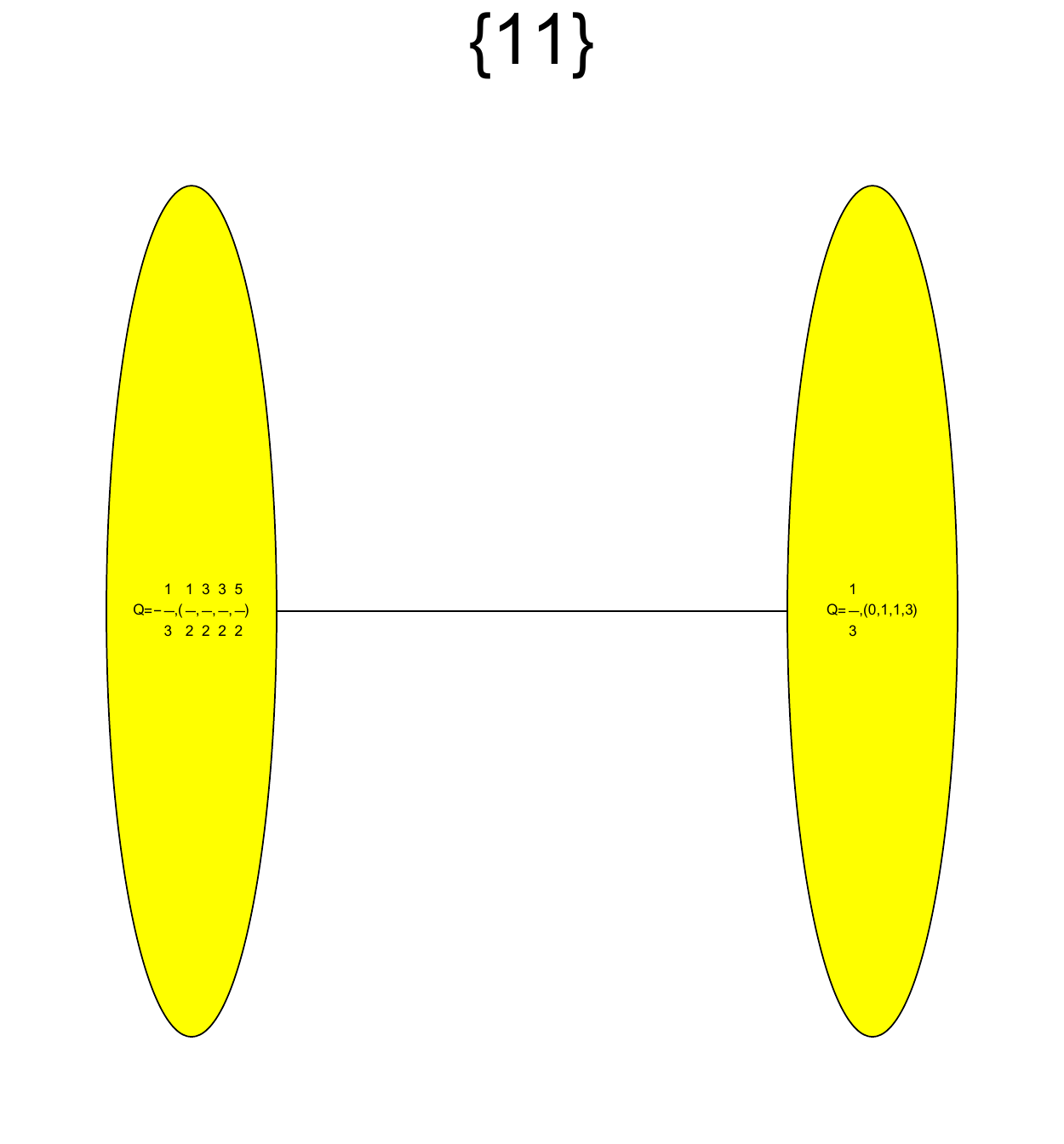}}
\scalebox{0.30}{\includegraphics{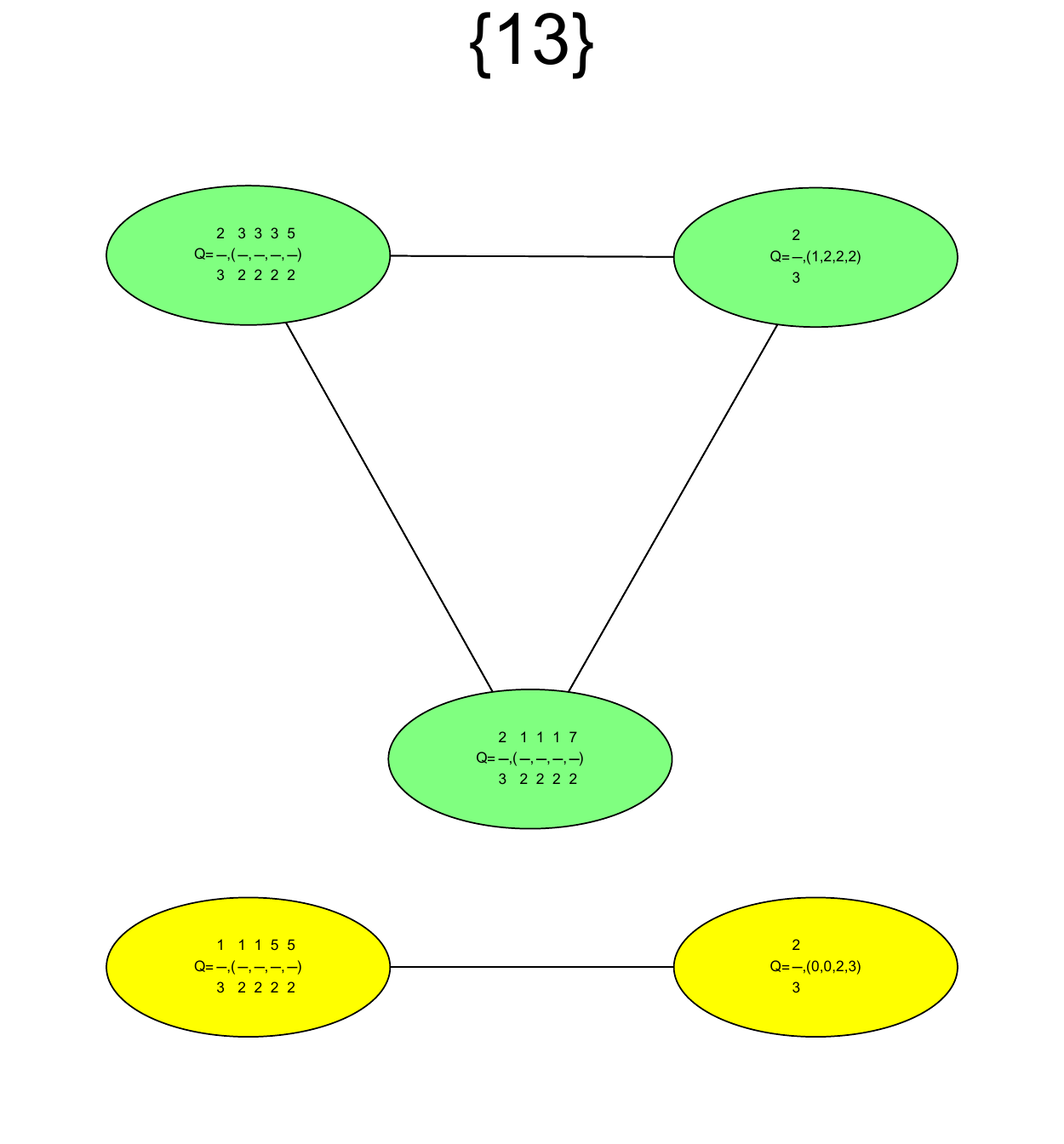}}
\scalebox{0.30}{\includegraphics{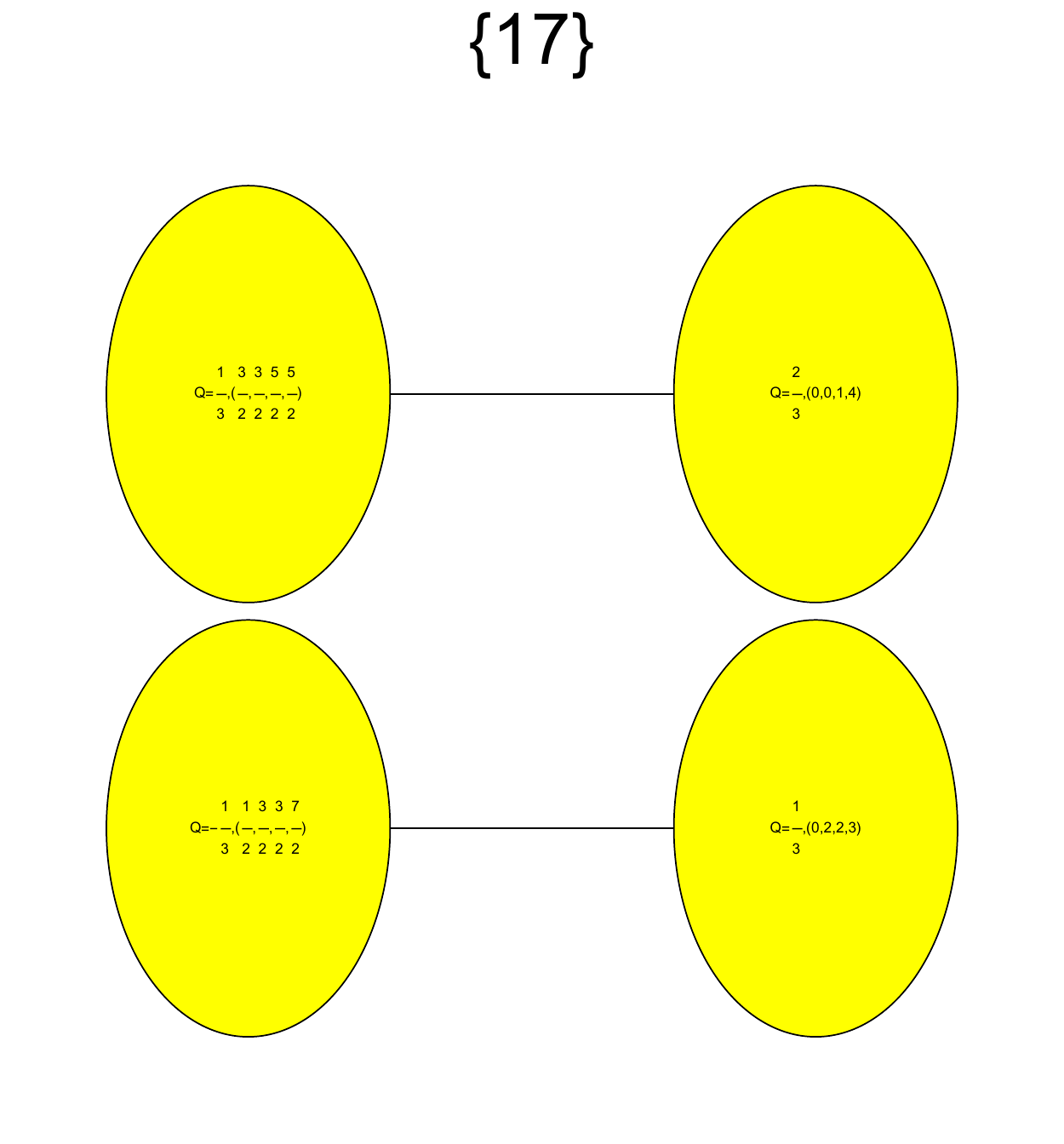}}
\scalebox{0.30}{\includegraphics{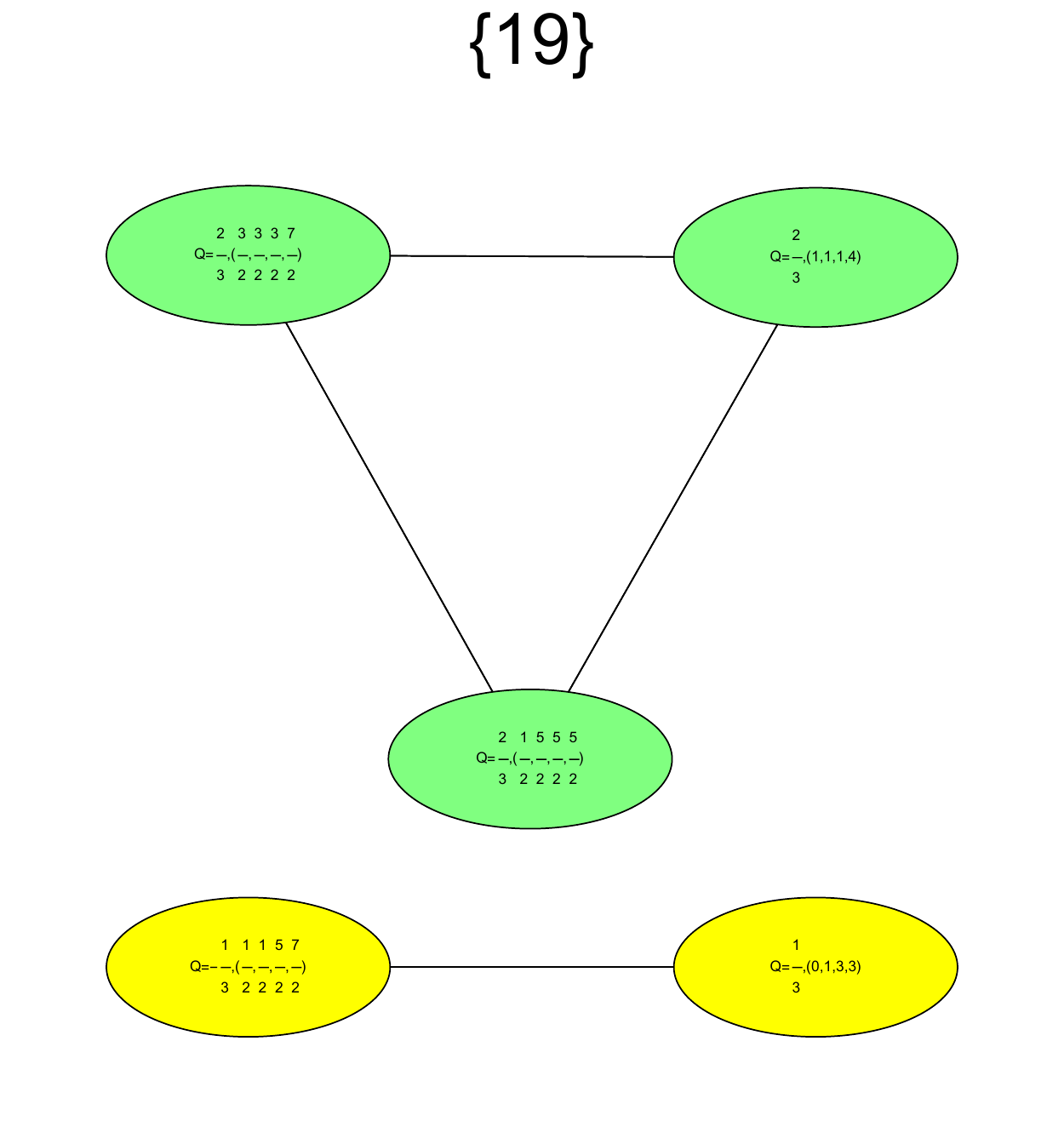}}
\scalebox{0.30}{\includegraphics{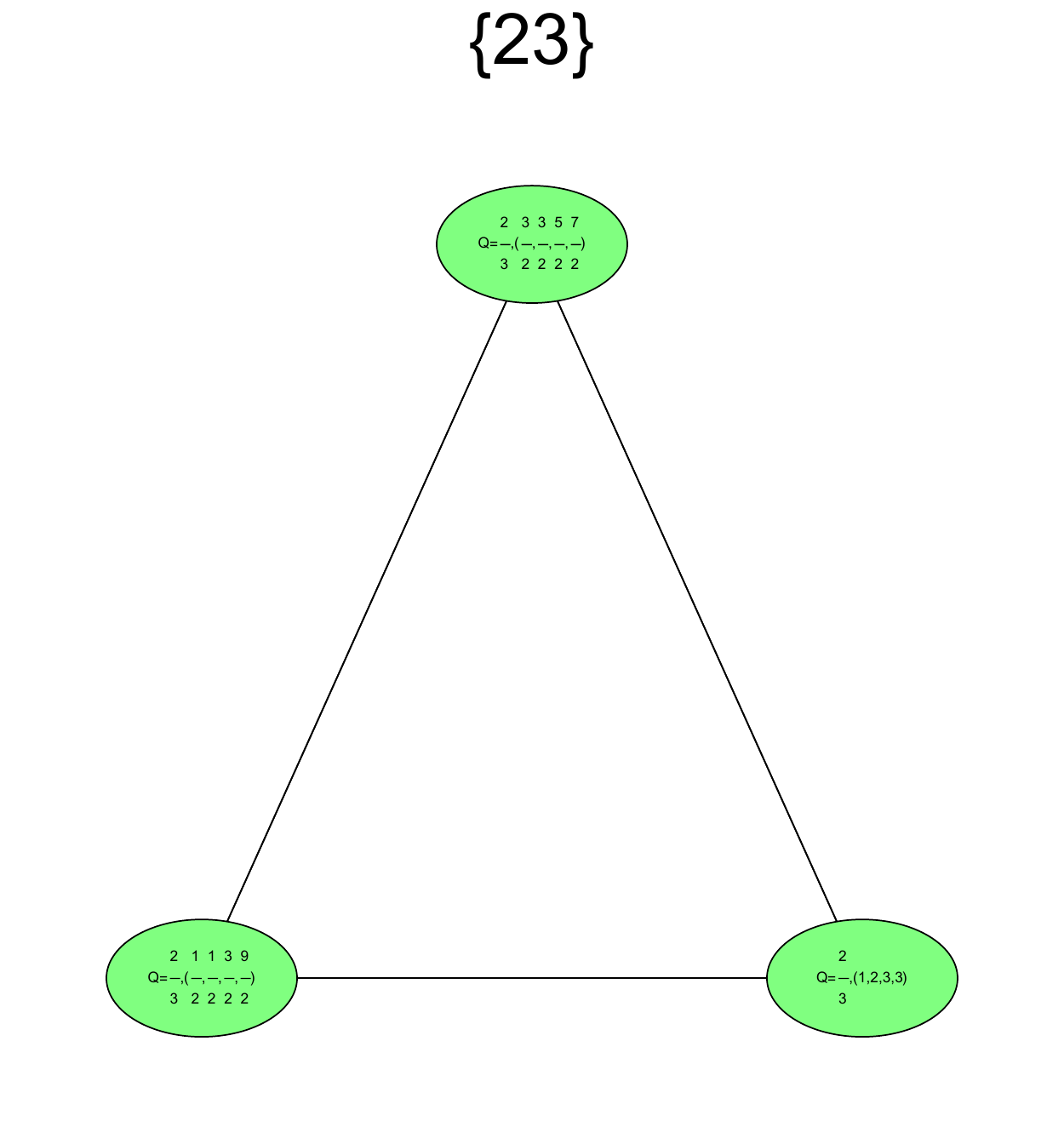}}
\scalebox{0.30}{\includegraphics{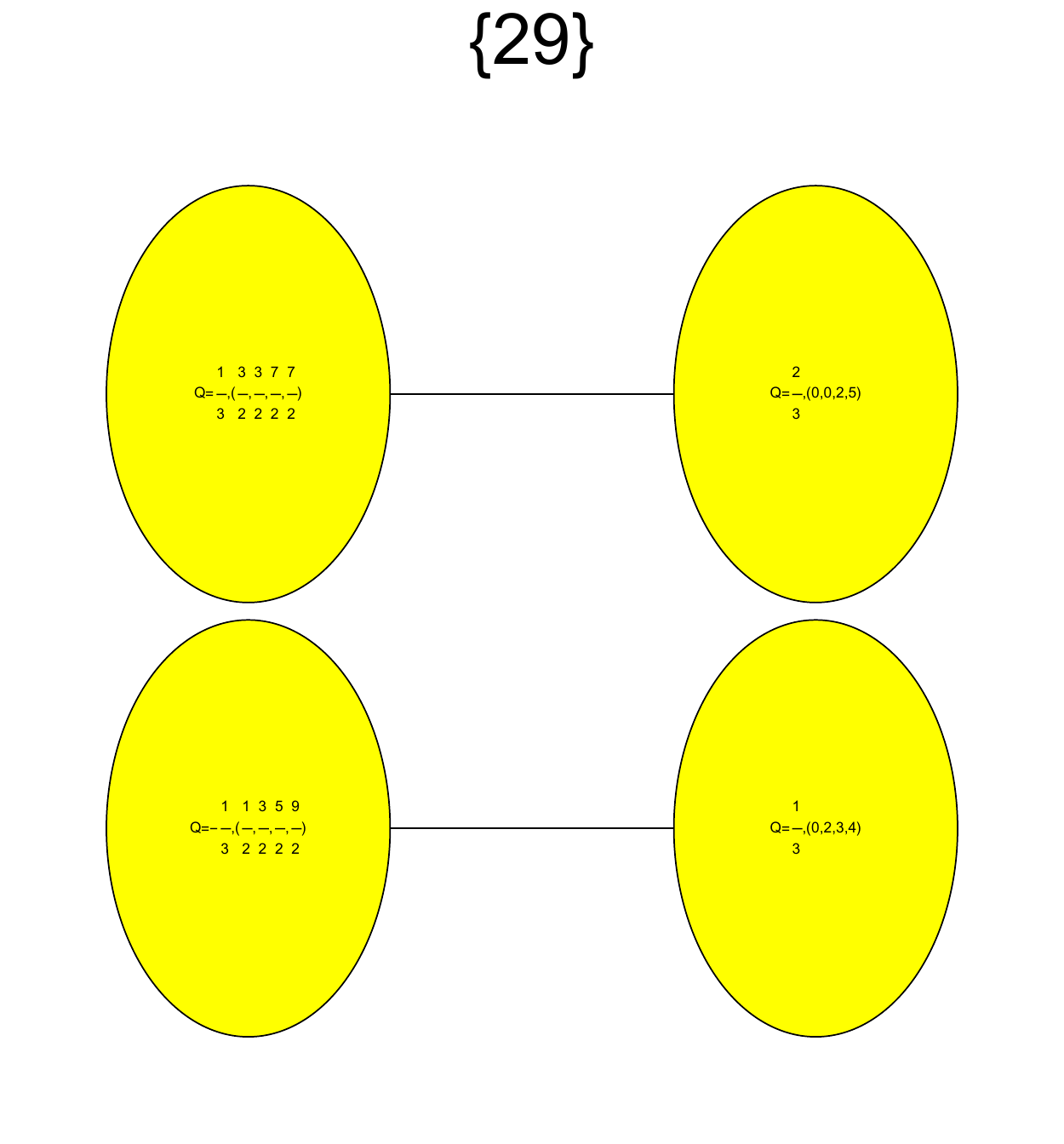}}
\scalebox{0.30}{\includegraphics{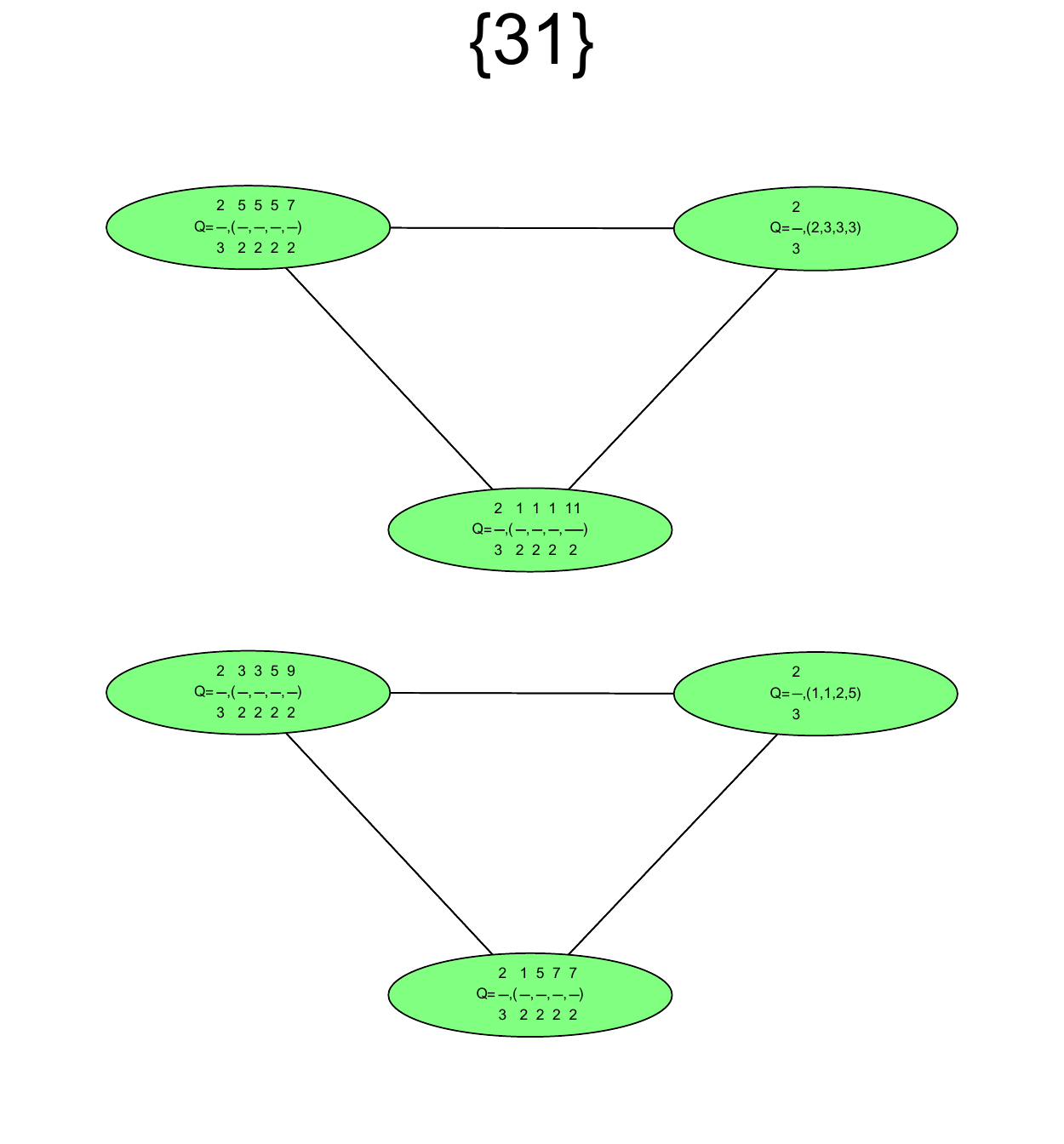}}
\scalebox{0.30}{\includegraphics{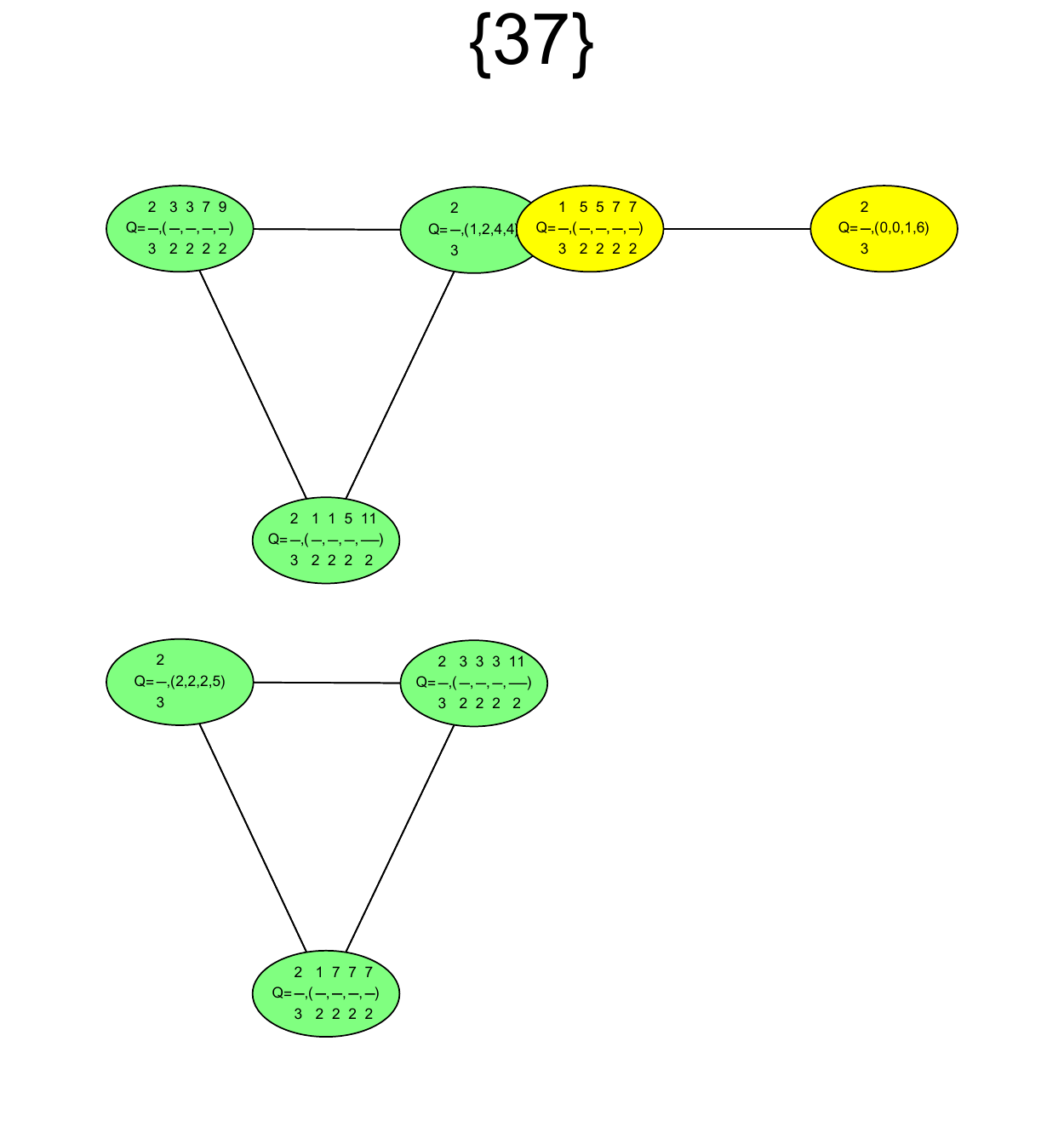}}
\scalebox{0.30}{\includegraphics{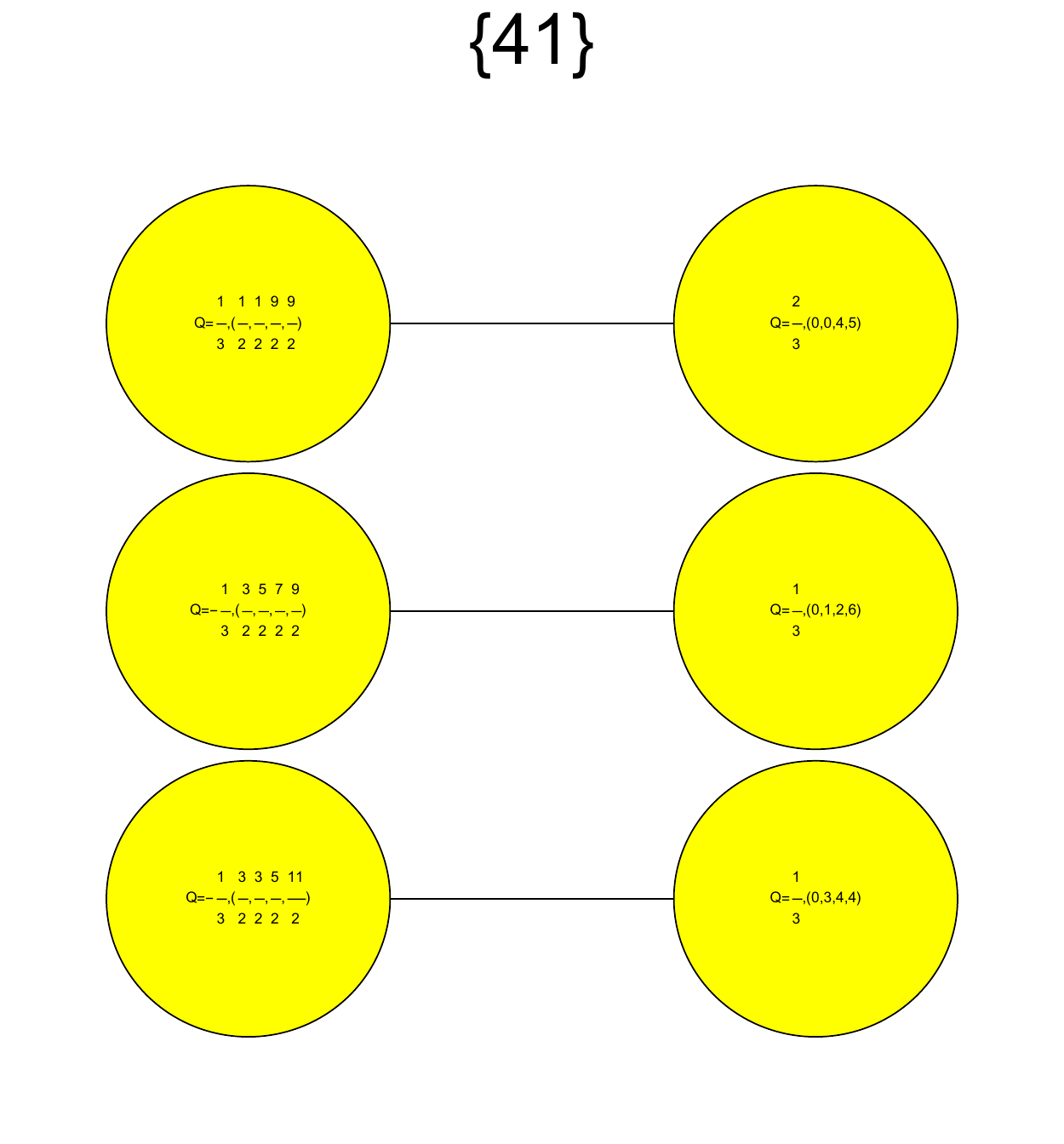}}
\scalebox{0.30}{\includegraphics{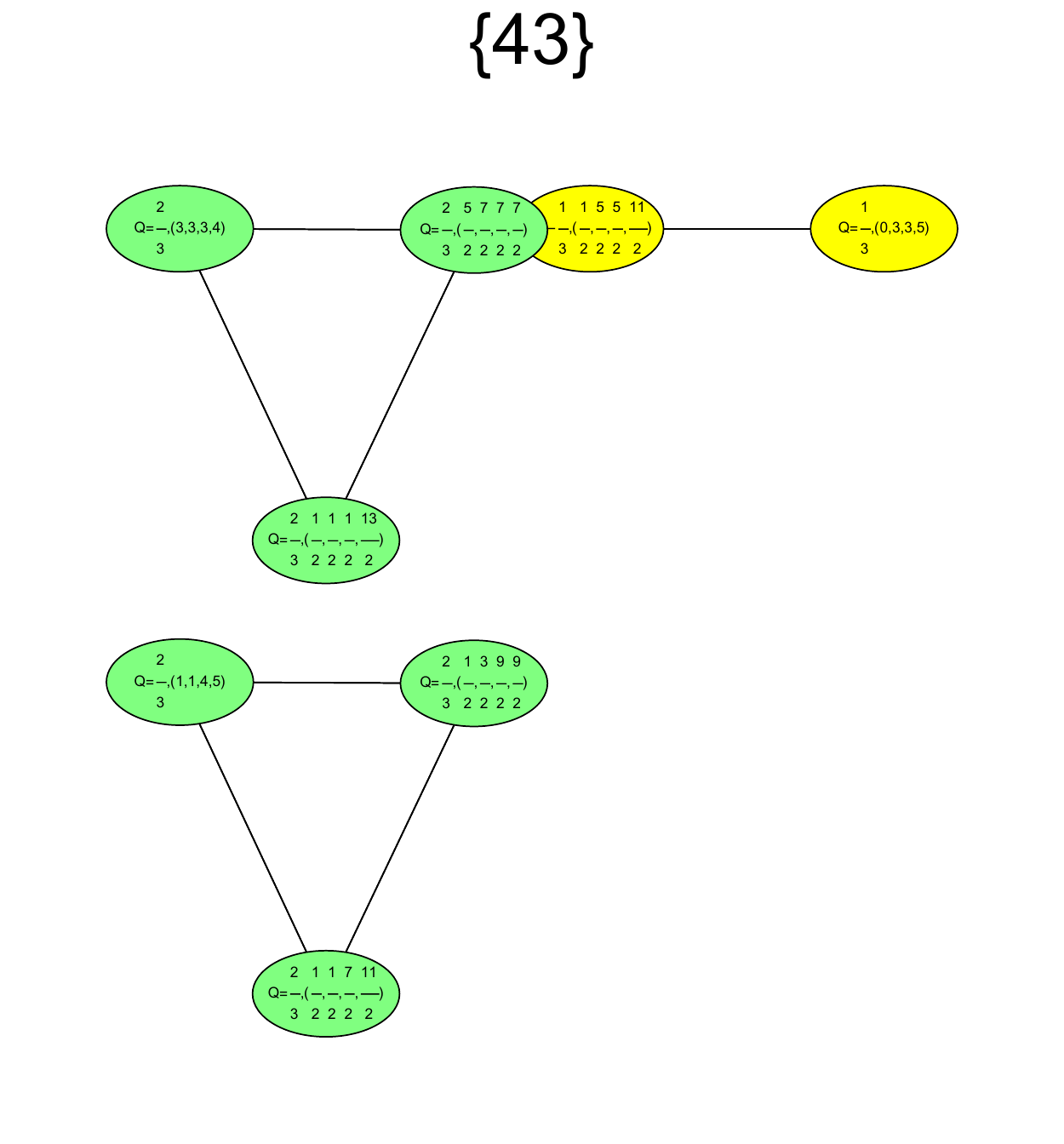}}
\scalebox{0.30}{\includegraphics{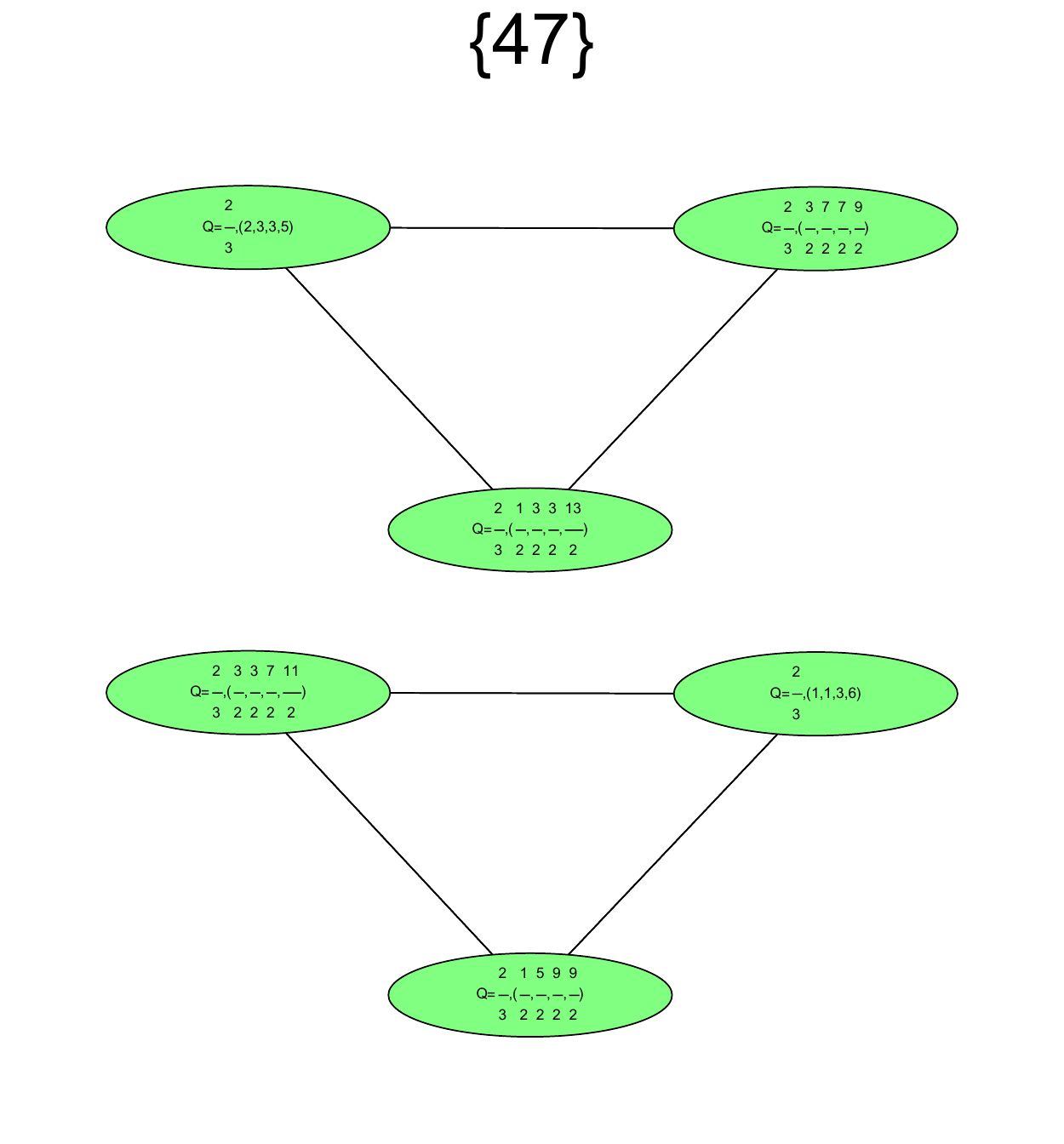}}
\scalebox{0.30}{\includegraphics{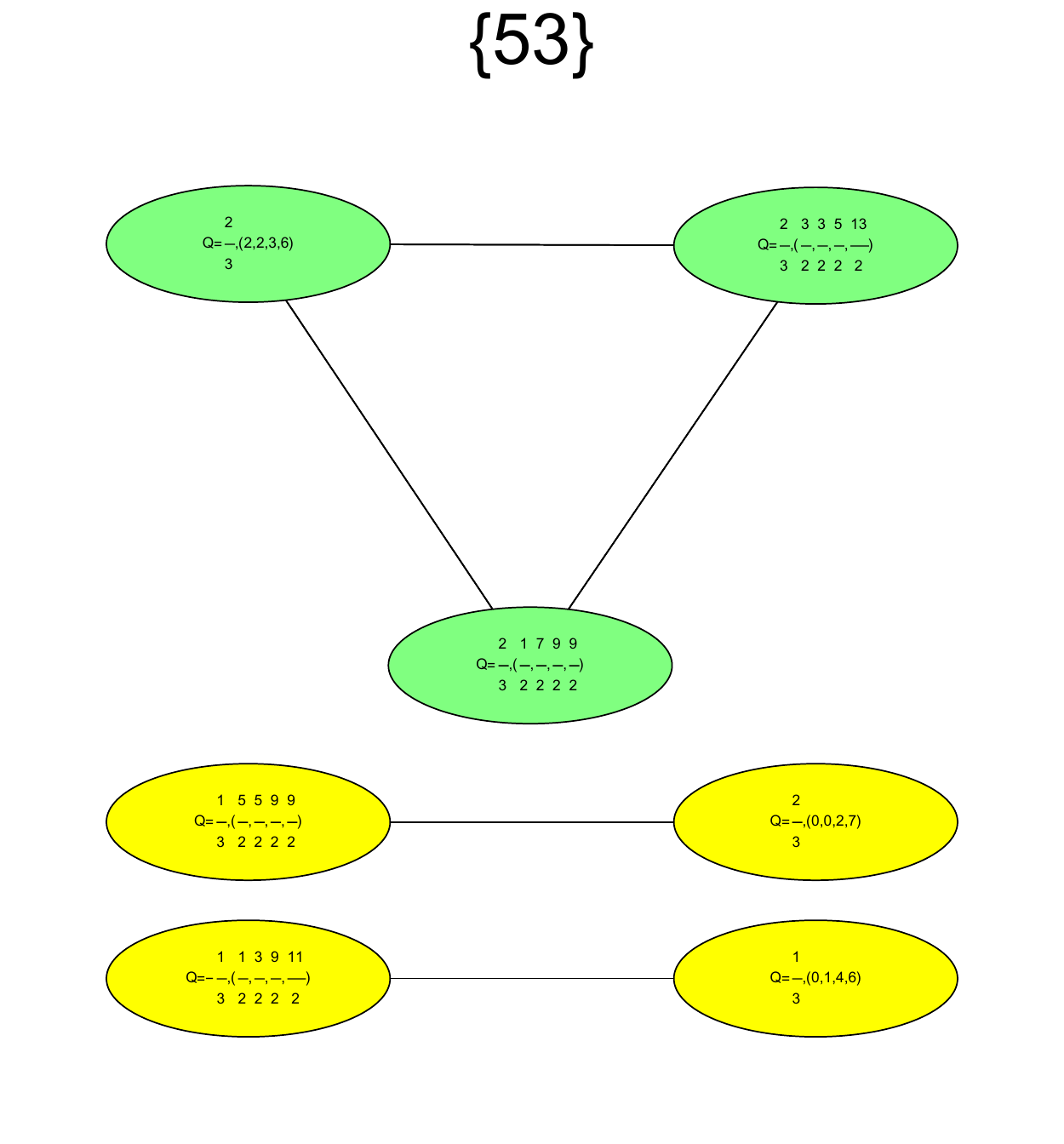}}
\scalebox{0.30}{\includegraphics{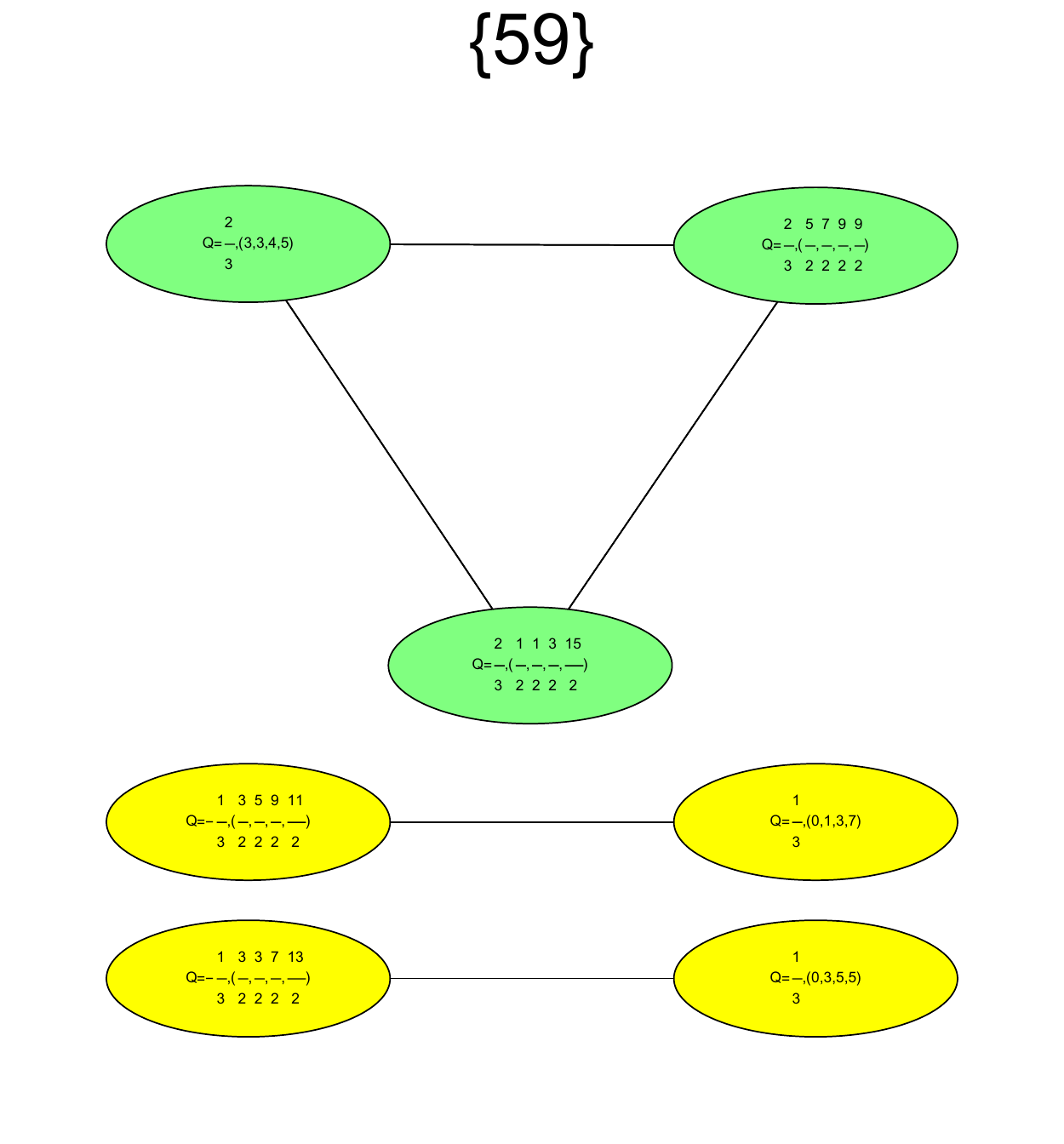}}
\scalebox{0.30}{\includegraphics{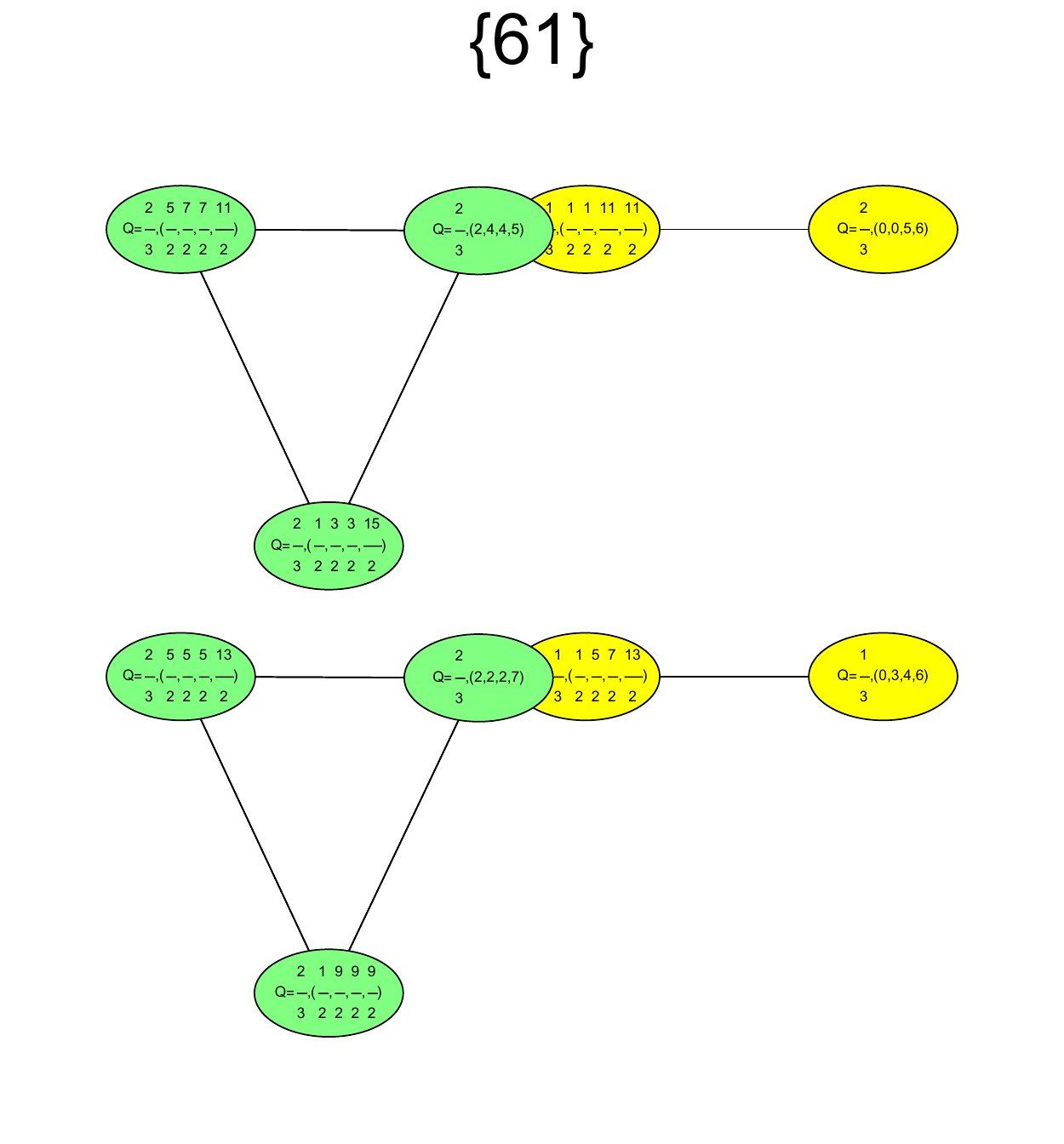}}
\scalebox{0.30}{\includegraphics{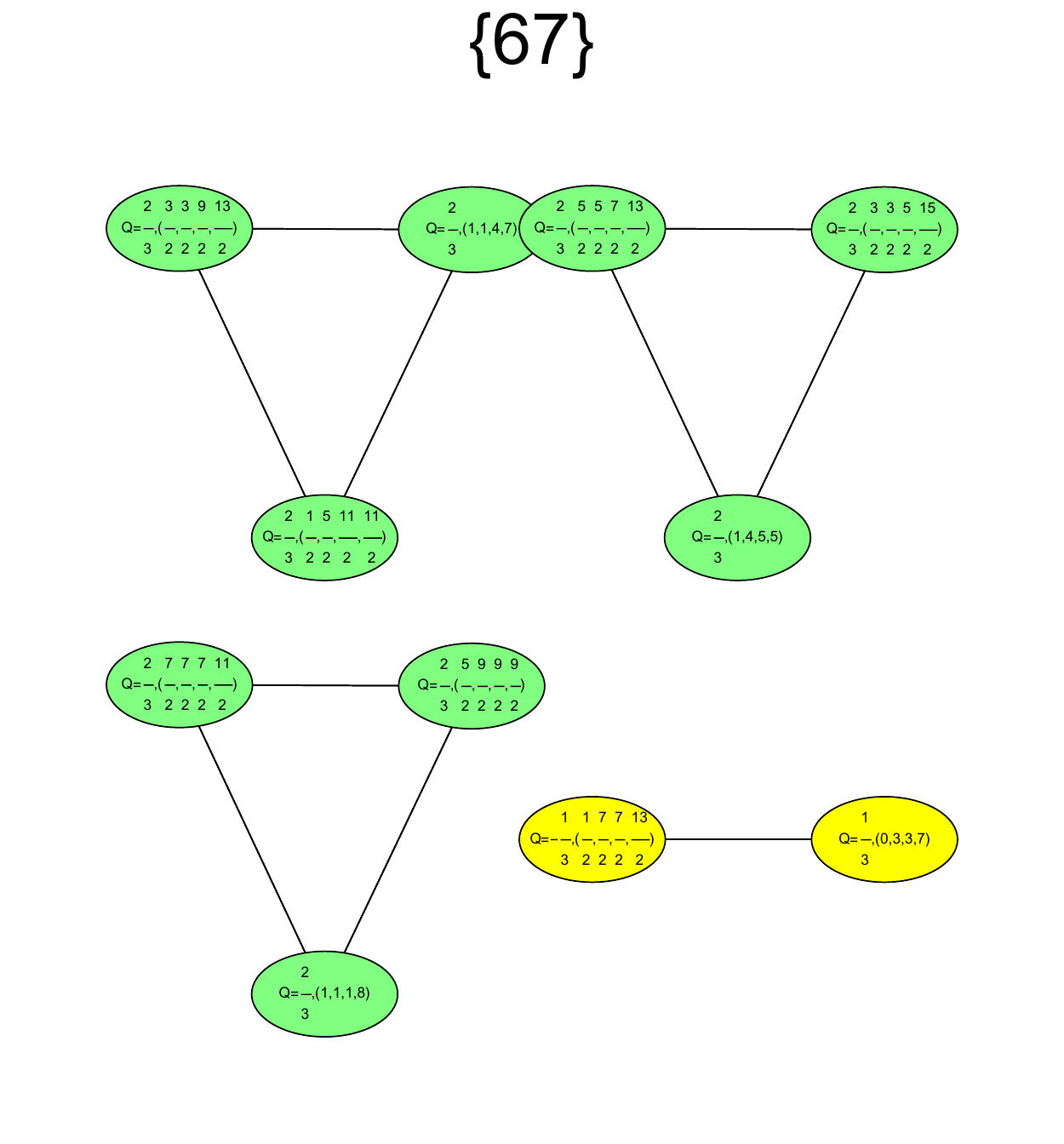}}
\scalebox{0.30}{\includegraphics{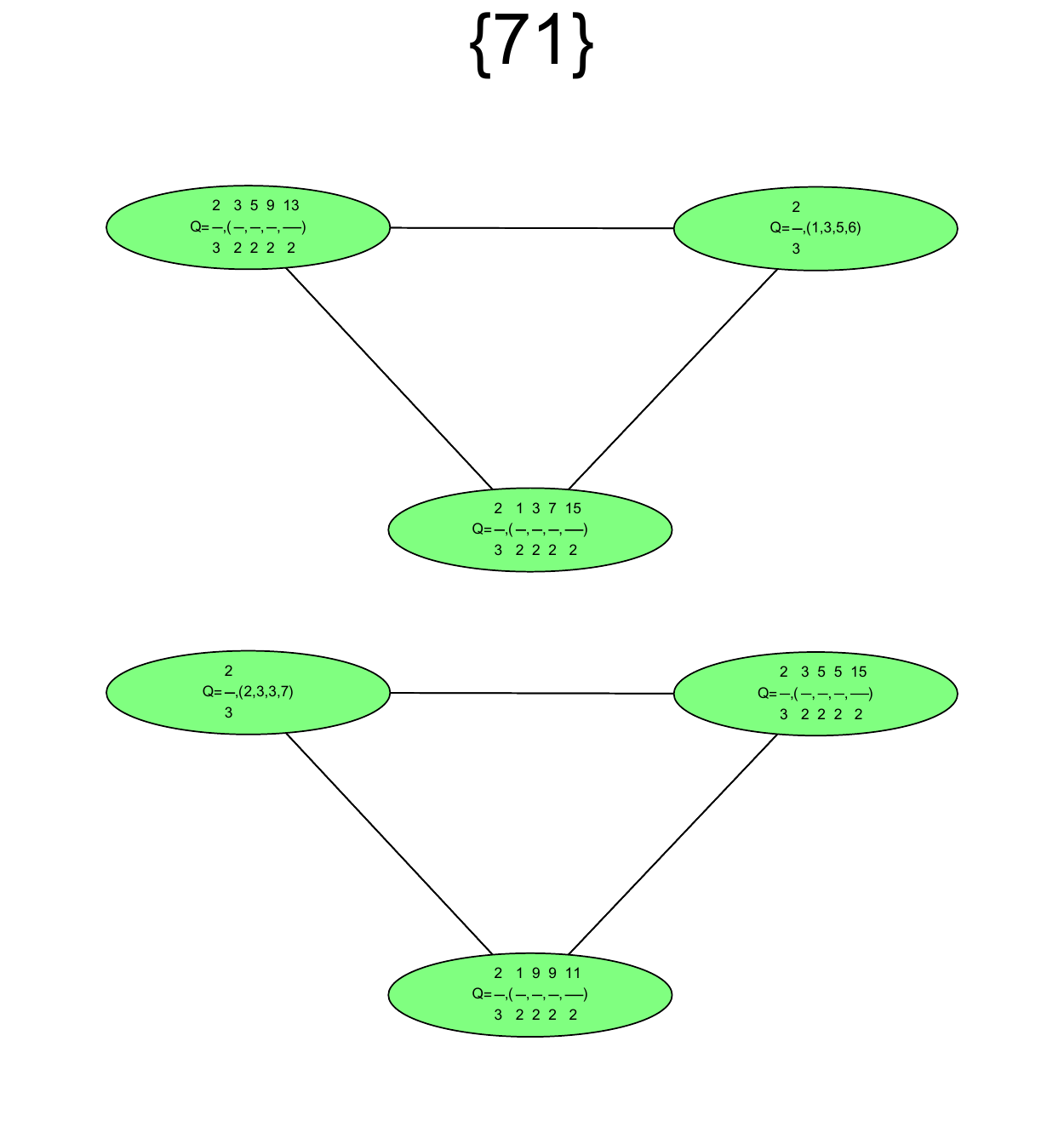}}
\caption{
\label{snowflake}
Hadrons with charges. 
}
\end{figure}

The $p=2$ is special $(1,1,0,0)$ and not included in the above remarks. It is neutral
and equivalent to its anti-particle. It is "Higgs like", as in the lepton case and
we can not place it yet. Also, we like to think about a situation, where space is actually 
a dyadic group of integers (a much more natural space as it is compact, features a smallest
translation), here scaling by a factor $2$ is a symmetry and multiplication by $2$ 
makes the grid finer (making the $2$-adic norm small). The mechanism of mass can anyway
only be understood when looking at dynamical setups, where particles travel. When looking
at wave equations in dyadic groups, the multiplication by 2 plays a special role and it can 
slow down particles, similar as mass does. \\

For $p=3$, $(0,1,1,1),(1,1,1,3)/2$ form a Meson of charge $1=2/3+1/3$.
For $p=5$, we have a Meson $(0,0,1,2),(1,1,3,3)/2$ of charge $1$. 
For $p=13$, we have a Baryon $(1,1,1,7)/2,(1,2,2,2),(3,3,3,5)/2$ of charge $0$ 
and a Meson $(0,0,2,3),(1,1,5,5)/2$ of charge $1$.
For $p=41$ we have a meson $(0,3,4,4),(3,3,5,11)/2$ of charge $\pm 1$. 
and a meson $(0,1,2,6),(3,5,7,9)/2)$ for which all coordinates are different. 

The integer units $i,j,k$ are the gluons while the $(1+i+j+k)/2$ etc are vector bosons. 
In the complex the $1,i$ generate photons. The $(1+i)$ is the Higgs in all cases.
Can not have Lipschitz prime $(a,a,b,b)$ as this is $2a^2+2b^2 = 2(a^2+b^2)$ which is 
$0,2$ modulo $4$.  You see in the figures some pictures. We are able to attach to any 
prime a collection of Baryons, for which the charges is determined. We believe it would
be interesting to study the combinatorics of this setup more. 

\bibliographystyle{plain}

\end{document}